\documentclass[12pt, footinclude=true,headinclude=true,  twoside=semi, headsepline=true, chapterprefix=false, appendixprefix=false, headings=big, numbers=enddot]{scrbook}
\usepackage[english, basque]{babel}
\usepackage{tocloft}

\usepackage{footnote}
\usepackage{empheq}
\allowdisplaybreaks[1]

\usepackage{hyphenat}

\usepackage{algpseudocode}
\usepackage{stmaryrd}

\usepackage{gensymb}

\usepackage{amsmath}
\usepackage{bm}

\usepackage{listings}
\usepackage{color} 

\usepackage{xcolor, colortbl}
 
\definecolor{forestgreen}{rgb}{0.13, 0.55, 0.13}
\DeclareOldFontCommand{\bf}{\normalfont\bfseries}{\mathbf}

\usepackage[colorlinks=true,  linkcolor={blue!50!black}, citecolor={red!50!black}, urlcolor={green!60!black}, menucolor=black,pdfhighlight=/I,draft=false,pagebackref=false, hypertexnames=false]{hyperref}
\usepackage[a4paper,textwidth=14.5cm]{geometry}
\usepackage{amsmath,amsthm,mathrsfs,amssymb,amsfonts,dsfont,nicefrac,stmaryrd,yhmath, verbatim, xcolor, graphicx, tikz, pgfplots, mathtools, mathrsfs, ifthen, import, float, subcaption, bbm, mathtools, csquotes, centernot,aliascnt}

\pgfplotsset{compat = newest}
\usetikzlibrary{cd}
\usetikzlibrary{babel}
\usetikzlibrary{3d, shapes.arrows}
\usepackage[inline]{enumitem}
\usepackage[titletoc]{appendix}
\usepackage[ruled,vlined]{algorithm2e}

\usepackage{adjustbox}
\def\layersep{2.5cm}

\usepackage{float}
\usepackage[section,below, above]{placeins}
\usepackage{framed} 

 \usepackage{pgfplotstable} 
 \usepackage{booktabs}      
  \usepackage{longtable}     
 \pgfplotstableset{
 begin table=\begin{longtable},
 end table=\end{longtable},
 }
 
\usepackage[automark]{scrlayer-scrpage}

\KOMAoptions{headsepline=true}
\AddToLayerPageStyleOptions{scrheadings}
  {oninit={
    \Ifstr{\headmark}{}{\KOMAoptions{headsepline=false}}{}%
  }}

\clearpairofpagestyles

\theoremstyle{definition}

\newtheorem{example}{Example}[chapter]
\usepackage[noabbrev,capitalise,nameinlink]{cleveref}
\Crefname{example}{Example}{Examples}

\newtheoremstyle{theoremdd}
  {\topsep}
  {\topsep}
  {}
  {0pt}
  {\bfseries}
  {.}
  { }
  {\thmnumber{#2}\autodot\thmname{ #1}\textnormal{\thmnote{ (#3)}}}
  
\theoremstyle{theoremdd}

\usepackage{listofitems} 

\usepackage[outline]{contour} 

\pgfplotsset{every axis legend/.append style={cells={anchor=west} }}

\usepgfplotslibrary{external}
\usepackage{multirow,booktabs}

\usepackage[labelfont=bf]{caption}
\usepackage{mfirstuc}

\usepackage{acro}

\acsetup{
make-links=true,
link-only-first=true,
list/name = Acronyms,
}
\DeclareAcronym{PDE}{
  short=PDE,
  long=Partial Differential Equation,
}

\DeclareAcronym{BVP}{
  short=BVP,
  long=Boundary-Value Problem,
}

\DeclareAcronym{FDM}{
  short=FDM,
  long=Finite Difference Method,
}

\DeclareAcronym{FEM}{
  short=FEM,
  long=Finite Element Method,
}

\DeclareAcronym{FFNN}{
  short=FFNN,
  long=Feed-Forward Neural Network,
}

\DeclareAcronym{ANN}{
  short=ANN,
  long=Artificial Neural Network,
}

\DeclareAcronym{NN}{
  short=NN,
  long=Neural Network,
}

\DeclareAcronym{DNN}{
  short=DNN,
  long=Deep Neural Network,
}

\DeclareAcronym{UAT}{
  short=UAT,
  long=Universal Approximation Theorem,
}

\DeclareAcronym{TF}{
  short=TF,
  long=TensorFlow,
}

\DeclareAcronym{TF2}{
  short=TF2,
  long=TensorFlow 2,
}

\DeclareAcronym{SGD}{
  short=SGD,
  long=Stochastic Gradient Descent,
}

\DeclareAcronym{SGA}{
  short=SGA,
  long=Stochastic Gradient Ascent,
}

\DeclareAcronym{SGDM}{
  short=SGDM,
  long=SGD with Momentum,
}

\DeclareAcronym{MC}{
  short=MC,
  long=Monte Carlo
}

\DeclareAcronym{AD}{
  short=AD,
  long=Automatic Differentiation,
}

\DeclareAcronym{PINN}{
  short=PINN,
  long=Physics-Informed Neural Network,
}

\DeclareAcronym{ODE}{
  short=ODE,
  long=Ordinary Differential Equation,
}

\DeclareAcronym{D2RM}{
  short=D\ensuremath{^2}RM,
  long=Deep Double Ritz Method,
}

\DeclareAcronym{DRM}{
  short=DRM,
  long=Deep Ritz Method,
}

\DeclareAcronym{DeepFEM}{
  short=DeepFEM,
  long=Deep Finite Element Method,
}

\DeclareAcronym{1D}{
  short=1D,
  long=one-dimensional,
}

\DeclareAcronym{2D}{
  short=2D,
  long=two-dimensional,
}

\DeclareAcronym{3D}{
  short=3D,
  long=three-dimensional,
}

\DeclareAcronym{WANs}{
  short=WANs,
  long=Weak Adversarial Networks,
}

\DeclareAcronym{GDRM}{
  short=GDRM,
  long=Generalized Deep Ritz Method,
}

\DeclareAcronym{GPU}{
  short=GPU,
  long=Graphics Processing Unit,
}

\begin{document}


\newcommand{\titleEN}{Solving Partial Differential Equations Using Artificial Neural Networks}
\renewcommand{\date}{2023}
\renewcommand{\author}{Carlos Uriarte}
\newcommand{\supervisors}{David Pardo, Elisabete Alberdi}

\renewcommand{\thefootnote}{\fnsymbol{footnote}}
\selectlanguage{english}

\begin{titlepage}

\vfill
\begin{center}
\includegraphics[height=2.5cm]{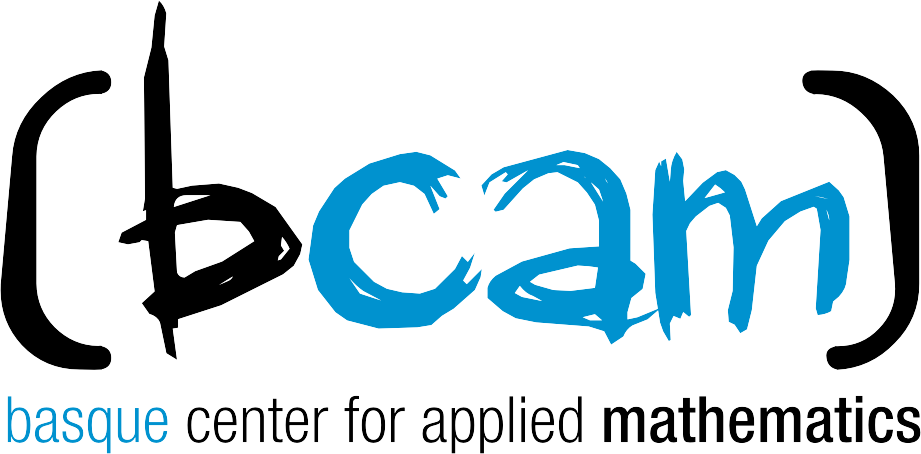}
\hfill
\includegraphics[height=2.5cm]{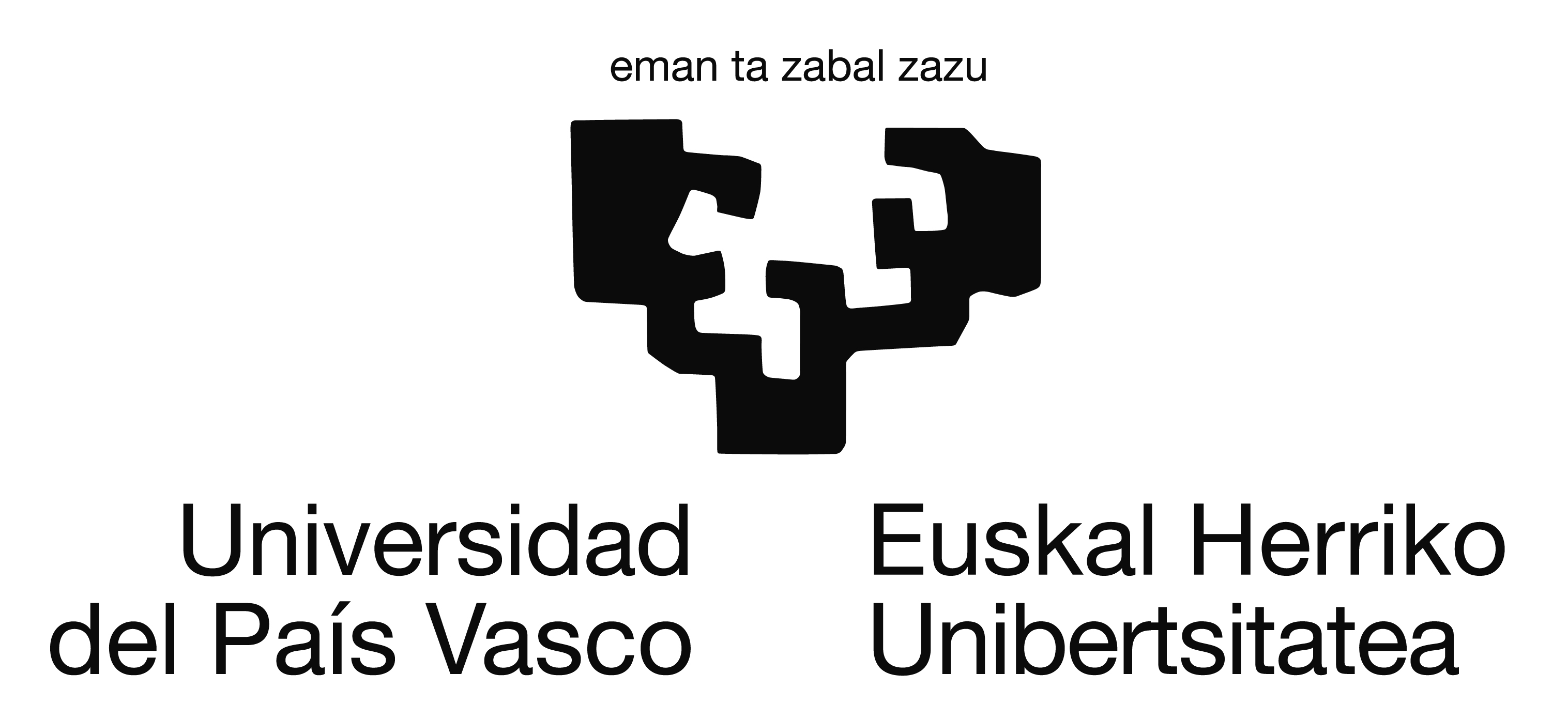}\\
\vspace{1.5cm}
{\Huge {\textsc{{\begin{minipage}{\textwidth}\hrulefill\begin{center}\vspace{0.3cm}\titleEN\end{center} \hrulefill \end{minipage}}}}}
\end{center}
\begin{center}
\vspace{0.5cm}
by\\
\vspace{1cm}
{\Huge{\textsc{\author}}}\\
\vspace{1cm}
\textit{A dissertation submitted to the University of the Basque Country (UPV/EHU) in fulfillment of the requirements for the degree of Doctor of Philosophy}\footnote{This document is an unofficial, revised, and English-only version of the original bilingual dissertation submitted to the University Repository prior to defense. The official version is available at \url{https://addi.ehu.es/handle/10810/68335}.}\\
\vspace{1cm}
{\large\textsc{Supervisors: \supervisors}\\
\vspace{1cm}
\date\footnote{Original deposit: December 2023. Last update: November 2024. Erratum reports sent to \texttt{carlos.uribar@gmail.com} are greatly appreciated.}
\vfill
}

\end{center}
\end{titlepage}

\renewcommand{\thefootnote}{\arabic{footnote}}


\frontmatter

\selectlanguage{english}

\clearpage
 
\selectlanguage{english}
\begingroup
\tableofcontents
\addtocontents{toc}{\protect\setcounter{tocdepth}{2}}


\selectlanguage{english}
\chapter{Acknowledgments}

Over the period of completion of this dissertation, I have met people to whom I will always be grateful for many and diverse reasons.

First of all, I want to express my deepest gratitude to David: a singular guy with an assertive and carefree attitude, and an unusual way of looking at life.  These, combined with his patience, incredible interpersonal skills, and great intuition, make him a wonderful person to work with, be advised by, and enjoy funny and not-so-funny moments.  I thank him for: his continuous teachings,  our long conversations, and that innate discipline of him to destroy works that had cost me so much to build---to tell the truth, he used to be right.  I will be eternally grateful to him for the successes and stumbles I have found along this exciting and unrepeatable PhDing journey. 

I would also like to thank my other supervisor, Elisabete, for bringing me the opportunity to develop the thesis under her supervision.  I would especially like to thank her for her impeccable follow-up during the organization and writing of this dissertation. Without her,  it would have been an infinitely more complex task. 

In addition to my supervisors, much of the credit for the sublime working environment goes to the rest of my research colleagues.  In particular,  Javier was a fundamental workmate during the first stage of my PhD, always available to help me virtually during the pandemic.  His virtuosity in tackling tedious tasks without tiring fascinates me, as well as in continuously proposing to have ``a last beer''---which rarely was the last one---after a day of hard work.  I have learned enormously from him and shared unforgettable moments. The same goes for Ana,  who accompanied me from the beginning to the end of my predoctoral career and with whom I have also shared enjoyable moments in every professional and not-so-professional appointment.  To highlight some of the not-so-professional ones: our adventures during the ``Camino de Santiago'' and the ``Aste Nagusiak'' full of glitter.

I am incredibly grateful to Judit for hosting me in such an endearing way during my research stay at the Oden Institute for Computational Sciences and Engineering of the University of Texas at Austin, USA.  It was an adventure I enjoyed a lot.  I also thank her for introducing me to Prof. Demkowicz,  from whom I learned a great deal of knowledge I have applied in my research works.  Similarly, I thank Ignacio,  Paulina, Keko, and Patrick for their incredible kindness and hospitality during my research stay at the Institute of Mathematics of the Pontifical Catholic University of Valpara\'iso, Chile. I have many good anecdotes engraved in my memory in which they were taking part, particularly at moments of ``carreteo''.

It is difficult to sufficiently acknowledge all the lived experiences with so many colleagues during these four years.  To avoid making this section too long,  I briefly thank Jamie, Jon Ander, Julen, Lena, Oscar, Felipe, Mahdi, Ali, Jes\'us, Manuela, and Tom\'as.  Although the greeting is brief, it has a lot of meaning. I thank them all a lot.

I would also like to thank all the BCAM staff for assisting me during the bureaucratic difficulties encountered, which have not been few. In this regard, I want to thank Miguel especially. He has been infinitely patient and lovely with me on many occasions.  A thousand thanks!

Besides,  during these four years, I have met people who, although not directly related to my research environment,  have become friends of mine due to their relationships with workmates and significantly contributed to the easygoingness of this adventure. I particularly want to thank Claudio and Mario for our many good times together.

Finally, I want to thank my family and friends for their support.  In particular, I want to offer my eternal gratitude to my father, Enrique; my brother, Luis; and my partner, Iratxe, who unconditionally accompanied me and firsthand experienced my moments of tension and joy.  Thank you for understanding the conflict of a researcher lifestyle---working at ungodly hours, leaving aside the usual healthy social habits,  simply because I felt I could not procrastinate my moments of research inspiration. Heartfeltly, thanks.
\selectlanguage{english}


\selectlanguage{english}

\chapter{Abstract}

Partial differential equations have a wide range of applications in modeling multiple physical, biological, or social phenomena. Therefore, we need to approximate the solutions of these equations in computationally feasible terms. Nowadays, among the most popular numerical methods for solving partial differential equations in engineering, we encounter the finite difference and finite element methods. An alternative numerical method that has recently gained popularity for numerically solving partial differential equations is the use of artificial neural networks.

Artificial neural networks, or neural networks for short, are mathematical structures with universal approximation properties. In addition, thanks to the extraordinary computational development of the last decade, neural networks have become accessible and powerful numerical methods for engineers and researchers. For example, imaging and language processing are applications of neural networks today that show sublime performance inconceivable years ago.

This dissertation contributes to the numerical solution of partial differential equations using neural networks with the following two-fold objective: investigate the behavior of neural networks as approximators of solutions of partial differential equations and propose neural-network-based methods for frameworks that are hardly addressable via traditional numerical methods. 

As novel neural-network-based proposals, we first present a method inspired by the finite element method when applying mesh refinements to solve parametric problems. Secondly, we propose a general residual minimization scheme based on a generalized version of the Ritz method. Finally, we develop a memory-based strategy to overcome a usual numerical integration limitation when using neural networks to solve partial differential equations.

\selectlanguage{english}

\cleardoublepage
\phantomsection
\addcontentsline{toc}{chapter}{Acronyms}
\printacronyms
\newpage
\listoffigures
\newpage
\listoftables
\endgroup

\mainmatter

\cleardoublepage
\phantomsection


\chapter{Introduction} \label{chapter1}

\section{Motivation}\label[section]{section1.1}

\emph{\acfp{PDE}} are of great value to society due to their broad applicability in modeling multiple biological, physical, or social phenomena \cite{gershenfeld1999nature,strauss2007partial,farlow1993partial,salsa2016partial}. Using derivatives in these equations allows us to describe complex relationships and their rates of change over time and/or space.  For example, they allow to model: heat transfer \cite{incropera1996fundamentals, hahn2012heat},  electromagnetic fields \cite{jackson1999classical, van2007electromagnetic}, fluid dynamics \cite{batchelor1967introduction,anderson1995computational, versteeg2007introduction},  population evolution \cite{holmes1994partial, murray2002mathematical, brauer2011mathematical},  and financial \cite{baxter1996financial, hull2015options} or health-care \cite{saltzman2009biomedical, lysaker2003noise, oden2010general} forecasts.

However, establishing an equation to describe a physical or social situation is only the first step. There arise different kinds of problems when modeling via PDEs, and the way we treat or approach them depends on the available computational resources.

In \emph{forward problems}, we are interested in determining the function that satisfies a PDE given some initial or boundary conditions.  In this way,  the solution to forward problems can provide insights into the modeled phenomena and can be used to make predictions, usually via post-processed measurements taken from the solution.

In \emph{inverse problems}, we have measurements of the solution, but we need to discover the value of specific parameters in the equation. This is an indirect problem, as we need to determine the parameters that provide a PDE solution that agrees with the given measurements \cite{tarantola2005inverse}.  \Cref{Chapter1_electromagnetic_fields} illustrates a possible framework of forward and inverse problems in electromagnetic fields \cite{feynman1986feynman}.

\begin{example}[Electromagnetic fields]\label{Chapter1_electromagnetic_fields}
Maxwell's equations model electromagnetic fields according to four well-known physical laws \cite{feynman1986feynman, jackson1999classical}:
\begin{equation}
\begin{cases}
\nabla\times\mathbf{E} = -j\omega\boldsymbol{\mu} \ \mathbf{H} - \mathbf{M}, &\text{Faraday's Law},\\
\nabla\times\mathbf{H} = (\boldsymbol{\sigma}+j\omega\boldsymbol{\varepsilon})\mathbf{E} + \mathbf{J},&\text{Amper\`e's Law},\\
\nabla\cdot (\boldsymbol{\varepsilon\mathbf{E}}) = \rho_f,&\text{Gauss' Law of Electricity},\\
\nabla\cdot (\boldsymbol{\mu\mathbf{H}}) = 0, &\text{Gauss' Law of Magnetism}.
\end{cases}
\end{equation} Here, $\mathbf{E}$ and $\mathbf{H}$ denote the electric and magnetic fields in the frequency domain, respectively; $\mathbf{J}$ and $\mathbf{M}$ are given source terms; $\boldsymbol{\sigma}$ stands for the electrical conductivity of the media, $\boldsymbol{\varepsilon}$ for the electrical permittivity,  $\boldsymbol{\mu}$ for the magnetic permeability,  $\rho_f$ for the electric charge density, $j$ for the imaginary unit,  and $\omega$ for the angular frequency. Then, a possible forward-inverse scenario on electromagnetic fields could consist of parameters $\{\boldsymbol{\sigma},\boldsymbol{\mu}, \boldsymbol{\varepsilon}\}$ and measurements $\{\mathbf{Z}\}$, where $\mathbf{Z}$ denotes the impedance tensor defined by $\mathbf{E}=\mathbf{Z}\mathbf{H}$.  For further details, see, e.g., \cite{alvarez2015hp}. \autoref{Forward_inverse_sketch} shows a graphic summary of this modeled framework. 

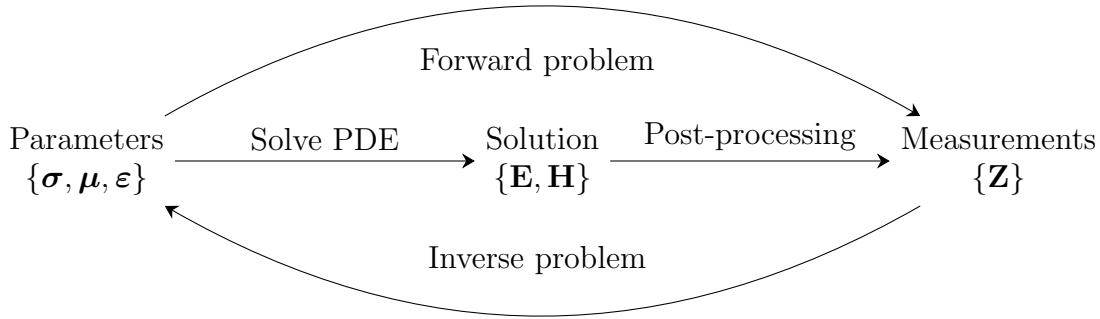
\begin{figure}[htbp]
\centering
\begin{tikzpicture}[node distance=6cm, every text node part/.style={align=center}]
    \node (params) {Parameters \\ $\{\boldsymbol{\sigma},\boldsymbol{\mu},\boldsymbol{\varepsilon}\}$};
    \node[right of=params] (solution) {Solution \\ $\{\mathbf{E}, \mathbf{H}\}$};
    \node[right of=solution] (measu) {Measurements \\ $\{\mathbf{Z}\}$};

    \draw[-{Stealth[width=2mm]}] (params) -- (solution) node[midway, above] {Solve PDE};
    \draw[-{Stealth[width=2mm]}] (solution) -- (measu) node[midway, above] {Post-processing};
    \draw[-{Stealth[width=2mm]},  bend left=30] (params) to node[midway, below=4mm] {Forward problem} (measu);
    \draw[-{Stealth[width=2mm]}, bend left=30] (measu) to node[midway, above=4mm] {Inverse problem} (params);
\end{tikzpicture}
\caption{Forward and inverse problems sketch in electromagnetic fields.}
\label{Forward_inverse_sketch}
\end{figure}
\end{example}

Solving forward problems requires, in general, fewer resources than inverse problems.  Given a choice of parameters, a forward problem typically consists of performing a single PDE simulation.  In contrast, an inverse problem usually requires the iterative resolution of forward problems, conveniently readjusting the choice of parameters at each iteration based on the outcome of previous forward simulations \cite{tarantola2005inverse, isakov2006inverse}. 

Connecting both types of problems,  we have \emph{parametric problems}, where the aim is to study the behavior of the solution/measurements as a function of the parameters. The solution can depend on one or several parameters, and we can study how changes in these parameters affect the behavior of the system.  Alternatively,  the parametric approach can be viewed as a forward problem where the parameters of the equation take the role of variables of the solution.  Typically,  the parametric forward problem is nonlinear with respect to the introduced new variables,  enormously escalating and hindering its resolvability.

In this dissertation, we will only address the task of solving parametric and non-parametric linear PDEs,  which is critical in many applications and frameworks.  We will indistinctly call PDEs both one-dimensional---also known as \emph{\acfp{ODE}}---and higher-dimensional differential equations.

\section{Traditional Numerical Methods}\label[section]{section1.2}

There are various methods for solving PDEs, which can be broadly classified into two categories: analytic and numerical methods.

Analytic methods involve using mathematical techniques such as separation of variables, integral transforms, or complex analysis to find an exact solution to the PDE \cite{sneddon2006elements, haberman2012applied,tikhonov2013equations}. However, these methods are often limited to relatively simple PDEs and geometries.  Additionally, even if the PDE is analytically solvable, the resulting solution may be expressed in a complicated form (e.g.,  as an infinite power series), which can be difficult to interpret and apply in practical scenarios.

On the other hand,  numerical methods involve approximating the solution to the PDE using a computational algorithm. Generally, these methods propose approximations that are easily interpretable or numerically recoverable (e.g.,  as a finite linear combination of simple prescribed functions).  Below, we briefly review the main aspects of three widely known and employed kinds of numerical methods for solving PDEs.

The \emph{\acf{FDM}} discretizes the spatial domain and/or time interval of the PDE into a finite number of subdomains \cite{smith1985numerical, press2007numerical, leveque2007finite}. Then, the evaluation of the derivatives at these discrete points is approximated by solving algebraic equations containing finite differences and values from nearby points.  The FDM converts a PDE into a system of linear equations that can be solved by matrix algebra techniques.  It is conceptually simple but challenging in design for complex problems or irregular domains. In particular, it suffers from the so-called \emph{curse of dimensionality} \cite{bellman1957dynamic,  bellman1961adaptive},  where the size of the involved matrix grows exponentially with the dimension of the problem.   

Similarly,  the \emph{\acf{FEM}} \cite{hughes1987finite, brenner2008mathematical, logan2010first, reddy2019introduction} proposes a mesh-based scheme where the approximated solution is in the form of a finite linear combination of some prescribed functions (typically,  piecewise polynomial) with local support.  The support is purposely designed on specific regions of the domain, called elements,  so the resulting matrix is sparse. The FEM is particularly suitable for dealing with PDEs in variational form but suffers from similar difficulties as the FDM.

\emph{Spectral methods} \cite{gottlieb1977numerical, boyd2001chebyshev, canuto2007spectral} are a class of numerical methods that approximate functions using a linear combination of orthogonal basis functions with global support, such as Chebyshev polynomials or Fourier series. These methods can achieve very high accuracy, and are particularly useful for smooth solutions with high oscillatory behavior, but may become computationally unfeasible for complex geometries.

In essence,  the numerical methods mentioned above are based on pre-establishing a finite-dimensional basis and parameterizing the approximation via the coefficients of the corresponding linear combination.  When using a mesh-based basis (in the FDM and the FEM),  it must be sufficiently fine to conveniently capture possible complex behaviors of the exact solution, but without falling into the computational disadvantage of over-refining. In spectral methods, the major challenge involves selecting orthogonal basis functions for arbitrary domains.

\section{Neural Networks are Universal Approximators}\label[section]{section1.3}

\emph{\acfp{ANN}}, or \emph{\acfp{NN}} for short, provide an alternative approximation approach to the traditional linear combination of prescribed functions.

The most basic form of NNs are \emph{\acfp{FFNN}} \cite{rosenblatt1958perceptron, ivakhnenko1967cybernetics, amari1967theory,goodfellow2016deep, schmidhuber2022annotated},  which are inspired by the structure and function of biological \acp{NN} in the brain. They consist of interconnected \emph{neurons} that are organized into layers and the information flows in one direction, from one layer to the following one.  Each neuron receives input from neurons in the previous layer and computes a weighted sum.  When such sum is calculated from \emph{all} neurons in the previous layer,  we say that the FFNN is \emph{fully connected}.  This sum is then passed through a non-linear activation function to later transmit the outcome to the neurons in the following layer. We defer the formal description of fully-connected FFNN to \Cref{chapter2}.

Fully-connected FFNNs are universal approximators.  Roughly speaking, this means that they are able to approximate almost any given function to a desired precision.

One of the first versions of the \emph{\acf{UAT}} entails the \emph{arbitrary-width} case, which states that a FFNN with a single hidden layer is a universal approximator of continuous functions as long as it is sufficiently wide. This result was parallelly and independently shown in 1989 for sigmoid activation functions \cite{cybenko1989approximation}, as well as for non-constant, bounded, and monotonically-increasing activation functions \cite{hornik1989multilayer}.  The principal author of the aforementioned last work showed two years later that it is not the specific choice of the activation function but the FFNN architecture itself that endows NNs the potential of being universal approximators \cite{hornik1991approximation}.  Related to this,  it was thereafter shown that the universal approximation property is equivalent to having a non-polynomial activation function \cite{leshno1993multilayer,  pinkus1999approximation}.

The ``dual'' version of the arbitrary-width case for the \ac{UAT} consists of fixing a bounded width and establishing an \emph{arbitrary depth}. Several authors have studied this reformulation since the early 2000s \cite{gripenberg2003approximation,kidger2020universal}, with a particular emphasis on ReLU activation functions \cite{yarotsky2017error, lu2017expressive,hanin2017approximating, hanin2019universal}.  Analogously, we encounter results entailing \emph{bounded-width} and \emph{bounded-depth} assumptions.  For example,  \cite{maiorov1999lower} claims the existence of a sigmoidal-type analytic activation function for a FFNN with depth two and bounded width being a universal approximator.  Additionally,  \cite{guliyev2018approximation} proved that FFNNs with a single hidden layer with bounded width are still universal approximators for univariate functions, which generally does not hold for multivariate functions.

Several extensions to the UAT have been developed over the past three decades dealing with discontinuous activation functions \cite{leshno1993multilayer},  non-compact domains \cite{kidger2020universal},  and alternative architectures and topologies \cite{kidger2020universal, lin2018resnet}.  We specifically highlight the quantitive version of the UAT \cite{barron1993universal}, which discusses the capability of FFNNs to be universal approximators overcoming the \emph{curse of dimensionality} compared to classical piecewise-linear function spaces.  For more extensive discussions about this, we refer to the works of \cite{poggio2017and, bach2017breaking, bauer2019deep, wojtowytsch2020can} and the references therein.

A proliferation of research initiatives and investigations have demonstrated through evidence and empirical tests that suitable combinations of different kinds of architectures, such as \emph{convolutional}, \emph{recurrent}, and \emph{residual} NNs (CNNs, RNNs, and ResNet, respectively), often surpass the approximation and computation capabilities of FFNNs in different scenarios and applications \cite{collobert2008unified, deng2009imagenet, li2015constructing, srivastava2015highway, goodfellow2016deep, he2016deep, vaswani2017attention,  abiodun2018state, dupond2019thorough, valueva2020application}.

In this dissertation, we will limit ourselves to fully-connected FFNN architectures. Therefore, the terms NNs and FFNNs can be interchangeably read from now on unless otherwise specified. Moreover, we emphasize that the study of the approximation capacities of the considered architectures is beyond the scope of this dissertation.  In other words, we will (naively) assume that our selected architectures during experimentation have a ``desired degree of approximation capacity''.

\section{Literature Review}\label[section]{section1.4}

In recent years,  a great deal of work has emerged and focused on employing NNs to solve PDEs. 

\subsection{The first works solving PDEs using NNs}\label[section]{section1.4.1}

The very first documented work addressing the PDE-solving task following the terms in \Cref{section1.3} is dated to 1994\footnote{Received by the publisher in August 1992,  revised in August 1993,  and first published in March 1994.} \cite{dissanayake1994neural}.  Although some earlier works already used NNs in various PDE-driven scenarios (see, e.g., \cite{lee1990neural, uchiyama1993solving, mcgee1994applications}),  proposal \cite{dissanayake1994neural} was the earliest to exploit the representation of a PDE solution by means of a NN that is supported on the universal approximation property.

Such work proposes a FFNN architecture with sigmoid activation functions to approximate the solution of a PDE with boundary conditions---also called \emph{\acf{BVP}}---described in the form of
\begin{equation}\label{first_PDEs_NN}
\begin{cases}
Au =f, \qquad &\text{in }\Omega,\\
Du =g, \qquad &\text{in }\partial\Omega,
\end{cases}
\end{equation} where $\Omega$ is a bounded open domain with boundary $\partial\Omega$,  $f$ is a given source function,  and $A$ and $D$ are some given (differential/boundary) operators.  The authors literally claimed\footnote{Original symbols and references have been adapted for suitability with our presentation.  [...] means that a portion of the original text has been omitted.}: ``A so-called strong solution of  \eqref{first_PDEs_NN} is an element of some underlying space $\mathbb{U}$ which satisfies both the PDE and the boundary conditions.  [...] However, if there is a `universal' approximator in this underlying space, that is, a function that has enough parameters and can be arbitrarily close to any element of the space by a judicious choice of the parameters, then this approximator can be a good choice for the approximate solution sought for. The closeness of two functions is described by the concept of denseness in $\mathbb{U}$. For a more mathematical definition, see Hornik et al. \cite{hornik1989multilayer}''.

Then,  the PDE-solving task is stated as a minimization of the \emph{objective function} given by
\begin{align}\label{first_PDEs_NN_loss}
\mathcal{F}(u_{\text{NN}}) &= \int_\Omega \Vert Au_{\text{NN}} - f \Vert^2 + \int_{\partial\Omega} \Vert Du_{\text{NN}} - g \Vert^2,
\end{align} where $u_\text{NN}$ is the considered FFNN model.  Originally,  $\Vert\cdot\Vert$ was unspecified but probably thought of as the usual discrete $2$-norm so that the minimization of $\mathcal{F}$ interprets within an $L^2$-minimization scheme (as suggested in \cite{cuomo2022scientific}).  After such presentation, the authors said\footnote{As before, original symbols and references have been adapted for suitability.  The portion of text between square brackets,  [\texttt{text}], indicates that we inserted or replaced the original text for convenience in the understanding.}: ``It is now possible to compute $Au$ and $D u$, in closed form, in terms of $x$.  We propose to use point collocation and enforce \eqref{first_PDEs_NN_loss} at selected locations in $\Omega$ and $\partial\Omega$. The resulting objective function may then be minimized to find values for [the learnable parameters of the NN]---hence the approximate solution to \eqref{first_PDEs_NN}''. They specified that the minimization of the objective function was carried out according to ``a quasi-Newton method''  and employing ``finite-difference gradients'' (see, e.g.,  \cite{dennis1977quasi} and \cite{milne2000calculus} for related references).

Since then, several works have emerged within the \ac{PDE}-solving framework using \acp{NN}.  Below, we review some of the most relevant ones from the late 1990s to the early 2010s.

In 1994,  \cite{meade1994numerical, meade1994solution, fernandez1994application} proposed the utilization of a FFNN with a single hidden layer activated by a ``hard limit'' function (so named by the authors) to solve linear and non-linear ODEs.  This choice of activation function allows interpreting the set of admissible NNs as a family of piecewise-linear functions, and thus approximate the ODE solution in those terms.  In 1998, \cite{lagaris1998artificial} introduced a scheme to strongly impose the boundary conditions as a composition of a (boundary-free) FFNN with a filter,  as opposed to \cite{dissanayake1994neural}, where the boundary conditions were imposed via the objective function.  In 1999,  \cite{nguyen1999approximation} reported a NN implementation for the numerical approximation of functions of several variables and their partial derivatives accompanied by several results.  Among those,  we encounter applications to PDEs as part of a boundary element method for analyzing viscoelastic flows.  In 2000,  \cite{he2000multilayer} proposed an extended backpropagation algorithm to solve a class of first-order PDEs with particular concern on nonlinear control systems, and \cite{lagaris2000neural} presented a NN-based model consisting of two different FFNN-type architectures to address PDEs with boundary conditions on irregular geometries.

Among the works from 2001 to the mid-2010s, we find multiple works which can be understood as continuations of those mentioned above.  Specifically,  \cite{lagaris1998artificial} and \cite{he2000multilayer} established significant precedents for the architecture designs and training setups, respectively,  of the following works: \cite{aarts2001neural} introduced a method combining NNs and evolutionary algorithms for solving PDEs and their boundary conditions; \cite{smaoui2004modelling} employed NNs to analyze the dynamics of Kuramato-Sivashinsky and Navier-Stokes nonlinear equations; \cite{malek2006numerical}
presented a hybrid method for solving high-order ODEs; \cite{shirvany2008numerical} addressed Schr\"odinger's equation by a \ac{FFNN} with an improved energy-function scheme via unsupervised training; \cite{beidokhti2009solving} proposed a combination of minimization techniques and collocation methods via NNs to establish approximate closed-analytic-form solutions to PDEs; \cite{tsoulos2009solving} used FFNNs in combination with grammatical evolution and local optimization to solve (systems of) ODEs and PDEs (see \cite{tsoulos2008neural} for a precedent in grammatical-evolution-based neural-network training); and \cite{mcfall2009artificial} presented a NN-based method for solving PDEs with irregular domains with mixed boundary conditions.  Additionally,  we encounter two major classes of extensions of NN-based PDE-solving frameworks: (Multiquadratic) Radial Basis Function (RBF) networks\footnote{See \cite{park1993approximation} for a UAT precedent on RBF networks.} \cite{mai2001numerical,  mai2002mesh, jianyu2002numerical, mai2003approximation,jianyu2003numerical, kansa2004volumetric, mai2005solving, golbabai2007radial, aminataei2008numerical, golbabai2009solving,chen2011numerical,kumar2011multilayer}; and NNs inspired in the FEM \cite{beltzer2003neural,  deng2003pillar, ziemianski2003hybrid,  ramuhalli2005finite, manevitz2005neural,  arndt2005approximating,  jilani2009adaptive, koroglu2010comparison}---sometimes also referred to as \emph{Finite Element Neural Networks}.  We refer to \cite{yadav2015introduction} for a detailed overview of many of the abovementioned works.

\subsection{Solving PDEs using NNs since TensorFlow}\label[section]{section1.4.2}

Prior to the first open-source release of \emph{\acf{TF}} \cite{tensorflow2015-whitepaper, abadi2016tensorflow} in 2015, the absence of efficient \emph{\acf{AD}} posed challenges in \ac{NN}-based research. AD is nowadays crucial to compute gradients of NNs during training \cite{griewank2008evaluating, baydin2018automatic, margossian2019review}. Without it, researchers and practitioners had to manually design the computation scheme of gradients, which typically was time-consuming and prone to errors, significantly as models grew in complexity.  Consequently, innovation and experimentation were hampered. The inability to quickly try new ideas and different types of models slowed down progress, making it difficult for newcomers to enter the field. 

With the introduction of TF \cite{tensorflow2015-whitepaper, abadi2016tensorflow} in 2015 and of PyTorch \cite{paszke2017automatic, paszke2019pytorch} in 2017, AD became more accessible for researchers, liberating them from manual computation designs and fueling faster development and innovation.  Thanks to the emergence of these NN-oriented, AD-based, and user-friendly platforms, there arose a specific interest in the PDE-solving community to investigate NN-based venues in their research from the late 2010s on.

Specifically, \emph{\acfp{PINN}} could be considered as one of the most significant precedents that have boosted the PDE-solving community to employ NNs in their works, especially during the last five years.  The first version of PINNs was released as a two-part treatise in 2017 \cite{raissi2017physics1, raissi2017physics2}.  There,  the authors proposed data-driven alternatives for either solving (in \cite{raissi2017physics1}) or discovering (in \cite{raissi2017physics2}) PDEs.  To distinguish the ``solving'' from the ``discovering'' frameworks, we think of a PDE with a differential operator governed by certain coefficients (parameters).  On the one hand,  if we know the values of all the coefficients, then we have complete access to the differential law involved in the PDE. Hence, we can reduce everything to the terms of the aforementioned ``first work'' \cite{dissanayake1994neural}\footnote{Interestingly, neither \cite{raissi2017physics1, raissi2017physics2} nor similar authorship reviews such as \cite{raissi2019physics, karniadakis2021physics} do cite \cite{dissanayake1994neural}.  We need to consult more general (and authorship-independent) reviews to find this relation. For example, \cite{cuomo2022scientific} states: ``The concept of incorporating prior knowledge into a machine learning algorithm is not entirely novel. In fact, Dissanayake and Phan-Thien \cite{dissanayake1994neural} can be considered one of the first \acp{PINN}''.} for approximating the PDE solution (recall the first part of \Cref{section1.4.1}). On the other hand,  if some of the coefficients are unknown, but we have access to some observations of the PDE solution at different spatio-temporal points,  then we can establish a supervised learning scheme where a NN (representing the PDE solution) together with some additional trainable parameters (representing the unknown PDE coefficients) are trained constrained to satisfy both the PDE and the accessible labeled data.

Continuing with the PINN terminology,  we encounter works addressing variational formulations,  named \emph{($hp$-)variational PINNs (($hp$-)VPINNs)} \cite{kharazmi2019variational, kharazmi2021hp, rojas2023robust}; conservation laws, named \emph{conservative PINNs (cPINNs)} \cite{ameya2020conservative}; fractional-order PDEs,  named \emph{fractional PINNs (fPINNs)} \cite{pang2019fpinns}; and Bayesian proposals for noisy data,  named \emph{Bayesian PINNs (B-PINNs)} \cite{yang2021b}.  Among PINN applications to real-world phenomena, we highlight: \cite{mao2020physics} for high-speed flows, \cite{misyris2020physics, huang2022applications} for power systems,  \cite{shukla2020physics} for ultrasound nondestructive quantification, \cite{cai2021physics1} for heat transfer,  and \cite{cai2021physics2,jin2021nsfnets} for fluid mechanics. We refer to \cite{cuomo2022scientific,lawal2022physics,blechschmidt2021three} for more extensive reviews in PINNs.

Alternatively,  related works adopting the \emph{deep} first name have also been proposed following the NN-based PDE-solving aim.  To mention a few: \cite{e2017deep, e2018deep} introduces the \emph{Deep Ritz Method (DRM)}, which minimizes the Ritz energy functional involved in variational formulations of PDEs---refer to \cite{lu2021priori, duan2022convergence, ming2021deep, uriarte2023deep} for related further analyses and extensions); \cite{sirignano2018dgm} presents the \emph{Deep Galerkin Method (DGM)}---refer to \cite{chen2020comparison, li2022deep,  al2022extensions,  shang2022deep} for further extensions, applications,  and comparisons; \cite{cai2020deep} proposes a \emph{Deep Least-Squares (DLS) method}---refer to \cite{cai2021least, liu2022adaptive, nareklishvili2023deep} for further extensions; and \cite{taylor2023deep1} presents a \emph{Deep Fourier Residual (DFR) method} with a recent extension to time-harmonic Maxwell's equations \cite{taylor2023deep2}. 

We also encounter several recent works regarding the resolution of parametric PDEs.  Among the most relevant: \cite{li2020multipole} proposes a multi-level graph NN;  \cite{li2021fourier} utilizes a parameterization in the Fourier space to achieve an expressive and efficient architecture suitable for parametric problems; \cite{bhattacharya2021model} combines NNs with model reduction; \cite{khoo2021solving} uses NNs to parameterize the physical quantity of interest as a function of input coefficients; \cite{kutyniok2022theoretical} analyzes theoretically the approximability of parametric maps by NNs for parametric PDEs; and \cite{penwarden2023metalearning} extends PINNs to a parametric version throughout a named ``metalearning'' approach.

\section{Major Contributions}\label[section]{section1.5}

This dissertation provides three major contributions to the resolution of PDEs using NNs: (i) a finite-element-based deep learning solver for parametric PDEs, (ii) a general residual minimization framework based on a double Ritz minimization scheme, and (iii) a numerical memory-based integration strategy for efficiently reducing the integration error during training.

\subsection{The Deep Finite Element Method}\label[section]{section1.5.1}

Solving inverse problems is of great value to our society \cite{yaman2013recent, LaTorre2015inverse} due to their broad applicability to multiple engineering disciplines such as imaging \cite{bertero2020introduction}, electromagnetics \cite{alvarez-aramberri2013inversion}, non-destructive evaluations \cite{omella2021sensitivity},  and geophysics \cite{rojas2016quadrature}, to mention a few.  Different methods exist to solve them, including gradient-based and statistical-based methods \cite{tarantola2005inverse, aster2019parameter}. These traditional methods typically evaluate the inverse solution pointwise (i.e., for a given set of measurements) but rarely provide a global representation of an inverse operator. To overcome this and approximate a full inverse operator, it is possible to use NNs (see, e.g., \cite{iturraran-viveros2014artificial, ye2018deep, ongie2020deep, shahriari2020deep, alyaev2021modeling, shahriari2021error}) due to their universal approximation property. 

Training NNs for solving inverse problems is typically performed in \emph{\acfp{GPU}}, especially when dealing with large databases or complex model architectures. In these situations, it is convenient to have an efficient parametric PDE solver implemented within the GPU in the form of a NN to rapidly iterate over the different parameter candidates during the inverse problem training. Under these conditions, there are two main approaches for training the inverse model: (i) pre-training beforehand the parametric forward model to later use it as a real-time solver/evaluator, or (ii) training both the forward and inverse models simultaneously in an (auto) encoder-decoder scheme  \cite{goodfellow2016deep}. 

As the first main contribution of this dissertation, \textbf{we propose a NN-based model that acts as a parametric forward solver inspired by the FEM}. Specifically, our NN architecture is mesh-dependent, mimicking the effect of the dynamics when applying mesh refinements \cite{babuvvska1978error, ainsworth1997posteriori, strouboulis2000design}.

Our presentation and implementation of the proposed method, called the \emph{\acf{DeepFEM}}, restricts to \emph{\acf{1D}} problems using piecewise-linear approximations and uniform refinements. The extension to higher-dimensional problems, higher-order polynomial approximations, and/or adaptive meshes is presumably straightforward but goes beyond the scope of this work. We test the DeepFEM on positive-definite and indefinite problems with constant and piecewise-constant parameters. In general, experiments show good approximations; however, on some occasions, the optimizer and convexity limitations prevent us from obtaining high-accuracy solutions.  

\subsection{The Deep Double Ritz Method}\label[section]{section1.5.2}

From an abstract point of view, variational formulations \cite{mikhlin1964variational} (on which, for example, the FEM is based) consist of translating BVPs into functional equations of the form
\begin{equation}\label{varform}
Bu=l\in\mathbb{V}',
\end{equation} where $\mathbb{U}$ and $\mathbb{V}$ are the so-called \emph{trial} and \emph{test} (normed) spaces, respectively; $\mathbb{V}'$ is the dual of $\mathbb{V}$; and $B:\mathbb{U}\longrightarrow\mathbb{V}'$ and $l\in\mathbb{V}'$ are the ``variational analogs'' of the original differential operator and source term, respectively. We refer to \cite[p. 150, eqs. (2.1)--(2.3)]{demkowicz2014overview} for details. Specifically, we highlight the claim: ``the nature of variational problems lies in the fact that the corresponding operator takes values in a dual space''.

Under this abstract formulation, \emph{residual minimization} naturally arises as a reformulation of \eqref{varform} as follows:
\begin{equation}\label{minres}
	\min_{u\in\mathbb{U}} \Vert Bu-l\Vert_{\mathbb{V}'},
\end{equation} where the use of the dual norm $\Vert \cdot\Vert_{\mathbb{V}'}$ is a must because $Bu-l$ takes values in the dual space $\mathbb{V}'$ \cite[p. 151, eq. (2.5)]{demkowicz2014overview}. 

Following the supremum definition of dual norms, we end up with a saddle-point (min-max\footnote{In Hilbert spaces, which is where we will limit ourselves later on, the supremum is attained and therefore, it can indistinctly be replaced by a maximum.}) problem over the trial and test spaces as follows:
\begin{equation}\label{min_max_approach}
\min_{u\in\mathbb{U}} \max_{v\in\mathbb{V}\setminus\{0\}} \frac{\langle Bu-l, v\rangle_{\mathbb{V}'\times\mathbb{V}}}{\Vert v\Vert_\mathbb{V}}.
\end{equation} Here, $\langle \cdot, \cdot\rangle_{\mathbb{V}'\times\mathbb{V}}$ stands for the duality pairing and $\Vert\cdot\Vert_{\mathbb{V}}$ for the norm in $\mathbb{V}$.

In the context of NNs, \cite{zang2020weak, bao2020numerical} proposed performing residual minimization employing two NNs to represent the trial and test functions during the min-max optimization strategy \eqref{min_max_approach}. They called the resulting method \emph{\acf{WANs}}.  Although \cite{zang2020weak, bao2020numerical} considered suboptimal norms during experimentation---they employed the $L^2$-norm instead of the proper $H^1$-norm, the major underlying drawback of such min-max approach relies on the fact that seeking the involved unitary test maximizer for each trial minimization iteration is numerically unstable (we will show this in \Cref{chapter4}).
 
To overcome the above instability issue, \textbf{we propose a general residual minimization method using NNs}, called the \emph{\acf{D2RM}}, that follows a similar spirit as in WANs and utilizes two \acp{NN} to nestedly iterate over trial and test functions, but relying on a stable Ritz minimization formulation \cite{ritz1909neue}. Indeed, \ac{D2RM} can be interpreted either as a stabilized version of WANs or as a generalization of the DRM \cite{e2018deep} to problems that are not necessarily symmetric and positive-definite (recall \Cref{section1.4.2}).

Following the interrelation between WANs, the DRM, and the \ac{D2RM}, we test and compare these three methods on several diffusion and convection problems to show the generalizability and stability of the \ac{D2RM} (compared to the other two) up to the approximation properties of the considered NNs and the training capacity of the optimizers.

\subsection{Memory-based Monte Carlo integration}\label[section]{section1.5.3}

Solving a BVP with NNs is typically described as a minimization of an objective function in the form of a definite integral that involves a NN in the integrand, e.g., recall \eqref{first_PDEs_NN_loss} or \eqref{min_max_approach}. Thus, the quality of the numerical approximation or quadrature is critical. If we were to consider a fixed quadrature grid for numerical integration during the entire training, it is highly probable that the NN acquires a greedy behavior that produces undesired NN predictions \cite{rivera2022quadrature} as a result of the well-known overfitting phenomenon \cite{hawkins2004problem, claeskens2008model,goodfellow2016deep}.

A usual choice for avoiding overfitting consists of considering \emph{\acf{MC} integration} \cite{newman1999monte, leobacher2014introduction} due to its stochastic nature. Moreover, its mesh-free structure allows a straightforward scalability to complex geometries or high-dimensional domains.   Unfortunately, \acs{MC} integration commits an error of order $\mathcal{O}(1/\sqrt{N})$, where $N$ is the number of integration points \cite{newman1999monte}. This low convergence rate of integration accuracy is typically overcome by sampling tens or hundreds of thousands of integration points (at each training iteration), drastically decreasing the training speed.

As the third main contribution of this dissertation, \textbf{we propose a memory-based MC integration scheme when using NNs for solving PDEs}. Here, the memory is aimed at taking advantage of integral estimations computed in previous iterations to contribute to the integral estimate at the current iteration. This scheme gradually decreases the integration error without incurring the computational overhead of dealing with disproportionate samples during training. Additionally, by conveniently translating this approach to the \emph{\acf{SGD}} optimizer, we find a favorable reinterpretation of the momentum method \cite{polyak1964some}. 

\section{Outline}\label[section]{section1.6}

The remainder of this dissertation is structured as follows. \Cref{chapter2} formally introduces and analyzes \acp{FFNN}. {\color{blue!50!black}Chapters} \ref{chapter3}, \ref{chapter4} and \ref{chapter5} are devoted to the detailed development of the main contributions of the dissertation briefly motivated in {\color{blue!50!black}Sections} \ref{section1.5.1}, \ref{section1.5.2} and \ref{section1.5.3}, respectively. \Cref{chapter6} exposes the main conclusions of the dissertation and sketches the challenges to be addressed in future work.  Finally,  \Cref{chapter7} describes the main achievements of this dissertation, including the published works,  scientific contributions to conferences and workshops, and other relevant scientific activities developed during the PhD program.



\chapter[Feed-Forward Neural Networks]{Feed-Forward Neural Networks} \label{chapter2}

\tikzset{arrowed/.style={decorate,
decoration={show path construction, 
moveto code={},
lineto code={
\draw[#1] (\tikzinputsegmentfirst) --  (\tikzinputsegmentlast);
},
curveto code={},
closepath code={},
}},arrowed/.default={-stealth}}
\pgfplotsset{gradient function/.initial=f,
dx/.initial=0.01,dy/.initial=0.01}
\pgfmathdeclarefunction{xgrad}{2}{%
\begingroup%
\pgfkeys{/pgf/fpu,/pgf/fpu/output format=fixed}%
\edef\myfun{\pgfkeysvalueof{/pgfplots/gradient function}}%
\pgfmathparse{(\myfun(#1+\pgfkeysvalueof{/pgfplots/dx},#2)%
-\myfun(#1,#2))/\pgfkeysvalueof{/pgfplots/dx}}%
\pgfmathsmuggle\pgfmathresult\endgroup%
}%
\pgfmathdeclarefunction{ygrad}{2}{%
\begingroup%
\pgfkeys{/pgf/fpu,/pgf/fpu/output format=fixed}%
\edef\myfun{\pgfkeysvalueof{/pgfplots/gradient function}}%
\pgfmathparse{(\myfun(#1,#2+\pgfkeysvalueof{/pgfplots/dy})%
-\myfun(#1,#2))/\pgfkeysvalueof{/pgfplots/dy}}%
\pgfmathsmuggle\pgfmathresult\endgroup%
}%

Herein, we formally introduce the \emph{\acfp{FFNN}} that we will exhaustively exploit along \Cref{chapter3}, \ref{chapter4}, and \ref{chapter5}. 

This chapter is organized as follows.  \Cref{section2.1} presents the architecture framework; \Cref{section2.2} discusses continuum, parameterized, and discretized setups; \Cref{section2.3} describes gradient-based training; and \Cref{section2.4} illustrates the functioning of \acp{FFNN} with an easy-to-analyze case of study.

\section{Architecture}\label{section2.1}

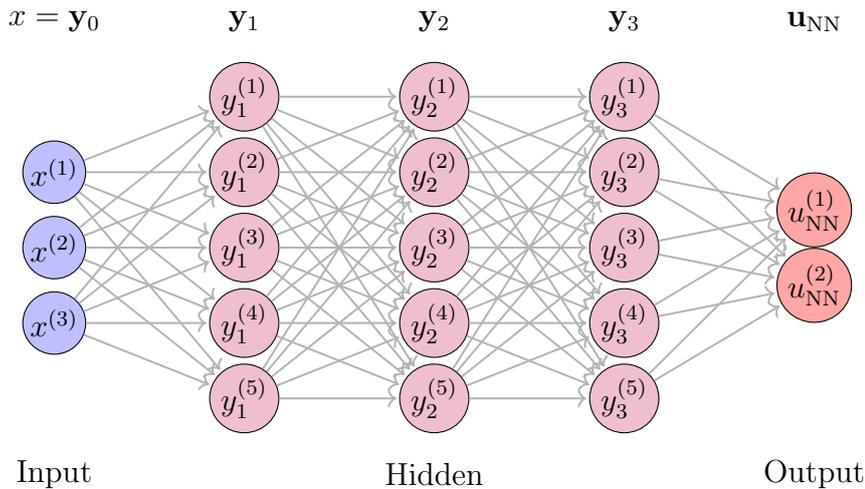
\begin{figure}
\centering
\begin{tikzpicture}[shorten >=1pt,->,
    draw=gray!60,
    node distance=\layersep,
    every pin edge/.style={<-,shorten <=1pt},
    neuron/.style={circle,draw=black,fill=black, inner sep=1pt},
    element/.style={rectangle,fill=black,minimum size=16pt,inner sep=0pt},
    input neuron/.style={neuron, fill=blue!25},
    output neuron/.style={neuron, fill=red!35},
    hidden neuron/.style={neuron, fill=purple!25},
    annot/.style={text width=4em, text centered},
    pil/.style={
           ->,
           thick,
           line width=2pt,
           shorten <=2pt,
           shorten >=2pt,}
]
    \def\layersep{1.5cm}

    \node[] at (0,1) {$x=\mathbf{y}_0$};
    \node[] at (2.5,1) {$\mathbf{y}_1$};
    \node[] at (5,1) {$\mathbf{y}_2$};
    \node[] at (7.5,1) {$\mathbf{y}_3$};
    \node[] at (10,1) {$\mathbf{u}_\text{NN}$};
    \node[] at (0,-5) {Input};
    \node[] at (5,-5) {Hidden};
    \node[] at (10,-5) {Output};
    \node[input neuron] (I-1) at (0,-1) {$x^{(1)}$};
    \node[input neuron] (I-2) at (0,-2) {$x^{(2)}$};
    \node[input neuron] (I-3) at (0,-3) {$x^{(3)}$};
    \node[hidden neuron] (H1-0) at (2.5,0) {$y^{(1)}_1$};
    \node[hidden neuron] (H1-1) at (2.5,-1) {$y^{(2)}_1$};
    \node[hidden neuron] (H1-2) at (2.5,-2) {$y^{(3)}_1$};
    \node[hidden neuron] (H1-3) at (2.5,-3) {$y^{(4)}_1$};
    \node[hidden neuron] (H1-4) at (2.5,-4){$y^{(5)}_1$};
    \node[hidden neuron] (H2-0) at (5,0) {$y^{(1)}_2$};
    \node[hidden neuron] (H2-1) at (5,-1) {$y^{(2)}_2$};
    \node[hidden neuron] (H2-2) at (5,-2) {$y^{(3)}_2$};
    \node[hidden neuron] (H2-3) at (5,-3) {$y^{(4)}_2$};
    \node[hidden neuron] (H2-4) at (5,-4) {$y^{(5)}_2$};
    \node[hidden neuron] (H3-0) at (7.5,0) {$y^{(1)}_3$};
    \node[hidden neuron] (H3-1) at (7.5,-1) {$y^{(2)}_3$};
    \node[hidden neuron] (H3-2) at (7.5,-2) {$y^{(3)}_3$};
    \node[hidden neuron] (H3-3) at (7.5,-3) {$y^{(4)}_3$};
    \node[hidden neuron] (H3-4) at (7.5,-4) {$y^{(5)}_3$};
    \node[output neuron] (O-1) at (10,-1.5) {$u_{\text{NN}}^{(1)}$};
    \node[output neuron] (O-2) at (10,-2.5) {$u_{\text{NN}}^{(2)}$};
    \path[line width=0.75pt] (I-1) edge (H1-0);
    \path[line width=0.75pt] (I-1) edge (H1-1);
    \path[line width=0.75pt] (I-1) edge (H1-2);
    \path[line width=0.75pt] (I-1) edge (H1-3);
    \path[line width=0.75pt] (I-1) edge (H1-4);
    \path[line width=0.75pt] (I-2) edge (H1-0);
    \path[line width=0.75pt] (I-2) edge (H1-1);
    \path[line width=0.75pt] (I-2) edge (H1-2);
    \path[line width=0.75pt] (I-2) edge (H1-3);
    \path[line width=0.75pt] (I-2) edge (H1-4);
    \path[line width=0.75pt] (I-3) edge (H1-0);
    \path[line width=0.75pt] (I-3) edge (H1-1);
    \path[line width=0.75pt] (I-3) edge (H1-2);
    \path[line width=0.75pt] (I-3) edge (H1-3);
    \path[line width=0.75pt] (I-3) edge (H1-4);
    \path[line width=0.75pt] (H1-0) edge (H2-0);
    \path[line width=0.75pt] (H1-1) edge (H2-0);
    \path[line width=0.75pt] (H1-2) edge (H2-0);
    \path[line width=0.75pt] (H1-3) edge (H2-0);
    \path[line width=0.75pt] (H1-4) edge (H2-0);
    \path[line width=0.75pt] (H1-0) edge (H2-1);
    \path[line width=0.75pt] (H1-1) edge (H2-1);
    \path[line width=0.75pt] (H1-2) edge (H2-1);
    \path[line width=0.75pt] (H1-3) edge (H2-1);
    \path[line width=0.75pt] (H1-4) edge (H2-1);
    \path[line width=0.75pt] (H1-0) edge (H2-2);
    \path[line width=0.75pt] (H1-1) edge (H2-2);
    \path[line width=0.75pt] (H1-2) edge (H2-2);
    \path[line width=0.75pt] (H1-3) edge (H2-2);
    \path[line width=0.75pt] (H1-4) edge (H2-2);
    \path[line width=0.75pt] (H1-0) edge (H2-3);
    \path[line width=0.75pt] (H1-1) edge (H2-3);
    \path[line width=0.75pt] (H1-2) edge (H2-3);
    \path[line width=0.75pt] (H1-3) edge (H2-3);
    \path[line width=0.75pt] (H1-4) edge (H2-3);
    \path[line width=0.75pt] (H1-0) edge (H2-4);
    \path[line width=0.75pt] (H1-1) edge (H2-4);
    \path[line width=0.75pt] (H1-2) edge (H2-4);
    \path[line width=0.75pt] (H1-3) edge (H2-4);
    \path[line width=0.75pt] (H1-4) edge (H2-4);
    \path[line width=0.75pt] (H2-0) edge (H3-0);
    \path[line width=0.75pt] (H2-1) edge (H3-0);
    \path[line width=0.75pt] (H2-2) edge (H3-0);
    \path[line width=0.75pt] (H2-3) edge (H3-0);
    \path[line width=0.75pt] (H2-4) edge (H3-0);
    \path[line width=0.75pt] (H2-0) edge (H3-1);
    \path[line width=0.75pt] (H2-1) edge (H3-1);
    \path[line width=0.75pt] (H2-2) edge (H3-1);
    \path[line width=0.75pt] (H2-3) edge (H3-1);
    \path[line width=0.75pt] (H2-4) edge (H3-1);
    \path[line width=0.75pt] (H2-0) edge (H3-2);
    \path[line width=0.75pt] (H2-1) edge (H3-2);
    \path[line width=0.75pt] (H2-2) edge (H3-2);
    \path[line width=0.75pt] (H2-3) edge (H3-2);
    \path[line width=0.75pt] (H2-4) edge (H3-2);
    \path[line width=0.75pt] (H2-0) edge (H3-3);
    \path[line width=0.75pt] (H2-1) edge (H3-3);
    \path[line width=0.75pt] (H2-2) edge (H3-3);
    \path[line width=0.75pt] (H2-3) edge (H3-3);
    \path[line width=0.75pt] (H2-4) edge (H3-3);
    \path[line width=0.75pt] (H2-0) edge (H3-4);
    \path[line width=0.75pt] (H2-1) edge (H3-4);
    \path[line width=0.75pt] (H2-2) edge (H3-4);
    \path[line width=0.75pt] (H2-3) edge (H3-4);
    \path[line width=0.75pt] (H2-4) edge (H3-4);
    \path[line width=0.75pt] (H3-0) edge (O-1);
    \path[line width=0.75pt] (H3-1) edge (O-1);
    \path[line width=0.75pt] (H3-2) edge (O-1);
    \path[line width=0.75pt] (H3-3) edge (O-1);
    \path[line width=0.75pt] (H3-4) edge (O-1);
    \path[line width=0.75pt] (H3-0) edge (O-2);
    \path[line width=0.75pt] (H3-1) edge (O-2);
    \path[line width=0.75pt] (H3-2) edge (O-2);
    \path[line width=0.75pt] (H3-3) edge (O-2);
    \path[line width=0.75pt] (H3-4) edge (O-2);
\end{tikzpicture}
\caption{Graph of a fully-connected FFNN with a three-dimensional input,  two-dimensional output,  depth three,  and hidden layers of width five.}
\label{FFNN_illustration}
\end{figure}

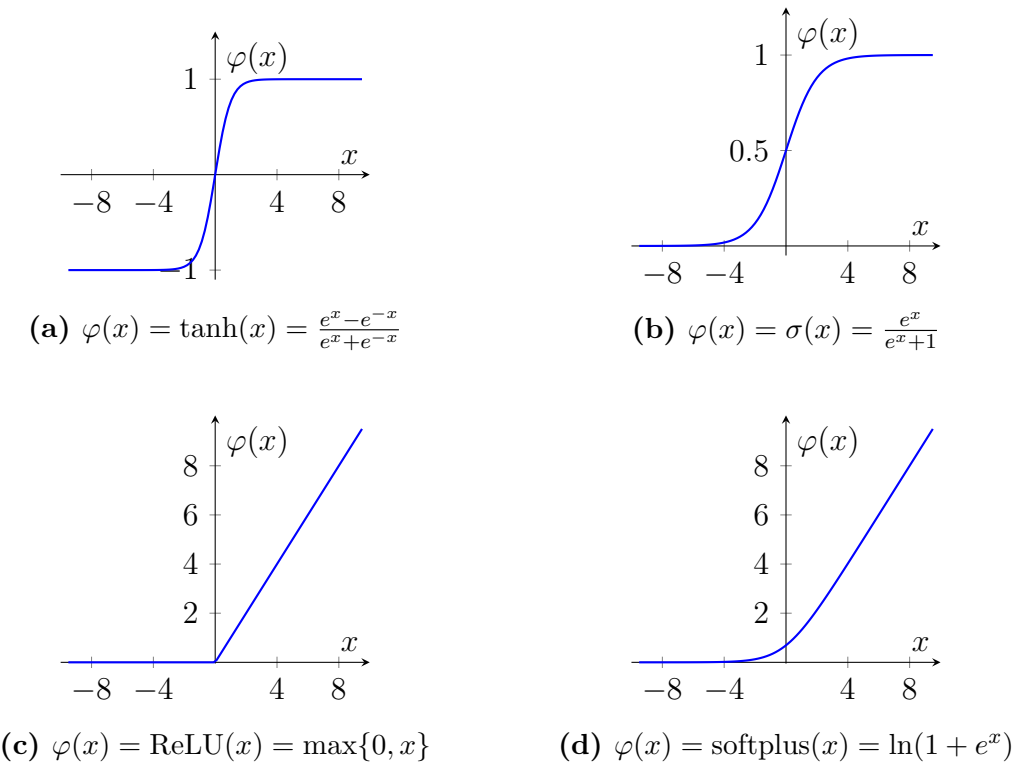
\begin{figure}
  \centering
    \begin{subfigure}[b]{0.48\textwidth}
    \centering
\begin{tikzpicture}
\begin{axis}[
    axis lines=middle,
    xmax=10,
    xmin=-10,
    xtick={-8,-4,4,8},
    ymin=-1.1,
    ymax=1.5,
    xlabel={$x$},
    ylabel={$\varphi(x)$},
    height=.7\textwidth]
\addplot [domain=-9.5:9.5, samples=100,
     thick, blue] {(exp(x) - exp(-x))/(exp(x) + exp(-x))};
\end{axis}
\end{tikzpicture}
    \caption{$\varphi(x)=\tanh(x)=\frac{e^x-e^{-x}}{e^x+e^{-x}}$}
  \end{subfigure}
    \hskip 1em
  \begin{subfigure}[b]{0.48\textwidth}
    \centering
\begin{tikzpicture}
\begin{axis}[
    axis lines=middle,
    xmax=10,
    xmin=-10,
    xtick={-8,-4,4,8},
    ymin=-0.05,
    ymax=1.25,
    xlabel={$x$},
    ylabel={$\varphi(x)$},
    height=.7\textwidth
]
\addplot [domain=-9.5:9.5, samples=100,
          thick, blue] {1/(1+exp(-x)};
\end{axis}
\end{tikzpicture}
    \caption{$\varphi(x)=\sigma(x)=\frac{e^x}{e^x+1}$}
  \end{subfigure}
  \vskip 2em
  \begin{subfigure}[b]{0.48\textwidth}
    \centering
\begin{tikzpicture}
\begin{axis}[
    axis lines=middle,
    xmax=10,
    xmin=-10,
    xtick={-8,-4,4,8},
    ymin=-0.05,
    ymax=10.05,
    ytick={2,4,6,8},
    xlabel={$x$},
    ylabel={$\varphi(x)$},
    height=.7\textwidth]
\addplot [domain=-9.5:9.5, samples=100, thick, blue] {max(0, x)};
\end{axis}
\end{tikzpicture}
    \caption{$\varphi(x)=\text{ReLU}(x)=\max\{0,x\}$}
  \end{subfigure}
  \hskip 1em
  \begin{subfigure}[b]{0.48\textwidth}
    \centering
\begin{tikzpicture}
\begin{axis}[
    axis lines=middle,
    xmax=10,
    xmin=-10,
    xtick={-8,-4,4,8},
    ymin=-0.05,
    ymax=10.05,
    ytick={2,4,6,8},
    xlabel={$x$},
    ylabel={$\varphi(x)$},
    height=.7\textwidth]
\addplot [domain=-9.5:9.5, samples=100, thick, blue] {ln(1+exp(x))};
\end{axis}
\end{tikzpicture}
    \caption{$\varphi(x)=\text{softplus}(x)=\ln(1+e^x)$}
  \end{subfigure}
  \caption{Some typical activation functions appearance: (a) hyperbolic tangent, (b) sigmoid, (c) rectified linear unit, and (d) softplus.}
  \label{Chapter1_activation_functions}
\end{figure}

In our case, a \ac{FFNN} is a (vector-valued) function $\mathbf{u}_{\text{NN}}: X\longrightarrow Y$ defined as follows:
\begin{subequations}\label{FFNN_architecture} 
\begin{align}
    \mathbf{y}_0(x) &:= x\in X\subset\mathbb{R}^{n_0},\label{input_layer}\\
    \mathbf{y}_j(x) &:= \varphi(\mathbf{W}_j \mathbf{y}_{j-1} + \mathbf{b}_j)\in \mathbb{R}^{n_j}, \qquad 1\leq j\leq K, \label{hidden_layer}\\
    \mathbf{u}_{\text{NN}}(x) &:=\mathbf{W} \mathbf{y}_K \in Y\subset\mathbb{R}^{n_{K+1}}\label{output_layer},
\end{align}
\end{subequations} where $\mathbf{W}_j \mathbf{y}_{j-1} + \mathbf{b}_j$ is an affine transformation determined by \emph{weights} $\mathbf{W}_j\in \mathbb{R}^{n_j\times n_{j-1}}$ and \emph{bias}  $\mathbf{b}_j\in \mathbb{R}^{n_j}$ for $1\leq j\leq K$,  $\mathbf{W} \mathbf{y}_K$ is a linear transformation via weights $\mathbf{W}\in \mathbb{R}^{n_{K+1}\times n_K}$, and $\varphi$ is a non-linear \emph{activation function} that acts componentwise,  i.e.,  for any $\mathbf{z}=(z^{(1)},  z^{(2)}, \ldots,z^{(k)})\in\mathbb{R}^k$ and any $k\in\mathbb{N}$,
\begin{equation}
\varphi\left(\mathbf{z}\right)=\left(\varphi\left(z^{(1)}\right),  \varphi\left(z^{(2)}\right), \ldots,\varphi\left(z^{(k)}\right)\right).
\end{equation} We call \emph{hidden layer} the mapping $\mathbf{y}_{j-1}\mapsto \mathbf{y}_j$ given in \eqref{hidden_layer}, and the dimension of its output vector, $n_j$, is known as its \emph{width}. Similarly, we have \emph{input} and \emph{output layers} in \eqref{input_layer} and \eqref{output_layer} with corresponding input and output dimensions $n_0$ and $n_{K+1}$,  respectively.  The \emph{depth} of $\textbf{u}_\text{NN}$ is its number of hidden layers, $K$,  that when reaching a significant value (e.g.,  more than three),  the \ac{NN} is typically referred to as a \emph{\acf{DNN}}.  \autoref{FFNN_illustration} shows the graph of a fully-connected \ac{FFNN} architecture and \autoref{Chapter1_activation_functions} displays different widely used activation functions.

The parameterization of $\mathbf{u}_\text{NN}$ is commonly given by all the elements involved in the affine transformations.  Hence,  to formally split the input domain versus the parameterization domain,  we redefine the \ac{FFNN} as follows:
\begin{subequations}\label{Chapter1_FFNN_definition}
\begin{equation}
\begin{array}{lccccccc}
\mathbf{u}_\text{NN}: & X & \times & \Theta & \longrightarrow & Y,\\
& x & , & \theta & \longmapsto & \mathbf{u}_\text{NN}(x; \theta),
\end{array}
\end{equation} where 
\begin{equation}
\theta := \{\mathbf{W}_1,\mathbf{b}_1,\mathbf{W}_2,\mathbf{b}_2,\ldots,\mathbf{W}_K,\mathbf{b}_K,\mathbf{W}\}
\end{equation} stands for the set of learnable parameters with domain 
\begin{equation}
\Theta:= \mathbb{R}^{n_{1} \times n_{0}}\times \mathbb{R}^{n_1} \times \mathbb{R}^{n_{2} \times n_{1}}\times \mathbb{R}^{n_2} \times \cdots \times \mathbb{R}^{n_{K} \times n_{K-1}}\times \mathbb{R}^{n_{K}} \times \mathbb{R}^{n_{K+1}\times n_{K}}.
\end{equation}
\end{subequations} The functions arising from this parameterization,  also known as \emph{realizations},  belong to a space that, in general,  is not a vector space but a \emph{manifold} \cite{petersen2021topological}.  Following the above notation, we denote it by
\begin{equation}
\mathbb{U}_\text{NN} :=  \{\mathbf{u}_\text{NN}( \ \cdot \ ; \theta): X\longrightarrow Y\}_{\theta\in\Theta},
\end{equation}  and by $\Phi_\text{NN}:\Theta\longrightarrow\mathbb{U}_\text{NN}$ the \emph{realization map} that relates each configuration of parameters with the corresponding realization for the given \ac{NN} architecture, i.e.,
\begin{equation}\label{realization_map}
\Phi_\text{NN}(\theta):=\mathbf{u}_\text{NN}(\ \cdot \ ;\theta):X\longrightarrow Y,\qquad \theta\in\Theta.
\end{equation} We will drop the boldface notation when the NN is specified as a scalar-valued function, i.e.,  we will write $u_\text{NN}$ instead of  $\mathbf{u}_\text{NN}$ whenever $n_{K+1}=1$.

The set of realizations produced by a given FFNN defined as in \eqref{Chapter1_FFNN_definition} could be interpreted as a family of finite-dimensional vector spaces. See \Cref{appendix2.a} for details.

\section{Continuum, parameterized, and discretized setups}\label{section2.2}

Although \acp{NN} are often venerated for their tremendous approximation capabilities, they also possess some noteworthy topological undesirable properties.  For example, given a \ac{NN} architecture, the set of realizations is generally non-closed and non-convex \cite{petersen2021topological, brevis2022neural}. This typically translates into possibly being unable to establish an ``optimal'' realization (that is better than any other) or that multiple ones may exist. This is a substantial departure from the usual approach where the optimal approximation candidate is typically well-defined as the projection onto a pre-established finite-dimensional vector subspace.  Aware of these ill-posed omens,  we formally present the approximation task via \acp{NN} as follows.

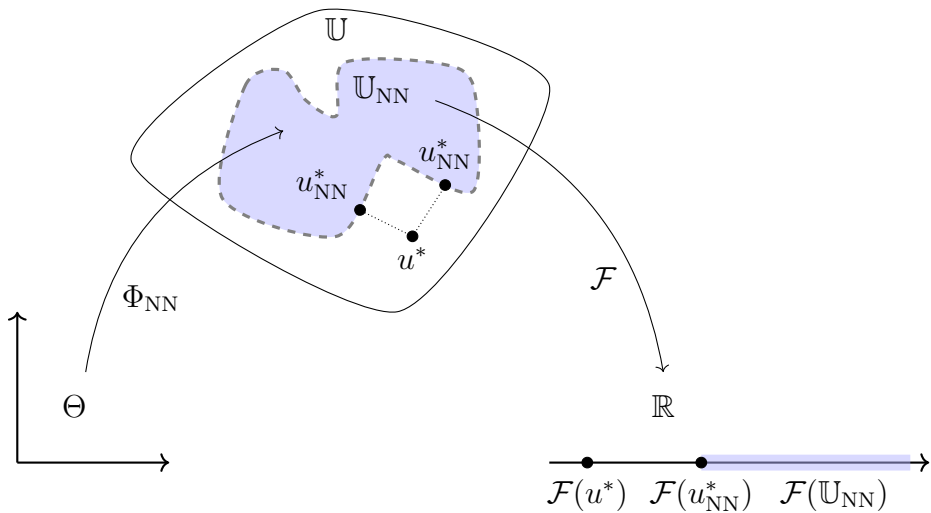
\begin{figure}[htbp]
\centering
\begin{tikzpicture}
    \draw[smooth cycle, tension=0.4] plot coordinates{(2,2) (-0.5,0) (3,-2) (5,1)} node at (2.2,1.75) {$\mathbb{U}$};
    \draw[dashed, line width=1.2,fill=blue!30!white, opacity=0.5] plot [smooth cycle] coordinates {(0.7,-0.5) (1, 1) (1.5, 1.2) (1.8, 0.8) (2.2, 0.6) (2.3, 1.3) (3.5,1.3) (4,1) (4, -0.4) (3,0) (2.8, 0) (2.2, -1)} ;
    \draw[] node at (2.8,0.95) {$\mathbb{U}_{\text{NN}}$};
    \filldraw (3.2,-1) circle (2pt) node[below=0.001] {$u^*$};
    \draw[densely dotted, line width=0.4pt] (3.2,-1) -- (2.46, -0.63);
    \filldraw (2.51, -0.65) circle (2pt) node[above left=0.0001cm] {$u^*_\text{NN}$};
    \draw[densely dotted, line width=0.4pt] (3.2,-1) -- (3.65, -0.3);
    \filldraw (3.63, -0.32) circle (2pt) node[above=0.07cm] {$u^*_\text{NN}$};

    \draw[thick, ->] (-2,-4) -- (-2, -2) node {};
    \draw[thick, ->] (-2,-4) -- (0, -4) node {};
    \draw[] node at (-1.25,-3.25) {$\Theta$};
    
    \path[->] (-1.1, -2.8) edge [bend left] node [right, xshift=-5mm, yshift=-10mm] {$\Phi_\text{NN}$} (1.5, 0.4);
    \path[->] (3.5, 0.8) edge [bend left] node[left, xshift=5mm, yshift=-10mm] {$\mathcal{F}$} (6.5, -2.8);

    \draw[thick, ->] (5, -4) -- (10, -4);
    \draw[color=blue!30!white, line width = 6pt, opacity=0.5] (7, -4) -- (9.75, -4);
    \draw[] (8.75, -4) node[below=0.07cm]  {$\mathcal{F}(\mathbb{U}_\text{NN})$};
   \filldraw (5.5,-4) circle (2pt) node[below=0.07cm]  {$\mathcal{F}(u^*)$};
   \filldraw (7,-4) circle (2pt) node[below=0.07cm]  {$\mathcal{F}(u^*_\text{NN})$};
   \draw[] node at (6.5,-3.25) {$\mathbb{R}$};
\end{tikzpicture}
\caption{Sketch of the parameterization of a NN and its objective-function minimization formulation. }
\label{Chapter1_NN_topology}
\end{figure}

Let  $u^*:X\longrightarrow Y$ be a function that belongs to a certain suitable space of functions $\mathbb{U}$, e.g., an infinite-dimensional vector space, and let $\mathcal{F}:\mathbb{U}\longrightarrow\mathbb{R}$ be an \emph{objective function} that characterizes $u^*$ as its well-defined minimizer,  i.e.,
\begin{equation}
u^*=\arg\min_{u\in \mathbb{U}}\mathcal{F}(u).
\end{equation}  Let $u_{\text{NN}}:X\times\Theta\longrightarrow Y$ be a \ac{NN} architecture according to definition \eqref{Chapter1_FFNN_definition}, and assume that the corresponding set of realizations $\mathbb{U}_{\text{NN}}$ is embedded in $\mathbb{U}$.  Then,  we define \emph{an optimal approximation of} $u^*$ \emph{via} $u_{\text{NN}}$ as
\begin{equation}\label{Chapter1_optiml_approximation_NNs}
u_{\text{NN}}^*\in\mathbb{U} \text{ such that } \mathcal{F}(u_\text{NN}^*)= \inf_{u_{\text{NN}}\in \mathbb{U}_{\text{NN}}} \mathcal{F}(u_{\text{NN}}).
\end{equation} Although $u_\text{NN}^*$ may not belong to $\mathbb{U}_\text{NN}$,  by assuming that both $\Phi_\text{NN}$ and $\mathcal{F}$ are continuous on $\Theta$ and $\mathbb{U}$, respectively, we obtain that both $\mathbb{U}_\text{NN}$ and $\mathcal{F}(\mathbb{U}_\text{NN})$ are (path-)connected.  Then,  there always exists a sequence of elements in $\mathbb{U}_\text{NN}$ that converges to the accumulation point $u_\text{NN}^*$.  In other words,  $u_\text{NN}^*$ can be approximated with arbitrary precision by means of $\mathbb{U}_\text{NN}$,  as well as $\mathcal{F}(u_\text{NN}^*)$ by means of $\mathcal{F}(\mathbb{U}_\text{NN})$.  \Cref{Chapter1_NN_topology} illustrates this setup.

So far,  we presented the \emph{parameterization} in \acp{NN} as an embedding of a (finite-dimensional) Euclidean space $\Theta$ into $\mathbb{U}$ via the realization mapping $\Phi_\text{NN}$.  Some authors might refer to this process as a \emph{discretization} due to the utilization of a finite number of parameters to represent the approximating functions.  However,  to avoid confusion,  we will avoid the latter term to distinguish it from the one concerning the input domain $X$ that we discuss below. 

Assume that the objective function possesses the following integral form:
\begin{equation}
\mathcal{F}(u):=\int_X I(u)(x) \ dx,\qquad u\in\mathbb{U},
\end{equation} where $I(u)$ stands for the integrand that consists of transformations of $u$.  Then,  the natural restriction of $\mathcal{F}$ to $\mathbb{U}_\text{NN}$ is by means of the parameterization $\Phi_\text{NN}$, i.e., 
\begin{equation}
\mathcal{F}(\Phi_\text{NN}(\theta))=\int_X I(\Phi_\text{NN}(\theta))(x) \ dx,\qquad \theta\in\Theta.
\end{equation} To address the challenges associated with the exact integration of \ac{NN},  we consider
\begin{subequations}
\label{Chapter1_discrete_stochastic_functional}
 \begin{equation}
\mathcal{L}:\Theta\times X^N \longrightarrow\mathbb{R}
\end{equation} defined by 
\begin{equation}
\mathcal{L}(\theta; x_1,x_2,\ldots,x_N)=\sum_{i=1}^N \omega_i \ I(\Phi_\text{NN}(\theta))(x_i),
\end{equation}
where $x_i\in X$ and $\omega_i>0$ are conveniently selected integration points and weights, respectively,  so as to resemble a suitable quadrature rule to approximate $\mathcal{F}\circ\Phi_\text{NN}$, i.e.,
\begin{equation}
\mathcal{F}(\Phi_\text{NN}(\theta))\approx\mathcal{L}(\theta;x_1,x_2,\ldots, x_N),\qquad \theta\in\Theta,  x_i\in X, 1\leq i\leq N.
\end{equation}
\end{subequations} Thus, we obtain a \emph{parameterized} and \emph{discretized} representative $\mathcal{L}$ of the objective function $\mathcal{F}\vert_{\mathbb{U}_\text{NN}}$ intended to be computable.  Moreover,  we may endow a \emph{stochastic} nature to \eqref{Chapter1_discrete_stochastic_functional} by sampling the integration points $\{x_i\}_{i=1}^N \subset X$ according to a random distribution every time we estimate the value of the integral.  \Cref{Chapter1_Monte_Carlo} shows a widely used stochastic integral approximation approach when using \acp{NN}.

\begin{example}[Monte Carlo integration] \label{Chapter1_Monte_Carlo}
For $\text{Vol}(X):=\int_X 1 dx<\infty$ and $\{x_i\}_{i=1}^N\subset X$ following a random uniform distribution,  \emph{\acf{MC} integration} consists of selecting $\omega_i=\text{Vol}(X)/N$ in \eqref{Chapter1_discrete_stochastic_functional} for all $1\leq i\leq N$, i.e.,
\begin{subequations}
\begin{alignat}{3}
&\mathcal{F}(\Phi_\text{NN}(\theta)) &&=\int_X I(\Phi_\text{NN}(\theta))(x) \ dx, \qquad &&\theta\in\Theta, \\ 
& \; &&\approx\frac{\text{Vol}(X)}{N}\sum_{i=1}^N I(\Phi_\text{NN}(\theta))(x_i)= \mathcal{L}(\theta; \{x_i\}_{i=1}^N),  \qquad &&\theta\in\Theta.
\end{alignat}
\end{subequations}
\end{example}

Henceforth, we will refer to the parameterized, discretized, and possibly stochastic rendition of the objective function as the \emph{loss function}. In this way, we distinguish the continuum-level functional that characterizes our minimization formulation, $\mathcal{F}$,  from its ``computationally feasible'' counterpart, $\mathcal{L}$. 

In data science, it is common to directly encounter the minimization formulation in a discretized and stochastic form.  \Cref{Chapter1_supervised_learning} presents a typical supervised learning approach with \acp{NN}.

\begin{example}[Supervised learning with \acp{NN}]\label{Chapter1_supervised_learning}
Let $D=\{x_i, y_i\}_{i \geq 1}\subset X\times Y$ be an available labeled large database, and let $u_\text{NN}$ be a given \ac{NN} model. Then,  the approximation task is typically presented as the minimization of
\begin{equation}
\mathcal{L}(\theta; \{x_i,y_i\}_{i=1}^N)=\frac{1}{N}\sum_{i=1}^N \Vert u_{\text{NN}}(x_i; \theta) - y_i\Vert, 
\end{equation} where $\Vert\cdot\Vert$ is a pre-established discrete norm and $\{x_i, y_i\}_{i=1}^N$ is a stochastically selected subsample from $D$, usually called \emph{batch}, such that $N<<\vert D\vert < \infty$.  

This loss-function-level presentation could be viewed as a Monte Carlo approximation of the continuum-level objective function
\begin{equation}
\mathcal{F}(u_{\text{NN}}) = \frac{1}{\text{Vol}(X)}\int_X \Vert u_{\text{NN}}(x;\theta)-u^*(x)\Vert \ dx, 
\end{equation} where $u^*$ stands for the continuum-level function from which we extracted the labeled database, i.e., $y_i=u^*(x_i)$.
\end{example}

\emph{Overfitting} \cite{hawkins2004problem, claeskens2008model,goodfellow2016deep} is a common issue during training, which makes the \ac{NN} learn to perform exceptionally well on the training data but struggles to generalize to new (unseen) data.  In our integration context, this usually means that the \ac{NN} focuses on satisfying the quadrature requirement only at the integration points (the training data) but possibly misbehaves elsewhere.  In \cite{rivera2022quadrature},  the authors show situations where integrating the \ac{NN} with few fixed points leads to greedy behavior during training, producing extremely undesirable solutions.  To overcome this, MC integration is typically utilized due to its stochastic nature.  Unfortunately, MC integration demands an immense sample size to avoid making large integration errors (recall \Cref{section1.5.3}).

\section{Gradient-based training}\label{section2.3}

To carry out the minimization,  we adopt first-order gradient-based methods.  Here,  the term ``first order'' specifies that we \emph{only} calculate gradients with respect to the learnable parameters (and not higher-order derivatives).

The cost of computing the gradient of a \emph{scalar-valued function} is known to be bounded by the cost of evaluating the function itself times a constant when using a usual backpropagation algorithm via \acf{AD} \cite{griewank2008evaluating, baydin2018automatic, margossian2019review}.  Indeed,  such a constant is known to be smaller than or equal to three under specific theoretical conditions (Baur-Strassen Theorem \cite{baur1983complexity}). Hence,  the computational cost of calculating the gradient of the objective/loss function is comparable to simply evaluating it. 

However, this fact cannot be extrapolated to \emph{vector-valued functions}. There, the cost of computing its gradient behaves as the cost of computing the componentwise gradients,  i.e.,  it is comparable to the cost of evaluating the vector function multiplied by the number of its (scalar-valued) output components. Then, second-order minimization methods for \acp{NN} have a cost proportional to evaluating the loss/objective function times the number of learnable parameters.  Since \acp{NN} typically follow an over-parameterized design to exploit their approximation capacities, second or higher-order optimization methods become extremely expensive.

In what follows, we present the gradient-descent method that serves as the foundation of the majority of currently employed optimizers when training \acp{NN}.  We first introduce this method at a continuum level to later deduce its parameterized, discretized, and stochastic counterpart.

\subsection{Continuum-level Gradient-Descent method}

\begin{figure}[htbp]
\centering
\begin{tikzpicture}
\begin{axis}[width=0.7*\textwidth, view = {45}{55},  zticklabels={,,},yticklabels={,,},xticklabels={,,},  ylabel = {$u\in \mathbb{U}$}, zlabel={$\mathcal{F}(u)\in\mathbb{R}$},%
    declare function={f(\x,\y)=cos(deg(\x)*0.8)*cos(deg(\y)*0.6)*exp(0.1*\x);}]
 \addplot3[surf, opacity=0.5, domain=-2.0:2.0, y domain=-7:-3
 ]{f(x,y)};
 \edef\myx{-2} 
 \edef\myy{-3.75} 
 \edef\mystep{-0.65}
 \pgfmathsetmacro{\myf}{f(\myx,\myy)}
 \edef\lstCoords{(\myx,\myy,\myf)}
 \pgfplotsforeachungrouped\X in{0,...,12}
 {
 \pgfmathsetmacro{\myx}{\myx+\mystep*xgrad(\myx,\myy)}
 \pgfmathsetmacro{\myy}{\myy+\mystep*ygrad(\myx,\myy)}
 \pgfmathsetmacro{\myf}{f(\myx,\myy)}
 \edef\lstCoords{\lstCoords\space (\myx,\myy,\myf)}
 }
 \addplot3[samples y=0,arrowed, line width=.6] coordinates \lstCoords;
\end{axis}
\end{tikzpicture}
\caption{Gradient-descent performance sketch in a convex scenario.}
\label{Chapter1_GD_sketch}
\end{figure}
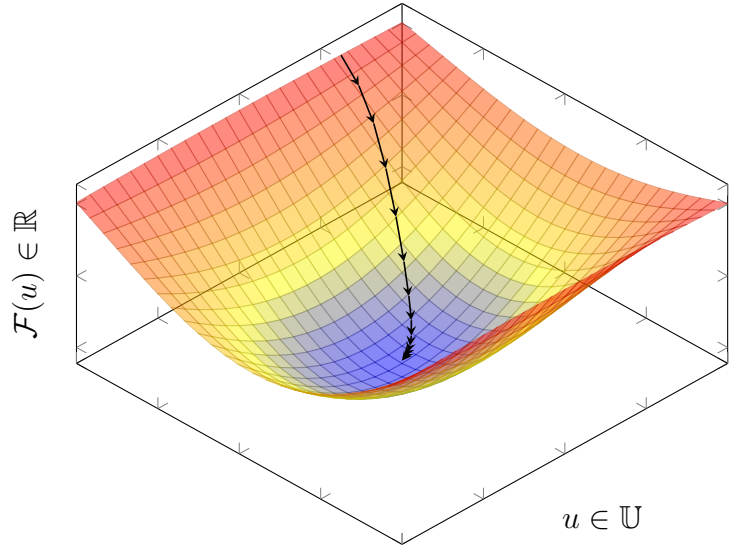

Let $\mathbb{U}$ be a Hilbert space and,  in addition to the well-posedness, assume that $\mathcal{F}:\mathbb{U}\longrightarrow\mathbb{R}$ is \emph{convex}, i.e.,  for $u, w\in\mathbb{U}$ and $0\leq \lambda\leq 1$,  it satisfies
\begin{equation}\label{Chapter1_convex_definition}
\mathcal{F}(\lambda u+(1-\lambda)w) \leq \lambda\mathcal{F}(u)+(1-\lambda)\mathcal{F}(w),
\end{equation} and \emph{sufficiently differentiable}.  Then,  the gradient\footnote{In Hilbert spaces,  we refer to the operator defined via the Riesz representatives of the Fréchet derivative at each point, i.e.,  at a Fréchet differentiable point $u\in\mathbb{U}$,  $\nabla\mathcal{F}(u)$ denotes the unique element in $\mathbb{U}$ that satisfies $D_u\mathcal{F}(w) = (w,\nabla\mathcal{F}(u))_\mathbb{U}$ for all $w\in\mathbb{U}$, where $D_u\mathcal{F}$ denotes the usual Fréchet derivative at $u\in\mathbb{U}$ and $(\cdot,\cdot)_\mathbb{U}$ stands for the inner product in $\mathbb{U}$. } 
operator $\nabla \mathcal{F}$ is monotonically increasing (Kachurovskii's Theorem \cite{kachurovskii1960monotone}),  and thus provides an ideal framework to carry out a continuum-level gradient-descent minimization \cite{smyrlis2004local} defined as follows:
\begin{equation}\label{Chapter1_gradient_descent_continuum}
u_t = u_{t-1} - \eta \nabla\mathcal{F}(u_{t-1}),
\end{equation} where $\eta>0$ is the \emph{learning rate} and $u_{t}\in\mathbb{U}$ is the function at the $t$$^\text{th}$ iteration.  Intuitively,  the gradient-descent method begins selecting an initial candidate $u^{(0)}\in\mathbb{U}$ and, iteratively,  follows steps against the gradient, which by an appropriate (adaptive) control of the size of the learning rate,  it is well known that we can approach $u^*$ as much as desired \cite{berger1977nonlinearity, smyrlis2004local}.   \autoref{Chapter1_GD_sketch} illustrates a graphical performance of the gradient-descent method on a convex scenario, and \Cref{Chapter1_example_gradient_descent} briefly analyzes a simple continuum-level case of study.

\begin{example}\label{Chapter1_example_gradient_descent}
Let $\Vert\cdot\Vert_\mathbb{U}$ denote the norm induced from the inner product of the Hilbert space $\mathbb{U}$.  Consider the objective function defined by $\mathcal{F}(u):=\Vert u\Vert_\mathbb{U}$. Then,  $\mathcal{F}$ is convex on $\mathbb{U}$ and Fréchet differentiable on $\mathbb{U}\setminus\{0\}$ with gradient $\nabla \mathcal{F}(u)=\frac{u}{\Vert u\Vert_\mathbb{U}}$.  Thus,  given an initial candidate $u_{0}\in\mathbb{U}\setminus\{0\}$, the continuum-level gradient-descent iterative method described in \eqref{Chapter1_gradient_descent_continuum} becomes
\begin{equation}
u_{t} = u_{t-1} - \eta \frac{u_{t-1}}{\Vert u_{t-1} \Vert_\mathbb{U}} = \left(1 - \frac{\lambda}{\Vert u_{t-1} \Vert_\mathbb{U}}\right) u_{t-1}.
\end{equation} Selecting $0<\eta=\eta(t)\leq \Vert u_{t-1} \Vert_\mathbb{U}$,  the iterative method converges to $0\in\mathbb{U}$,  which is the global minimum of $\mathcal{F}$ on $\mathbb{U}$.  Indeed, if at any iteration we select $\eta(t) = \Vert u_{t-1} \Vert_\mathbb{U}$,  then we reach the minimum.
\end{example}

\subsection{Parameterized-level Gradient-Descent method}

Because the continuum-level gradient-descent method is computationally intract\-able, we obtain a practical method by composing the objective function with the realization mapping,
\begin{equation}
\mathcal{F}\circ \Phi_\text{NN}:\Theta\longrightarrow\mathbb{R},
\end{equation} and apply the corresponding gradient operator $\nabla(\mathcal{F}\circ\Phi_\text{NN})$.  Because we are now in a (finite-dimensional) Euclidean domain, we identify the gradient with the vector of partial derivatives given by
\begin{equation}
\frac{\partial(\mathcal{F}\circ\Phi_\text{NN})}{\partial\theta}:=\\\left[\frac{\partial(\mathcal{F}\circ\Phi_\text{NN})}{\partial\theta^{(i)}}\right]_{i=1}^{\text{dim}(\Theta)},
\end{equation} where $\theta^{(i)}$ denotes the $i$$^\text{th}$ learnable parameter of $\theta$.  Consequently,  the gradient-descent iterative method translates into the domain of learnable parameters as follows:
\begin{equation}\label{Chapter1_GD}
\theta_{t} = \theta_{t-1} - \eta \frac{\partial(\mathcal{F}\circ\Phi_\text{NN})}{\partial\theta}(\theta_{t-1}), 
\end{equation} where $\theta_{t}\in\Theta$ stands for the set of learnable parameters at the $t^\text{th}$ iteration.

While \eqref{Chapter1_GD} can be seen as a parameterized-level version of \eqref{Chapter1_gradient_descent_continuum}, it is essential to note that some of the ideal conditions observed at the continuum level in the gradient-descent method, such as the convexity, do not necessarily satisfy at the parameterized level when dealing with \acp{NN}.  Indeed,  \eqref{Chapter1_GD} merely aims to converge to parameters in $\Theta$ whose partial derivatives tend to vanish,  which may lead to undesirable results (e.g., when falling in local minima or saddle points).  This supposes another departure from the usual vector-space-based approach where the convexity is preserved under linear parameterizations (see \Cref{Chapter1_convexity_linear_combination}). 

\begin{example}[Convexity preservation under linear transformation]\label{Chapter1_convexity_linear_combination}
For simplicity, we first present the case of linear combinations.  Let $\mathcal{F}$ be convex on $\mathbb{U}$ and consider $\Phi$ defined by
\begin{equation}
\Phi(\theta) := \sum_{i=1}^{\text{dim}(\Theta)} \theta^{(i)} u^{(i)},\qquad \theta=\{\theta^{(i)}\}_{i=1}^{\text{dim}(\Theta)},
\end{equation} where $\{u^{(i)}\}_{i=1}^{\text{dim}(\Theta)}\subset\mathbb{U}$ is a pre-established set of functions.  Then,  $\Phi(\Theta)$ is a vector space of dimension at most $\text{dim}(\Theta)$ and  $\mathcal{F}\circ\Phi$ is convex on $\Theta$. 
\begin{proof}
The dimension of $\Phi(\Theta)$ is the number of linearly independent functions in $\{u^{(i)}\}_{i=1}^{\text{dim}(\Theta)}\subset\mathbb{U}$. To prove the convexity of $\mathcal{F}\circ\Phi$,  let $0\leq \lambda\leq 1$ and $\theta,\vartheta\in\Theta$.  By linearity,  $\Phi(\lambda \theta + (1-\lambda)\vartheta) = \lambda \Phi(\theta) + (1-\lambda)\Phi(\vartheta)$. We obtain the desired result by applying $\mathcal{F}$ and the convexity definition.
\end{proof}

Following a similar reasoning,  it is straightforward to check that the convexity condition is also preserved when replacing the linear combination with a general linear transformation.  In particular,  the restriction of $\mathcal{F}\circ\Phi_{\text{NN}}$ to the learnable parameters of the output layer in a \ac{FFNN} according to \Cref{section2.1} is convex whenever $\mathcal{F}$ is convex.
\end{example}

\subsection{Stochastic Gradient-Descent method}

To avoid the computational challenges associated with the integral form of $\mathcal{F}$,  we resort to the loss-function approximation introduced in \eqref{Chapter1_discrete_stochastic_functional}.  Then,  the so-called \emph{\acf{SGD} optimizer} arises as
\begin{equation}\label{Chapter1_SGD}
\theta_{t} = \theta_{t-1} - \eta \frac{\partial\mathcal{L}}{\partial\theta}(\theta_{t-1}; \mathbf{x}_{t}),
\end{equation} where $\mathbf{x}_{t}\in X^N$ denotes the stochastically sampled set of integration points at the $t$$^\text{th}$ iteration.  Note that \eqref{Chapter1_SGD} introduces discrete dependencies and uncertainty during training compared to \eqref{Chapter1_GD}.

Multiple variants of stochastic gradient-based optimizers have been designed and exploited in recent years with the aim of improving many of the aforementioned poor conditions in which \acp{NN} are immersed.  Here are some of the most popular today: SGD with momentum \cite{qian1999momentum} and/or Nesterov acceleration \cite{nesterov1983method},  Adagrad \cite{duchi2011adaptive}, Adadelta \cite{zeiler2012adadelta},  RMSprop \cite{hinton2012neural}, Adam/AdaMax \cite{kingma2014adam},  and  Nadam \cite{dozat2016incorporating}. We refer to \cite{ruder2016overview} for an overview.

\section{A case of study}\label{section2.4}

To illustrate many of the (undesirable) properties pointed out so far, we discuss the following easy-to-analyze case of study.

At the continuum level,  let $X:=(-1,1)$,  $Y:=\mathbb{R}$,  and
\begin{equation}
u^*(x):=\text{sgn}(x)=\begin{cases}-1,&\text{if } x<0,\\ \phantom{-}0,&\text{if } x=0,\\ \phantom{-}1,&\text{if } x>0.\end{cases}
\end{equation} Consider $\mathbb{U}:=L^2(X)=\{u:X\longrightarrow \mathbb{R} : \int_X u^2 <\infty\}$ and let 
\begin{equation}\label{Chapter1_case_study}
\mathcal{F}(u):=\Vert u-u^*\Vert_{L^2(X)}^2=\int_X (u(x)-u^*(x))^2 \ dx,
\end{equation} whose minimum in $\mathbb{U}$ is zero and is attained only at $u^*$.

\begin{figure}[htbp]
\centering
\begin{subfigure}{\textwidth}
\centering
\begin{tikzpicture} 
\begin{axis}[width=0.7*\textwidth,  view = {0}{90},  colorbar,  colorbar style={fill opacity=0.5, tick label style={opacity=1}}, xlabel = {$a$}, ylabel={$c$}, title={$\displaystyle \frac{\partial(\mathcal{F}\circ\Phi_\text{NN})}{\partial a}(a,c)$}]
    \addplot3 [surf, shader = interp, opacity=0.4, domain=-5:5, samples=50,
 ]{(2*y/x^2)*(2*ln(cosh(x))+ (y-2*x)*tanh(x)- x*y*(1-(tanh(x))^2))};
     \draw[color=black, line width=0.8, dashed] plot[variable=\x,  domain=-5:5,samples=5,smooth] (\x,  0,1);
     \draw[color=black, line width=0.8, dashed] plot[variable=\x,  domain=-5:-0.1,samples=50,smooth] (\x,  {2*( ln(cosh(\x)) - \x * tanh(\x) ) / ( \x * ( 1 - (tanh(\x))^2 ) - tanh(\x) )},  1);
     \draw[color=black, line width=0.8, dashed] plot[variable=\x,  domain=0.1:5,samples=50,smooth] (\x,  {2*( ln(cosh(\x)) - \x * tanh(\x) ) / ( \x * ( 1 - (tanh(\x))^2 ) - tanh(\x) )},  1);
     \node (a) at (1.9,-1.3) {$c=\frac{2(\log(\cosh(a)) - a\tanh(a))}{a \text{ sech}^2(a) - \tanh(a))}$};
     \node (b) at (-2,0.5) {$c=0$};
     \draw[->] (a) -- (3, {2*( ln(cosh(3)) - 3 * tanh(3) ) / ( 3 * ( 1 - (tanh(3))^2 ) - tanh(3) )-0.2});
     \draw[->] (a) -- (-0.25,{2*( ln(cosh(-0.4)) - (-0.4) * tanh(-0.4) ) / ( -0.4 * ( 1 - (tanh(-0.4))^2 ) - tanh(-0.4) )}); 
\end{axis}
\end{tikzpicture}
\end{subfigure}
\vskip 1em
\begin{subfigure}{\textwidth}
\centering
\begin{tikzpicture} 
\begin{axis}[width=0.7*\textwidth,  view = {0}{90},  colorbar,  colorbar style={fill opacity=0.5, tick label style={opacity=1}}, xlabel = {$a$}, ylabel={$c$}, title={$\displaystyle \frac{\partial(\mathcal{F}\circ\Phi_\text{NN})}{\partial c}(a,c)$}]
    \addplot3 [surf, shader = interp, opacity=0.4, domain=-5:5, samples=50]{(4/x)*(x*y-ln(cosh(x))-y*tanh(x))};
     \draw[color=black, line width=0.8, dashed] plot[variable=\x,  domain=-5:5,samples=5,smooth] (0, \x,1);
     \draw[color=black, line width=0.8, dashed] plot[variable=\x,  domain=-5:-0.1,samples=50,smooth] (\x,  {ln(cosh(\x))/(\x-tanh(\x))},  1);
     \draw[color=black, line width=0.8, dashed] plot[variable=\x,  domain=0.1:5,samples=50,smooth] (\x,  {ln(cosh(\x))/(\x-tanh(\x))},  1);
     \node (a) at (2.5,-2) {$c=\frac{\log(\cosh(a))}{a-\tanh(a)}$};
     \node (b) at (-1,2) {$a=0$};
     \draw[->] (a) -- (3, {ln(cosh(3))/(3-tanh(3))-0.2});
     \draw[->] (a) -- (-0.25,{ln(cosh(-0.4))/(-0.4-tanh(-0.4))}); 
\end{axis}
\end{tikzpicture}
\end{subfigure}
\caption{Graph of the partial derivatives of $\mathcal{F}\circ\Phi_\text{NN}$ for \eqref{Chapter1_case_study}. Dashed curves indicate vanishing regions.}
\label{Chapter1_analytic_gradF}
\end{figure}
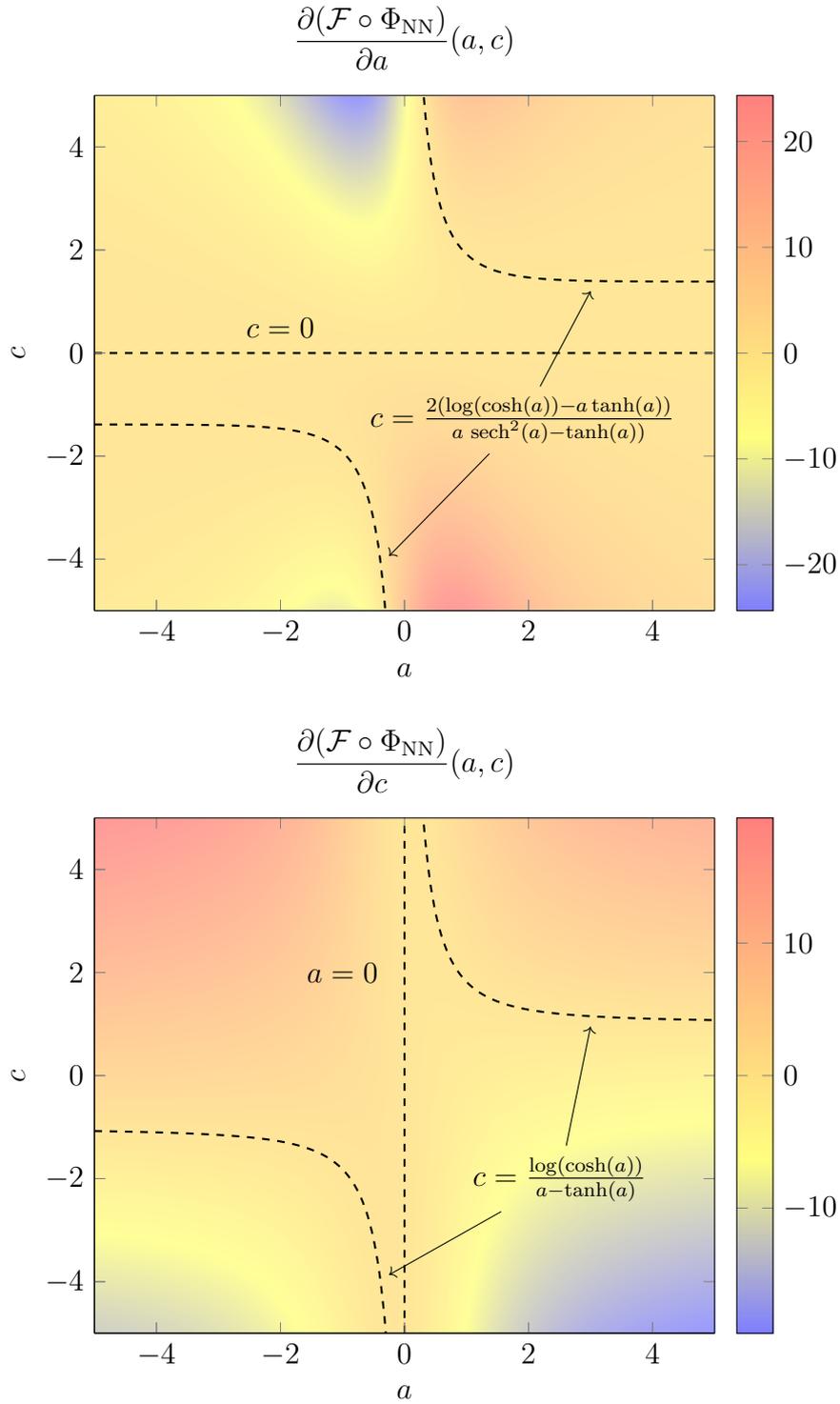

At the parameterization level,  let $\Theta:=\mathbb{R}^2$ and let $u_\text{NN}:X\times\Theta\longrightarrow\mathbb{R}$ be defined by
\begin{equation}
u_\text{NN}(x;\theta):=c\tanh(ax),
\end{equation} where $\theta=\{a,c\}$ denotes the set of learnable parameters.  Then,  
\begin{equation}
\mathbb{U}_\text{NN}=\Phi_\text{NN}(\Theta)=\{\Phi_\text{NN}(a,c) = c\tanh(ax):X\longrightarrow Y\}_{a,c\in\mathbb{R}}, 
\end{equation} and routine calculations yield
\begin{align}
(\mathcal{F}\circ\Phi_\text{NN})(a,c)&= 2(1+c^2) - \frac{4c\log(\cosh(a))+2c^2\tanh(a)}{a},  \label{Chapter1_objetive_function}\\
\frac{\partial(\mathcal{F}\circ\Phi_\text{NN})}{\partial a}(a,c)&= \frac{2c\left( 2\log(\cosh(a))+ (c-2a)\tanh(a)- ac(\text{sech}^2(a)\right)}{a^2}, \\
\frac{\partial(\mathcal{F}\circ\Phi_\text{NN})}{\partial c}(a,c)&= \frac{4\left(ac-\log(\cosh(a))-c\tanh(a)\right)}{a},
\end{align} for $a\neq 0$,  $(\mathcal{F}\circ\Phi_\text{NN})(0,c)=2$,  $\frac{\partial(\mathcal{F}\circ\Phi_\text{NN})}{\partial c}(0,c)=-2c$,  and $\frac{\partial(\mathcal{F}\circ\Phi_\text{NN})}{\partial a}(0,c)=0$.  \Cref{Chapter1_analytic_F} shows the graph of $\mathcal{F}\circ\Phi_\text{NN}$,  and \Cref{Chapter1_analytic_gradF} shows the graphs of its partial derivatives.  

\begin{figure}[htbp]
\centering
\begin{subfigure}{\textwidth}
\centering
\begin{tikzpicture}
\begin{axis}[width=0.7*\textwidth,  view = {35}{55},  colorbar,  colorbar style={fill opacity=0.4, tick label style={opacity=1}}, xlabel = {$a$}, ylabel={$c$}, ytick={-4,-2,0,2,4},zlabel={$(\mathcal{F}\circ\Phi_\text{NN})(a,c)$},ztick={0,25,50}]
 \addplot3[surf, opacity=0.4, domain=-5:5,  samples=50,
 ]{2*(1+y^2)-(4*y*ln(cosh(x))+2*y^2*tanh(x))/x};
\end{axis}
\end{tikzpicture}
\end{subfigure}
\vskip 1em
\begin{subfigure}{\textwidth}
\centering
\begin{tikzpicture} 
\begin{axis}[width=0.7*\textwidth,  view = {0}{90},  colorbar,  colorbar style={fill opacity=0.5, tick label style={opacity=1}}, xlabel = {$a$}, ylabel={$c$}, title={$(\mathcal{F}\circ\Phi_\text{NN})(a,c)$}]
    \addplot3 [surf, shader = interp, opacity=0.4, domain=-5:5, samples=50,
 ]{2*(1+y^2)-(4*y*ln(cosh(x))+2*y^2*tanh(x))/x};
\end{axis}
\end{tikzpicture}
\end{subfigure}
\caption{Graph of $\mathcal{F}\circ \Phi_\text{NN}$ for \eqref{Chapter1_case_study}. Top panel: 3D view. Bottom panel: Top view.}
\label{Chapter1_analytic_F}
\end{figure}
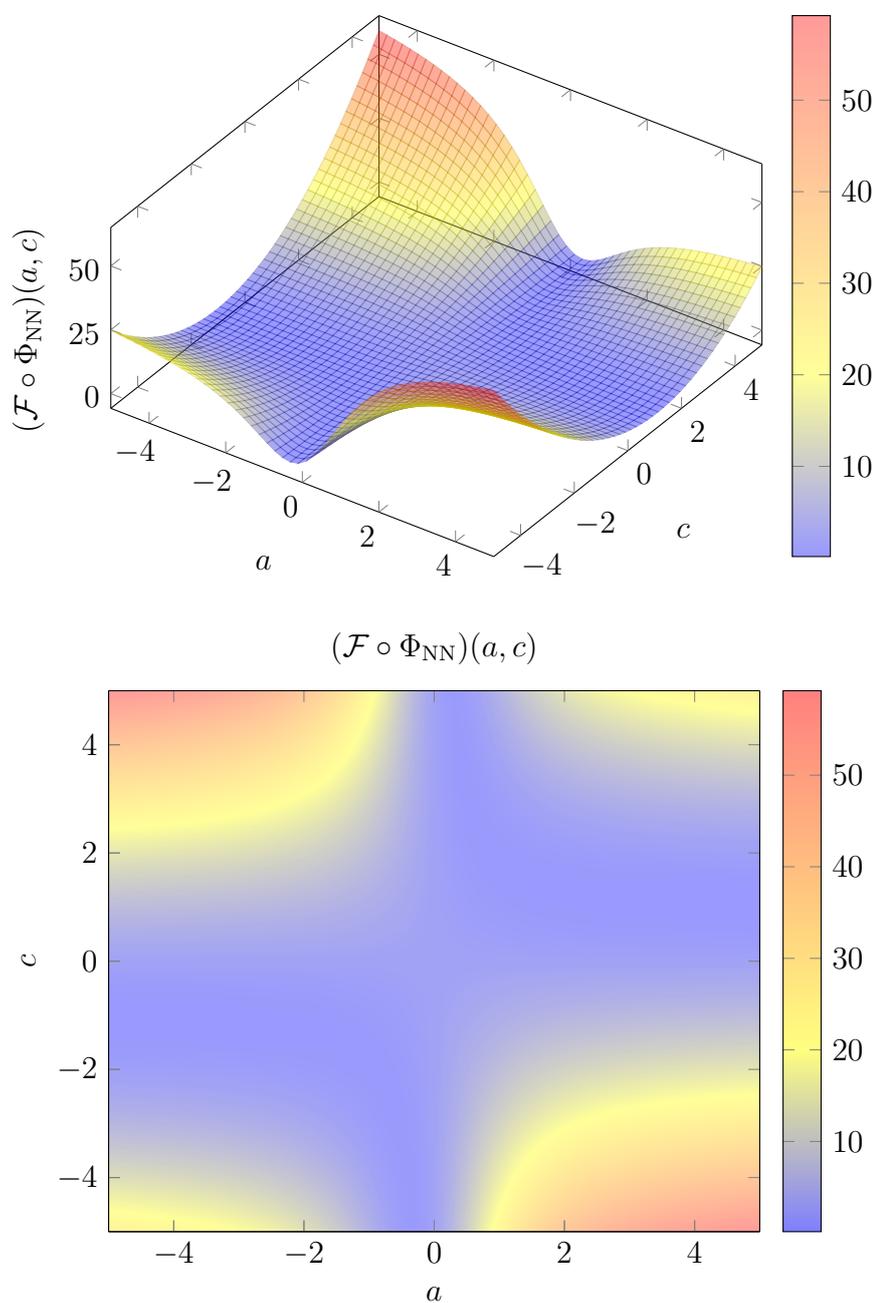

The following statements expose a series of properties of this case of study:

\begin{enumerate}
\item The architecture $u_\text{NN}$ posseses \emph{arbitrary precision approximation capacity},  which means that $u^*_\text{NN}=u^*$.
\begin{proof}
Let $c=1$.  Then,  it is straightforward to check that $(\mathcal{F}\circ\Phi_\text{NN})(a,1)$ converges to zero as $a\to \infty$,  which implies that $\Phi_\text{NN}(1,a)$ converges to $u^*$ as $a\to \infty$.
\end{proof}
\item The set of realizations $\mathbb{U}_\text{NN}$ is \emph{non-closed}.
\begin{proof}
The optimal realization $u_\text{NN}^*\in L^2(X)\setminus C^\infty(X)$ is unattainable because $\Phi_\text{NN}(a,c)\in C^\infty(X)$ for all $a,c\in\mathbb{R}$.
\end{proof}
\item The realization map $\Phi_\text{NN}$ is non-injective.  In particular,  the sequence of learnable parameters whose images through $\Phi_\text{NN}$ converge to $u^*_\text{NN}$ is non-unique, although the limit realization $u^*_\text{NN}$ is unique.
\begin{proof}
The odd symmetry of $\tanh(x)$ implies $\Phi_\text{NN}(a,c)=\Phi_\text{NN}(-a,-c)$ for all $(a,c)\in\Theta$,  which shows the non-injectivity.  Then, for $c= 1$, we find distinct sequences $\{(a,1)\}_{a\geq 1}$ and $\{(-a,-1)\}_{a\geq 1}$ of learnable parameters with distinct limits but whose images through $\Phi_\text{NN}$ converge to $u^*_\text{NN}$ as $a\to\infty$. The uniqueness of $u^*_\text{NN}$ follows from item 1.  and the well-posedness of the objective-function minimization. 
\end{proof}
\item The set of realizations $\mathbb{U}_\text{NN}$ is \emph{non-convex}.
\begin{proof}
We need to find $0<\lambda<1$ and $a_1,a_2,c_1,c_2\in\mathbb{R}$ such that $f(x)=\lambda c_1\tanh(a_1 x) + (1-\lambda) c_2\tanh(a_2 x)$ is unrealizable via $\Phi_\text{NN}$, i.e.,  it is not of the form $\Phi_\text{NN}(a,c)=u_\text{NN}(x;a,c) = c\tanh(a x)$ for any $a,c\in\mathbb{R}$.  Indeed,  let $\lambda=1/2$,  $c_1=c_2=2$, $a_1=2$, and $a_2=-1$.  Then,  $f(x)=\tanh(2x)-\tanh(x)$,  which is neither constant nor strict monotonic on $X$.  Because $u_\text{NN}(x;a,c)$ is either strictly monotonic or constant on $X$,  $f(x)$ is unrealizable via $\Phi_\text{NN}$.
\end{proof}
\item $\mathcal{F}$ is \emph{convex} on $\mathbb{U}$,  but  $\mathcal{F}\circ\Phi_\text{NN}$ is \emph{non-convex} on $\Theta$.
\begin{proof}
Let $0\leq \lambda\leq 1$ and $u_1,u_2\in\mathbb{U}$. Then,  the triangular inequality $\Vert \lambda(u_1-u^*) + (1-\lambda)(u_2-u^*) \Vert_{L^2(X)}\leq \lambda\Vert u_1-u^*\Vert_{L^2(X)} + (1-\lambda)\Vert u_2-u^*\Vert_{L^2(X)}$ implies $\mathcal{F}(\lambda u_1+(1-\lambda)u_2)\leq \lambda \mathcal{F}(u_1)+(1-\lambda)\mathcal{F}(u_2)$, which shows the convexity of $\mathcal{F}$.  To show the non-convexity of $\mathcal{F}\circ\Phi_\text{NN}$,  fix any $c\neq 0$ and study the curvature of $f_c(a):=(\mathcal{F}\circ\Phi_\text{NN})(a,c)$---recall \Cref{Chapter1_analytic_F}.  For example,  $f''_{1}(a)$ is positive when $a>0$ and negative when $a<-1$.
\end{proof}
\item The negative gradient field possesses undesirable attractive basin regions that may lead to poor results when performing the gradient-descent method.
\begin{proof}
The only parameter configurations that lead to desirable realizations occur when $c=\pm 1$ and $a\to\pm\infty$,  which implies $\frac{\partial(\mathcal{F}\circ\Phi_\text{NN})}{\partial\{a,c\}}(a,\pm 1)\to(0,0)$ as $a\to\pm\infty$.  However, we can also obtain near-zero gradients in undesirable parameter configurations. \Cref{Chapter1_analytic_F_flow} shows poor gradient-descent performances that get trapped in near-zero-gradient regions. 
\end{proof}

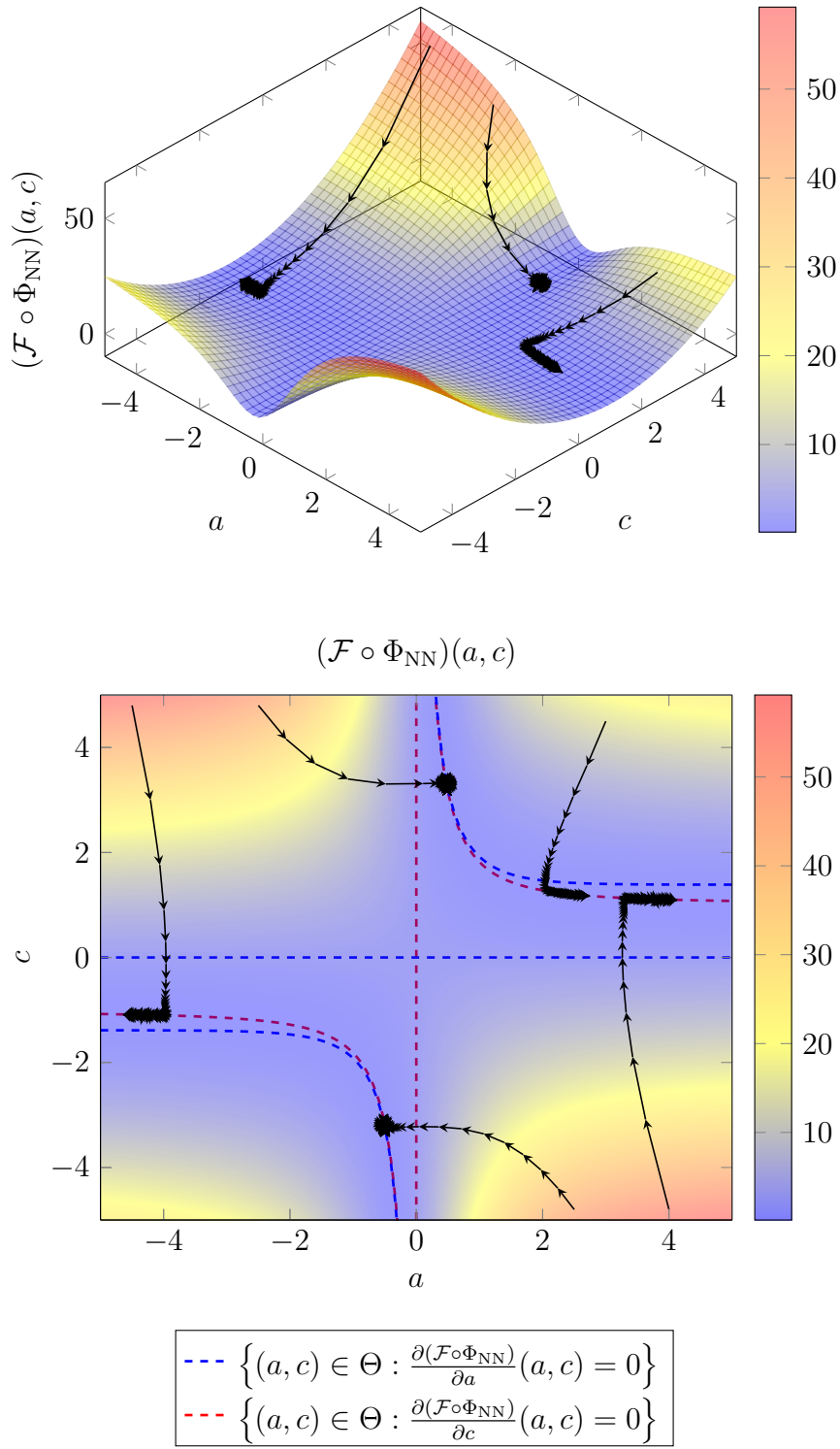
\begin{figure}[htbp]
\centering
\begin{subfigure}{\textwidth}
\centering
\begin{tikzpicture}
\begin{axis}[width=0.7*\textwidth, view = {45}{55},  colorbar,  colorbar style={fill opacity=0.4, tick label style={opacity=1}}, xlabel={$a$},  ylabel = {$c$}, zlabel={$(\mathcal{F}\circ\Phi_\text{NN})(a,c)$},%
    declare function={f(\x,\y)= 2*(1+\y^2)-(4*\y*ln(cosh(\x))+2*\y^2*tanh(\x))/\x;}]

 \addplot3[surf, opacity=0.4, domain=-5:5,samples=50,
 ]{f(x,y)};
 \edef\myx{-2.5} 
 \edef\myy{4.8} 
 \edef\mystep{-0.05}
 \pgfmathsetmacro{\myf}{f(\myx,\myy)}
 \edef\lstCoords{(\myx,\myy,\myf)}
 \pgfplotsforeachungrouped\X in{0,...,200}
 {
 \pgfmathsetmacro{\myx}{\myx+\mystep*xgrad(\myx,\myy)}
 \pgfmathsetmacro{\myy}{\myy+\mystep*ygrad(\myx,\myy)}
 \pgfmathsetmacro{\myf}{f(\myx,\myy)}
 \edef\lstCoords{\lstCoords\space (\myx,\myy,\myf)}
 }
 \addplot3[samples y=0,arrowed, line width=.6, shader=interp] coordinates \lstCoords;
 \edef\myx{-4.5} 
 \edef\myy{4.8} 
 \edef\mystep{-0.1}
 \pgfmathsetmacro{\myf}{f(\myx,\myy)}
 \edef\lstCoords{(\myx,\myy,\myf)}
 \pgfplotsforeachungrouped\X in{0,...,200}
 {
 \pgfmathsetmacro{\myx}{\myx+\mystep*xgrad(\myx,\myy)}
 \pgfmathsetmacro{\myy}{\myy+\mystep*ygrad(\myx,\myy)}
 \pgfmathsetmacro{\myf}{f(\myx,\myy)}
 \edef\lstCoords{\lstCoords\space (\myx,\myy,\myf)}
 }
 \addplot3[samples y=0,arrowed, line width=.6, shader=interp] coordinates \lstCoords;
  \edef\myx{3} 
 \edef\myy{4.5} 
 \edef\mystep{-0.1}
 \pgfmathsetmacro{\myf}{f(\myx,\myy)}
 \edef\lstCoords{(\myx,\myy,\myf)}
 \pgfplotsforeachungrouped\X in{0,...,200}
 {
 \pgfmathsetmacro{\myx}{\myx+\mystep*xgrad(\myx,\myy)}
 \pgfmathsetmacro{\myy}{\myy+\mystep*ygrad(\myx,\myy)}
 \pgfmathsetmacro{\myf}{f(\myx,\myy)}
 \edef\lstCoords{\lstCoords\space (\myx,\myy,\myf)}
 }
 \addplot3[samples y=0,arrowed, line width=.6, shader=interp] coordinates \lstCoords;
\end{axis}
\end{tikzpicture}
\end{subfigure}
\vskip 2em
\begin{subfigure}{\textwidth}
\centering
\begin{tikzpicture}
\begin{axis}[width=0.7*\textwidth, view = {0}{90},  xmin=-5,xmax=5,ymin=-5,ymax=5, colorbar,  colorbar style={fill opacity=0.5,  tick label style={opacity=1}},  xlabel={$a$},  ylabel = {$c$}, zlabel={$(\mathcal{F}\circ\Phi_\text{NN})(a,c)$},title={$(\mathcal{F}\circ\Phi_\text{NN})(a,c)$},%
declare function={f(\x,\y)= 2*(1+\y^2)-(4*\y*ln(cosh(\x))+2*\y^2*tanh(\x))/\x;}]

\draw[color=red, line width=1, dashed] plot[variable=\x,  domain=-5:-0.1,samples=50,smooth] (\x,  {ln(cosh(\x))/(\x-tanh(\x))},  1);
\draw[color=red, line width=1, dashed] plot[variable=\x,  domain=0.1:5,samples=50,smooth] (\x,  {ln(cosh(\x))/(\x-tanh(\x))},  1);
\draw[color=red, line width=1, dashed] plot[variable=\x,  domain=-5:5,samples=50,smooth] (0,\x, 1);
\draw[color=blue, line width=1, dashed] plot[variable=\x,  domain=-5:-0.1,samples=50,smooth] (\x,  {2*( ln(cosh(\x)) - \x * tanh(\x) ) / ( \x * ( 1 - (tanh(\x))^2 ) - tanh(\x) )},  1);
\draw[color=blue, line width=1, dashed] plot[variable=\x,  domain=0.1:5,samples=50,smooth] (\x,  {2*( ln(cosh(\x)) - \x * tanh(\x) ) / ( \x * ( 1 - (tanh(\x))^2 ) - tanh(\x) )},  1);
\draw[color=blue, line width=1, dashed] plot[variable=\x,  domain=-5:5,samples=50,smooth] (\x, 0, 1);

 \addplot3[surf, opacity=0.4, domain=-5:5,samples=50,shader=interp
 ]{f(x,y)};
 \edef\myx{-2.5} 
 \edef\myy{4.8} 
 \edef\mystep{-0.05}
 \pgfmathsetmacro{\myf}{f(\myx,\myy)}
 \edef\lstCoords{(\myx,\myy,\myf)}
 \pgfplotsforeachungrouped\X in{0,...,200}
 {
 \pgfmathsetmacro{\myx}{\myx+\mystep*xgrad(\myx,\myy)}
 \pgfmathsetmacro{\myy}{\myy+\mystep*ygrad(\myx,\myy)}
 \pgfmathsetmacro{\myf}{f(\myx,\myy)}
 \edef\lstCoords{\lstCoords\space (\myx,\myy,\myf)}
 }
 \addplot3[samples y=0,arrowed, line width=.6, shader=interp] coordinates \lstCoords;
 \edef\myx{-4.5} 
 \edef\myy{4.8} 
 \edef\mystep{-0.1}
 \pgfmathsetmacro{\myf}{f(\myx,\myy)}
 \edef\lstCoords{(\myx,\myy,\myf)}
 \pgfplotsforeachungrouped\X in{0,...,200}
 {
 \pgfmathsetmacro{\myx}{\myx+\mystep*xgrad(\myx,\myy)}
 \pgfmathsetmacro{\myy}{\myy+\mystep*ygrad(\myx,\myy)}
 \pgfmathsetmacro{\myf}{f(\myx,\myy)}
 \edef\lstCoords{\lstCoords\space (\myx,\myy,\myf)}
 }
 \addplot3[samples y=0,arrowed, line width=.6, shader=interp] coordinates \lstCoords;
  \edef\myx{3} 
 \edef\myy{4.5} 
 \edef\mystep{-0.1}
 \pgfmathsetmacro{\myf}{f(\myx,\myy)}
 \edef\lstCoords{(\myx,\myy,\myf)}
 \pgfplotsforeachungrouped\X in{0,...,100}
 {
 \pgfmathsetmacro{\myx}{\myx+\mystep*xgrad(\myx,\myy)}
 \pgfmathsetmacro{\myy}{\myy+\mystep*ygrad(\myx,\myy)}
 \pgfmathsetmacro{\myf}{f(\myx,\myy)}
 \edef\lstCoords{\lstCoords\space (\myx,\myy,\myf)}
 }
 \addplot3[samples y=0,arrowed, line width=.6, shader=interp] coordinates \lstCoords;
 \edef\myx{2.5} 
 \edef\myy{-4.8} 
 \edef\mystep{-0.03}
 \pgfmathsetmacro{\myf}{f(\myx,\myy)}
 \edef\lstCoords{(\myx,\myy,\myf)}
 \pgfplotsforeachungrouped\X in{0,...,200}
 {
 \pgfmathsetmacro{\myx}{\myx+\mystep*xgrad(\myx,\myy)}
 \pgfmathsetmacro{\myy}{\myy+\mystep*ygrad(\myx,\myy)}
 \pgfmathsetmacro{\myf}{f(\myx,\myy)}
 \edef\lstCoords{\lstCoords\space (\myx,\myy,\myf)}
 }
 \addplot3[samples y=0,arrowed, line width=.6, shader=interp] coordinates \lstCoords;
  \edef\myx{4} 
 \edef\myy{-4.8} 
 \edef\mystep{-0.1}
 \pgfmathsetmacro{\myf}{f(\myx,\myy)}
 \edef\lstCoords{(\myx,\myy,\myf)}
 \pgfplotsforeachungrouped\X in{0,...,200}
 {
 \pgfmathsetmacro{\myx}{\myx+\mystep*xgrad(\myx,\myy)}
 \pgfmathsetmacro{\myy}{\myy+\mystep*ygrad(\myx,\myy)}
 \pgfmathsetmacro{\myf}{f(\myx,\myy)}
 \edef\lstCoords{\lstCoords\space (\myx,\myy,\myf)}
 }
 \addplot3[samples y=0,arrowed, line width=.6, shader=interp] coordinates \lstCoords;
\end{axis}
\end{tikzpicture}
\end{subfigure}
\vskip 1em
\begin{subfigure}{\textwidth}
\centering
\begin{tikzpicture}
\begin{axis}[hide axis,
	    xmin=-0.1,
	    xmax=1.1,
	    ymin=-0.1,
	    ymax=1.1,
	    legend style={fill=none, row sep=3pt, inner sep = 3pt},
	    legend columns = 1,
            ]

    	\addlegendimage{color=blue, line width=1, dashed};
	\addlegendentry{$\left\lbrace(a,c)\in\Theta:\frac{\partial(\mathcal{F}\circ\Phi_\text{NN})}{\partial a}(a,c)=0\right\rbrace$};
	\addlegendimage{color=red, line width=1, dashed};
	\addlegendentry{$\left\lbrace(a,c)\in\Theta:\frac{\partial(\mathcal{F}\circ\Phi_\text{NN})}{\partial c}(a,c)=0\right\rbrace$};
    \addlegendimage{color=white, line width=1, dashed};
\end{axis}
\end{tikzpicture}
\end{subfigure}
\caption{Different executions of the gradient-descent method for \eqref{Chapter1_case_study}. We maintain constant learning rates in all cases and carry out $200$ iterations. Top panel: 3D view. Bottom panel: Top view.}
\label{Chapter1_analytic_F_flow}
\end{figure}
\end{enumerate}

At the discretization level,  we approximate $\mathcal{F}\circ\Phi_\text{NN}$ via a loss function that we expect to inherit similar conditions to above.  Indeed, let $\mathbf{x}=\{x_i\}_{i=1}^N\subset X$. Then, 
\begin{subequations}
\begin{align}
(\mathcal{F}\circ\Phi_\text{NN})(a,c)\approx \mathcal{L}(a,c;\mathbf{x}) &= \sum_{i=1}^N \omega_i (c\tanh(ax_i)\pm 1)^2,\\
\frac{\partial(\mathcal{F}\circ\Phi_\text{NN})}{\partial a}(a,c)\approx \frac{\partial \mathcal{L}}{\partial a}(a,c;\mathbf{x}) &= 2ac \sum_{i=1}^N \omega_i  (c\tanh(a x_i)\pm 1)\text{sech}^2(a x_i),\\
\frac{\partial(\mathcal{F}\circ\Phi_\text{NN})}{\partial c}(a,c)\approx \frac{\partial \mathcal{L}}{\partial c}(a,c;\mathbf{x}) &= 2 \sum_{i=1}^N \omega_i (c\tanh(a x_i)\pm 1) \tanh(a x_i),
\end{align}
\end{subequations} where $\omega_i$ denotes the integration weight related to the integration point $x_i$, and where the positive or negative sign depends on whether $x_i$ is negative or positive, respectively.

When $c=\pm 1$ and $a\to \pm\infty$, we obtain the desired realizations convergence $\Phi_\text{NN}(a,c)\to u^*$ with corresponding vanishing of the gradient of the loss function
\begin{equation}
\frac{\partial \mathcal{L} }{\partial\{a,c\}}(a,\pm 1;\mathbf{x}) \to (0,0) \text{ as } a\to\pm\infty, \quad \text{for any } \mathbf{x}\in X^N  \text{ and } N>0.
\end{equation} However,  $\frac{\partial \mathcal{L} }{\partial\{a,c\}}(0,c;\mathbf{x})=(0,0)$ for any $c\in\mathbb{R},  \mathbf{x}\in X^N  \text{ and } N>0$,  but $\Phi_\text{NN}(0,c) = 0\neq u^*$ for all $c\in\mathbb{R}$.  Hence,  parameters of the form $(0,c)\in\Theta$ are undesirable critical points susceptible to trapping the learnable parameters during stochastic gradient-descent training.

\begin{subappendices}

\section{Interpretation of the set of realizations as a family of vector spaces}\label[appendix]{appendix2.a}
The set of realizations produced by a given FFNN defined as in \eqref{FFNN_architecture} could be interpreted as a family of finite-dimensional vector spaces (up to a prescribed dimension). Thus,  we argue guided under the following intuition: the collection of hidden layers is in charge of generating different finite-dimensional vector spaces,  while the output layer is in charge of producing linear combinations over such spaces. 

For simplicity, we first consider the scalar-valued case ($n_{K+1}=1$), and then we extend it to the vector-valued case ($n_{K+1}>1$).

\subsubsection*{Scalar-valued FFNN (\boldmath{$n_{K+1}=1$})}

Let $u_\text{NN}$ be a scalar-valued FFNN and let us split its corresponding set of learnable parameters $\theta$ as follows: 
\begin{subequations}
\begin{equation}
\theta:=(\eta,\mathbf{W})\in\Theta=\Lambda \times \mathbb{R}^{1\times n_{K}},
\end{equation} where
\begin{equation}
\eta:=\{\mathbf{W}_j,\mathbf{b}_j\}_{j=1}^K\in\Lambda = \{\mathbb{R}^{n_j\times n_{j-1}} \times \mathbb{R}^{n_j}\}_{j=1}^K.
\end{equation} 
\end{subequations}

Then, given any  $\eta\in\Lambda$,  we obtain a corresponding vector space $\mathbb{U}_\eta$ spanned by the component functions of $\mathbf{y}_K(\cdot\; ;\eta)$,
\begin{subequations}
\begin{align}
    \mathbb{U}_\eta := \text{span}\big\{\mathbf{y}_K(\cdot\; ;\eta)\big\} = \text{span}\big\{ y_K^{(i)}(\cdot,\eta) : 1\leq i\leq n_K\big\}.
\end{align} Here, $y_K^{(i)}(\cdot\; ;\eta):X\longrightarrow \mathbb{R}$ denotes the (scalar-valued) $i^\text{th}$ component function of
\begin{equation}
\mathbf{y}_K(\cdot \ ;\eta)=\begin{bmatrix} y_K^{(1)}(\cdot \ ;\eta) \\[0.2cm] y_K^{(2)}(\cdot \ ;\eta) \\[0.2cm] \vdots \\[0.2cm] y_K^{(n_K)}(\cdot \ ;\eta)\end{bmatrix}: X\longrightarrow \mathbb{R}^{n_K}.
\end{equation}
\end{subequations}
Then, 
\begin{equation}
\mathbb{U}_\text{NN} = \bigcup_{\eta\in\Lambda} \mathbb{U}_\eta.
\end{equation}

\subsubsection*{Vector-valued FFNN (\boldmath{$n_{K+1} > 1$})}

All of the above is scalable to a vector-valued case by constructing the Cartesian product vector space $\mathbb{U}_\eta^{n_{K+1}}$ consisting of $n_{K+1}$ copies of $\mathbb{U}_\eta$, i.e.,
\begin{equation}
\mathbb{U}_\eta^{n_{K+1}} = \mathbb{U}_\eta \times \mathbb{U}_\eta \times \cdots \times \mathbb{U}_\eta, \qquad (n_{K+1} \text{ times}).
\end{equation} Then, 
\begin{subequations}
\begin{equation}
\mathbb{U}_\text{NN} = \bigcup_{\eta\in\Lambda} \mathbb{U}_\eta^{n_{K+1}},
\end{equation} where each element in $\mathbb{U}_\eta^{n_{K+1}}$ is of the form
\begin{align}
\textbf{u}_\text{NN}(\cdot\; ; \eta,\mathbf{W}) = 
\begin{bmatrix}
        \displaystyle \sum_{r=1}^{n_K} W^{(1,r)} \ y_{K}^{(r)}(\cdot\; ; \eta) \\[0.6cm] 
        \displaystyle \sum_{r=1}^{n_K} W^{(2,r)} \ y_{K}^{(r)}(\cdot\; ; \eta) \\[0.6cm] 
        \vdots \\[0.6cm] 
        \displaystyle \sum_{r=1}^{n_K} W^{(n_{K+1},r)} \ y_{K}^{(r)}(\cdot\; ; \eta)
\end{bmatrix}\in \mathbb{U}^{n_{K+1}}_\eta = \begin{bmatrix}
        \mathbb{U}_\eta \\[0.9cm] 
        \mathbb{U}_\eta \\[0.9cm] 
        \vdots \\[0.9cm] 
        \mathbb{U}_\eta
\end{bmatrix}, 
\end{align}
\end{subequations} with $W^{(j,r)}$ denoting the $(j,r)^{\text{th}}$ entry of $\mathbf{W}$.

By construction,  for any $\eta\in\Lambda$,
\begin{equation}
\text{dim}(\mathbb{U}_\eta)\leq n_K,\qquad \text{dim}(\mathbb{U}_\eta^{n_{K+1}}) = n_{K+1} \text{dim}(\mathbb{U}_\eta)\leq n_{K+1} n_K.
\end{equation} Note that the linear independence among the component functions of  $\mathbf{y}_K(\cdot\; ;\eta)$ is generally not guaranteed.
\end{subappendices}

\chapter{The Deep Finite Element Method} \label{chapter3}

\begin{quote}
\textbf{Summary.} We introduce a dynamic Deep Learning architecture based on the Finite Element Method to solve linear parametric Partial Differential Equations. The connections between neurons in the architecture mimic the Finite Element connectivity graph when applying mesh refinements. We select and discuss several loss functions employing preconditioners and different norms to enhance convergence. For simplicity, we implement the resulting method in one spatial domain (1D). Numerical experiments show, in general, good approximations for symmetric problems that are either positive-definite or indefinite in parametric and non-parametric equations. However, lack of convexity prevents us from obtaining high-accuracy solutions on some occasions.  \emph{Refer to \cite{uriarte2022finite} for the published version}.
\end{quote}

\section{Introduction}\label[section]{section3.1}

Many of the works discussed in the literature review (recall \Cref{section1.4}) address the idea of finding a continuous function (the \ac{NN}) that approximates the solution of a (parametric) \ac{PDE}.  Generally, those designs allow evaluating the \ac{NN} at any point in the domain, i.e., they have a mesh-free structure.  However, they also present some limitations. We highlight the following two:  (a) the resulting deep-learning-based architectures lack explainability \cite{samek2017explainable, barredoarrieta2020explainable}, and (b) numerical integration rules are challenging to design within the loss, as mentioned in \cite{kharazmi2019variational} and further developed in \cite{rivera2022quadrature}.

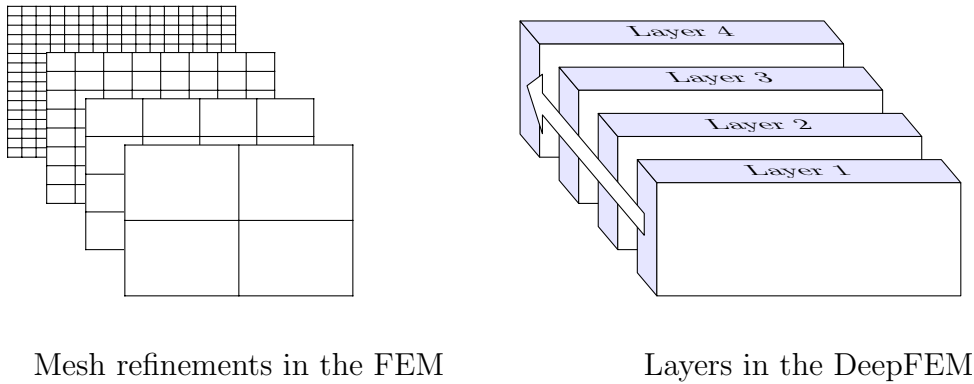
\begin{figure}
    \centering
 \begin{tikzpicture}[
  line join=round,
  z=(130:.8cm), sh/.style={shift={(0,0,#1)}},
  sides/.style={fill=blue!10!white}]
\foreach \i in {4, 3, 2, 1}{
  \tikzset{sh=\i}
  \fill[white] (0,0) rectangle coordinate(@) (3,2);
  \draw[line cap=rect, xstep=3/(2^\i), ystep=2/(2^\i)] (0,0) grid (3,2);
}
\node[below] at (@|-0,0) {Mesh refinements in the FEM};

\tikzset{shift=(right:7)}

\foreach \i in {4, 3, 2, 1}{
  \tikzset{sh=\i}
  \draw[sides, canvas is yz plane at x=0]
    (0,0) rectangle (1.5,.5) coordinate (@);
  \draw[sides, canvas is xz plane at y=1.5]
    (@) rectangle node[transform shape, yscale=.75] {Layer \i}
                  (4,0) coordinate (@);
  \draw[fill=white] (0,0) rectangle coordinate(@) (@);
}
\node[below] at (@|-0,0) {Layers in the DeepFEM};

\node[single arrow, draw, fill=white, minimum height=3cm,
      canvas is yz plane at x=0, rotate=90] at (0, .75, 2.75) {};
\end{tikzpicture}
\caption{Relation of mesh refinements in FEM with layers in the DeepFEM.}
    \label{fig:Deep-FEM_diagram}
\end{figure}

To overcome both limitations, we propose a deep-learning-based method for solving parametric \acp{PDE} that resembles the \ac{FEM}, called the \emph{\acf{DeepFEM}}.  The proposed \ac{NN} architecture aims to act as a solver of the parametric system of linear equations arising in the \ac{FEM} and to mimic the Finite Element connectivity graph when applying mesh refinements---we associate each NN layer with a mesh refinement. Each NN layer has a ResNet \cite{he2016deep} design and extends coarse solutions to finer meshes.  \Cref{fig:Deep-FEM_diagram} illustrates the relation between \ac{FEM} refinements and \ac{NN} layers. In this way, the architecture provides a certain degree of explainability, being the output the vector of nodal values corresponding to the finest mesh. In addition, this discrete approach enables exact numerical integration during training because the \ac{NN} prediction lies on a piecewise-polynomial space.

\ac{DeepFEM} first sets an initial architecture that produces coarse solutions after training. Then, we iteratively and dynamically insert layers into the architecture, maintaining the previously trained learnable parameters and adding new ones. Subsequently, we retrain the learnable variables\footnote{Throughout this chapter, we change the terminology used in \Cref{chapter2} from \emph{parameters} to \emph{variables}.  We do so to distinguish the parameters associated with the \ac{PDE} coefficients from those related to the weights and biases of the \ac{NN}.} of the new model and we repeat this process until we achieve a desired degree of accuracy in the solution.  In this way, the proposed \ac{NN} allows us to perform a mesh-convergence study.

For simplicity,  our implementation restricts to \ac{1D} problems with piecewise-linear approximations and uniform refinements. The extension to higher dimensional problems, higher order polynomial approximations, and/or adaptive meshes is straightforward; but it requires a more elaborated implementation to deal with geometric criteria and node numbering between refinements, which we do not delve into in this work. We show numerical results of model problems with constant and piecewise-constant parameters.

The main contribution of the presented technology lies in the ability of the \ac{NN} to solve parametric problems. For illustration, we first introduce it for the non-parametric case before considering the parametric problem. We do the same when describing the numerical results since the comprehensibility and limitations of the method for the non-parametric setting are extrapolable to the parametric one.

The remainder of this chapter is as follows.  \Cref{section3.2} introduces our problem of interest and the corresponding variational formulation.  \Cref{section3.3} describes our selected \ac{DeepFEM} solver architecture and \Cref{section3.4} defines the loss when employing several preconditioners.  \Cref{section3.5} shows implementation details and the limitations we encountered when using the library \ac{TF} \cite{tensorflow2015-whitepaper, abadi2016tensorflow}. Finally, \Cref{section3.6} discusses numerical results.

\section{Model problem and Finite Element formulation}\label[section]{section3.2}

We focus on a particular parametric \acp{BVP}, although the presented approach applies to other problems that can be solved using the \ac{FEM}.

Let $\Omega$ be a smooth domain and consider the following parametric linear \ac{BVP}:
\begin{align}
\begin{cases}
- \nabla \cdot \sigma\nabla u + \alpha u = f, &\text{in } \Omega,\\
\phantom{- \nabla \cdot \sigma\nabla u + \alpha}u = 0, &\text{in } \Gamma_D,\\
\phantom{\nabla \cdot \sigma\nabla u} -\sigma \ \frac{\partial u}{\partial n} = g, &\text{in } \Gamma_N,
\end{cases}
\label{original problem}
\end{align} where parameters $\sigma>0$ and $\alpha\in\mathbb{R}$ are piecewise-constant functions, $f$ is the source,  and $g$ is the Neumann data.  $\Gamma_D$ and $\Gamma_N$ are the Dirichlet and Neumann boundaries, respectively. $n$ denotes the outer normal vector at each point of $\Gamma_N$ and $\partial u/\partial n = \nabla u \cdot n$.  Note that this model problem encompasses Poisson's ($\alpha=0$) and Helmholtz's ($\alpha <0; \sigma=1$) equations.

A variational formulation of the above \ac{BVP} reads as follows: \begin{equation}
\label{eq:variational_formulation}
\displaystyle{\left\{
\begin{tabular}{l}
Find $u^*\in H^1_0(\Omega)$ such that\\
$(\sigma\nabla u^*, \nabla v)_{\Omega} + (\alpha u^*, v)_{\Omega} = (f, v)_{\Omega} - (g, v)_{\Gamma_N}, \forall v\in H^1_0(\Omega)$,\\
\end{tabular}
    \right.}
\end{equation} with
\begin{equation}
(u,v)_\Omega := \int_\Omega u\cdot v,
\end{equation} and where we have the following underlying spaces:
\begin{subequations}
\begin{align}
    L^2(\Omega) &= \{ u:\Omega\longrightarrow\mathbb{R} :  (u,u)_{\Omega} < \infty\},\\
    H^1(\Omega)&=\{u\in L^2(\Omega): (\nabla u,\nabla u)_\Omega < \infty\},\\
    H^1_0(\Omega)&=\{u\in H^1(\Omega): u\vert_{\Gamma_D} = 0 \}.
\end{align}
\end{subequations} We assume that both $f$ and $g$ are sufficiently regular along the entire process, namely, $f\in L^2(\Omega)$ and $g\in H^{1/2}(\partial\Omega)$---we refer to \cite{girault2012finite} for a detailed discussion about boundary/fractional Sobolev spaces.

\begin{figure}
    \centering
    \begin{tikzpicture} 
        \draw[dotted, thick][-] (0,0) -- (10.5,0) ; 
		\draw (0,-0.35) node{$x_0$};
		\draw[color=red, thick][-]  (1.5,0) -- (0,1.5);
		\draw (0,1.9) node{\color{red}$\psi_{0}(x)$};
		\draw (1.5,-0.35) node{$x_1$};
		\draw (3,-0.35) node{$x_2$};
        \draw (4.5,-0.35) node{$x_3$};
        \draw (6,-0.35) node{$x_4$};
        \draw (7.5,-0.35) node{$\ldots$};
        \draw (9,-0.35) node{$x_{J-1}$};
		\draw (10.5,-0.35) node{$x_{J}$};
		\draw[thick][-]  (0,0) -- (1.5,1.5) node[anchor=north] {};
		\draw[thick][-]  (1.5,1.5) -- (3,0) node[anchor=north] {};
		\draw (1.5,1.9) node{$\psi_{1}(x)$};
		\draw[color=green!40!black, thick][-]  (1.5,0) -- (3,1.5) node[anchor=north] {};
		\draw[color=green!40!black, thick][-]  (3,1.5) -- (4.5,0) node[anchor=north] {};
		\draw (3,1.9) node{\color{green!40!black}$\psi_2(x)$};		
		\draw[color=blue, thick][-]  (3,0) -- (4.5,1.5) node[anchor=north] {};
		\draw[color=blue, thick][-]  (4.5,1.5) -- (6,0) node[anchor=north] {};
		\draw (4.5,1.9) node{\color{blue}$\psi_3(x)$};
		\draw[color=purple, thick][-]  (9,0) -- (10.5,1.5);
		\draw (10.5,1.9) node{\color{purple}$\psi_{J}(x)$}; node[anchor=north] {};
	\foreach \x in {0,1.5,3,4.5,6,9,10.5}{
         \node[circle, draw=black, fill=black,anchor=north, inner sep=1.5pt, minimum size=0.5pt] (a) at (\x, 0.08){};
         }
\end{tikzpicture}
    \caption{Support of the tent-shape piecewise-linear functions $\psi_j$.}
    \label{fig:tent_functions}
\end{figure}
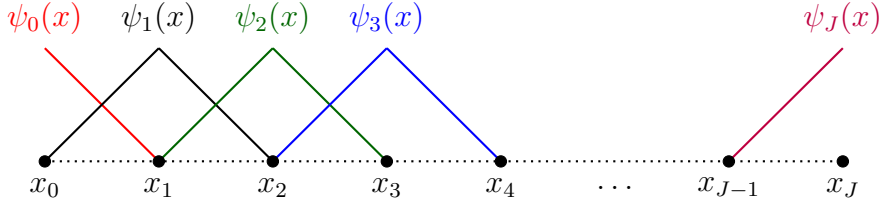

Following a FEM formulation in \ac{1D}, we look for a solution of the form 
\begin{equation}\label{FEM_form_solution}
u_\text{FEM}(x; \sigma,\alpha) := \sum_{j=0}^J u_{\text{FEM}, j}(\sigma,\alpha) \ \psi_j(x), 
\end{equation} where $u_{\text{FEM}, j}(\sigma,\alpha)$ and $\psi_j(x)$ are the unknown coefficients and pre-established tent-shape piecewise-linear basis functions (see \Cref{fig:tent_functions}), respectively.  

Inserting \eqref{FEM_form_solution} in the variational formulation \eqref{eq:variational_formulation} and testing against all the basis functions $\psi_j$, we arrive at the following system of linear equations:
\begin{equation}
\label{eq:system_FEM}
\mathbf{A}\mathbf{u} = \mathbf{f},
\end{equation} where $\mathbf{A}=\mathbf{A}(\sigma,\alpha):=[(\sigma \psi_j',  \psi_i')_\Omega+(\alpha \psi_j, \psi_i)_\Omega]_{i,j=0}^J$,  $\mathbf{u}=\mathbf{u}_{\text{FEM}}(\sigma,\alpha):=[u_{\text{FEM}, j}]_{j=0}^J$ is the vector of unknown coefficients, and $\mathbf{f}:=[(f,\psi_j)_\Omega-(g,\psi_j)_{\Gamma_N}]_{j=0}^J$ is the load vector.

\section{Dynamic architecture}
\label[section]{section3.3}
We first describe our proposed architecture for the non-parametric case and then extend it to the parametric case. Finally, we consider both constant and piecewise-constant parameter alternatives. 
 
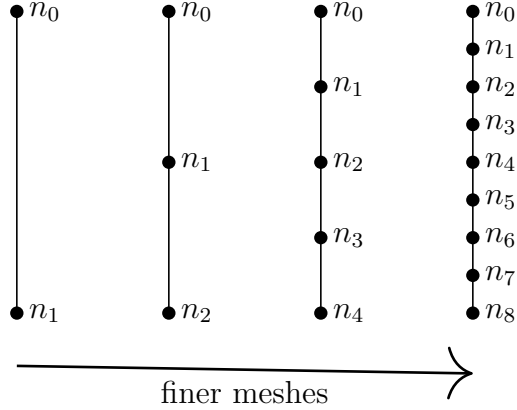
\begin{figure}[htb]
\centering
\begin{tikzpicture}
\usetikzlibrary{arrows.meta}

\foreach \l in {0,1,2,3} {
    \draw[black, line width = 0.6] (2*\l,0) -- (2*\l,-4);
    }
    
\foreach \x in {0,1} {
    	\path[black, fill=black] (2*0, -4*\x) circle[black, radius=2.5pt] node[right=0.5, black] {$n_{\x}$};
    }

\foreach \x in {0,1,2} {
    	\path[black, fill=black] (2*1, -2*\x) circle[black, radius=2.5pt] node[right=0.5, black] {$n_{\x}$};
    }
    
\foreach \x in {0,1,2,3,4} {
    	\path[black, fill=black] (2*2, -1*\x) circle[black, radius=2.5pt] node[right=0.5, black] {$n_{\x}$};
    }
    
\foreach \x in {0,1,2,3,4,...,8} {
    	\path[black, fill=black] (2*3, -0.5*\x) circle[black, radius=2.5pt] node[right=0.5, black] {$n_{\x}$};
    }
    
\draw[black, line width = 1, arrows={->[scale=2.5]}] (0,-4.7) -- (6,-4.8) node[midway, below] {finer meshes};;

\end{tikzpicture}
\caption{Considered node numbering criteria for uniform mesh refinements in 1D.}
\label{fig:coarse_fine_numbering}
\end{figure}

\Cref{fig:coarse_fine_numbering} shows our selected node numbering when performing uniform mesh refinements. Accordingly, the extension operator $\mathbf{E}$ is given by a sparse matrix filled with ones and halves depending on the contribution with which each node propagates from the coarse to the fine mesh. Note that similar extension operators exist for \ac{2D} and \ac{3D} problems as well as for higher-order elements, and for $H(\text{div})$, $H(\text{curl})$, and $L^2$ discretizations \cite{demkowicz2000rham}.

\subsection{Non-parametric scheme for constant PDE coefficients}

Let $\sigma$ and $\alpha$ be real-valued constant coefficients. For a one-element mesh, we propose the following architecture: 
\begin{align}
    \mathbf{u}_{\text{NN}}^{(1)}(\sigma,\alpha) = \mathbf{W}^{(1)}  \ \varphi\left( \mathbf{W}_{\sigma,\alpha}^{(1)} \begin{bmatrix} \sigma \\ \alpha \end{bmatrix} + \mathbf{b}_{\sigma,\alpha}^{(1)}\right) \in\mathbb{R}^2, 
\end{align} where $\varphi$ is an activation function, and the output is a vector of size two that aims to approximate the FEM solution in the considered (single-element) coarse mesh (step $s=1$), i.e.,  
\begin{align}
\mathbf{u}_{\text{FEM}}^{(1)}(\sigma,\alpha) = 
\begin{bmatrix}
 u_{\text{FEM},0}^{(1)}(\sigma,\alpha) \\ u_{\text{FEM},1}^{(1)}(\sigma,\alpha)
\end{bmatrix} \approx \begin{bmatrix}
 u_{\text{NN},0}^{(1)}(\sigma,\alpha) \\ u_{\text{NN},1}^{(1)}(\sigma,\alpha)
\end{bmatrix} = \mathbf{u}_{\text{NN}}^{(1)}(\sigma,\alpha).
\end{align} The set of learnable variables\ consists of
\begin{equation}
\theta^{(1)} = \{\mathbf{W}_{\sigma,\alpha}^{(1)}, \mathbf{b}_{\sigma,\alpha}^{(1)}, \mathbf{W}^{(1)}\}.
\end{equation}

Then,  we refine and obtain a two-element mesh (step $s=2$).  Its counterpart in the architecture consists of adding an input-dependent ResNet \cite{he2016deep} to $\mathbf{u}_{\text{NN}}^{(1)}$ to define $ \mathbf{u}_{\text{NN}}^{(2)}$ as follows:
\begin{subequations}
\begin{align}
    \mathbf{r}_{\text{NN}}^{(2)}(\sigma,\alpha) &= \mathbf{W}^{(2)}  \ \varphi\left( \mathbf{W}_{\sigma,\alpha}^{(2)} \begin{bmatrix} \sigma \\ \alpha \end{bmatrix} + \mathbf{b}_{\sigma,\alpha}^{(2)}\right) \in\mathbb{R}^3, \\
    \mathbf{u}_\text{NN}^{(2)}(\sigma,\alpha) &= \mathbf{E}_{1}^{2}\mathbf{u}_{\text{NN}}^{(1)}(\sigma,\alpha) + \mathbf{r}_\text{NN}^{(2)}(\sigma,\alpha)\in\mathbb{R}^3,
\end{align}
\end{subequations} with $\mathbf{E}_1^2$ denoting the extension matrix of $\mathbf{u}_{\text{NN}}^{(1)}$ on the fine (two-element) mesh.  The set of learnable variables $\theta^{(2)}$ is now 
\begin{equation}
\{\mathbf{W}_{\sigma,\alpha}^{(2)}, \mathbf{b}_{\sigma,\alpha}^{(2)}, \mathbf{W}^{(2)}\}\qquad \text{or}\qquad \theta^{(1)} \cup \{\mathbf{W}_{\sigma,\alpha}^{(2)}, \mathbf{b}_{\sigma,\alpha}^{(2)}, \mathbf{W}^{(2)}\}
\end{equation} depending on whether we perform a \emph{layer-by-layer} or an \emph{end-to-end} training, respectively.  Note that the new inserted learnable variables correspond only to $\mathbf{r}_\text{NN}^{(2)}$.  We initialize the variables of $\mathbf{r}_\text{NN}^{(2)}$ so that it is zero at the beginning of the retraining (e.g., initializing $\mathbf{W}^{(2)},  \mathbf{b}^{(2)} = 0$),  and we maintain the learned values in $\theta^{(1)}$ in the previous step at the beginning of the retraining at (the current) step $s=2$.  We retrain $\mathbf{u}_\text{NN}^{(2)}$ (via layer-by-layer or end-to-end training regimes) so it approximates the FEM solution on the (two-element) fine mesh, i.e.,
\begin{align}
\mathbf{u}_{\text{FEM}}^{(2)}(\sigma,\alpha) = 
\begin{bmatrix}
 u_{\text{FEM}, 0}^{(2)}(\sigma,\alpha) \\ u_{\text{FEM}, 1}^{(2)}(\sigma,\alpha) \\ u_{\text{FEM}, 2}^{(2)}(\sigma,\alpha)
\end{bmatrix} \approx \begin{bmatrix}
 u_{\text{NN}, 0}^{(2)}(\sigma,\alpha) \\ u_{\text{NN}, 1}^{(2)}(\sigma,\alpha) \\ u_{\text{NN}, 2}^{(2)}(\sigma,\alpha)
\end{bmatrix} = \mathbf{u}_{\text{NN}}^{(2)}(\sigma,\alpha).
\end{align}

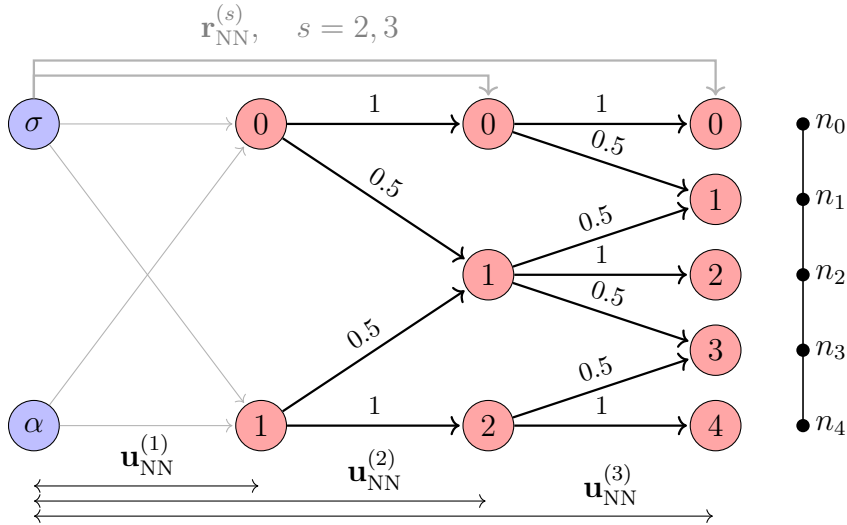
\begin{figure}[htb]
\centering
\begin{tikzpicture}[
    shorten >=1pt,->,
    draw=gray!60,
    node distance=\layersep,
    every pin edge/.style={<-,shorten <=1pt},
    neuron/.style={circle, fill=black, draw=black, minimum size=19pt,inner sep=2pt},
    element/.style={rectangle,fill=black,minimum size=16pt,inner sep=0pt},
    input element/.style={neuron, fill=blue!25},
    output neuron/.style={neuron, fill=red!35},
    hidden neuron/.style={neuron, fill=purple!50},
    annot/.style={text width=4em, text centered},
    punkt/.style={
           rectangle,
           rounded corners,
           draw=black, very thick,
           text width=4em,
           minimum height=2em,
           text centered},
    pil/.style={
           ->,
           thick,
           shorten <=2pt,
           shorten >=2pt,}
]
    \def\layersep{2.3cm}
    
    \draw[black,<->] (0,-4.8) -- (1.3*\layersep,-4.8);
    \path[draw, fill=black] (0.65*\layersep, -4.8) circle[radius=0pt] node [above=0.01, black] {$\mathbf{u}^{(1)}_{\text{NN}}$};
    \foreach \name / \y in {0,...,0} {
        \node[input element] (IS-\name) at (0,-\y) {$\sigma$};
        }
        
    \foreach \name / \y in {0,...,0} {
        \node[input element] (IA-\name) at (0,-\y-4) {$\alpha$};
        }

    \foreach \name / \y in {0,...,1} {
        \node[output neuron] (O1-\name) at (1.3*\layersep,-4*\y) {$\y$};
        }
        
    \foreach \source in {0,...,0}
        \foreach \dest in {0,1}
            \path (IS-\source) edge (O1-\dest);
            
    \foreach \source in {0,...,0}
        \foreach \dest in {0,1}
            \path (IA-\source) edge (O1-\dest);
    
    \draw[black,<->] (0,-5) -- (2.6*\layersep,-5);
    \path[draw, fill=black] (1.95*\layersep, -5) circle[radius=0pt] node [above=0.01, black] {$\mathbf{u}^{(2)}_{\text{NN}}$};
        
    \foreach \name / \y in {0,...,2} {
        \node[output neuron] (O2-\name) at (2.6*\layersep,-2*\y) {$\y$};
        }
    
    \path[black, line width=0.85pt] (O1-0) edge node[sloped, anchor=south, auto=false] {\footnotesize $1$}(O2-0);
    \path[black, line width=0.85pt] (O1-0) edge node[sloped, anchor=south, auto=false] {\footnotesize $0.5$} (O2-1);
    \path[black, line width=0.85pt] (O1-1) edge  node[sloped, anchor=south, auto=false] {\footnotesize $0.5$} (O2-1);
    \path[black, line width=0.85pt] (O1-1) edge  node[sloped, anchor=south, auto=false] {\footnotesize $1$} (O2-2); 
    
    \path[gray!60, to path={-- ++(0,.3) -| (\tikztotarget)}, line width=0.85] (IS-0.north) edge (O2-0.north);
    
    \draw[black,<->] (0,-5.2) -- (3.9*\layersep,-5.2);
    \path[draw, fill=black] (3.30*\layersep, -5.2) circle[radius=0pt] node [above=0.01, black] {$\mathbf{u}^{(3)}_{\text{NN}}$};
    
    \draw[black,-, line width = 0.6] (4.4*\layersep,0) -- (4.4*\layersep,-4);
    \foreach \x in {0,1,2,3,4} {
        \pgfmathparse{ \x/4 };
        \pgfmathresult;
        \let\y\pgfmathresult;
    	\path[black, fill=black] (4.4*\layersep, -\x) circle[black, radius=2.5pt] node[right=0.5, black] {$n_{\x}$};
    }
        
    \foreach \name / \y in {0,...,4} {
        \node[output neuron] (O3-\name) at (3.9*\layersep,-\y) {$\y$};
        }
    
    \path[black,  draw=black, line width=0.85pt] (O2-0) edge  node[sloped, anchor=south, auto=false] {\footnotesize $1$} (O3-0);
    \path[black, line width=0.85pt] (O2-0) edge  node[sloped, anchor=south, auto=false] {\footnotesize $0.5$} (O3-1);
    \path[black, line width=0.85pt] (O2-1) edge  node[sloped, anchor=south, auto=false] {\footnotesize $0.5$} (O3-1);
    \path[black, line width=0.85pt] (O2-1) edge  node[sloped, anchor=south, auto=false] {\footnotesize $1$} (O3-2);
    \path[black, line width=0.85pt] (O2-1) edge  node[sloped, anchor=south, auto=false] {\footnotesize $0.5$} (O3-3);
    \path[black, line width=0.85pt] (O2-2) edge  node[sloped, anchor=south, auto=false] {\footnotesize $0.5$} (O3-3);
    \path[black, line width=0.85pt] (O2-2) edge  node[sloped, anchor=south, auto=false] {\footnotesize $1$} (O3-4);
    
    \path[gray!60, to path={-- ++(0,.5) -| (\tikztotarget)}, line width=0.85] (IS-0.north) edge (O3-0.north);
    \path[draw, fill=gray] (3.5,0.9) circle[radius=0pt] node [above=0.01, gray] {$\mathbf{r}^{(s)}_{\text{NN}}, \quad s=2,3$};
    
\end{tikzpicture}
\caption{Sketch of a non-parametric \acs{DeepFEM} dynamic architecture.  Gray arrows within the architecture indicate learnable variables,  while black arrows refer to (non-trainable) extension operations.}
\label{fig:DL_onesample}
\end{figure}

We repeat this process iteratively (step by step), increasing the depth of our \ac{NN} architecture until the number of elements is sufficient to accurately approximate the exact solution.  \Cref{fig:DL_onesample} sketches this dynamic architecture.

\subsection{Parametric scheme for constant PDE coefficients}
\label{section:Parametric Deep-FEM architecture for constant parameters}

The non-parametric scheme sets weights and biases of very low dimensionality (e.g.,$\mathbf{W}^{(1)}, \mathbf{W}^{(1)}_{\sigma,\alpha}\in\mathbb{R}^{2\times 2}$ and $\mathbf{W}^{(2)}\in\mathbb{R}^{3\times 2}, \mathbf{W}^{(2)}_{\sigma,\alpha}\in\mathbb{R}^{2\times 2}$) because we were only interested in training over a single sample of coefficients.

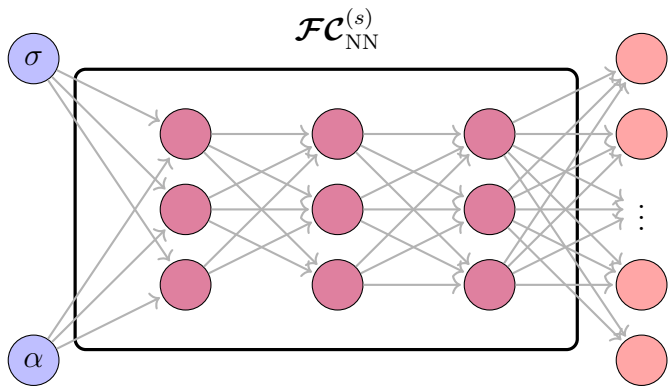
\begin{figure}[htb]
\centering
\begin{tikzpicture}[shorten >=1pt,->,
    draw=gray!60,
    node distance=\layersep,
    every pin edge/.style={<-,shorten <=1pt},
    neuron/.style={circle, draw=black, fill=black,minimum size=19pt,inner sep=2pt},
    element/.style={rectangle,fill=black,minimum size=16pt,inner sep=0pt},
    input neuron/.style={neuron, fill=blue!25},
    output neuron/.style={neuron, fill=red!35},
    hidden neuron/.style={neuron, fill=purple!50},
    annot/.style={text width=4em, text centered},
    punkt/.style={
           rectangle,
           rounded corners,
           draw=black, very thick,
           text width=4em,
           minimum height=2em,
           text centered},
    pil/.style={
           ->,
           thick,
           shorten <=2pt,
           shorten >=2pt,}
]
    \def\layersep{1.5cm}

    \node[input neuron] (I-1) at (0,0) {$\sigma$};
    \node[input neuron] (I-2) at (0,-4) {$\alpha$};
    \node[punkt, minimum width = 16em, minimum height = 9em] (block) at (3.85,-2) {};
    \node[hidden neuron] (H1-1) at (2,-1) {};
    \node[hidden neuron] (H1-2) at (2,-2) {};
    \node[hidden neuron] (H1-3) at (2,-3) {};
    \node[hidden neuron] (H2-1) at (4,-1) {};
    \node[hidden neuron] (H2-2) at (4,-2) {};
    \node[hidden neuron] (H2-3) at (4,-3) {};
    \node[hidden neuron] (H3-1) at (6,-1) {};
    \node[hidden neuron] (H3-2) at (6,-2) {};
    \node[hidden neuron] (H3-3) at (6,-3) {};
    \node[output neuron] (O-1) at (8,0) {};
    \node[output neuron] (O-2) at (8,-1) {};
    \node[] (O-3) at (8,-2) {$\vdots$};
    \node[output neuron] (O-4) at (8,-3) {};
    \node[output neuron] (O-5) at (8,-4) {};
    
    \node[] (FC) at (4,0.4) {$\boldsymbol{\mathcal{FC}}_\text{NN}^{(s)}$};
    
    \path[line width=0.8] (I-1) edge (H1-1);
    \path[line width=0.8] (I-1) edge (H1-2);
    \path[line width=0.8] (I-1) edge (H1-3);
    \path[line width=0.8] (I-2) edge (H1-1);
    \path[line width=0.8] (I-2) edge (H1-2);
    \path[line width=0.8] (I-2) edge (H1-3);
    \path[line width=0.8] (H1-1) edge (H2-1);
    \path[line width=0.8] (H1-1) edge (H2-2);    
    \path[line width=0.8] (H1-1) edge (H2-3);
    \path[line width=0.8] (H1-2) edge (H2-1);
    \path[line width=0.8] (H1-2) edge (H2-2);    
    \path[line width=0.8] (H1-2) edge (H2-3);
    \path[line width=0.8] (H1-3) edge (H2-1);
    \path[line width=0.8] (H1-3) edge (H2-2);    
    \path[line width=0.8] (H2-1) edge (H3-1);
    \path[line width=0.8] (H2-1) edge (H3-2);    
    \path[line width=0.8] (H2-1) edge (H3-3);
    \path[line width=0.8] (H2-2) edge (H3-1);
    \path[line width=0.8] (H2-2) edge (H3-2);    
    \path[line width=0.8] (H2-2) edge (H3-3);
    \path[line width=0.8] (H2-3) edge (H3-1);
    \path[line width=0.8] (H2-3) edge (H3-2);    
    \path[line width=0.8] (H2-3) edge (H3-3);
    
    \path[line width=0.8] (H3-1) edge (O-1);
    \path[line width=0.8] (H3-1) edge (O-2);    
    \path[line width=0.8] (H3-1) edge (O-3);
    \path[line width=0.8] (H3-1) edge (O-4);
    \path[line width=0.8] (H3-1) edge (O-5);
    \path[line width=0.8] (H3-2) edge (O-1);
    \path[line width=0.8] (H3-2) edge (O-2);    
    \path[line width=0.8] (H3-2) edge (O-3);
    \path[line width=0.8] (H3-2) edge (O-4);
    \path[line width=0.8] (H3-2) edge (O-5);
    \path[line width=0.8] (H3-3) edge (O-1);
    \path[line width=0.8] (H3-3) edge (O-2);    
    \path[line width=0.8](H3-3) edge (O-3);
    \path[line width=0.8] (H3-3) edge (O-4);
    
\end{tikzpicture}
\caption{Sketch of the architecture of a trainable block at step $s$.  Gray arrows indicate learnable variables.}
\label{fig:trainable_block_architecture}
\end{figure}

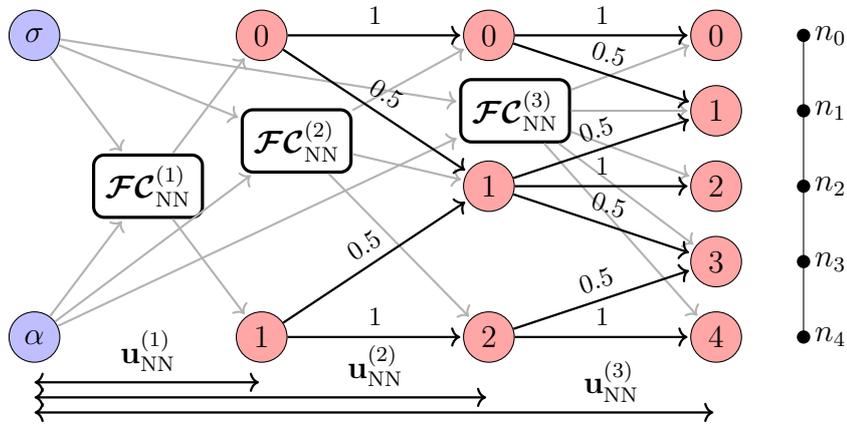
\begin{figure}[htbp]
\centering
\begin{tikzpicture}[
    shorten >=1pt,->,
    draw=gray!60,
    node distance=\layersep,
    every pin edge/.style={<-,shorten <=1pt},
    neuron/.style={circle, draw=black, fill=black!25,minimum size=19pt,inner sep=0pt},
    element/.style={rectangle,fill=black!25,minimum size=16pt,inner sep=0pt},
    input element/.style={neuron, fill=blue!25},
    output neuron/.style={neuron, fill=red!35},
    hidden neuron/.style={neuron, fill=purple!50},
    annot/.style={text width=4em, text centered},
    punkt/.style={
           rectangle,
           rounded corners,
           draw=black, very thick,
           text width=4em,
           minimum height=2em,
           text centered},
    pil/.style={
           ->,
           thick,
           shorten <=2pt,
           shorten >=2pt,}
]
    \def\layersep{2.3cm}
    
    \draw[black,<->, line width = 0.85] (0,-4.6) -- (1.3*\layersep,-4.6);
    \path[draw, fill=black] (0.65*\layersep, -4.6) circle[radius=0pt] node [above=0.01, black] {$\mathbf{u}_\text{NN}^{(1)}$};
    
    \foreach \name / \y in {0,...,0} {
        \node[input element] (IS-\name) at (0,-\y) {$\sigma$};
        }
        
    \foreach \name / \y in {0,...,0} {
        \node[input element] (IA-\name) at (0,-\y-4) {$\alpha$};
        }
    
    \foreach \name / \y in {0,...,4} {
        \node[output neuron] (O3-\name) at (3.9*\layersep,-\y) {$\y$};
        }
        
    \foreach \name / \y in {0,...,2} {
        \node[output neuron] (O2-\name) at (2.6*\layersep,-2*\y) {$\y$};
        }
        
    \foreach \name / \y in {0,...,1} {
        \node[output neuron] (O1-\name) at (1.3*\layersep,-4*\y) {$\y$};
        }

    \foreach \name / \y in {0,...,0} {
        \node[punkt, text width = 2.8em, fill=none] (H1-\name) at (0.65*\layersep,-\y-2) {$\boldsymbol{\mathcal{FC}}_\text{NN}^{(1)}$};
    }
    
    \foreach \name / \y in {0,...,0} {
        \node[punkt, text width = 2.8em, fill=none] (H3-\name) at (2.75*\layersep,-\y-1) {$\boldsymbol{\mathcal{FC}}_\text{NN}^{(3)}$};
    }
    
    \foreach \source in {0,...,0}
        \foreach \dest in {0,...,0}
            \path[line width=0.8] (IS-\source) edge (H3-\dest);
            
    \foreach \source in {0,...,0}
        \foreach \dest in {0,...,0}
            \path[line width=0.8] (IA-\source) edge (H3-\dest);
            
    \foreach \source in {0,...,0}
        \foreach \dest in {0,...,4}
            \path[line width=0.8] (H3-\source) edge (O3-\dest);
        
    \foreach \source in {0,...,0}
        \foreach \dest in {0,...,0}
            \path[line width=0.8] (IS-\source) edge (H1-\dest);
            
    \foreach \source in {0,...,0}
        \foreach \dest in {0,...,0}
            \path[line width=0.8] (IA-\source) edge (H1-\dest);
            
    \foreach \source in {0,...,0}
        \foreach \dest in {0,...,1}
            \path[line width=0.8] (H1-\source) edge (O1-\dest);
    
    \foreach \name / \y in {0,...,0} {
        \node[punkt, text width = 2.8em] (H2-\name) at (1.5*\layersep,-\y-1.4) {$\boldsymbol{\mathcal{FC}}_\text{NN}^{(2)}$};
          
    \foreach \source in {0,...,0}
        \foreach \dest in {0,...,0}
            \path[line width=0.8] (IS-\source) edge (H2-\dest);
            
    \foreach \source in {0,...,0}
        \foreach \dest in {0,...,0}
            \path[line width=0.8] (IA-\source) edge (H2-\dest);
            
    \foreach \source in {0,...,0}
        \foreach \dest in {0,...,2}
            \path[line width=0.8] (H2-\source) edge (O2-\dest);
    }
    
    \draw[black,<->, line width = 0.85] (0,-4.8) -- (2.6*\layersep,-4.8);
    \path[draw, fill=black] (1.95*\layersep, -4.8) circle[radius=0pt] node [above=0.01, black] {$\mathbf{u}_\text{NN}^{(2)}$};
    
    \path[black, line width=0.85pt] (O1-0) edge node[sloped, anchor=south, auto=false] {\footnotesize $1$}(O2-0);
    \path[black, line width=0.85pt] (O1-0) edge node[sloped, anchor=south, auto=false] {\footnotesize $0.5$} (O2-1);
    \path[black, line width=0.85pt] (O1-1) edge  node[sloped, anchor=south, auto=false] {\footnotesize $0.5$} (O2-1);
    \path[black, line width=0.85pt] (O1-1) edge  node[sloped, anchor=south, auto=false] {\footnotesize $1$} (O2-2); 
    
    \draw[black,<->, line width = 0.85] (0,-5) -- (3.9*\layersep,-5);
    \path[draw, fill=black] (3.30*\layersep, -5) circle[radius=0pt] node [above=0.01, black] {$\mathbf{u}_\text{NN}^{(3)}$};
    
    \draw[black,-] (4.4*\layersep,0) -- (4.4*\layersep,-4);
    \foreach \x in {0,1,2,3,4} {
        \pgfmathparse{ \x/4 };
        \pgfmathresult;
        \let\y\pgfmathresult;
    	\path[black, fill=black] (4.4*\layersep, -\x) circle[radius=2.5pt] node [right=0.1, black] {$n_{\x}$};
    }   
    \path[black, line width=0.85pt] (O2-0) edge  node[sloped, anchor=south, auto=false] {\footnotesize $1$} (O3-0);
    \path[black, line width=0.85pt] (O2-0) edge  node[sloped, anchor=south, auto=false] {\footnotesize $0.5$} (O3-1);
    \path[black, line width=0.85pt] (O2-1) edge  node[sloped, anchor=south, auto=false] {\footnotesize $0.5$} (O3-1);
    \path[black, line width=0.85pt] (O2-1) edge  node[sloped, anchor=south, auto=false] {\footnotesize $1$} (O3-2);
    \path[black, line width=0.85pt] (O2-1) edge  node[sloped, anchor=south, auto=false] {\footnotesize $0.5$} (O3-3);
    \path[black, line width=0.85pt] (O2-2) edge  node[sloped, anchor=south, auto=false] {\footnotesize $0.5$} (O3-3);
    \path[black, line width=0.85pt] (O2-2) edge  node[sloped, anchor=south, auto=false] {\footnotesize $1$} (O3-4);
\end{tikzpicture}
\caption{Sketch of a parametric \acs{DeepFEM} dynamic architecture.  Gray arrows within the architecture indicate learnable variables,  while black arrows refer to (non-trainable) extension operations. }
\label{fig:DL_parametric}
\end{figure}

For the parametric problem (i.e., to train over data of coefficient samples), we add width and depth to the trainable pieces of the dynamic architecture as follows:
\begin{subequations}
\begin{alignat}{3}
    &\mathbf{u}_{\text{NN}}^{(1)} =  \boldsymbol{\mathcal{FC}}^{(1)}_\text{NN}(\sigma, \alpha),\\
    &\mathbf{r}_{\text{NN}}^{(s)} = \boldsymbol{\mathcal{FC}}^{(s)}_\text{NN}(\sigma, \alpha), \qquad &s\geq 2,\\
    &\mathbf{u}_{\text{NN}}^{(s)} = \mathbf{E}_{s-1}^{(s)}\mathbf{u}_{\text{NN}}^{(s-1)} + \mathbf{r}_{\text{NN}}^{(s)},\qquad &s\geq 2.
\end{alignat}
\end{subequations} Here,  $\boldsymbol{\mathcal{FC}}_{\text{NN}}^{(s)}$ denotes a fully-connected \ac{FFNN} with non-activated last layer (recall \Cref{section2.1}) and whose output dimension is specified to match the dimensionality of the underlying FEM discretization.  For simplicity,  we call \emph{trainable blocks} to the $\boldsymbol{\mathcal{FC}}_{\text{NN}}^{(s)}$ models.  \Cref{fig:trainable_block_architecture} shows the architecture of the trainable blocks and \Cref{fig:DL_parametric} shows them inside the parametric DeepFEM.

\subsection{Parametric scheme for piecewise-constant PDE coefficients}
\begin{figure}[htbp]
\centering
\begin{tikzpicture}[
    shorten >=1pt,->,
    draw=gray!60,
    node distance=\layersep,
    every pin edge/.style={<-,shorten <=1pt},
    neuron/.style={circle, draw=black, fill=black!25,minimum size=19pt,inner sep=0pt},
    element/.style={rectangle,fill=black!25,minimum size=16pt,inner sep=0pt},
    input element/.style={neuron, fill=blue!25},
    output neuron/.style={neuron, fill=red!35},
    output neuron2/.style={neuron2, fill=red!35},
    hidden neuron/.style={neuron, fill=purple!50},
    annot/.style={text width=4em, text centered},
    punkt/.style={
           rectangle,
           rounded corners,
           draw=black, very thick,
           text width=4em,
           minimum height=2em,
           text centered},
    pil/.style={
           ->,
           thick,
           shorten <=2pt,
           shorten >=2pt,}
]
    \def\layersep{2.3cm}
    
    \draw[black,<->, line width = 0.85] (0,-4.6) -- (1.3*\layersep,-4.6);
    \path[draw, fill=black] (0.65*\layersep, -4.6) circle[radius=0pt] node [above=0.01, black] {$\mathbf{u}_{\text{NN}}^{(1)}$};
    \foreach \name / \y in {0,...,1} {
        \pgfmathparse{ int(\y+1) };
        \pgfmathresult;
        \let\z\pgfmathresult;
        \node[input element] (IS-\name) at (0,-\y) {$\sigma_{\z}$};
        }
        
    \foreach \name / \y in {0,...,1} {
        \pgfmathparse{ int(\y+1) };
        \pgfmathresult;
        \let\z\pgfmathresult;
        \node[input element] (IA-\name) at (0,-\y-3) {$\alpha_{\z}$};
        }

    \foreach \name / \y in {0,...,0} {
        \node[punkt, text width = 2.8em] (H1-\name) at (0.65*\layersep,-\y-1.5) {$\boldsymbol{\mathcal{FC}}_\text{NN}^{(1)}$};
    }

    \foreach \name / \y in {0,...,2} {
        \node[output neuron] (O1-\name) at (1.3*\layersep,-2*\y) {$\y$};
        }
        
    \foreach \source in {0,...,1}
        \foreach \dest in {0,...,0}
            \path[line width=0.8pt] (IS-\source) edge (H1-\dest);
            
    \foreach \source in {0,...,1}
        \foreach \dest in {0,...,0}
            \path[line width=0.8pt] (IA-\source) edge (H1-\dest);
            
    \foreach \source in {0,...,0}
        \foreach \dest in {0,...,2}
            \path[line width=0.8pt] (H1-\source) edge (O1-\dest);
    
    \foreach \name / \y in {0,...,0} {
        \node[punkt, text width = 2.8em] (H2-\name) at (1.5*\layersep,-\y-0.9) {$\boldsymbol{\mathcal{FC}}_\text{NN}^{(2)}$};
      
    \foreach \source in {0,...,1}
        \foreach \dest in {0,...,0}
            \path[line width=0.8pt] (IS-\source) edge (H2-\dest);
            
    \foreach \source in {0,...,1}
        \foreach \dest in {0,...,0}
            \path[line width=0.8pt] (IA-\source) edge (H2-\dest);
            
    \foreach \name / \y in {0,...,4} {
        \node[output neuron] (O2-\name) at (2.6*\layersep,-\y) {$\y$};
        }
            
    \foreach \source in {0,...,0}
        \foreach \dest in {0,...,4}
            \path[line width=0.8pt] (H2-\source) edge (O2-\dest);
    }
    
    \draw[black,<->, line width = 0.85] (0,-4.8) -- (2.6*\layersep,-4.8);
    \path[draw, fill=black] (1.95*\layersep, -4.8) circle[radius=0pt] node [above=0.01, black] {$\mathbf{u}_{\text{NN}}^{(2)}$};
            
    \path[black!70, line width=0.85pt] (O1-0) edge  node[sloped, anchor=south, auto=false] {\footnotesize $1$} (O2-0);
    \path[black!70, line width=0.85pt] (O1-0) edge  node[sloped, anchor=south, auto=false] {\footnotesize $0.5$} (O2-1);
    \path[black!70, line width=0.85pt] (O1-1) edge  node[sloped, anchor=south, auto=false] {\footnotesize $0.5$} (O2-1);
    \path[black!70, line width=0.85pt] (O1-1) edge  node[sloped, anchor=south, auto=false] {\footnotesize $1$} (O2-2);
    \path[black!70, line width=0.85pt] (O1-1) edge  node[sloped, anchor=south, auto=false] {\footnotesize $0.5$} (O2-3);
    \path[black!70, line width=0.85pt] (O1-2) edge  node[sloped, anchor=south, auto=false] {\footnotesize $0.5$} (O2-3);
    \path[black!70, line width=0.85pt] (O1-2) edge  node[sloped, anchor=south, auto=false] {\footnotesize $1$} (O2-4);
    
    \draw[black,<->, line width = 0.85] (0,-5) -- (3.9*\layersep,-5);
    \path[draw, fill=black] (3.30*\layersep, -5) circle[radius=0pt] node [above=0.01, black] {$\mathbf{u}_{\text{NN}}^{(3)}$};
    
    \draw[black,-, line width = 0.6] (4.4*\layersep,0) -- (4.4*\layersep,-4);
    \foreach \x in {0,1,2,3,4,...,8} {
        \pgfmathparse{ \x/8 };
        \pgfmathresult;
        \let\y\pgfmathresult;
    	\path[black, fill=black] (4.4*\layersep, -0.5*\x) circle[radius=2.5pt] node [right=0.1, black] {$n_{\x}$};
    }
        
    \foreach \name / \y in {0,...,8} {
        \node[output neuron] (O3-\name) at (3.9*\layersep,-0.5*\y) {$\y$};
        }

    \path[black!70, line width=0.85pt] (O2-0) edge  node[sloped, anchor=south, auto=false] {} (O3-0);
    \path[black!70, line width=0.85pt] (O2-0) edge  node[sloped, anchor=south, auto=false] {} (O3-1);
    \path[black!70, line width=0.85pt] (O2-1) edge  node[sloped, anchor=south, auto=false] {} (O3-1);
    \path[black!70, line width=0.85pt] (O2-1) edge  node[sloped, anchor=south, auto=false] {} (O3-2);
    \path[black!70, line width=0.85pt] (O2-1) edge  node[sloped, anchor=south, auto=false] {} (O3-3);
    \path[black!70, line width=0.85pt] (O2-2) edge  node[sloped, anchor=south, auto=false] {} (O3-3);
    \path[black!70, line width=0.85pt] (O2-2) edge  node[sloped, anchor=south, auto=false] {} (O3-4);
    \path[black!70, line width=0.85pt] (O2-2) edge  node[sloped, anchor=south, auto=false] {} (O3-5);
    \path[black!70, line width=0.85pt] (O2-3) edge  node[sloped, anchor=south, auto=false] {} (O3-5);
    \path[black!70, line width=0.85pt] (O2-3) edge  node[sloped, anchor=south, auto=false] {} (O3-6);
    \path[black!70, line width=0.85pt] (O2-3) edge  node[sloped, anchor=south, auto=false] {} (O3-7);
    \path[black!70, line width=0.85pt] (O2-4) edge  node[sloped, anchor=south, auto=false] {} (O3-7);
    \path[black!70, line width=0.85pt] (O2-4) edge  node[sloped, anchor=south, auto=false] {} (O3-8);
    
    \path[gray!60, to path={-- ++(0,.5) -| (\tikztotarget)}, line width=1.25] (IS-0.north) edge (O3-0.north);
    
    \node[punkt, text width = 2.8em, fill=white] (H3-0) at (3.2*\layersep,0.9) {$\boldsymbol{\mathcal{FC}}_\text{NN}^{(3)}$};
    
\end{tikzpicture}
\caption{Sketch of a parametric \acs{DeepFEM} dynamic architecture for piecewise-constant PDE coefficients.  Gray arrows within the architecture indicate learnable variables,  while black arrows refer to (non-trainable) extension operations.}
\label{fig:DL_parametric_piecewiseconstant}
\end{figure}
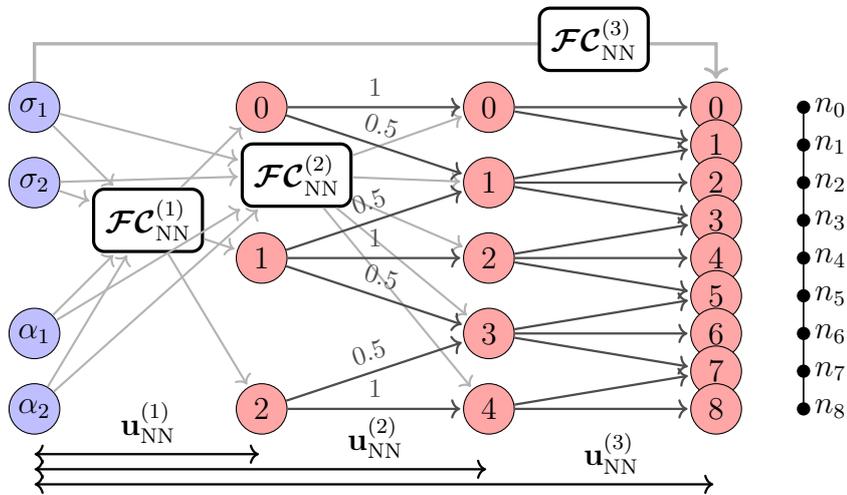

Finally, we add the piecewise-constant behavior to the parameters that we assume taking constant values on each element of the initial mesh.  Hence, the architecture is modified only in the input layer by replacing the $\sigma$ and $\alpha$ values with vectors of sizes equal to the number of elements of the initial mesh. \Cref{fig:DL_parametric_piecewiseconstant} illustrates this.

%
%

\section{Loss function and training}
\label[section]{section3.4}

Within the step-by-step methodology described in \Cref{section3.3}, we drop the superscript ``$(s)$'' for simplicity, although the learnable variables,  loss functions, and considered norms are step-dependent.

To make the \ac{NN} approximate the solution to the parametric system of linear equations arising in the FEM,  i.e., 
\begin{equation}
 \mathbf{u}_\text{NN}\approx\mathbf{u}_\text{FEM},
\end{equation}
we consider the loss function given by
\begin{equation} \label{loss_residual}
      \mathcal{L}(\theta,  \{\sigma_i,\alpha_i\}_{i=1}^N) := \frac{1}{N} \sum_{i=1}^N \Vert \mathbf{A}(\sigma_i,\alpha_i) \mathbf{u}_{\text{NN}}(\sigma_i,\alpha_i;\theta) - \mathbf{f}\Vert,
\end{equation} where  $\mathbf{u}_{\text{NN}}(\sigma_i,\alpha_i;\theta)$ is the prediction of the \ac{DeepFEM} model for the $i$$^\text{th}$ sample, $\mathbf{A} (\sigma_i,\alpha_i)$ is the corresponding FEM matrix---recall \eqref{eq:system_FEM}, $\mathbf{A}(\sigma_i,\alpha_i) \mathbf{u}_{\text{NN}}(\sigma_i,\alpha_i;\theta) - \mathbf{f}$ is the residual vector, and $\lVert\cdot\rVert$ is a pre-established vector norm (e.g., the $2$-norm).

\subsection{Gradient-based training} 
To minimize the loss function,  at the beginning of each step, we apply the Adam optimizer \cite{kingma2014adam}.  When the number of iterations attains an established maximum, we stop its execution.  Alternatively, to control stagnation, we monitor and compare the value of the loss function obtained at any given iteration with the lowest loss obtained few iterations before.  When the loss improvement is insignificant after a prescribed number of iterations, we stop the execution.  Moreover,  since the loss may increase or decrease during the Adam performance, we monitor the loss so as to return the best configuration of trainable variables at the end of its execution.  Thereafter, we apply a customized gradient-descent-based optimizer with a loss-dependent adaptive learning rate (called \emph{Adalr}). We maintain the same stopping and stagnation criteria as in Adam. The Adalr optimizer readjusts the trainable variables to prevent a loss increase.  \Cref{GD_loss-lr} shows the learning rate adaptivity and trainable variables acceptance-rejection criteria.  In the numerical results, we will compare in more detail both optimizers. As a general overview: Adam rapidly decreases the loss initially but tends to plateau at a suboptimal level, causing oscillations without significant improvement.  On the other hand, Adalr monotonically reduces the loss,  usually with a slower convergence speed compared to Adam.

\begin{algorithm}
\SetKwInput{KwInput}{Input}  
\SetKwInput{KwOutput}{Output}
\SetKw{KwBy}{by}
\tcc{Superscripts denote training iterations.}
\KwInput{$\theta^{(0)}, \eta^{(0)}, D$ \tcp*{Initial variables, learning rate,  data}} 
\KwOutput{$\theta^*$, $\mathcal{L}^*$ \tcp*{Final variables and loss values}}
 $\mathcal{L}^{(0)} = \mathcal{L}(\theta^{(0)}; D); \; \theta^{(1)} =  \theta^{(0)}-\lambda^{(0)}\frac{\partial\mathcal{L}}{\partial\theta}(\theta^{(0)}; D)$\;
 $\theta^{*} = \theta^{(0)} ; \; \mathcal{L}^{*} = \mathcal{L}^{(0)}; \; t = 1$\;
 \While{\upshape not STOP}{
   $\mathcal{L}^{(t)} = \mathcal{L}(\theta^{(t)}; D); \; \theta^{(t+1)} =  \theta^{(t)}-\lambda^{(t)}\frac{\partial\mathcal{L}}{\partial\theta}(\theta^{(t)};D)$ \;
     \tcc{If the new loss is worse}
  	\eIf{$\mathcal{L}^{(t)} > \mathcal{L}^*$}{
   	        $\lambda^{(t+1)} = \texttt{Decrease}(\lambda^{(t}); \; \theta^{(t+1)} = \theta^{*}$\;
   	}{
        $\theta^{*} = \theta^{(t)}$\;
    \tcc{If convergence is slow}
  	  \eIf{\upshape convergence is slow \textbf{\upshape and} \upshape no variable rejection in $t-1$}{
   	      $\lambda^{(t+1)} = \texttt{Increase}(\lambda^{(t)})$\; 
   	  }{$\lambda^{(t+1)} = \lambda^{(t)}$\;}     
   	  $\mathcal{L}^{*} = \mathcal{L}^{(t)}$\;    
   }
 $t = t+1$\;
 }
 \Return{$\theta^{*}$, $\mathcal{L}^{*}$}

 \caption{Adalr optimizer}
 \label[algorithm]{GD_loss-lr}
\end{algorithm}

\subsection{Norm selection and preconditioning}

We want to select a discrete norm in the loss function---recall \eqref{loss_residual}---so that its behavior is similar to the error in the \emph{energy norm}.  In a symmetric and positive-definite problem, the energy norm is given by
\begin{equation}
    \lVert\mathbf{v}\rVert_{\mathbf{A}} := \sqrt{\mathbf{v}^T \mathbf{A} \mathbf{v}},
\end{equation} where $\mathbf{v}^T$ stands for the transpose of the column vector $\mathbf{v}$, and $\mathbf{A}$ denotes the symmetric and positive-definite matrix of the system of linear equations. For each $(\sigma,\alpha)$ sample, by writing the error vector as 
\begin{equation}
\mathbf{e}(\sigma, \alpha) := \mathbf{u}_\text{NN}(\sigma, \alpha) -  \mathbf{A}^{-1}(\sigma, \alpha)\mathbf{f} = \mathbf{u}_\text{NN}(\sigma, \alpha) - \mathbf{u}_\text{FEM}(\sigma, \alpha),
\end{equation} we arrive at the following norm relation with the residual:
\begin{equation}
    \lVert\mathbf{e}(\sigma, \alpha)\rVert_{\mathbf{A}(\sigma, \alpha)} = \lVert \mathbf{A}(\sigma, \alpha)\mathbf{u}_\text{NN}(\sigma, \alpha) -  \mathbf{f} \rVert_{\mathbf{A}^{-1}(\sigma, \alpha)}.
\end{equation} As in iterative methods for solving a system of linear equations, we consider preconditioners to mimic the inverse operator and thus decrease the condition number of the system \cite{pardo2004integration}. Thus, we define the discrete norm $\Vert\cdot\Vert$ in \eqref{loss_residual} as
\begin{equation} \label{loss_preconditioned_residual}
    \Vert\mathbf{v}\Vert := \lVert \mathbf{v} \rVert_{\mathbf{P}} = \sqrt{\mathbf{v}^T \mathbf{P} \mathbf{v}},
\end{equation} where $\mathbf{P}=\mathbf{P}(\sigma, \alpha)$ is a preconditioner for $\mathbf{A}(\sigma, \alpha)$.  In particular, we consider block-Jacobi preconditioners of different block sizes and with one-element overlap \cite{pardo2004integration}.  Note that when $\mathbf{P}=\mathbf{I}$ is the identity matrix, the corresponding norm is the discrete $2$-norm, i.e., $\Vert\cdot\Vert_\mathbf{I}=\Vert\cdot\Vert_2$.

In case the problem is indefinite, we consider a positive definite operator and corresponding matrix $\mathbf{B}$. The relation between the error and the residual in this new norm is given by
\begin{align}
    \lVert\mathbf{e}(\sigma, \alpha)\rVert_{\mathbf{B}} &=  \lVert\mathbf{A}^{-1}(\sigma, \alpha)\left\{\mathbf{A}(\sigma, \alpha)\mathbf{u}_\text{NN}(\sigma, \alpha)-\mathbf{f}\right\}\rVert_{\mathbf{B}}.
\end{align} Following the above reasoning leads us to consider, in this occasion,  the norm $\Vert\cdot\Vert$ in \eqref{loss_residual} as
\begin{equation} \label{loss_preconditioned_residual_2}
   \Vert\mathbf{v}\Vert :=  \lVert \mathbf{P}\mathbf{v} \rVert_{\mathbf{B}} = \sqrt{\mathbf{v}^T\mathbf{P}^T\mathbf{B}\mathbf{P}\mathbf{v}} = \lVert\mathbf{v} \rVert_{\mathbf{P}^T\mathbf{B}\mathbf{P}}.
\end{equation} Again,  when $\mathbf{P}=\mathbf{B}=\mathbf{I}$ are the identity operators, the above equation reduces to the discrete $2$-norm of the residual.

In the numerical results, in addition to using the energy norm for positive definite problems, we employ the $L^2$ and $H^1$ norms for testing and monitoring the (preconditioned) residual and the error. These continuum-level norms have their discretized counterparts as follows:
\begin{align}
&\lVert u \rVert_{L^2} = \sqrt{\mathbf{u}^T \mathbf{M} \mathbf{u}} = \lVert \mathbf{u}\rVert_{\mathbf{M}}, \label{eq:continuous vs discrete norms L2}\\ 
&\lVert u \rVert_{H^1} = \sqrt{\mathbf{u}^T \mathbf{M} \mathbf{u} + \mathbf{u}^T \mathbf{K} \mathbf{u}} = \lVert \mathbf{u}\rVert_{\mathbf{K}+\mathbf{M}}, \label{eq:continuous vs discrete norms H1}
\end{align} where $\mathbf{u}=[u_{j}]_{j=0}^J$ is the vector of evaluations at the nodal points of the mesh, $\psi_{j}$ is the basis function related to the $j$$^\text{th}$ node, and $\mathbf{M}$ and $\mathbf{K}$ are the mass and stiffness matrices defined by $(\psi_{s}, \psi_{r})_{\Omega}$ and $(\psi'_{s}, \psi'_{r})_{\Omega}$ in the $(r,s)$$^\text{th}$ entry, respectively. For convenience, we maintain the naming of the discrete norms \eqref{eq:continuous vs discrete norms L2} and \eqref{eq:continuous vs discrete norms H1} by the names inherited at their continuum level, namely, $L^2$ and $H^1$ norms, respectively.

Similarly, we drop the boldface notation when referring to the scalar-valued functions of the corresponding vector functions. For the NN predictions and FEM solutions,
\begin{subequations}
\begin{align}
u_\text{NN}(x;\sigma,\alpha) &= \sum_{j=0}^J u_{\text{NN}, j}(\sigma,\alpha) \ \psi_j(x),\\
u_\text{FEM}(x;\sigma,\alpha) &= \sum_{j=0}^J u_{\text{FEM}, j}(\sigma,\alpha) \ \psi_j(x),
\end{align}
\end{subequations} where $u_{\text{NN}, j}(\sigma,\alpha)$ and $u_{\text{FEM}, j}(\sigma,\alpha)$ are the $j^\text{th}$ components of the DeepFEM prediction and FEM solution vectors, $\mathbf{u}_{\text{NN}}(\sigma,\alpha)$ and $\mathbf{u}_{\text{FEM}}(\sigma,\alpha)$, respectively.

\section{Implementation}
\label[section]{section3.5}

We implement the described framework in the Python programming language, and we used the library \ac{TF2} \cite{tensorflow2015-whitepaper, abadi2016tensorflow}, NumPy \cite{harris2020array}, and SciPy \cite{virtanen2020scipy} to build the \ac{NN} models and to generate and manage the \ac{FEM} data. We create the layers by redefining the corresponding base classes in Keras inside TF2 (\texttt{tf.keras})\footnote{During the TF1.X era,  Keras was an external library dedicated to facilitating the use of TF in Python.  Since the TF2 release in September 2019,  Keras became a sublibrary of TF2 until version TF2.12 in March 2023, where Keras became (again) an external framework built on top of TF that additionally supports JAX and PyTorch.  See \url{https://keras.io} for further details (last accessed: September 4, 2023).}. We use a dedicated sparse library within TF2 (\texttt{tf.sparse}) to handle the extension operations and the FEM matrices to save memory and achieve high performance. We use double precision (float64) instead of the default single precision (float32). We calculate the loss function in an auxiliary non-trainable layer applied after the main model $\mathbf{u}_\text{NN}$. Below, we describe the three main difficulties encountered during our implementation.

\subsection{Reshape of the batch sparse tensors}
Tensors flowing through  Keras are prepared to maintain their first dimension (a.k.a axis) for the batch of samples. Consequently, matrix-vector multiplications in the loss function transform into an operation between a 3D sparse tensor (batch of sparse matrices) and a 2D tensor (batch of vectors) in the Keras model. We could define our operation via Einstein summation if both tensors were dense. However, there is still not an equivalent TF2 function supporting sparse tensors\footnote{We opened an issue concerning this regard in September 2020 in Github (\url{https://github.com/tensorflow/tensorflow/issues/43497}); however, we did not receive any robust support response yet (last accessed: September 4, 2023).}. We overcome this by reshaping the 3D sparse tensor as a sparse non-overlapping block 2D matrix, and the 2D tensor as a long 1D vector. This leads us to miss the batch flow behavior of the model and be forced to avoid using the Keras training function. We thus perform a low-level optimization via AD in \emph{eager execution}, which is significantly slower than training via \emph{graph execution} within Keras\footnote{See \url{https://www.tensorflow.org/guide/intro_to_graphs} for further details (last accessed: September 4, 2023).}.

\subsection{FEM data generation}
We utilize SciPy as the FEM environment to build the extension operators and the matrices and preconditioners of the system of linear equations. We also use it to solve the sparse systems via its built-in solver to compare model predictions. For more complicated geometries, there exist other more efficient software platforms to assemble the matrices, preconditioners, and extension operators (e.g., FEniCS \cite{alnaes2015fenics}). We recommend their use for more complex FE systems (e.g., in 2D and 3D problems). In any case, these offline operations are calculated prior to the training of the NNs to not affect the optimization time of the model.

\subsection{Preconditioners assembly and action}
Ideally, we should manage preconditioners as LU block decompositions of the matrices of the system, and calculate their actions at each loss evaluation using a forward-backward substitution algorithm \cite{abur1988parallel,moshfegh2017direct}.  Again, TF2 lacks an equivalent built-in function to evaluate these actions with sparse tensors. For this reason, we manage the preconditioners as already assembled tensors.

\section{Numerical results}
\label[section]{section3.6}

We conduct a series of experiments to analyze the \ac{DeepFEM} performance.  \Cref{section3.6.1}, \Cref{section3.6.2}, and \Cref{section3.6.3} contemplate non-parametric problems, while \Cref{section3.6.4} addresses parametric problems. All experiments are carried out in the spatial domain $(0,1)$ with Dirichlet and Neumann boundary conditions at $0$ and $1$, respectively.

\subsection{A single PDE example}
\label[section]{section3.6.1}

Let
\begin{equation} \label{eq:poisson_x5_10refinements}
\begin{cases}
-u'' = -20x^{3}, \\
u(0)=0, u'(1)=5.
\end{cases}
\end{equation} Its exact solution is $u^*(x) = x^5$. We solve it employing the \ac{DeepFEM} starting from a one-element mesh ($s=1$). We perform ten mesh refinements (eleven steps, i.e.,  $1\leq s\leq 11$) to finish with a $1024$-element overkill mesh. We analyze our proposed method in this first experiment when performing multiple refinements. All the training blocks of the model consist of one-neuron width and one-layer depth (recall \Cref{section:Parametric Deep-FEM architecture for constant parameters}) employing ReLU activation functions on all trainable blocks.

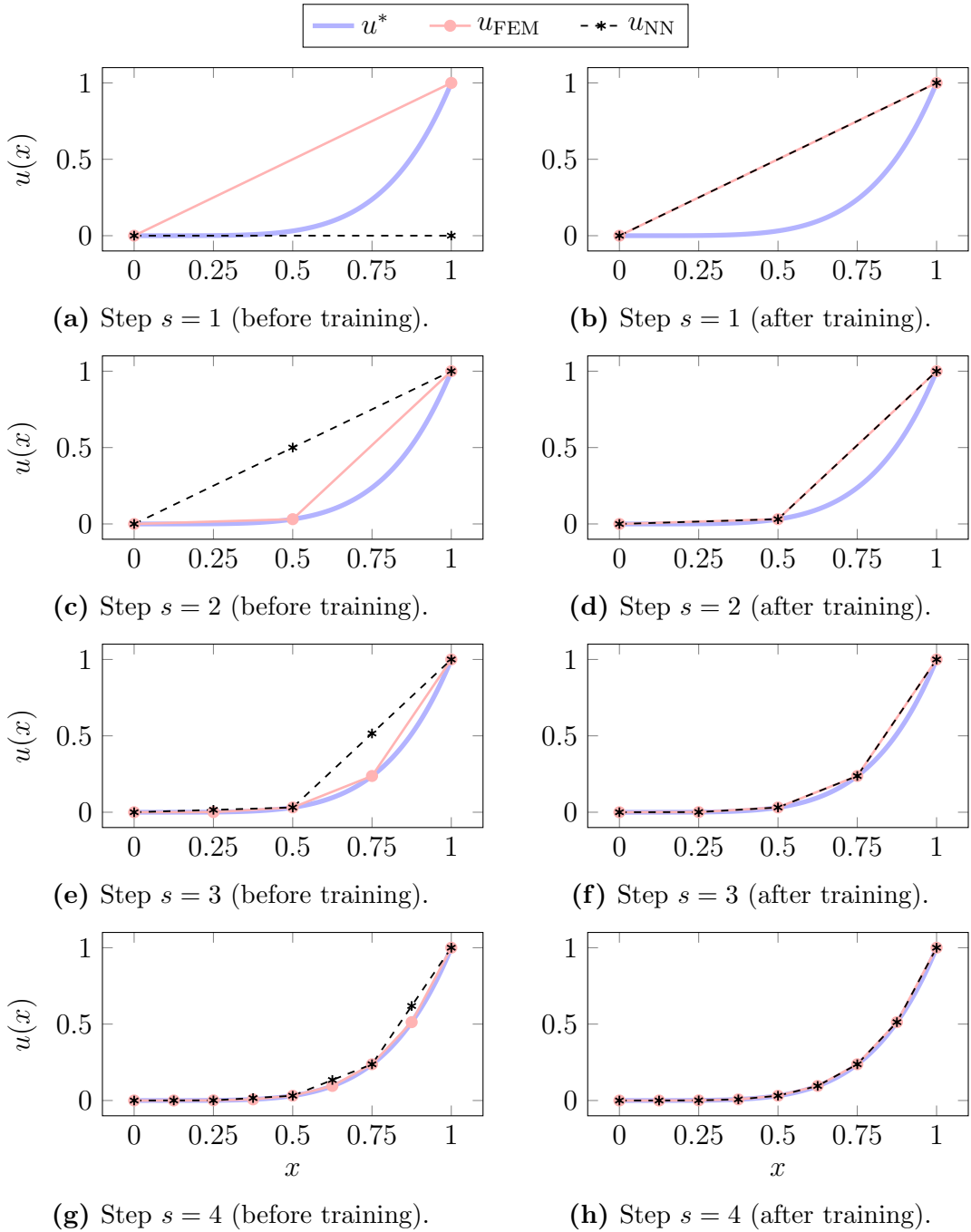
\begin{figure}[htbp]
\centering
\begin{subfigure}[b]{\textwidth}
\centering
 \begin{tikzpicture}
\pgfplotsset{A/.append style={
		hide axis,
	    xmin=-0.1,
	    xmax=1.1,
	    ymin=-0.1,
	    ymax=1.1,
	 	legend style={fill=none, anchor=north, /tikz/every even column/.append style={column sep=0.5cm}},
	 	legend columns = -1
            }
}
\begin{axis}[A]
\addlegendimage{color=blue!30!white, line width=2.};
\addlegendentry{$u^*$};
\addlegendimage{color=red!30!white, line width=1., mark=*, mark options={scale=1, solid}};
\addlegendentry{$u_{\text{FEM}}$};
\addlegendimage{color=black, style=dashed, line width=.7 ,mark=asterisk, mark options = {scale=1, solid}};
\addlegendentry{$u_{\text{NN}}$};
\end{axis}
\end{tikzpicture}
\end{subfigure}
\vskip 0.5em
\begin{subfigure}[b]{0.49\textwidth}
\centering
\begin{tikzpicture}
\newcommand{\fileDL}{Chapters/3.Chapter/figures/numerical_results/DATA/results_poisson_x5_Sinv_3refinements/start/test_prediction_2nodes.csv}
\newcommand{\fileFEM}{Chapters/3.Chapter/figures/numerical_results/DATA/results_poisson_x5_Sinv_3refinements/fem/test_FEM_2nodes.csv}
\pgfplotsset{A/.append style={
	    xmin=-0.1,
	    xmax=1.1,
	    ymin=-0.1,
	    ymax=1.1,
	    ylabel = {$u(x)$},
	    height=0.6*\textwidth,
	    width=\textwidth,
		xtick={0,0.25,...,1},
	 	legend style={fill=none, at={(1.18,1.3)},anchor=north, /tikz/every even column/.append style={column sep=0.5cm}},
	 	legend columns = -1
            }
}
\begin{axis}[A]
\addplot[domain = 0:1, samples = 200, smooth, line width=2., blue!30] {x^5};
	\addplot[color=red!30!white, line width=1., mark=*, mark options={scale=1, solid}] table[x expr=\thisrow{x},y=u_FEM0]{\fileFEM};
	\addplot[color=black, style=dashed, line width=.7 ,mark=asterisk, mark options = {scale=1, solid}] table[x expr=\thisrow{x},y=u_pred0]{\fileDL};
\end{axis}
\end{tikzpicture}
\caption{Step $s=1$ (before training).}
\label{fig:step1_start_poisson_x5}
\end{subfigure}\hfill
\begin{subfigure}[b]{0.49\textwidth}
\centering
\begin{tikzpicture}
\newcommand{\fileDL}{Chapters/3.Chapter/figures/numerical_results/DATA/results_poisson_x5_Sinv_3refinements/end/test_prediction_2nodes.csv}
\newcommand{\fileFEM}{Chapters/3.Chapter/figures/numerical_results/DATA/results_poisson_x5_Sinv_3refinements/fem/test_FEM_2nodes.csv}
\pgfplotsset{A/.append style={
	    xmin=-0.1,
	    xmax=1.1,
	    ymin=-0.1,
	    ymax=1.1,
	    height=0.6*\textwidth,
	    width=\textwidth,
		xtick={0,0.25,...,1},
        }
}
\begin{axis}[A]
\addplot[domain = 0:1, samples = 200, smooth, line width=2., blue!30] {x^5};
	\addplot[color=red!30!white, line width=1., mark=*, mark options={scale=1, solid}] table[x expr=\thisrow{x},y=u_FEM0]{\fileFEM};
	\addplot[color=black, style=dashed, line width=.7 ,mark=asterisk, mark options = {scale=1, solid}] table[x expr=\thisrow{x},y=u_pred0]{\fileDL};
\end{axis}
\end{tikzpicture}
\caption{Step $s=1$ (after training).}
\label{fig:step1_end_poisson_x5}
\end{subfigure}
\vskip 0.5em
\begin{subfigure}[b]{0.49\textwidth}
\centering
\begin{tikzpicture}
\newcommand{\fileDL}{Chapters/3.Chapter/figures/numerical_results/DATA/results_poisson_x5_Sinv_3refinements/start/test_prediction_3nodes.csv}
\newcommand{\fileFEM}{Chapters/3.Chapter/figures/numerical_results/DATA/results_poisson_x5_Sinv_3refinements/fem/test_FEM_3nodes.csv}
\pgfplotsset{A/.append style={
	    xmin=-0.1,
	    xmax=1.1,
	    ymin=-0.1,
	    ymax=1.1,
	    ylabel = {$u(x)$},
	    height=0.6*\textwidth,
	    width=\textwidth,
		xtick={0,0.25,...,1},
       }
}
\begin{axis}[A]
\addplot[domain = 0:1, samples = 200, smooth, line width=2., blue!30] {x^5};
\addplot[color=red!30!white, line width=1., mark=*, mark options={scale=1, solid}] table[x expr=\thisrow{x},y=u_FEM0]{\fileFEM};
\addplot[color=black, style=dashed, line width=.7 ,mark=asterisk, mark options = {scale=1, solid}] table[x expr=\thisrow{x},y=u_pred0]{\fileDL};
\end{axis}
\end{tikzpicture}
\caption{Step $s=2$ (before training).}
\label{fig:step2_start_poisson_x5}
\end{subfigure}\hfill
\begin{subfigure}[b]{0.49\textwidth}
\centering
\begin{tikzpicture}
\newcommand{\fileDL}{Chapters/3.Chapter/figures/numerical_results/DATA/results_poisson_x5_Sinv_3refinements/end/test_prediction_3nodes.csv}
\newcommand{\fileFEM}{Chapters/3.Chapter/figures/numerical_results/DATA/results_poisson_x5_Sinv_3refinements/fem/test_FEM_3nodes.csv}
\pgfplotsset{A/.append style={
	    xmin=-0.1,
	    xmax=1.1,
	    ymin=-0.1,
	    ymax=1.1,
	    height=0.6*\textwidth,
	    width=\textwidth,
		xtick={0,0.25,...,1},
       }
}
\begin{axis}[A]
\addplot[domain = 0:1, samples = 200, smooth, line width=2., blue!30] {x^5};
\addplot[color=red!30!white, line width=1., mark=*, mark options={scale=1, solid}] table[x expr=\thisrow{x},y=u_FEM0]{\fileFEM};
\addplot[color=black, style=dashed, line width=.7 ,mark=asterisk, mark options = {scale=1, solid}] table[x expr=\thisrow{x},y=u_pred0]{\fileDL};
\end{axis}
\end{tikzpicture}
\caption{Step $s=2$ (after training).}
\label{fig:step2_end_poisson_x5}
\end{subfigure}
\vskip 0.5em    
\begin{subfigure}[b]{0.49\textwidth}
\centering
\begin{tikzpicture}
\newcommand{\fileDL}{Chapters/3.Chapter/figures/numerical_results/DATA/results_poisson_x5_Sinv_3refinements/start/test_prediction_5nodes.csv}
\newcommand{\fileFEM}{Chapters/3.Chapter/figures/numerical_results/DATA/results_poisson_x5_Sinv_3refinements/fem/test_FEM_5nodes.csv}
\pgfplotsset{A/.append style={
	    xmin=-0.1,
	    xmax=1.1,
	    ymin=-0.1,
	    ymax=1.1,
	    ylabel = {$u(x)$},
	    height=0.6*\textwidth,
	    width=\textwidth,
		xtick={0,0.25,...,1},
       }
}
\begin{axis}[A]
\addplot[domain = 0:1, samples = 200, smooth, line width=2., blue!30] {x^5};
\addplot[color=red!30!white, line width=1., mark=*, mark options={scale=1, solid}] table[x expr=\thisrow{x},y=u_FEM0]{\fileFEM};
\addplot[color=black, style=dashed, line width=.7 ,mark=asterisk, mark options = {scale=1, solid}] table[x expr=\thisrow{x},y=u_pred0]{\fileDL};
\end{axis}
\end{tikzpicture}
\caption{Step $s=3$ (before training).}
\label{fig:step3_start_poisson_x5}
\end{subfigure}\hfill
\begin{subfigure}[b]{0.49\textwidth}
\centering
\begin{tikzpicture}
\newcommand{\fileDL}{Chapters/3.Chapter/figures/numerical_results/DATA/results_poisson_x5_Sinv_3refinements/end/test_prediction_5nodes.csv}
\newcommand{\fileFEM}{Chapters/3.Chapter/figures/numerical_results/DATA/results_poisson_x5_Sinv_3refinements/fem/test_FEM_5nodes.csv}
\pgfplotsset{A/.append style={
	    xmin=-0.1,
	    xmax=1.1,
	    ymin=-0.1,
	    ymax=1.1,
	    height=0.6*\textwidth,
	    width=\textwidth,
		xtick={0,0.25,...,1},
       }
}
\begin{axis}[A]
\addplot[domain = 0:1, samples = 200, smooth, line width=2., blue!30] {x^5};
\addplot[color=red!30!white, line width=1., mark=*, mark options={scale=1, solid}] table[x expr=\thisrow{x},y=u_FEM0]{\fileFEM};
\addplot[color=black, style=dashed, line width=.7 ,mark=asterisk, mark options = {scale=1, solid}] table[x expr=\thisrow{x},y=u_pred0]{\fileDL};
\end{axis}
\end{tikzpicture}
\caption{Step $s=3$ (after training). }
\label{fig:step3_end_poisson_x5}
\end{subfigure}
\vskip 0.5em
\begin{subfigure}[b]{0.49\textwidth}
\centering
\begin{tikzpicture}
\newcommand{\fileDL}{Chapters/3.Chapter/figures/numerical_results/DATA/results_poisson_x5_Sinv_3refinements/start/test_prediction_9nodes.csv}
\newcommand{\fileFEM}{Chapters/3.Chapter/figures/numerical_results/DATA/results_poisson_x5_Sinv_3refinements/fem/test_FEM_9nodes.csv}
\pgfplotsset{A/.append style={
	    xmin=-0.1,
	    xmax=1.1,
	   	xlabel = {$x$},
	    ymin=-0.1,
	    ymax=1.1,
	    ylabel = {$u(x)$},
	    height=0.6*\textwidth,
	    width=\textwidth,
		xtick={0,0.25,...,1},
       }
}
\begin{axis}[A]
\addplot[domain = 0:1, samples = 200, smooth, line width=2., blue!30] {x^5};
\addplot[color=red!30!white, line width=1., mark=*, mark options={scale=1, solid}] table[x expr=\thisrow{x},y=u_FEM0]{\fileFEM};
\addplot[color=black, style=dashed, line width=.7 ,mark=asterisk, mark options = {scale=1, solid}] table[x expr=\thisrow{x},y=u_pred0]{\fileDL};
\end{axis}
\end{tikzpicture}
\caption{Step $s=4$ (before training).}
\label{fig:step4_start_poisson_x5}
\end{subfigure}\hfill
\begin{subfigure}[b]{0.49\textwidth}
\centering
\begin{tikzpicture}
\newcommand{\fileDL}{Chapters/3.Chapter/figures/numerical_results/DATA/results_poisson_x5_Sinv_3refinements/end/test_prediction_9nodes.csv}
\newcommand{\fileFEM}{Chapters/3.Chapter/figures/numerical_results/DATA/results_poisson_x5_Sinv_3refinements/fem/test_FEM_9nodes.csv}
\pgfplotsset{A/.append style={
	    xmin=-0.1,
	    xmax=1.1,
	    xlabel = {$x$},
	    ymin=-0.1,
	    ymax=1.1,
	    height=0.6*\textwidth,
	    width=\textwidth,
		xtick={0,0.25,...,1},
       }
}
\begin{axis}[A]
\addplot[domain = 0:1, samples = 200, smooth, line width=2., blue!30] {x^5};
\addplot[color=red!30!white, line width=1., mark=*, mark options={scale=1, solid}] table[x expr=\thisrow{x},y=u_FEM0]{\fileFEM};
\addplot[color=black, style=dashed, line width=.7 ,mark=asterisk, mark options = {scale=1, solid}] table[x expr=\thisrow{x},y=u_pred0]{\fileDL};
\end{axis}
\end{tikzpicture}
\caption{Step $s=4$ (after training).}
\label{fig:step4_end_poisson_x5}
\end{subfigure}
    
\caption{First four steps of the \acs{DeepFEM} for problem \eqref{eq:poisson_x5_10refinements}.  $u^*$ is the exact solution, $u_{\text{FEM}}$ is the finite element solution, and $u_{\text{NN}}$ is the NN prediction at each training step.}
    \label{fig:poisson_x5}
\end{figure}

\Cref{fig:poisson_x5} shows the predictions of \ac{DeepFEM} over the first four steps before and after training. At each step,  we observe that the model extends the coarse prediction to the fine mesh, and then (re)trains the resulting NN to obtain a proper fine grid approximation.  Results show superb accuracy.

\begin{figure}[htbp]
\centering
\begin{subfigure}[b]{.98\textwidth}
\begin{tikzpicture}
\pgfplotsset{loss/.append style={ 
     xmin=0,         
     xmax=8e3,      
     xlabel = {iteration}, 
     ymin=1e-13,     
     ymax=10,      
     ylabel = {loss},
     xtick={0,1000,2000,3000,4000,5000,6000,7000},
    scaled x ticks=false,
     height=0.3*\textwidth,    
     width=0.98*\textwidth,    
     legend style={fill=none, at={(0.5,1.45)}, anchor=north, minimum height = 0.6cm, /tikz/every even column/.append style={column sep=0.5cm}},
	 legend columns = -1    
      }
}
\begin{semilogyaxis}[loss]
\addplot[color=blue!50, line width=1.5 ,] table[x expr=\thisrow{iteration},y=loss_history]{Chapters/3.Chapter/figures/numerical_results/DATA/results_poisson_x5_errorenergy_10refinements_endtoend/loss_history_reduced3.csv};
\addlegendentry{ loss $=\lVert\mathbf{A}\mathbf{u}_\text{NN}-\mathbf{f}\rVert_{\mathbf{A}^{-1}}$ = energy-norm error};
\draw[black,<->, line width = 0.85] (0, 5e-10) -- (1400, 5e-10) node[midway, below] {$s=1$};
\draw[black,<->, line width = 0.85] (1401, 1e-9) -- (2937, 1e-9) node[midway, below] {$s=2$};
\draw[black,<->, line width = 0.85] (2938, 3e-9) -- (4523, 3e-9) node[midway, below] {$s=3$};
\draw[black,<->, line width = 0.85] (4524, 1e-8) -- (6313, 1e-8) node[midway, below] {$s=4$};
\draw[red,<->, line width = 0.65] (0, 5e-6) -- (798, 5e-6) node[midway, below] {\footnotesize Adam};
\draw[red,<->, line width = 0.65] (5082, 3e-4) -- (6313, 3e-4) node[midway, above] {\footnotesize Adalr};
\end{semilogyaxis}
\end{tikzpicture}
\caption{Loss function evolution during the first four steps.}
\label{fig:loss_poisson_x5_errorenergy_4refinements_endtoend}
\end{subfigure}
\begin{subfigure}[b]{\textwidth}
\centering
\begin{tikzpicture}
\pgfplotsset{loss/.append style={ 
     xlabel = {iteration}, 
     ymin=1e-13,     
     ymax=10,      
     ylabel = {loss},
    scaled x ticks=false,
     height=0.3*\textwidth,    
     width=0.98*\textwidth,    
     }
}
\begin{semilogyaxis}[loss]
\addplot[color=blue!50, line width=1.5 ,] table[x expr=\thisrow{iteration},y=loss_history]{Chapters/3.Chapter/figures/numerical_results/DATA/results_poisson_x5_errorenergy_10refinements_endtoend/loss_history_reduced3.csv};
\end{semilogyaxis}
\end{tikzpicture}
\caption{Loss function evolution during eleven steps.}
\label{fig:loss_poisson_x5_errorenergy_10refinements_endtoend}
\end{subfigure}
\vskip 0.5em
\centering
\begin{subfigure}[t]{.98\textwidth}
\centering
\begin{tikzpicture}
\pgfplotsset{A/.append style={
	    xmin=-0.1,
	    xmax=1.1,
	    xlabel = {$x$},
	    ylabel = {$u_{\text{NN}}-u_\text{FEM}$},
	    height=0.3*\textwidth,
	    width=0.5*\textwidth,
		xtick={0,0.25,...,1},
        }
}
\begin{axis}[A]
\addplot[color=red, line width=1.3] table[x expr=\thisrow{x},y=error0]{Chapters/3.Chapter/figures/numerical_results/DATA/results_poisson_x5_errorenergy_10refinements_endtoend/error/test_error_1025nodes.csv};
\end{axis}
\end{tikzpicture}
\caption{Error function at the end of training.}
\label{fig:poisson_x5_errorfunction_errorenergy_1025nodes_endtoend}
\end{subfigure}
\caption{End-to-end training of \acs{DeepFEM} for problem \eqref{eq:poisson_x5_10refinements}. The loss function coincides with the energy-norm error.}
\label{fig:loss_poisson_x5_errorenergy_4-10refinements_endtoend}
\end{figure}
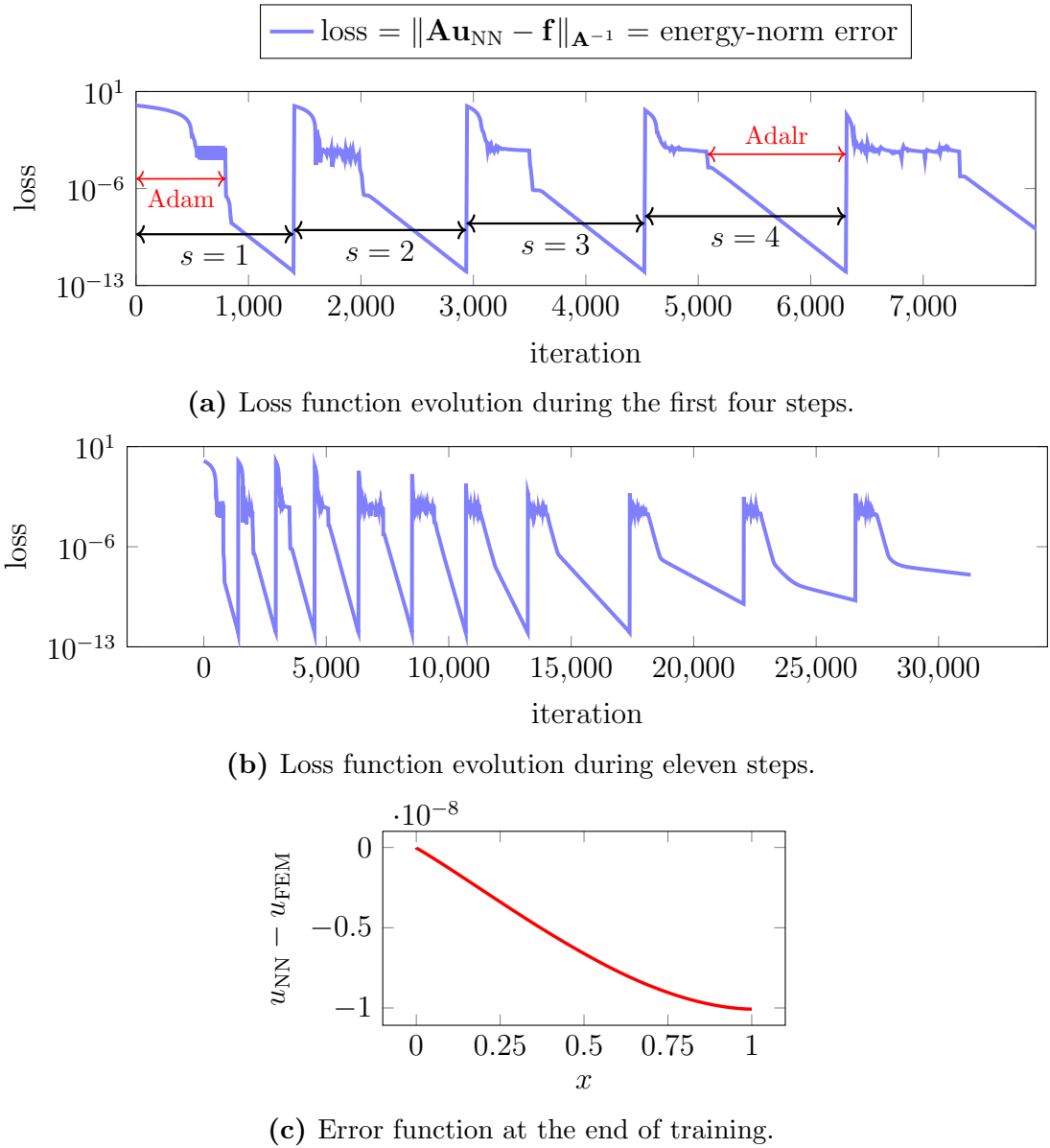

\Cref{fig:loss_poisson_x5_errorenergy_4refinements_endtoend} shows the loss function evolution along the first four steps when it coincides with the energy-norm error,  while \Cref{fig:loss_poisson_x5_errorenergy_10refinements_endtoend} shows all eleven steps of the training process.  This example illustrates the convergence behavior and limitations of DL under optimal conditions, i.e., when the inverse of the matrix is already part of the loss. The vertical jumps observed in the loss function correspond to the extension from one grid to the next one. After each jump, we observe a noisy convergence of the loss, which corresponds to the use of the Adam optimizer that works as an initial aggressive loss descender that gets stunned quickly. We then switch to the Adalr optimizer, which exhibits a monotonic loss decrease. At each phase, we select an initial learning rate equal to the first loss evaluation multiplied by $10^{-3}$ (for Adam) or by $10^{-2}$ (for Adalr). We set a maximum of $2\mathord{,}000$ (Adam) and $4\mathord{,}000$ (Adalr) iterations for each optimizer performance at each step. Moreover, if we attain a loss below $10^{-12}$, we stop the optimization at each step. We carry out an end-to-end training of the \ac{DeepFEM} model. We observe a convergence deterioration as we increase the step number (and grid size) that remains below $10^{-8}$ in all occasions (see \Cref{fig:poisson_x5_errorfunction_errorenergy_1025nodes_endtoend}).

\subsubsection{End-to-end vs. layer-by-layer training}

\begin{figure}[htbp]
\centering
\begin{subfigure}[t]{.98\textwidth}
\centering
\begin{tikzpicture}
\pgfplotsset{loss/.append style={ 
     xlabel = {iteration},     
     ylabel = {loss},
     scaled x ticks=false,
     height=0.3*\textwidth,    
     width=0.98*\textwidth,    
     legend style={fill=none, at={(0.5,1.5)},anchor=north, minimum height = 0.6cm, /tikz/every even column/.append style={column sep=0.5cm}},
	 legend columns = -1
      }
}
\begin{semilogyaxis}[loss]
\addplot[color=blue!50, line width=1.5 ,] table[x expr=\thisrow{iteration},y=loss_history]{Chapters/3.Chapter/figures/numerical_results/DATA/results_poisson_x5_errorenergy_10refinements_onlylast/loss_history_reduced3.csv};
    \addlegendentry{loss = $\lVert\mathbf{A}\mathbf{u}_\text{NN}-\mathbf{f}\rVert_{\mathbf{A}^{-1}}$ = energy-norm error};
\end{semilogyaxis}
\end{tikzpicture}
\caption{Loss function evolution}
\label{fig:loss_poisson_x5_errorenergy_10refinements_onlylast}
\end{subfigure}
\vskip 0.5em
\centering
\begin{subfigure}[t]{\textwidth}
\centering
\begin{tikzpicture}
\pgfplotsset{A/.append style={
		compat=1.17,
	    xmin=-0.1,
	    xmax=1.1,
	   	ymin=-0.0015,
	    ymax=0.0001,
	    xlabel = {$x$},
	    ylabel = {$u_{\text{NN}}-u_{\text{FEM}}$},
	    height=0.3*\textwidth,
	    width=0.5*\textwidth,
		xtick={0,0.25,...,1},
        }
}
\begin{axis}[A]
	\addplot[color=red, line width=1.3] table[x expr=\thisrow{x},y=error0]{Chapters/3.Chapter/figures/numerical_results/DATA/results_poisson_x5_errorenergy_10refinements_onlylast/error/new_error3.csv};
\end{axis}
\end{tikzpicture}
\caption{Error function at the end of the training.}
\label{fig:L2error_poisson_x5_errorenergy_10refinements_onlylast}
\end{subfigure}
\caption{Layer-by-layer training of the \acs{DeepFEM} model for problem \eqref{eq:poisson_x5_10refinements} employing the energy norm in the loss function. }
\label{fig:training_poisson_x5_errorenergy_10refinements_onlylast}
\end{figure}
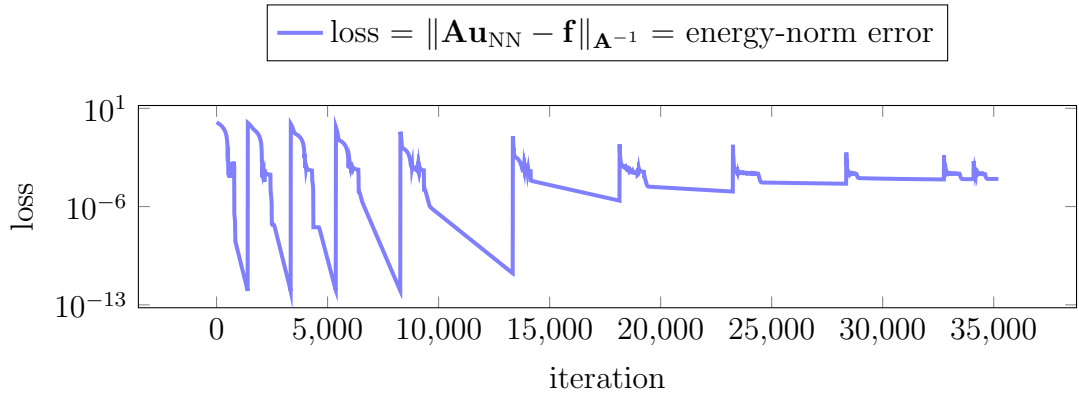
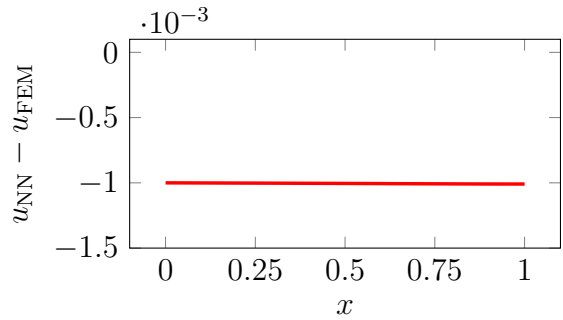

By construction, each layer of \ac{DeepFEM} produces the coefficients related to the basis functions associated with each Finite Element mesh. The finer the mesh, the more local the supports of the basis functions associated with the coefficients. While end-to-end training allows adjustments in the entire hierarchy of the coefficients related to basis functions of coarse and fine meshes, layer-by-layer training only adjusts coefficients associated with the finest mesh \Cref{fig:loss_poisson_x5_errorenergy_10refinements_onlylast} shows the loss convergence under the same conditions as above but performing a layer-by-layer training. We observe that convergence deteriorates as the number of iterations increases compared to performing an end-to-end training (recall \Cref{fig:loss_poisson_x5_errorenergy_10refinements_endtoend}). This occurs because the loss involves the gradient of the error when employing the energy-norm. \Cref{fig:L2error_poisson_x5_errorenergy_10refinements_onlylast} shows the error function, which is almost constant. Thus, the derivative of the error is nearly zero, so it is challenging to minimize the energy-norm error by adjusting individual coefficients associated with local-support basis functions since this often implies an increase in the derivative of the error. This convergence accelerates using both global-support and local-support basis functions, as in multigrid methods \cite{bramble2019multigrid}.

\begin{figure}[htbp]
\centering
\begin{subfigure}[t]{.98\textwidth}
\centering
\begin{tikzpicture}
\pgfplotsset{loss/.append style={ 
     xlabel = {iteration}, 
     ymin=1e-13,     
     ymax=10,      
     ylabel = {loss},
     scaled x ticks=false,
     height=0.3*\textwidth,    
     width=0.98*\textwidth,    
     legend style={fill=none, at={(0.5,1.5)},anchor=north, minimum height = 0.6cm, /tikz/every even column/.append style={column sep=0.5cm}},
	 legend columns = -1
      }    
}
\begin{semilogyaxis}[loss]
\addplot[color=blue!50, line width=1.5 ,] table[x expr=\thisrow{iteration},y=loss_history]{Chapters/3.Chapter/figures/numerical_results/DATA/results_poisson_x5_errorL2_10refinements_onlylast/loss_history_reduced3.csv};
\addlegendentry{loss = $\lVert\mathbf{A}^{-1}(\mathbf{A}\mathbf{u}_\text{NN}-\mathbf{f})\rVert_{\mathbf{M}}$ = $L^2$-norm error};
\end{semilogyaxis}
\end{tikzpicture}
\caption{Loss function evolution when perfoming a layer-by-layer training.}
\label{fig:loss_poisson_x5_errorL2_10refinements_onlylast}
\end{subfigure}
\begin{subfigure}[t]{.98\textwidth}
\begin{tikzpicture}
\pgfplotsset{loss/.append style={ 
     xtick={0,10000,20000,30000,40000,50000},     
     xlabel = {iteration}, 
     ymin=1e-13,     
     ymax=10,      
     ylabel = {loss},
     scaled x ticks=false,
     height=0.3*\textwidth,    
     width=0.98*\textwidth,    
      }    
}
\begin{semilogyaxis}[loss]
\addplot[color=blue!50, line width=1.5 ,] table[x expr=\thisrow{iteration},y=loss_history]{Chapters/3.Chapter/figures/numerical_results/DATA/results_poisson_x5_errorL2_10refinements_endtoend/loss_history_reduced3.csv};
\end{semilogyaxis}
\end{tikzpicture}
\caption{Loss function evolution when performing an end-to-end training.}
\label{fig:loss_poisson_x5_errorL2_10refinements_endtoend}
\end{subfigure}
\caption{Training of the \acs{DeepFEM} for model problem \eqref{eq:poisson_x5_10refinements} using the $L^2$-norm for the loss function.}
\label{fig:loss_poisson_x5_errorL2_10refinements}
\end{figure}
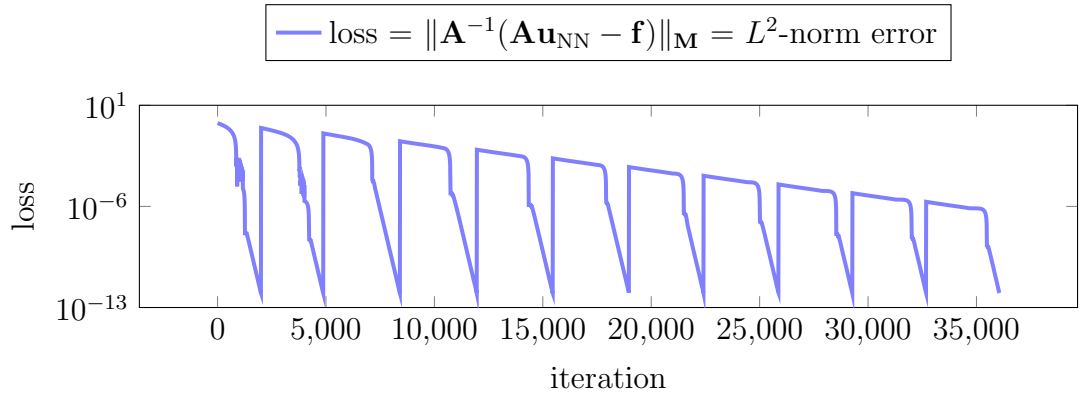
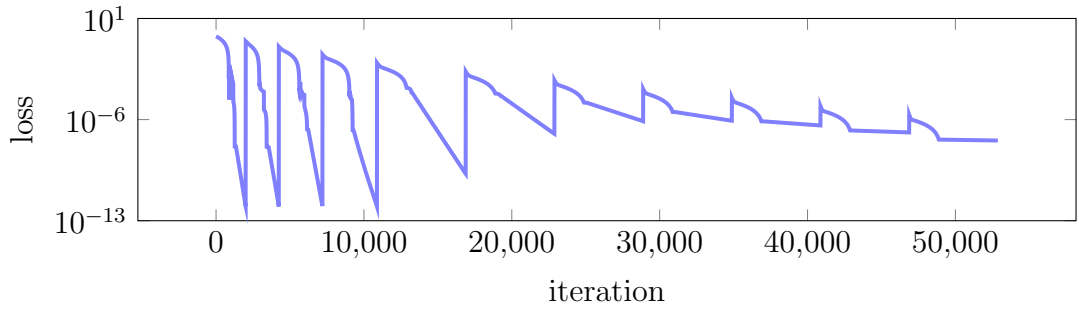

If we select a norm for the loss that ignores gradients, e.g., the $L^2$-norm, training of local-support basis functions provides outstanding results, as shown in \Cref{fig:loss_poisson_x5_errorL2_10refinements_onlylast}. However, optimizing with respect to the $L^2$-norm is discouraged for solving differential equations.

Even if the end-to-end training is the best alternative when dealing with losses involving the gradient of the residual/error, we observe a convergence deterioration of the loss function as iterations move forward (recall \Cref{fig:loss_poisson_x5_errorenergy_10refinements_endtoend}). We suspect this is independent of the norm selection, but it occurs because of the conflicting coexistence of many trainable variables. To illustrate this, we consider the above $L^2$-norm case where layer-by-layer training suffices to achieve an outstanding convergence, and we perform an end-to-end training. \Cref{fig:loss_poisson_x5_errorL2_10refinements_endtoend} shows the convergence deterioration of the loss function. 

In the following, we only consider end-to-end training cases of study.

\subsubsection{Preconditioners action}

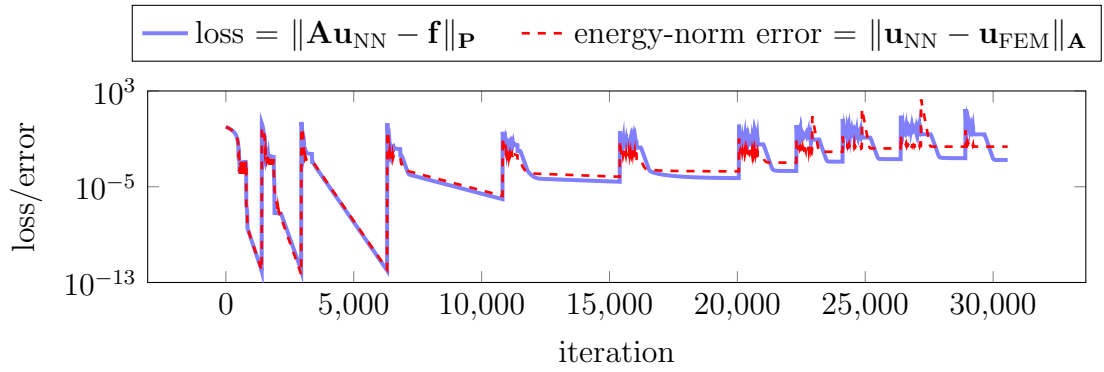
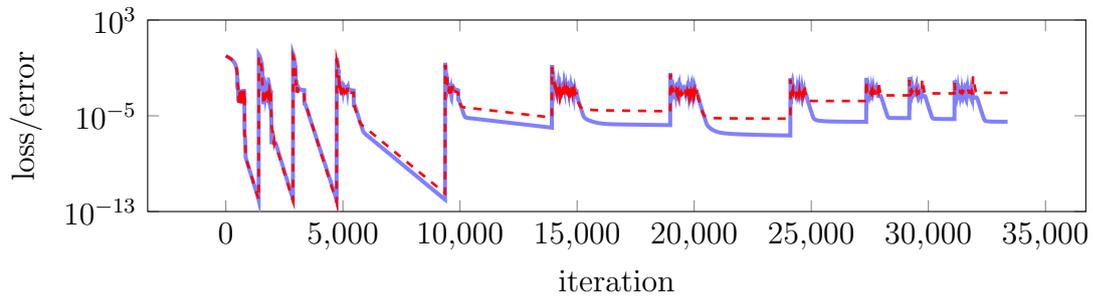
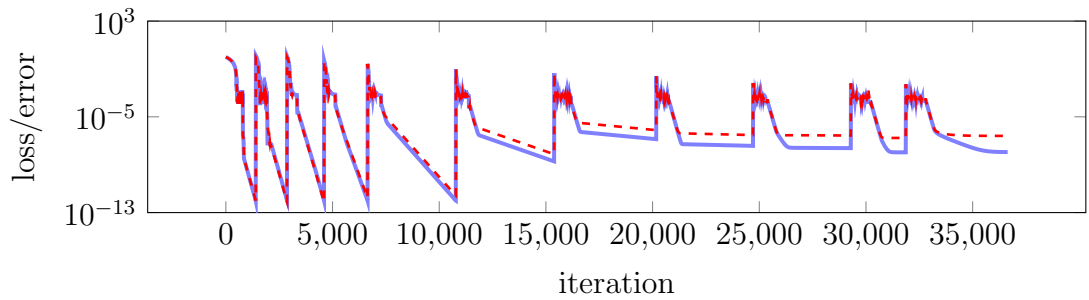
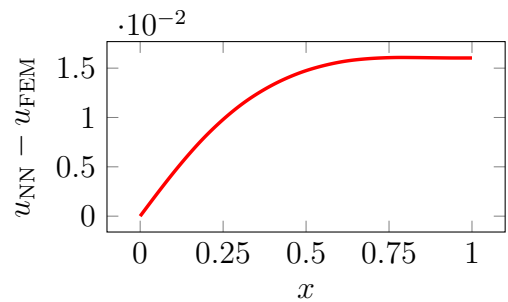
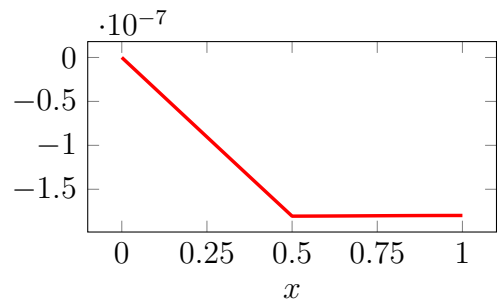
\begin{figure}[htbp]
\centering
\begin{subfigure}[b]{.98\textwidth}
\centering
\begin{tikzpicture}
\pgfplotsset{loss/.append style={ 
     xlabel = {iteration}, 
     ymin=1e-13,     
     ymax=10e2,      
     ylabel = {loss/error},
    scaled x ticks=false,
     height=0.29*\textwidth,    
     width=0.98*\textwidth,    
     legend style={fill=none, at={(0.5,1.45)},anchor=north, minimum height = 0.6cm, /tikz/every even column/.append style={column sep=0.5cm}},
	 legend columns = -1
      }    
}
\begin{semilogyaxis}[loss]
\addplot[color=blue!50, line width=1.5 ,] table[x expr=\thisrow{iteration},y=loss_history]{Chapters/3.Chapter/figures/numerical_results/DATA/results_poisson_x5_residuall2_10refinements_endtoend/loss_history_reduced3.csv};
\addlegendentry{loss = $\lVert\mathbf{A}\mathbf{u}_\text{NN}-\mathbf{f}\rVert_{\mathbf{P}}$};
\addplot[color=red!, line width=1 , style=dashed] table[x expr=\thisrow{iteration},y=error_history]{Chapters/3.Chapter/figures/numerical_results/DATA/results_poisson_x5_residuall2_10refinements_endtoend/error_history_reduced3.csv};
\addlegendentry{energy-norm error = $\lVert\mathbf{u}_\text{NN}-\mathbf{u}_\text{FEM}\rVert_{\mathbf{A}}$};
\end{semilogyaxis}
\end{tikzpicture}
\caption{Employing no preconditioner, i.e.,  $\mathbf{P}=\mathbf{I}$.}
\label{fig:loss_poisson_x5_residuall2_10refinements_endtoend}
\end{subfigure}
\vskip 0.5em
\begin{subfigure}[b]{.98\textwidth}
\centering
\begin{tikzpicture}
\pgfplotsset{loss/.append style={ 
     xlabel = {iteration}, 
     ymin=1e-13,     
     ymax=10e2,      
     ylabel = {loss/error},
     scaled x ticks=false,
     height=0.29*\textwidth,    
     width=0.98*\textwidth,    
     }    
}
\begin{semilogyaxis}[loss]
\addplot[color=blue!50, line width=1.5 ,] table[x expr=\thisrow{iteration},y=loss_history]{Chapters/3.Chapter/figures/numerical_results/DATA/results_poisson_x5_residualS2_10refinements_endtoend/loss_history_reduced3.csv};
\addplot[color=red!, line width=1 , style=dashed] table[x expr=\thisrow{iteration},y=error_history]{Chapters/3.Chapter/figures/numerical_results/DATA/results_poisson_x5_residualS2_10refinements_endtoend/error_history_reduced3.csv};
\end{semilogyaxis}
\end{tikzpicture}
\caption{Employing blocks of size two.}
\label{fig:loss_poisson_x5_residualS2_10refinements_endtoend}
\end{subfigure}
\vskip 0.5em
\begin{subfigure}[t]{.98\textwidth}
\centering
\begin{tikzpicture}
\pgfplotsset{loss/.append style={ 
     xlabel = {iteration}, 
     ymin=1e-13,     
     ymax=10e2,      
     ylabel = {loss/error},
     xtick={0,5000,10000,15000,20000,25000,30000,35000},
     scaled x ticks=false,
     height=0.29*\textwidth,    
     width=0.98*\textwidth,      
      }
}  
\begin{semilogyaxis}[loss]
\addplot[color=blue!50, line width=1.5 ,] table[x expr=\thisrow{iteration},y=loss_history]{Chapters/3.Chapter/figures/numerical_results/DATA/results_poisson_x5_residualShalves_10refinements_endtoend/loss_history_reduced3.csv};
\addplot[color=red!, line width=1 , style=dashed] table[x expr=\thisrow{iteration},y=error_history]{Chapters/3.Chapter/figures/numerical_results/DATA/results_poisson_x5_residualShalves_10refinements_endtoend/error_history_reduced3.csv};
\end{semilogyaxis}
\end{tikzpicture}
\caption{Employing blocks of size equal to half the size of the mesh.}
\label{fig:loss_poisson_x5_residualShalves_10refinements_endtoend} 
\end{subfigure}
\vskip 0.5em
\begin{subfigure}[t]{.48\textwidth}
\centering
\begin{tikzpicture}
\pgfplotsset{A/.append style={
	    xmin=-0.1,
	    xmax=1.1,
	    xlabel = {$x$},
	    ylabel = {$u_\text{NN}-u_\text{FEM}$},
	    height=0.59*\textwidth,
	    width=0.98*\textwidth,
		xtick={0,0.25,...,1},
        }
}
\begin{axis}[A]
\addplot[color=red, line width=1.3] table[x expr=\thisrow{x},y=error0]{Chapters/3.Chapter/figures/numerical_results/DATA/results_poisson_x5_residuall2_10refinements_endtoend/error/test_error_1025nodes.csv};
\end{axis}
\end{tikzpicture}
\caption{Error function of case (a).}\label{fig:poisson_x5_errorfunction_residuall2_1025nodes_endtoend}\end{subfigure} \hfill
\begin{subfigure}[t]{0.48\textwidth}
\centering
\begin{tikzpicture}
\pgfplotsset{A/.append style={
	    xmin=-0.1,
	    xmax=1.1,
	    xlabel = {$x$},
	    height=0.59*\textwidth,
	    width=\textwidth,
		xtick={0,0.25,...,1},
        }
}
\begin{axis}[A]
\addplot[color=red, line width=1.3] table[x expr=\thisrow{x},y=error0]{Chapters/3.Chapter/figures/numerical_results/DATA/results_poisson_x5_residualShalves_10refinements_endtoend/error/test_error_1025nodes.csv};
\end{axis}
\end{tikzpicture}
\caption{Error function of case (c).}\label{fig:poisson_x5_errorfunction_residualShalves_1025nodes_endtoend}
\end{subfigure}
\caption{End-to-end training of the \acs{DeepFEM} for problem \eqref{eq:poisson_x5_10refinements} using block-Jacobi preconditioners of different sizes.}
    \label{fig:loss_poisson_x5_residual_10refinements_endtoend}
\end{figure}

We now consider the loss given in \eqref{loss_preconditioned_residual} with three different preconditioners: (a) $\mathbf{P}$ being the identity matrix (\Cref{fig:loss_poisson_x5_residuall2_10refinements_endtoend}); (b) $\mathbf{P}$ being a block-Jacobi preconditioner with blocks of size two (\Cref{fig:loss_poisson_x5_residualS2_10refinements_endtoend}); and (c) $\mathbf{P}$ being a block-Jacobi preconditioner with blocks of size equal to half the number of elements in the mesh (\Cref{fig:loss_poisson_x5_residualShalves_10refinements_endtoend}). In all the cases, we show the loss function evolution along with the energy-norm error evolution.

We observe some differences between the energy-norm error and the loss that increases with the mesh size, as expected. We also observe that the larger the size of the block-Jacobi, the discrepancy between the loss and the norm error reduces. In addition, the further the loss is from the energy-norm error, the earlier the loss stagnates convergence---compare errors at \Cref{fig:poisson_x5_errorfunction_residuall2_1025nodes_endtoend} and \Cref{fig:poisson_x5_errorfunction_residualShalves_1025nodes_endtoend}. This suggests that the loss induced by the energy-norm error is more convex with respect to the variables than other simplified variants, such as the loss caused by the $2$-norm of the residual vector.

\subsubsection{Norm interchange during training}

Depending on the norm we select, we deal with different convexity shapes of the loss with respect to the variables of the NN. Within certain portion of the learnable variables domain, it often happens that a given loss is more convex than others, which has a direct impact on the optimizer convergence. With the aim of avoiding stagnation, we propose to change the norm during optimization, expecting to improve convexity. In this way, we consider
\begin{equation}
\label{loss_norm_change}
\mathcal{L}(\theta;\sigma, \alpha) = C_E\lVert\mathbf{A}\mathbf{u}_\text{NN}-\mathbf{f}\rVert_{\mathbf{P}} + \ C_{L^2}\lVert\mathbf{P}(\mathbf{A}\mathbf{u}_{\text{NN}} - \mathbf{f})\rVert_{\mathbf{M}},
\end{equation} where $C_E, C_{L^2}\in\{0,1\}$ are distinct values that interchange when the convergence stagnates. 

\begin{figure}[htb]
\centering
\begin{tikzpicture}
\pgfplotsset{loss/.append style={ 
     xmin=8e3,         
     xtick={10000,12000,...,22000},     
     xlabel = {iteration}, 
     ylabel = {loss},
    scaled x ticks=false,
     height=0.3*\textwidth,    
     width=0.98*\textwidth,    
     legend style={fill=none, at={(0.5,1.6)}, anchor=north, minimum height = 0.6cm, /tikz/every even column/.append style={column sep=0.5cm}},
	 legend cell align={left}    
     }
}
\begin{semilogyaxis}[loss]
\addplot[color=blue!40, line width=2.] table[x expr=\thisrow{iteration},y=loss_history]{Chapters/3.Chapter/figures/numerical_results/DATA/results_poisson_x5_residualS2_4refinements_moreiterations_endtoend/loss_history_reduced2.csv};
\addlegendentry{loss = $\lVert\mathbf{A}\mathbf{u}_\text{NN} - \mathbf{f}\rVert_{\mathbf{P}}$};
\addplot[color=green!20!black, line width=1.25 , style=dashed] table[x expr=\thisrow{iteration},y=loss_history]{Chapters/3.Chapter/figures/numerical_results/DATA/results_poisson_x5_residualS2_4refinements_lossinterchange_endtoend/loss_history_reduced2.csv};
\addlegendentry{loss = $C_E \lVert\mathbf{A}\mathbf{u}_\text{NN} - \mathbf{f}\rVert_{\mathbf{P}} + C_{L^2} \lVert\mathbf{P}(\mathbf{A}\mathbf{u}_\text{NN} - \mathbf{f})\rVert_{\mathbf{M}}$};
\node[green!20!black] at (11000,1e-8) {$C_E = 1$};
\node[green!20!black] at (16000,6e-12) {$C_{L^2} = 1$};
\node[green!20!black] at (21000,3e-12) {$C_E = 1$};
\end{semilogyaxis}
\end{tikzpicture}
\caption{Energy-norm vs. energy- and $L^2$-norm interchange training at step $s=5$ in \Cref{fig:loss_poisson_x5_residualS2_10refinements_endtoend}.}
\label{fig:loss_poisson_x5_residualS2_4refinements_lossinterchange_endtoend}
\end{figure}
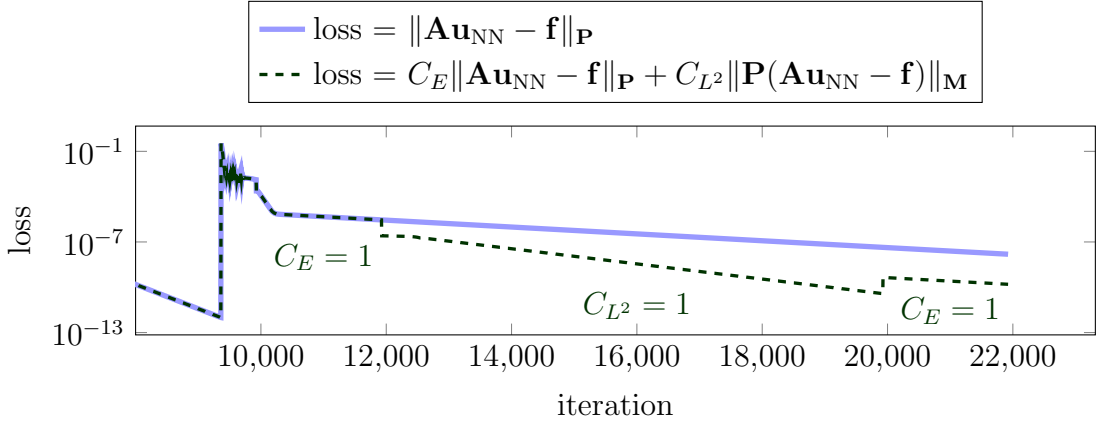

To illustrate the above idea, we consider the case of \Cref{fig:loss_poisson_x5_residualS2_10refinements_endtoend}. We maintain the loss at equation \eqref{loss_preconditioned_residual} for the first four steps. Then, we consider two variants for the loss in the fifth step when employing the Adalr optimizer: (a) same loss as in the previous steps---maintaining $C_E=1$ in \eqref{loss_norm_change}---but with a maximum of $12\mathord{,}000$ iterations; and (b) the loss at equation \eqref{loss_norm_change} with $C_E = 1$ for $2\mathord{,}000$ iterations, changing to $C_{L^2}=1$ for another $8\mathord{,}000$ iterations, and returning to $C_E = 1$ for $2\mathord{,}000$ additional iterations. Although the total number of iterations is the same in all situations, we obtain a loss at the end of the fifth step that is lower when performing the norm change (around $10^{-10}$) than when maintaining the energy-norm (around $10^{-8}$)---see \Cref{fig:loss_poisson_x5_residualS2_4refinements_lossinterchange_endtoend}. We observe that the slope of the loss function evolution is higher when employing the $L^2$-norm, allowing us to start from a lower loss value when returning to the energy norm than maintaining it during the entire training step.

\newpage

\subsection{Sine solution in Poisson and Helmholtz equations}
\label[section]{section3.6.2}

We consider two different \acp{BVP} whose exact solutions are $u^*(x)=\sin(10\pi x)$:
\begin{equation} \label{eq:poisson_helmholtz_sin5_3refinements}
\begin{cases}
-u'' = 100\pi^2\sin(10\pi x), &\\
u(0)=0, u'(1)=10\pi,
\end{cases} \qquad
\begin{cases}
-u''-100\pi^2 u = 0, &\\
u(0)=0, u'(1)=10\pi.
\end{cases}
\end{equation} The former is a symmetric and positive-definite problem, while the second is indefinite. We solve both employing the \ac{DeepFEM} selecting the loss as the $H^1$-norm of the preconditioned residual,
\begin{equation}
\mathcal{L}(\theta; \sigma, \alpha) = \lVert\mathbf{P}(\mathbf{A}\mathbf{u}_\text{NN}-\mathbf{f})\rVert_{\mathbf{K}+\mathbf{M}}.
\end{equation} We start from $32$ elements and perform three mesh refinements. We select $\mathbf{P}$ with blocks of size $32$ for $1\leq s\leq 4$. The optimizers and NN architecture are those of \Cref{section3.6.1}. 

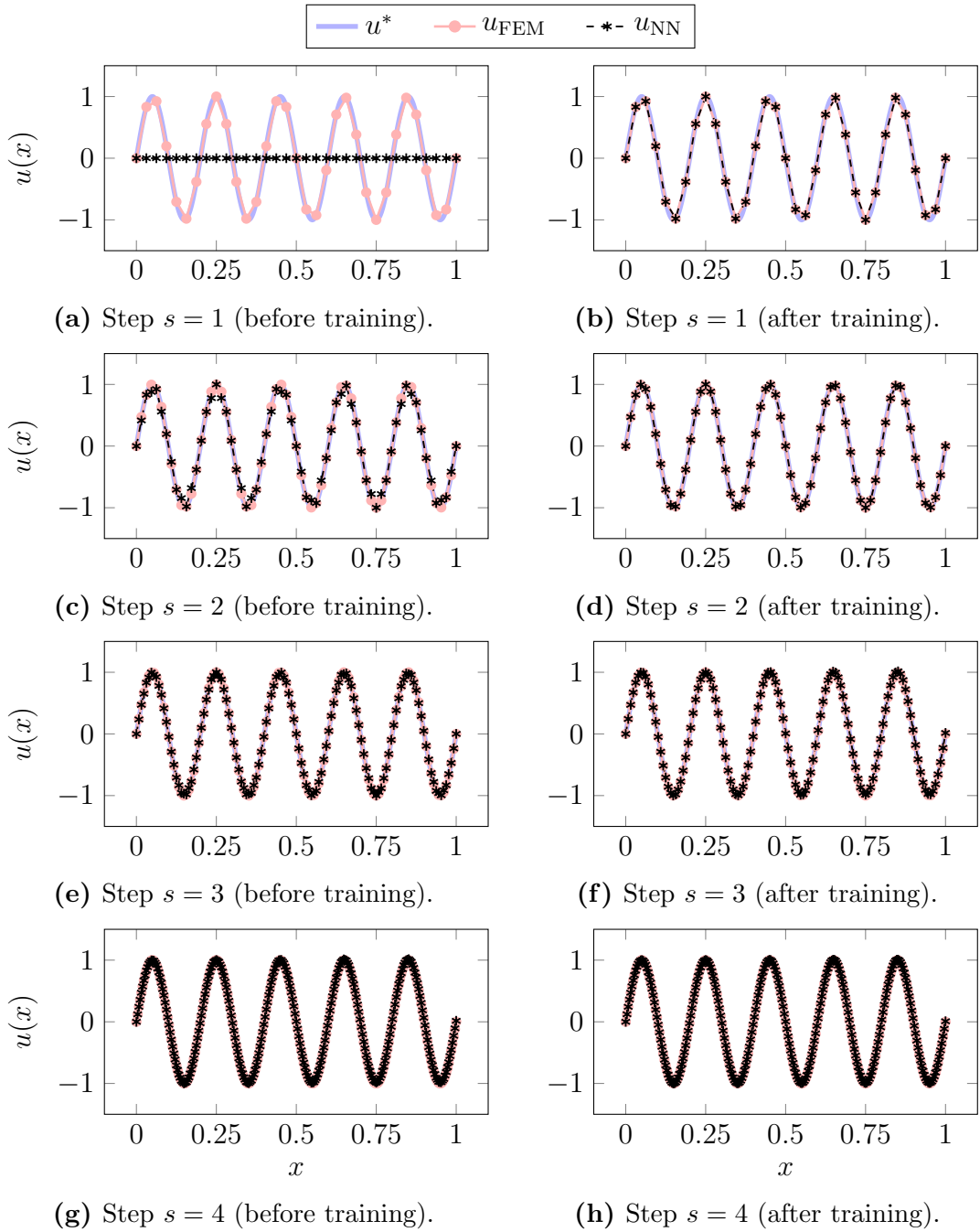
\begin{figure}[htbp]
\centering
\begin{subfigure}[b]{\textwidth}
\centering
 \begin{tikzpicture}
\pgfplotsset{A/.append style={
		hide axis,
	    xmin=-0.1,
	    xmax=1.1,
	    ymin=-0.1,
	    ymax=1.1,
	 	legend style={fill=none, anchor=north, /tikz/every even column/.append style={column sep=0.5cm}},
	 	legend columns = -1
            }
}
\begin{axis}[A]
\addlegendimage{color=blue!30!white, line width=2.};
\addlegendentry{$u^*$};
\addlegendimage{color=red!30!white, line width=1., mark=*, mark options={scale=1, solid}};
\addlegendentry{$u_{\text{FEM}}$};
\addlegendimage{color=black, style=dashed, line width=.7 ,mark=asterisk, mark options = {scale=1, solid}};
\addlegendentry{$u_{\text{NN}}$};
\end{axis}
\end{tikzpicture}
\end{subfigure}
\vskip 0.5em
\begin{subfigure}[b]{0.49\textwidth}
\centering
\begin{tikzpicture}
\newcommand{\fileDL}{Chapters/3.Chapter/figures/numerical_results/DATA/results_poisson_sin5_residualS33_3refinements_endtoend/start/test_prediction_33nodes.csv}
\newcommand{\fileFEM}{Chapters/3.Chapter/figures/numerical_results/DATA/results_poisson_sin5_residualS33_3refinements_endtoend/fem/test_FEM_33nodes.csv}
\pgfplotsset{A/.append style={
	    xmin=-0.1,
	    xmax=1.1,
	    ymin=-1.5,
	    ymax=1.5,
	    ylabel = {$u(x)$},
	    height=0.6*\textwidth,
	    width=\textwidth,
		xtick={0,0.25,...,1},
            }
}
\begin{axis}[A]
\addplot[domain = 0:1, samples = 200, smooth, line width=2, blue!30] {sin(deg(31.4159265*x))};
\addplot[color=red!30!white, line width=1., mark=*, mark options={scale=.8, solid}] table[x expr=\thisrow{x},y=u_FEM0]{\fileFEM};
\addplot[color=black, style=dashed, line width=0.7 ,mark=asterisk, mark options = {scale=1, solid}] table[x expr=\thisrow{x},y=u_pred0]{\fileDL};
\end{axis}
\end{tikzpicture}
\caption{Step $s=1$ (before training).}
\label{fig:step1_start_poisson_sin5}
\end{subfigure} \hfill
\begin{subfigure}[b]{0.49\textwidth}
\centering
\begin{tikzpicture}
\newcommand{\fileDL}{Chapters/3.Chapter/figures/numerical_results/DATA/results_poisson_sin5_residualS33_3refinements_endtoend/end/test_prediction_33nodes.csv}
\newcommand{\fileFEM}{Chapters/3.Chapter/figures/numerical_results/DATA/results_poisson_sin5_residualS33_3refinements_endtoend/fem/test_FEM_33nodes.csv}
\pgfplotsset{A/.append style={
	    xmin=-0.1,
	    xmax=1.1,
	    ymin=-1.5,
	    ymax=1.5,
	    height=0.6*\textwidth,
	    width=\textwidth,
		xtick={0,0.25,...,1},
            }
}
\begin{axis}[A]
\addplot[domain = 0:1, samples = 200, smooth, line width=2, blue!30] {sin(deg(31.4159265*x))};
\addplot[color=red!30!white, line width=1., mark=*, mark options={scale=.8, solid}] table[x expr=\thisrow{x},y=u_FEM0]{\fileFEM};
\addplot[color=black, style=dashed, line width=0.7 ,mark=asterisk, mark options = {scale=1, solid}] table[x expr=\thisrow{x},y=u_pred0]{\fileDL};
\end{axis}
\end{tikzpicture}
\caption{Step $s=1$ (after training).}
\label{fig:step1_end_poisson_sin5}
\end{subfigure}
\vskip 0.5em
\begin{subfigure}[b]{0.49\textwidth}
\centering
\begin{tikzpicture}
\newcommand{\fileDL}{Chapters/3.Chapter/figures/numerical_results/DATA/results_poisson_sin5_residualS33_3refinements_endtoend/start/test_prediction_65nodes.csv}
\newcommand{\fileFEM}{Chapters/3.Chapter/figures/numerical_results/DATA/results_poisson_sin5_residualS33_3refinements_endtoend/fem/test_FEM_65nodes.csv}
\pgfplotsset{A/.append style={
	    xmin=-0.1,
	    xmax=1.1,
	    ymin=-1.5,
	    ymax=1.5,
	    ylabel = {$u(x)$},
	    height=0.6*\textwidth,
	    width=\textwidth,
		xtick={0,0.25,...,1},
        }
}
\begin{axis}[A]
\addplot[domain = 0:1, samples = 200, smooth, line width=2, blue!30] {sin(deg(31.4159265*x))};
\addplot[color=red!30!white, line width=1., mark=*, mark options={scale=.8, solid}] table[x expr=\thisrow{x},y=u_FEM0]{\fileFEM};
\addplot[color=black, style=dashed, line width=0.7 ,mark=asterisk, mark options = {scale=1, solid}] table[x expr=\thisrow{x},y=u_pred0]{\fileDL};
\end{axis}
\end{tikzpicture}
\caption{Step $s=2$ (before training).}
\label{fig:step2_start_poisson_sin5}
\end{subfigure} \hfill
\begin{subfigure}[b]{0.49\textwidth}
\centering
\begin{tikzpicture}
\newcommand{\fileDL}{Chapters/3.Chapter/figures/numerical_results/DATA/results_poisson_sin5_residualS33_3refinements_endtoend/end/test_prediction_65nodes.csv}
\newcommand{\fileFEM}{Chapters/3.Chapter/figures/numerical_results/DATA/results_poisson_sin5_residualS33_3refinements_endtoend/fem/test_FEM_65nodes.csv}
\pgfplotsset{A/.append style={
	    xmin=-0.1,
	    xmax=1.1,
	    ymin=-1.5,
	    ymax=1.5,
	    height=0.6*\textwidth,
	    width=\textwidth,
		xtick={0,0.25,...,1},
            }
}
\begin{axis}[A]
\addplot[domain = 0:1, samples = 200, smooth, line width=2, blue!30] {sin(deg(31.4159265*x))};
\addplot[color=red!30!white, line width=1., mark=*, mark options={scale=.8, solid}] table[x expr=\thisrow{x},y=u_FEM0]{\fileFEM};
\addplot[color=black, style=dashed, line width=0.7 ,mark=asterisk, mark options = {scale=1, solid}] table[x expr=\thisrow{x},y=u_pred0]{\fileDL};
\end{axis}
\end{tikzpicture}
\caption{Step $s=2$ (after training).}
\label{fig:step2_end_poisson_sin5}
\end{subfigure}
\vskip 0.5em
\begin{subfigure}[b]{0.49\textwidth}
\centering
\begin{tikzpicture}
\newcommand{\fileDL}{Chapters/3.Chapter/figures/numerical_results/DATA/results_poisson_sin5_residualS33_3refinements_endtoend/start/test_prediction_129nodes.csv}
\newcommand{\fileFEM}{Chapters/3.Chapter/figures/numerical_results/DATA/results_poisson_sin5_residualS33_3refinements_endtoend/fem/test_FEM_129nodes.csv}
\pgfplotsset{A/.append style={
	    xmin=-0.1,
	    xmax=1.1,
	    ymin=-1.5,
	    ymax=1.5,
	    ylabel = {$u(x)$},
	    height=0.6*\textwidth,
	    width=\textwidth,
		xtick={0,0.25,...,1},
            }
}
\begin{axis}[A]
\addplot[domain = 0:1, samples = 200, smooth, line width=2, blue!30] {sin(deg(31.4159265*x))};
\addplot[color=red!30!white, line width=1., mark=*, mark options={scale=.8, solid}] table[x expr=\thisrow{x},y=u_FEM0]{\fileFEM};
\addplot[color=black, style=dashed, line width=0.7 ,mark=asterisk, mark options = {scale=1, solid}] table[x expr=\thisrow{x},y=u_pred0]{\fileDL};
\end{axis}
\end{tikzpicture}
\caption{Step $s=3$ (before training).}
\label{fig:step3_start_poisson_sin5}
\end{subfigure}\hfill
\begin{subfigure}[b]{0.49\textwidth}
\centering
\begin{tikzpicture}
\newcommand{\fileDL}{Chapters/3.Chapter/figures/numerical_results/DATA/results_poisson_sin5_residualS33_3refinements_endtoend/end/test_prediction_129nodes.csv}
\newcommand{\fileFEM}{Chapters/3.Chapter/figures/numerical_results/DATA/results_poisson_sin5_residualS33_3refinements_endtoend/fem/test_FEM_129nodes.csv}
\pgfplotsset{A/.append style={
	    xmin=-0.1,
	    xmax=1.1,
	    ymin=-1.5,
	    ymax=1.5,
	    height=0.6*\textwidth,
	    width=\textwidth,
		xtick={0,0.25,...,1},
            }
}
\begin{axis}[A]
\addplot[domain = 0:1, samples = 200, smooth, line width=2, blue!30] {sin(deg(31.4159265*x))};
\addplot[color=red!30!white, line width=1., mark=*, mark options={scale=.8, solid}] table[x expr=\thisrow{x},y=u_FEM0]{\fileFEM};
\addplot[color=black, style=dashed, line width=0.7 ,mark=asterisk, mark options = {scale=1, solid}] table[x expr=\thisrow{x},y=u_pred0]{\fileDL};
\end{axis}
\end{tikzpicture}
\caption{Step $s=3$ (after training).}
\label{fig:step3_end_poisson_sin5}
\end{subfigure}
\vskip 0.5em
\begin{subfigure}[b]{0.49\textwidth}
\centering
\begin{tikzpicture}
\newcommand{\fileDL}{Chapters/3.Chapter/figures/numerical_results/DATA/results_poisson_sin5_residualS33_3refinements_endtoend/start/test_prediction_257nodes.csv}
\newcommand{\fileFEM}{Chapters/3.Chapter/figures/numerical_results/DATA/results_poisson_sin5_residualS33_3refinements_endtoend/fem/test_FEM_257nodes.csv}
\pgfplotsset{A/.append style={
	    xmin=-0.1,
	    xmax=1.1,
	    xlabel = {$x$},
	    ymin=-1.5,
	    ymax=1.5,
	    ylabel = {$u(x)$},
	    height=0.6*\textwidth,
	    width=\textwidth,
		xtick={0,0.25,...,1},
            }
}
\begin{axis}[A]
\addplot[domain = 0:1, samples = 200, smooth, line width=2, blue!30] {sin(deg(31.4159265*x))};
\addplot[color=red!30!white, line width=1., mark=*, mark options={scale=.8, solid}] table[x expr=\thisrow{x},y=u_FEM0]{\fileFEM};
\addplot[color=black, style=dashed, line width=0.7 ,mark=asterisk, mark options = {scale=1, solid}] table[x expr=\thisrow{x},y=u_pred0]{\fileDL};
\end{axis}
\end{tikzpicture}
\caption{Step $s=4$ (before training).}
\label{fig:step4_start_poisson_sin5}
\end{subfigure}\hfill
\begin{subfigure}[b]{0.49\textwidth}
\centering
\begin{tikzpicture}
\newcommand{\fileDL}{Chapters/3.Chapter/figures/numerical_results/DATA/results_poisson_sin5_residualS33_3refinements_endtoend/end/test_prediction_257nodes.csv}
\newcommand{\fileFEM}{Chapters/3.Chapter/figures/numerical_results/DATA/results_poisson_sin5_residualS33_3refinements_endtoend/fem/test_FEM_257nodes.csv}
\pgfplotsset{A/.append style={
	    xmin=-0.1,
	    xmax=1.1,
	    xlabel = {$x$},
	    ymin=-1.5,
	    ymax=1.5,
	    height=0.6*\textwidth,
	    width=\textwidth,
		xtick={0,0.25,...,1},
        }
}
\begin{axis}[A]
\addplot[domain = 0:1, samples = 200, smooth, line width=2, blue!30] {sin(deg(31.4159265*x))};
\addplot[color=red!30!white, line width=1., mark=*, mark options={scale=.8, solid}] table[x expr=\thisrow{x},y=u_FEM0]{\fileFEM};
\addplot[color=black, style=dashed, line width=0.7 ,mark=asterisk, mark options = {scale=1, solid}] table[x expr=\thisrow{x},y=u_pred0]{\fileDL};
\end{axis}
\end{tikzpicture}
\caption{Step $s=4$ (after training).}
\label{fig:step4_end_poisson_sin5}
\end{subfigure}
\caption{Four steps for the \acs{DeepFEM} in Poisson problem \eqref{eq:poisson_helmholtz_sin5_3refinements}. $u^*$ is the exact solution, $u_{\text{FEM}}$ is the finite element solution, and $u_{\text{NN}}$ is the \acs{DeepFEM} prediction.}
\label{fig:poisson_sin5}
\end{figure}

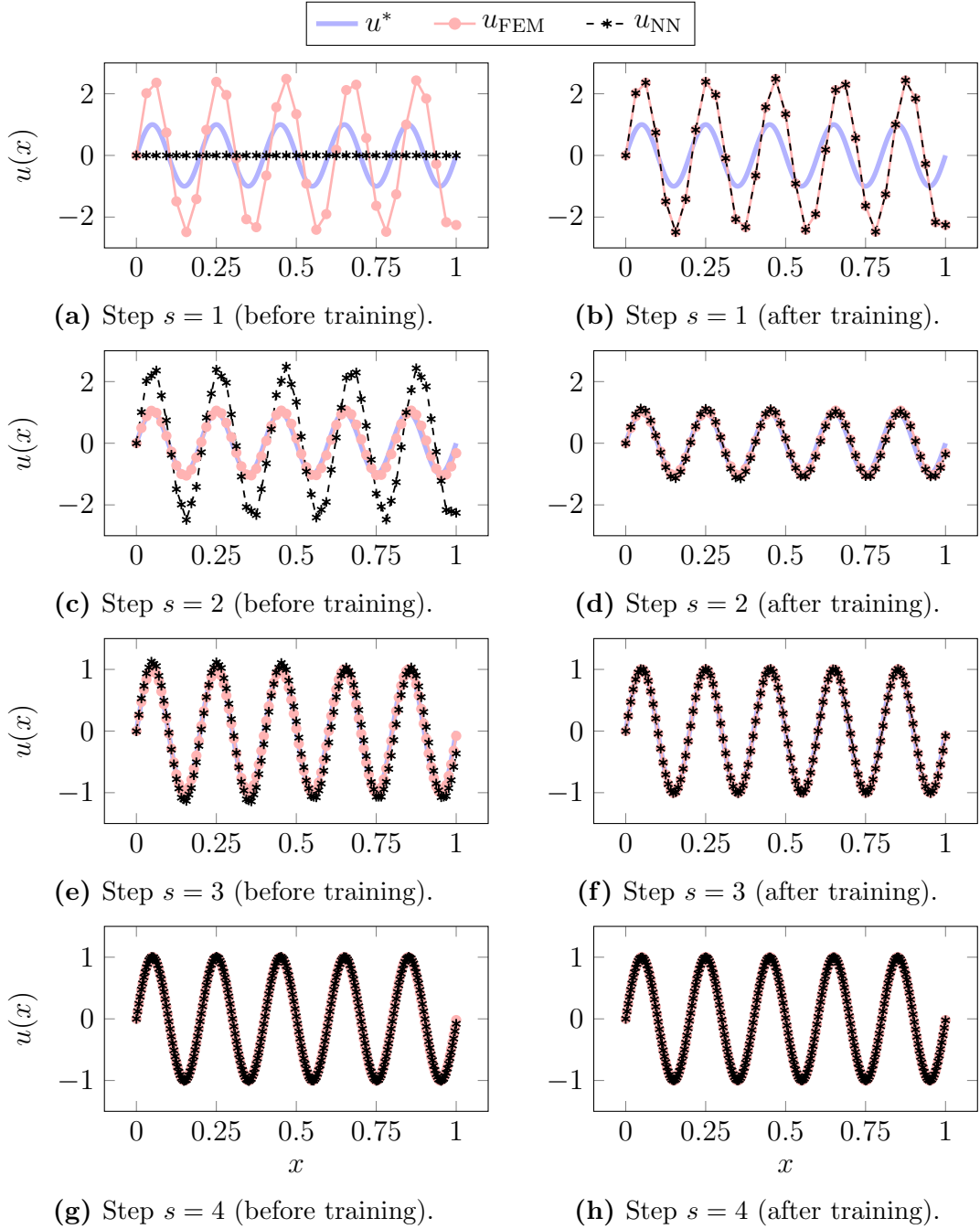
\begin{figure}[htbp]
\centering
\begin{subfigure}[b]{\textwidth}
\centering
 \begin{tikzpicture}
\pgfplotsset{A/.append style={
		hide axis,
	    xmin=-0.1,
	    xmax=1.1,
	    ymin=-0.1,
	    ymax=1.1,
	 	legend style={fill=none, anchor=north, /tikz/every even column/.append style={column sep=0.5cm}},
	 	legend columns = -1
            }
}
\begin{axis}[A]
\addlegendimage{color=blue!30!white, line width=2.};
\addlegendentry{$u^*$};
\addlegendimage{color=red!30!white, line width=1., mark=*, mark options={scale=1, solid}};
\addlegendentry{$u_{\text{FEM}}$};
\addlegendimage{color=black, style=dashed, line width=.7 ,mark=asterisk, mark options = {scale=1, solid}};
\addlegendentry{$u_{\text{NN}}$};
\end{axis}
\end{tikzpicture}
\end{subfigure}
\vskip 0.5em
\begin{subfigure}[b]{0.49\textwidth}
\centering
\begin{tikzpicture}
\newcommand{\fileDL}{Chapters/3.Chapter/figures/numerical_results/DATA/results_helmholtz_sin5_residualS33_3refinements_endtoend/start/test_prediction_33nodes.csv}
\newcommand{\fileFEM}{Chapters/3.Chapter/figures/numerical_results/DATA/results_helmholtz_sin5_residualS33_3refinements_endtoend/fem/test_FEM_33nodes.csv}
\pgfplotsset{A/.append style={
	    xmin=-0.1,
	    xmax=1.1,
	    ymin=-3,
	    ymax=3,
	    ylabel = {$u(x)$},
	    height=0.6*\textwidth,
	    width=\textwidth,
		xtick={0,0.25,...,1},
        }
}
\begin{axis}[A]
\addplot[domain = 0:1, samples = 200, smooth, line width=2, blue!30] {sin(deg(31.4159265*x))};
\addplot[color=red!30!white, line width=1., mark=*, mark options={scale=.8, solid}] table[x expr=\thisrow{x},y=u_FEM0]{\fileFEM};
\addplot[color=black, style=dashed, line width=0.7 ,mark=asterisk, mark options = {scale=1, solid}] table[x expr=\thisrow{x},y=u_pred0]{\fileDL};
\end{axis}
\end{tikzpicture}
\caption{Step $s=1$ (before training).}
\label{fig:step1_start_helmholtz_sin5}
\end{subfigure}\hfill
\begin{subfigure}[b]{0.49\textwidth}
\centering
\begin{tikzpicture}
\newcommand{\fileDL}{Chapters/3.Chapter/figures/numerical_results/DATA/results_helmholtz_sin5_residualS33_3refinements_endtoend/end/test_prediction_33nodes.csv}
\newcommand{\fileFEM}{Chapters/3.Chapter/figures/numerical_results/DATA/results_helmholtz_sin5_residualS33_3refinements_endtoend/fem/test_FEM_33nodes.csv}
\pgfplotsset{A/.append style={
	    xmin=-0.1,
	    xmax=1.1,
	    ymin=-3,
	    ymax=3,
	    height=0.6*\textwidth,
	    width=\textwidth,
		xtick={0,0.25,...,1},
        }
}
\begin{axis}[A]
\addplot[domain = 0:1, samples = 200, smooth, line width=2, blue!30] {sin(deg(31.4159265*x))};
\addplot[color=red!30!white, line width=1., mark=*, mark options={scale=.8, solid}] table[x expr=\thisrow{x},y=u_FEM0]{\fileFEM};
\addplot[color=black, style=dashed, line width=0.7 ,mark=asterisk, mark options = {scale=1, solid}] table[x expr=\thisrow{x},y=u_pred0]{\fileDL};
\end{axis}
\end{tikzpicture}
\caption{Step $s=1$ (after training).}
\label{fig:step1_end_helmholtz_sin5}
\end{subfigure}
\vskip 0.5em   
\begin{subfigure}[b]{0.49\textwidth}
\centering
\begin{tikzpicture}
\newcommand{\fileDL}{Chapters/3.Chapter/figures/numerical_results/DATA/results_helmholtz_sin5_residualS33_3refinements_endtoend/start/test_prediction_65nodes.csv}
\newcommand{\fileFEM}{Chapters/3.Chapter/figures/numerical_results/DATA/results_helmholtz_sin5_residualS33_3refinements_endtoend/fem/test_FEM_65nodes.csv}
\pgfplotsset{A/.append style={
	    xmin=-0.1,
	    xmax=1.1,
	    ymin=-3,
	    ymax=3,
	    ylabel = {$u(x)$},
	    height=0.6*\textwidth,
	    width=\textwidth,
		xtick={0,0.25,...,1},
        }
}
\begin{axis}[A]
\addplot[domain = 0:1, samples = 200, smooth, line width=2, blue!30] {sin(deg(31.4159265*x))};
\addplot[color=red!30!white, line width=1., mark=*, mark options={scale=.8, solid}] table[x expr=\thisrow{x},y=u_FEM0]{\fileFEM};
\addplot[color=black, style=dashed, line width=0.7 ,mark=asterisk, mark options = {scale=1, solid}] table[x expr=\thisrow{x},y=u_pred0]{\fileDL};
\end{axis}
\end{tikzpicture}
\caption{Step $s=2$ (before training).}
\label{fig:step2_start_helmholtz_sin5}
\end{subfigure}\hfill
\begin{subfigure}[b]{0.49\textwidth}
\centering
\begin{tikzpicture}
\newcommand{\fileDL}{Chapters/3.Chapter/figures/numerical_results/DATA/results_helmholtz_sin5_residualS33_3refinements_endtoend/end/test_prediction_65nodes.csv}
\newcommand{\fileFEM}{Chapters/3.Chapter/figures/numerical_results/DATA/results_helmholtz_sin5_residualS33_3refinements_endtoend/fem/test_FEM_65nodes.csv}
\pgfplotsset{A/.append style={
	    xmin=-0.1,
	    xmax=1.1,
	    ymin=-3,
	    ymax=3,
	    height=0.6*\textwidth,
	    width=\textwidth,
		xtick={0,0.25,...,1},
        }
}
\begin{axis}[A]
\addplot[domain = 0:1, samples = 200, smooth, line width=2, blue!30] {sin(deg(31.4159265*x))};
\addplot[color=red!30!white, line width=1., mark=*, mark options={scale=.8, solid}] table[x expr=\thisrow{x},y=u_FEM0]{\fileFEM};
\addplot[color=black, style=dashed, line width=0.7 ,mark=asterisk, mark options = {scale=1, solid}] table[x expr=\thisrow{x},y=u_pred0]{\fileDL};
\end{axis}
\end{tikzpicture}
\caption{Step $s=2$ (after training).}
\label{fig:step2_end_helmholtz_sin5}
\end{subfigure}
\vskip 0.5em
\begin{subfigure}[b]{0.49\textwidth}
\centering
\begin{tikzpicture}
\newcommand{\fileDL}{Chapters/3.Chapter/figures/numerical_results/DATA/results_helmholtz_sin5_residualS33_3refinements_endtoend/start/test_prediction_129nodes.csv}
\newcommand{\fileFEM}{Chapters/3.Chapter/figures/numerical_results/DATA/results_helmholtz_sin5_residualS33_3refinements_endtoend/fem/test_FEM_129nodes.csv}
\pgfplotsset{A/.append style={
	    xmin=-0.1,
	    xmax=1.1,
	    ymin=-1.5,
	    ymax=1.5,
	    ylabel = {$u(x)$},
	    height=0.6*\textwidth,
	    width=\textwidth,
		xtick={0,0.25,...,1},
		}
}
\begin{axis}[A]
\addplot[domain = 0:1, samples = 200, smooth, line width=2, blue!30] {sin(deg(31.4159265*x))};
\addplot[color=red!30!white, line width=1., mark=*, mark options={scale=.8, solid}] table[x expr=\thisrow{x},y=u_FEM0]{\fileFEM};
\addplot[color=black, style=dashed, line width=0.7 ,mark=asterisk, mark options = {scale=1, solid}] table[x expr=\thisrow{x},y=u_pred0]{\fileDL};
\end{axis}
\end{tikzpicture}
\caption{Step $s=3$ (before training).}
\label{fig:step3_start_helmholtz_sin5}
\end{subfigure}\hfill
\begin{subfigure}[b]{0.49\textwidth}
\centering
\begin{tikzpicture}
\newcommand{\fileDL}{Chapters/3.Chapter/figures/numerical_results/DATA/results_helmholtz_sin5_residualS33_3refinements_endtoend/end/test_prediction_129nodes.csv}
\newcommand{\fileFEM}{Chapters/3.Chapter/figures/numerical_results/DATA/results_helmholtz_sin5_residualS33_3refinements_endtoend/fem/test_FEM_129nodes.csv}
\pgfplotsset{A/.append style={
	    xmin=-0.1,
	    xmax=1.1,
	    ymin=-1.5,
	    ymax=1.5,
	    height=0.6*\textwidth,
	    width=\textwidth,
		xtick={0,0.25,...,1},
        }
}
\begin{axis}[A]
\addplot[domain = 0:1, samples = 200, smooth, line width=2, blue!30] {sin(deg(31.4159265*x))};
\addplot[color=red!30!white, line width=1., mark=*, mark options={scale=.8, solid}] table[x expr=\thisrow{x},y=u_FEM0]{\fileFEM};
\addplot[color=black, style=dashed, line width=0.7 ,mark=asterisk, mark options = {scale=1, solid}] table[x expr=\thisrow{x},y=u_pred0]{\fileDL};
\end{axis}
\end{tikzpicture}
\caption{Step $s=3$ (after training).}
\label{fig:step3_end_helmholtz_sin5}
\end{subfigure}
\vskip 0.5em
\begin{subfigure}[b]{0.49\textwidth}
\centering
\begin{tikzpicture}
\newcommand{\fileDL}{Chapters/3.Chapter/figures/numerical_results/DATA/results_helmholtz_sin5_residualS33_3refinements_endtoend/start/test_prediction_257nodes.csv}
\newcommand{\fileFEM}{Chapters/3.Chapter/figures/numerical_results/DATA/results_helmholtz_sin5_residualS33_3refinements_endtoend/fem/test_FEM_257nodes.csv}
\pgfplotsset{A/.append style={
	    xmin=-0.1,
	    xmax=1.1,
	    xlabel = {$x$},
	    ymin=-1.5,
	    ymax=1.5,
	    ylabel = {$u(x)$},
	    height=0.6*\textwidth,
	    width=\textwidth,
		xtick={0,0.25,...,1},
        }
}
\begin{axis}[A]
\addplot[domain = 0:1, samples = 200, smooth, line width=2, blue!30] {sin(deg(31.4159265*x))};
\addplot[color=red!30!white, line width=1., mark=*, mark options={scale=.8, solid}] table[x expr=\thisrow{x},y=u_FEM0]{\fileFEM};
\addplot[color=black, style=dashed, line width=0.7 ,mark=asterisk, mark options = {scale=1, solid}] table[x expr=\thisrow{x},y=u_pred0]{\fileDL};
\end{axis}
\end{tikzpicture}
\caption{Step $s=4$ (before training).}
\label{fig:step4_start_helmholtz_sin5}
\end{subfigure}\hfill
\begin{subfigure}[b]{0.49\textwidth}
\centering
\begin{tikzpicture}
\newcommand{\fileDL}{Chapters/3.Chapter/figures/numerical_results/DATA/results_helmholtz_sin5_residualS33_3refinements_endtoend/end/test_prediction_257nodes.csv}
\newcommand{\fileFEM}{Chapters/3.Chapter/figures/numerical_results/DATA/results_helmholtz_sin5_residualS33_3refinements_endtoend/fem/test_FEM_257nodes.csv}
\pgfplotsset{A/.append style={
	    xmin=-0.1,
	    xmax=1.1,
	    xlabel = {$x$},
	    ymin=-1.5,
	    ymax=1.5,
	    height=0.6*\textwidth,
	    width=\textwidth,
		xtick={0,0.25,...,1},
        }
}
\begin{axis}[A]
\addplot[domain = 0:1, samples = 200, smooth, line width=2, blue!30] {sin(deg(31.4159265*x))};
\addplot[color=red!30!white, line width=1., mark=*, mark options={scale=.8, solid}] table[x expr=\thisrow{x},y=u_FEM0]{\fileFEM};
\addplot[color=black, style=dashed, line width=0.7 ,mark=asterisk, mark options = {scale=1, solid}] table[x expr=\thisrow{x},y=u_pred0]{\fileDL};
\end{axis}
\end{tikzpicture}
\caption{Step $s=4$ (after training).}
\label{fig:step4_end_helmholtz_sin5}
\end{subfigure} 
\caption{Four steps of the \acs{DeepFEM} in Helmholtz problem \eqref{eq:poisson_helmholtz_sin5_3refinements}. $u^*$ is the exact solution, $u_{\text{FEM}}$ is the finite element solution, and $u_{\text{NN}}$ is the \acs{DeepFEM} prediction.}
\label{fig:helmholtz_sin5}
\end{figure}

\Cref{fig:poisson_sin5} and \Cref{fig:helmholtz_sin5} show the NN predictions of Poisson and Helmholtz problems, respectively.  In Poisson,  FEM solutions coincide with $u^*$ at the nodal points, but this is not true for Helmholtz's equation.

\begin{figure}[htbp]
\centering
\begin{subfigure}[b]{\textwidth}
\centering
 \begin{tikzpicture}
\pgfplotsset{A/.append style={
		hide axis,
	    xmin=-0.1,
	    xmax=1.1,
	    ymin=-0.1,
	    ymax=1.1,
	 	legend style={fill=none, anchor=north, /tikz/every even column/.append style={column sep=0.25cm}},
	 	legend columns = 1
            }
}
\begin{axis}[A]
\addlegendimage{color=blue!50, line width=1.5};
\addlegendentry{loss = $\lVert\mathbf{P}(\mathbf{A}\mathbf{u}_\text{NN}-\mathbf{f})\rVert_{\mathbf{K}+\mathbf{M}}$};
\addlegendimage{color=red!, line width=1, dashed};
\addlegendentry{$H^1$-norm error = $\lVert\mathbf{u}_\text{NN}-\mathbf{u}_\text{FEM}\rVert_{\mathbf{K}+\mathbf{M}}$};
\end{axis}
\end{tikzpicture}
\end{subfigure}
\begin{subfigure}[b]{\textwidth}
\begin{tikzpicture}
\pgfplotsset{loss1/.append style={ 
     xmax=3000,             
     xlabel = {iteration},
     xtick={0,2000},
     ymin=1e-13,     
     ymax=10e2,      
     ylabel = {loss/error},
    scaled x ticks=false,
     height=4cm,    
     width=4cm,    
     legend style={fill=none, at={(1.9,1.55)},anchor=north, minimum height = 0.6cm, /tikz/every even column/.append style={column sep=0.5cm}},
	 legend columns = -1
      }    
}
\begin{semilogyaxis}[loss1, at={(0,0)}]
\addplot[color=blue!50, line width=1.5 ,] table[x expr=\thisrow{iteration},y=loss_history]{Chapters/3.Chapter/figures/numerical_results/DATA/results_poisson_sin5_residualS33_3refinements_endtoend/loss_history_reduced2.csv};
\addplot[color=red!, line width=1 , style=dashed] table[x expr=\thisrow{iteration},y=error_history]{Chapters/3.Chapter/figures/numerical_results/DATA/results_poisson_sin5_residualS33_3refinements_endtoend/error_history_reduced2.csv};
\end{semilogyaxis}
\pgfplotsset{loss2/.append style={ 
     xmin=3000,
     xmax=13000,          
     xlabel = {iteration},
     xtick={3000,6000,...,12000},
     ymin=1e-5,     
     ymax=1,      
    scaled x ticks=false,
     height=4cm,    
     width=9cm,    
      }    
}
\begin{semilogyaxis}[loss2, at={(60,0)}]
\addplot[color=blue!50, line width=1.5 ,] table[x expr=\thisrow{iteration},y=loss_history]{Chapters/3.Chapter/figures/numerical_results/DATA/results_poisson_sin5_residualS33_3refinements_endtoend/loss_history_reduced2.csv};
\addplot[color=red!, line width=1 , style=dashed] table[x expr=\thisrow{iteration},y=error_history]{Chapters/3.Chapter/figures/numerical_results/DATA/results_poisson_sin5_residualS33_3refinements_endtoend/error_history_reduced2.csv};
\end{semilogyaxis}
\end{tikzpicture}
\caption{Loss function evolution along with the $H^1$-norm error for Poisson problem \eqref{eq:poisson_helmholtz_sin5_3refinements}.}
\label{fig:loss_poisson_sin5_residualS33_3refinements_endtoend}
\end{subfigure}
\vskip 0.5em
\begin{subfigure}[b]{\textwidth}
\begin{tikzpicture}
\pgfplotsset{loss1/.append style={ 
     xmax=3000,             
     xlabel = {iteration},
     xtick={0,2000},
     ymin=1e-13,     
     ymax=10e2,      
     ylabel = {loss/error},
     scaled x ticks=false,
     height=4cm,    
     width=4cm,    
      }    
}
\begin{semilogyaxis}[loss1, at={(0,0)}]
\addplot[color=blue!50, line width=1.5 ,] table[x expr=\thisrow{iteration},y=loss_history]{Chapters/3.Chapter/figures/numerical_results/DATA/results_helmholtz_sin5_residualS33_3refinements_endtoend/loss_history_reduced2.csv};
\addplot[color=red!, line width=1 , style=dashed] table[x expr=\thisrow{iteration},y=error_history]{Chapters/3.Chapter/figures/numerical_results/DATA/results_helmholtz_sin5_residualS33_3refinements_endtoend/error_history_reduced2.csv};
\end{semilogyaxis}
\pgfplotsset{loss2/.append style={ 
     xmin=3000,
     xmax=13000,          
     xlabel = {iteration},
     xtick={3000,6000,...,12000},
     ymin=1e-4,     
     ymax=1e2,
     ytick={1e2, 1e0, 1e-2},
     scaled x ticks=false,
     height=4cm,    
     width=9cm,    
      }    
}
\begin{semilogyaxis}[loss2, at={(60,0)}]
\addplot[color=blue!50, line width=1.5 ,] table[x expr=\thisrow{iteration},y=loss_history]{Chapters/3.Chapter/figures/numerical_results/DATA/results_helmholtz_sin5_residualS33_3refinements_endtoend/loss_history_reduced2.csv};
\addplot[color=red!, line width=1 , style=dashed] table[x expr=\thisrow{iteration},y=error_history]{Chapters/3.Chapter/figures/numerical_results/DATA/results_helmholtz_sin5_residualS33_3refinements_endtoend/error_history_reduced2.csv};
\end{semilogyaxis}
\end{tikzpicture}
\caption{Loss function evolution along with the $H^1$-norm error for Helmholtz problem \eqref{eq:poisson_helmholtz_sin5_3refinements}.}
\label{fig:loss_helmholtz_sin5_residualS33_3refinements_endtoend}
\end{subfigure}
\vskip 0.5em
\begin{subfigure}[t]{0.45\textwidth}
\centering
\begin{tikzpicture}
\pgfplotsset{A/.append style={
	    xmin=-0.1,
	    xmax=1.1,
	    xlabel = {$x$},
	    ylabel = {$u_\text{NN}-u_\text{FEM}$},
	    height=0.6*\textwidth,
	    width=\textwidth,
	    xtick={0,0.25,...,1},
        }
}
\begin{axis}[A]
\addplot[color=red, line width=1.3] table[x expr=\thisrow{x},y=error0]{Chapters/3.Chapter/figures/numerical_results/DATA/results_poisson_sin5_residualS33_3refinements_endtoend/error/test_error_257nodes.csv};
\end{axis}
\end{tikzpicture}
\caption{Error function at the end of the training (Poisson's problem).}\label{fig:poisson_sin5_errorfunction_residualS33_257nodes_endtoend}\end{subfigure} \hfill
\begin{subfigure}[t]{0.45\textwidth}
\centering
\begin{tikzpicture}
\pgfplotsset{A/.append style={
	    xmin=-0.1,
	    xmax=1.1,
	    xlabel = {$x$},
	    height=0.6*\textwidth,
	    width=\textwidth,
		xtick={0,0.25,...,1},
        }
}
\begin{axis}[A]
\addplot[color=red, line width=1.3] table[x expr=\thisrow{x},y=error0]{Chapters/3.Chapter/figures/numerical_results/DATA/results_helmholtz_sin5_residualS33_3refinements_endtoend/error/test_error_257nodes.csv};
\end{axis}
\end{tikzpicture}
\caption{Error function at the end of the training (Helmholtz's problem).}\label{fig:helmholtz_sin5_errorfunction_residualS33_257nodes_endtoend}
\end{subfigure}
\caption{End-to-end training of the \acs{DeepFEM} for problem \eqref{eq:poisson_helmholtz_sin5_3refinements} along three steps with $32$-size block-Jacobi preconditioners.}
    \label{fig:loss_poisson_helmholtz_sin5_residualS33_3refinements_endtoend}
\end{figure}
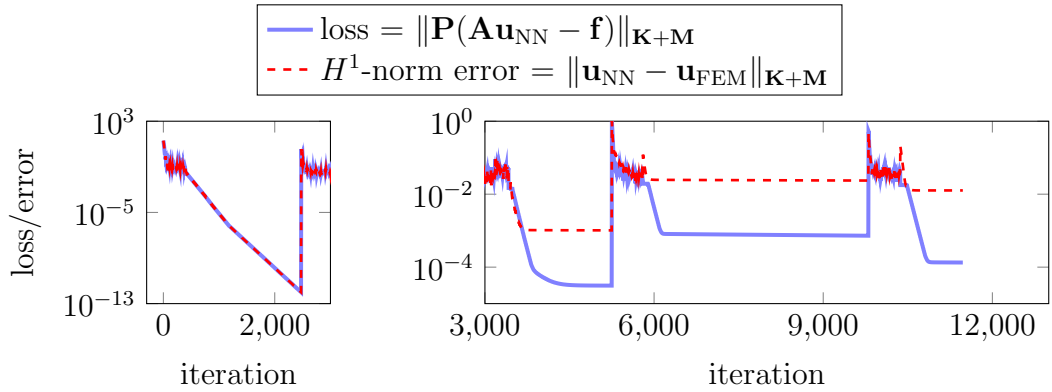
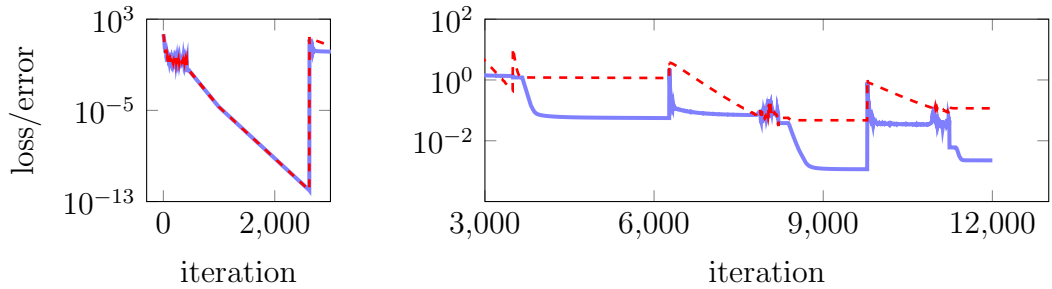
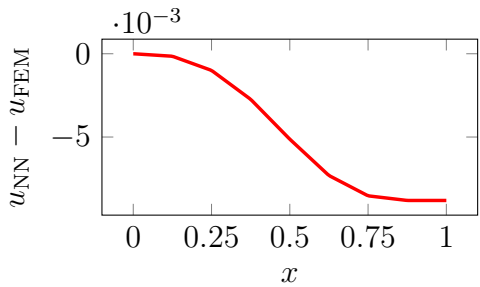
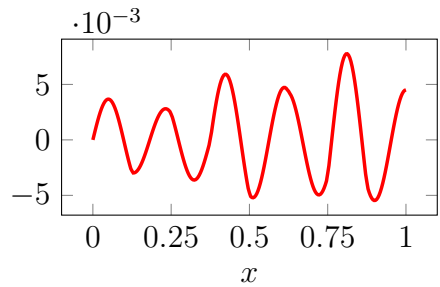

When training with the preconditioner being the inverse matrix, the NN converges perfectly with a monotonic decrease. However, when considering non-inverse preconditioning, convergence stagnates after decreasing a couple of orders of magnitude (see \Cref{fig:loss_poisson_sin5_residualS33_3refinements_endtoend} and \Cref{fig:loss_helmholtz_sin5_residualS33_3refinements_endtoend}). In the Adalr optimizer training phases, the loss decreases abruptly in the early iterations but then remains flat, producing a corresponding $H^1$-norm error reduction that is insignificant. In both experiments, we obtain errors around $10^{-3}$ (see \Cref{fig:poisson_sin5_errorfunction_residualS33_257nodes_endtoend} and \Cref{fig:helmholtz_sin5_errorfunction_residualS33_257nodes_endtoend}). These results illustrate that it is possible to obtain approximate solutions with a certain (but not high) degree of accuracy, probably due to the non-convexity of the losses---recall \Cref{section2.3}.

\subsection{Sinusoidal solution with piecewise-constant coefficients}
\label[section]{section3.6.3}

Following Helmhotz's equation, we add varying frequencies along the propagation domain by considering piecewise-constant coefficients as follows:
\begin{equation} \label{eq:helmholtz_piecewiseconstant_3refinements}
\begin{cases}
-(\sigma u')' + \alpha u = 0, &\\
u(0)=0, u'(1)=10\pi,
\end{cases}
\end{equation} with
\begin{equation}
\sigma = \begin{cases}
1, &\text{if } 0<x<1/3,\\
2, &\text{if } 1/3<x<2/3,\\
3, &\text{if } 2/3<x<1,\\
\end{cases}\qquad 
\alpha = \begin{cases}
-3000, &\text{if } 0<x<1/3,\\
-2000, &\text{if } 1/3<x<2/3,\\
-1000, &\text{if } 2/3<x<1.\\
\end{cases}
\end{equation} We solve it employing the $H^1$-norm for the loss. We start with $48$ elements and perform three uniform refinements. We select $\mathbf{P}$ with $48$-size blocks in the first two mesh refinements and $96$-size blocks in the last two. The NN architecture is the same as the one described in \Cref{section3.6.1}.

\begin{figure}[htbp]
\centering
\begin{subfigure}[b]{\textwidth}
\centering
 \begin{tikzpicture}
\pgfplotsset{A/.append style={
		hide axis,
	    xmin=-0.1,
	    xmax=1.1,
	    ymin=-0.1,
	    ymax=1.1,
	 	legend style={fill=none, anchor=north, /tikz/every even column/.append style={column sep=0.5cm}},
	 	legend columns = -1
            }
}
\begin{axis}[A]
\addlegendimage{color=red!30!white, line width=1., mark=*, mark options={scale=1, solid}};
\addlegendentry{$u_{\text{FEM}}$};
\addlegendimage{color=black, style=dashed, line width=.7 ,mark=asterisk, mark options = {scale=1, solid}};
\addlegendentry{$u_{\text{NN}}$};
\end{axis}
\end{tikzpicture}
\end{subfigure}
\vskip 0.5em
\begin{subfigure}[b]{0.49\textwidth}
\centering
\begin{tikzpicture}
\newcommand{\fileDL}{Chapters/3.Chapter/figures/numerical_results/DATA/results_helmholtz_piecewiseconstant_3refinements_endtoend/start/test_prediction_49nodes.csv}
\newcommand{\fileFEM}{Chapters/3.Chapter/figures/numerical_results/DATA/results_helmholtz_piecewiseconstant_3refinements_endtoend/fem/test_FEM_49nodes.csv}
\pgfplotsset{A/.append style={
	    xmin=-0.1,
	    xmax=1.1,
	    ymin=-1,
	    ymax=1,
	    ylabel = {$u(x)$},
	    height=0.6*\textwidth,
	    width=\textwidth,
		xtick={0,0.25,...,1},
	 	legend style={fill=none, at={(1.18,1.3)},anchor=north, /tikz/every even column/.append style={column sep=0.5cm}},
	 	legend columns = -1
            }
}
\begin{axis}[A]
\addplot[color=red!30!white, line width=1., mark=*, mark options={scale=.8, solid}] table[x expr=\thisrow{x},y=u_FEM0]{\fileFEM};
\addplot[color=black, style=dashed, line width=0.7 ,mark=asterisk, mark options = {scale=1, solid}] table[x expr=\thisrow{x},y=u_pred0]{\fileDL};
\end{axis}
\end{tikzpicture}
\caption{Step $s=1$ (before training).}
\label{fig:step1_start_helmholtz_piecewiseconstant}
\end{subfigure} \hfill
\begin{subfigure}[b]{0.49\textwidth}
\centering
\begin{tikzpicture}
\newcommand{\fileDL}{Chapters/3.Chapter/figures/numerical_results/DATA/results_helmholtz_piecewiseconstant_3refinements_endtoend/end/test_prediction_49nodes.csv}
\newcommand{\fileFEM}{Chapters/3.Chapter/figures/numerical_results/DATA/results_helmholtz_piecewiseconstant_3refinements_endtoend/fem/test_FEM_49nodes.csv}
\pgfplotsset{A/.append style={
	    xmin=-0.1,
	    xmax=1.1,
	    ymin=-1,
	    ymax=1,
	    height=0.6*\textwidth,
	    width=\textwidth,
		xtick={0,0.25,...,1},
            }
}
\begin{axis}[A]
\addplot[color=red!30!white, line width=1., mark=*, mark options={scale=.8, solid}] table[x expr=\thisrow{x},y=u_FEM0]{\fileFEM};
\addplot[color=black, style=dashed, line width=0.7 ,mark=asterisk, mark options = {scale=1, solid}] table[x expr=\thisrow{x},y=u_pred0]{\fileDL};
\end{axis}
\end{tikzpicture}
\caption{Step $s=1$ (after training).}
\label{fig:step1_end_helmholtz_piecewiseconstant}
\end{subfigure}
\vskip 0.5em
\begin{subfigure}[b]{0.49\textwidth}
\centering
\begin{tikzpicture}
\newcommand{\fileDL}{Chapters/3.Chapter/figures/numerical_results/DATA/results_helmholtz_piecewiseconstant_3refinements_endtoend/start/test_prediction_97nodes.csv}
\newcommand{\fileFEM}{Chapters/3.Chapter/figures/numerical_results/DATA/results_helmholtz_piecewiseconstant_3refinements_endtoend/fem/test_FEM_97nodes.csv}
\pgfplotsset{A/.append style={
	    xmin=-0.1,
	    xmax=1.1,
	    ymin=-1,
	    ymax=1,
	    ylabel = {$u(x)$},
	    height=0.6*\textwidth,
	    width=\textwidth,
		xtick={0,0.25,...,1},
        }
}
\begin{axis}[A]
\addplot[color=red!30!white, line width=1., mark=*, mark options={scale=.8, solid}] table[x expr=\thisrow{x},y=u_FEM0]{\fileFEM};
\addplot[color=black, style=dashed, line width=0.7 ,mark=asterisk, mark options = {scale=1, solid}] table[x expr=\thisrow{x},y=u_pred0]{\fileDL};
\end{axis}
\end{tikzpicture}
\caption{Step $s=2$ (before training).}
\label{fig:step2_start_helmholtz_piecewiseconstant}
\end{subfigure} \hfill
\begin{subfigure}[b]{0.49\textwidth}
\centering
\begin{tikzpicture}
\newcommand{\fileDL}{Chapters/3.Chapter/figures/numerical_results/DATA/results_helmholtz_piecewiseconstant_3refinements_endtoend/end/test_prediction_97nodes.csv}
\newcommand{\fileFEM}{Chapters/3.Chapter/figures/numerical_results/DATA/results_helmholtz_piecewiseconstant_3refinements_endtoend/fem/test_FEM_97nodes.csv}
\pgfplotsset{A/.append style={
	    xmin=-0.1,
	    xmax=1.1,
	    ymin=-1,
	    ymax=1,
	    height=0.6*\textwidth,
	    width=\textwidth,
		xtick={0,0.25,...,1},
            }
}
\begin{axis}[A]
\addplot[color=red!30!white, line width=1., mark=*, mark options={scale=.8, solid}] table[x expr=\thisrow{x},y=u_FEM0]{\fileFEM};
\addplot[color=black, style=dashed, line width=0.7 ,mark=asterisk, mark options = {scale=1, solid}] table[x expr=\thisrow{x},y=u_pred0]{\fileDL};
\end{axis}
\end{tikzpicture}
\caption{Step $s=2$ (after training).}
\label{fig:step2_end_helmholtz_piecewiseconstant}
\end{subfigure}
\vskip 0.5em
\begin{subfigure}[b]{0.49\textwidth}
\centering
\begin{tikzpicture}
\newcommand{\fileDL}{Chapters/3.Chapter/figures/numerical_results/DATA/results_helmholtz_piecewiseconstant_3refinements_endtoend/start/test_prediction_193nodes.csv}
\newcommand{\fileFEM}{Chapters/3.Chapter/figures/numerical_results/DATA/results_helmholtz_piecewiseconstant_3refinements_endtoend/fem/test_FEM_193nodes.csv}
\pgfplotsset{A/.append style={
	    xmin=-0.1,
	    xmax=1.1,
	    ymin=-1,
	    ymax=1,
	    ylabel = {$u(x)$},
	    height=0.6*\textwidth,
	    width=\textwidth,
		xtick={0,0.25,...,1},
            }
}
\begin{axis}[A]
\addplot[color=red!30!white, line width=1., mark=*, mark options={scale=.8, solid}] table[x expr=\thisrow{x},y=u_FEM0]{\fileFEM};
\addplot[color=black, style=dashed, line width=0.7 ,mark=asterisk, mark options = {scale=1, solid}] table[x expr=\thisrow{x},y=u_pred0]{\fileDL};
\end{axis}
\end{tikzpicture}
\caption{Step $s=3$ (before training).}
\label{fig:step3_start_helmholtz_piecewiseconstant}
\end{subfigure} \hfill
\begin{subfigure}[b]{0.49\textwidth}
\centering
\begin{tikzpicture}
\newcommand{\fileDL}{Chapters/3.Chapter/figures/numerical_results/DATA/results_helmholtz_piecewiseconstant_3refinements_endtoend/end/test_prediction_193nodes.csv}
\newcommand{\fileFEM}{Chapters/3.Chapter/figures/numerical_results/DATA/results_helmholtz_piecewiseconstant_3refinements_endtoend/fem/test_FEM_193nodes.csv}
\pgfplotsset{A/.append style={
	    xmin=-0.1,
	    xmax=1.1,
	    ymin=-1,
	    ymax=1,
	    height=0.6*\textwidth,
	    width=\textwidth,
		xtick={0,0.25,...,1},
            }
}
\begin{axis}[A]
\addplot[color=red!30!white, line width=1., mark=*, mark options={scale=.8, solid}] table[x expr=\thisrow{x},y=u_FEM0]{\fileFEM};
\addplot[color=black, style=dashed, line width=0.7 ,mark=asterisk, mark options = {scale=1, solid}] table[x expr=\thisrow{x},y=u_pred0]{\fileDL};
\end{axis}
\end{tikzpicture}
\caption{Step $s=3$ (after training).}
\label{fig:step3_end_helmholtz_piecewiseconstant}
\end{subfigure}
\vskip 0.5em    
\begin{subfigure}[b]{0.49\textwidth}
\centering
\begin{tikzpicture}
\newcommand{\fileDL}{Chapters/3.Chapter/figures/numerical_results/DATA/results_helmholtz_piecewiseconstant_3refinements_endtoend/start/test_prediction_385nodes.csv}
\newcommand{\fileFEM}{Chapters/3.Chapter/figures/numerical_results/DATA/results_helmholtz_piecewiseconstant_3refinements_endtoend/fem/test_FEM_385nodes.csv}
\pgfplotsset{A/.append style={
	    xmin=-0.1,
	    xmax=1.1,
	    xlabel={$x$},
	    ymin=-1,
	    ymax=1,
	    ylabel = {$u(x)$},
	    height=0.6*\textwidth,
	    width=\textwidth,
		xtick={0,0.25,...,1},
        }
}
\begin{axis}[A]
\addplot[color=red!30!white, line width=1., mark=*, mark options={scale=.8, solid}] table[x expr=\thisrow{x},y=u_FEM0]{\fileFEM};
\addplot[color=black, style=dashed, line width=0.7 ,mark=asterisk, mark options = {scale=1, solid}] table[x expr=\thisrow{x},y=u_pred0]{\fileDL};
\end{axis}
\end{tikzpicture}
\caption{Step $s=4$ (before training).}
\label{fig:step4_start_helmholtz_piecewiseconstant}
\end{subfigure} \hfill
\begin{subfigure}[b]{0.49\textwidth}
\centering
\begin{tikzpicture}
\newcommand{\fileDL}{Chapters/3.Chapter/figures/numerical_results/DATA/results_helmholtz_piecewiseconstant_3refinements_endtoend/end/test_prediction_385nodes.csv}
\newcommand{\fileFEM}{Chapters/3.Chapter/figures/numerical_results/DATA/results_helmholtz_piecewiseconstant_3refinements_endtoend/fem/test_FEM_385nodes.csv}
\pgfplotsset{A/.append style={
	    xmin=-0.1,
	    xmax=1.1,
	   	xlabel={$x$},
	    ymin=-1,
	    ymax=1,
	    height=0.6*\textwidth,
	    width=\textwidth,
		xtick={0,0.25,...,1},
            }
}
\begin{axis}[A]
\addplot[color=red!30!white, line width=1., mark=*, mark options={scale=.8, solid}] table[x expr=\thisrow{x},y=u_FEM0]{\fileFEM};
\addplot[color=black, style=dashed, line width=0.7 ,mark=asterisk, mark options = {scale=1, solid}] table[x expr=\thisrow{x},y=u_pred0]{\fileDL};
\end{axis}
\end{tikzpicture}
\caption{Step $s=4$ (after training).}
\label{fig:step4_end_helmholtz_piecewiseconstant}
\end{subfigure}
\caption{Four steps of the \acs{DeepFEM} for problem \eqref{eq:helmholtz_piecewiseconstant_3refinements}. $u_{\text{FEM}}$ is the finite element solution and $u_{\text{NN}}$ is the \acs{DeepFEM} prediction.}
\label{fig:helmholtz_piecewiseconstant}
\end{figure}
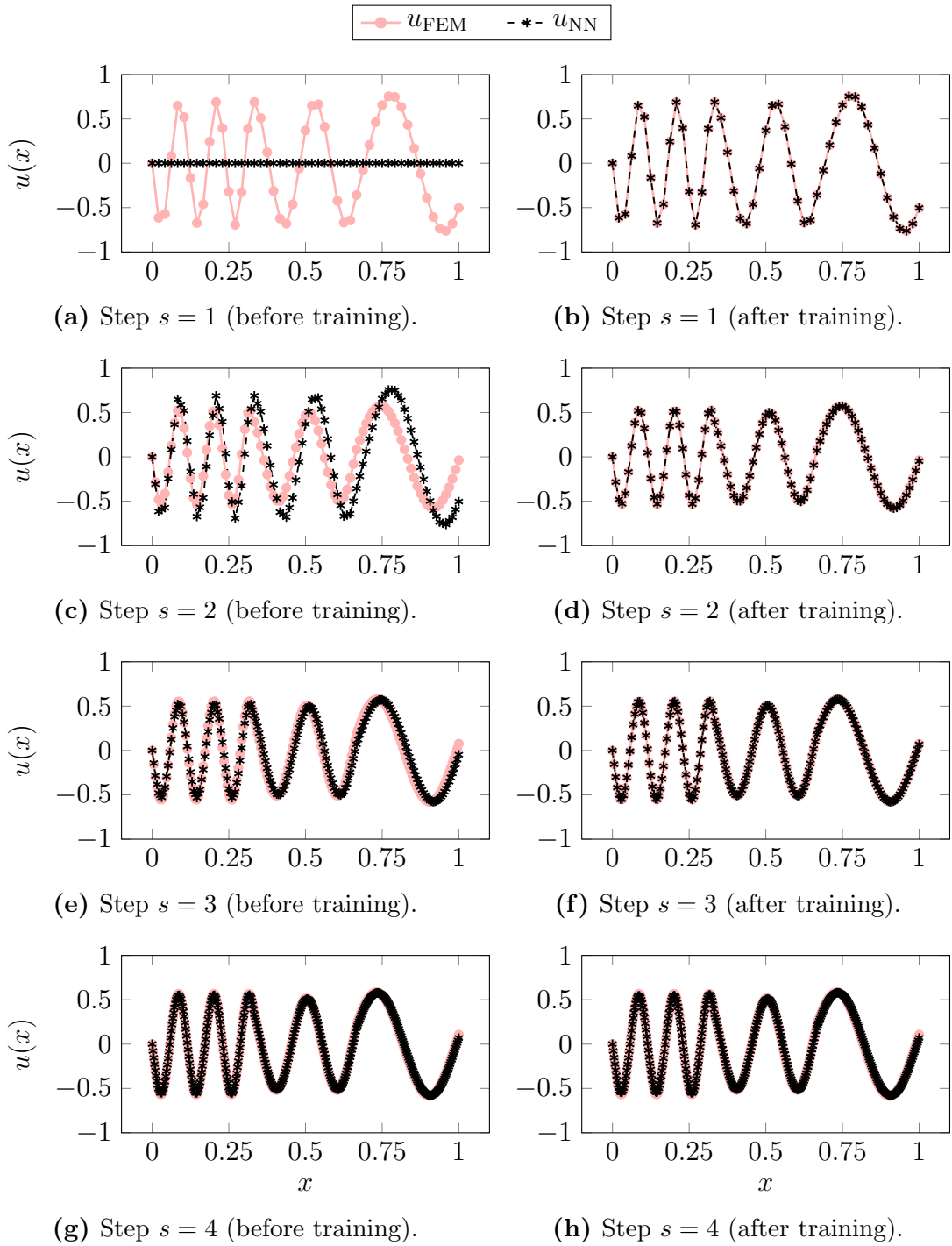

\begin{figure}[htbp]
\centering
\begin{subfigure}[b]{\textwidth}
\centering
 \begin{tikzpicture}
\pgfplotsset{A/.append style={
		hide axis,
	    xmin=-0.1,
	    xmax=1.1,
	    ymin=-0.1,
	    ymax=1.1,
	 	legend style={fill=none, anchor=north, /tikz/every even column/.append style={column sep=0.25cm}},
	 	legend columns = 1
            }
}
\begin{axis}[A]
\addlegendimage{color=blue!50, line width=1.5};
\addlegendentry{loss = $\lVert\mathbf{P}(\mathbf{A}\mathbf{u}_\text{NN}-\mathbf{f})\rVert_{\mathbf{K}+\mathbf{M}}$};
\addlegendimage{color=red!, line width=1, dashed};
\addlegendentry{$H^1$-norm error = $\lVert\mathbf{u}_\text{NN}-\mathbf{u}_\text{FEM}\rVert_{\mathbf{K}+\mathbf{M}}$};
\end{axis}
\end{tikzpicture}
\end{subfigure}
\begin{subfigure}[b]{\textwidth}
\begin{tikzpicture}
\pgfplotsset{loss1/.append style={ 
     xmax=3300,             
     xlabel = {iteration},
     xtick={0,2000},
     ymin=1e-13,     
     ymax=10e2,      
     ylabel = {loss/error},
    scaled x ticks=false,
     height=4cm,    
     width=4cm,    
     legend style={fill=none, at={(1.9,1.55)},anchor=north, minimum height = 0.6cm, /tikz/every even column/.append style={column sep=0.5cm}},
	 legend columns = -1
      }    
}
\begin{semilogyaxis}[loss1, at={(0,0)}]
\addplot[color=blue!50, line width=1.5 ,] table[x expr=\thisrow{iteration},y=loss_history]{Chapters/3.Chapter/figures/numerical_results/DATA/results_helmholtz_piecewiseconstant_3refinements_endtoend/loss_history_reduced2.csv};
\addplot[color=red!, line width=1 , style=dashed] table[x expr=\thisrow{iteration},y=error_history]{Chapters/3.Chapter/figures/numerical_results/DATA/results_helmholtz_piecewiseconstant_3refinements_endtoend/error_history_reduced2.csv};
\end{semilogyaxis}

\pgfplotsset{loss2/.append style={ 
     xmin=3000,
     xlabel = {iteration},
     xtick={4000,8000,...,12000},
     ymin=1e-3,     
     ymax=5,      
    scaled x ticks=false,
     height=4cm,    
     width=9cm,    
      }    
}
\begin{semilogyaxis}[loss2, at={(60,0)}]
\addplot[color=blue!50, line width=1.5 ,] table[x expr=\thisrow{iteration},y=loss_history]{Chapters/3.Chapter/figures/numerical_results/DATA/results_helmholtz_piecewiseconstant_3refinements_endtoend/loss_history_reduced2.csv};
\addplot[color=red!, line width=1 , style=dashed] table[x expr=\thisrow{iteration},y=error_history]{Chapters/3.Chapter/figures/numerical_results/DATA/results_helmholtz_piecewiseconstant_3refinements_endtoend/error_history_reduced2.csv};
\end{semilogyaxis}
\end{tikzpicture}
\caption{Loss function evolution along with the $H^1$-norm error for Helmholtz problem \eqref{eq:poisson_helmholtz_sin5_3refinements}}
\label{fig:loss_helmholtz_piecewiseconstant_3refinements_endtoend}
\end{subfigure}
\vskip 0.5em
\begin{subfigure}[t]{0.48\textwidth}
\centering
\begin{tikzpicture}
\pgfplotsset{A/.append style={
	    xmin=-0.1,
	    xmax=1.1,
	    xlabel = {$x$},
	    ylabel = {$u_\text{NN}-u_\text{FEM}$},
	    height=0.6*\textwidth,
	    width=\textwidth,
		xtick={0,0.25,...,1},
        }
}
\begin{axis}[A]
\addplot[color=red, line width=1.3] table[x expr=\thisrow{x},y=error0]{Chapters/3.Chapter/figures/numerical_results/DATA/results_helmholtz_piecewiseconstant_3refinements_endtoend/error/test_error_193nodes.csv};
\end{axis}
\end{tikzpicture}
\caption{Error function at the end of $s=3$.}\label{fig:helmholtz_piecewiseconstant_errorfunction_193nodes_endtoend}
\end{subfigure}
\hfill
\begin{subfigure}[t]{0.48\textwidth}
\centering
\begin{tikzpicture}
\pgfplotsset{A/.append style={
	    xmin=-0.1,
	    xmax=1.1,
	    xlabel = {$x$},
	    height=0.6*\textwidth,
	    width=\textwidth,
		xtick={0,0.25,...,1},
        }
}
\begin{axis}[A]
\addplot[color=red, line width=1.3] table[x expr=\thisrow{x},y=error0]{Chapters/3.Chapter/figures/numerical_results/DATA/results_helmholtz_piecewiseconstant_3refinements_endtoend/error/test_error_385nodes.csv};
\end{axis}
\end{tikzpicture}
\caption{Error function at the end of $s=4$.}\label{fig:helmholtz_piecewiseconstant_errorfunction_385nodes_endtoend}
\end{subfigure}
\caption{End-to-end training of the \acs{DeepFEM} for problem \eqref{eq:helmholtz_piecewiseconstant_3refinements} along four steps with block-Jacobi preconditioners of sizes $48$, $48$, $96$, and $96$ at steps $s=1,2,3,4$, respectively.}\label{fig:helmholtz_piecewiseconstant_3refinements_endtoend}
\end{figure}
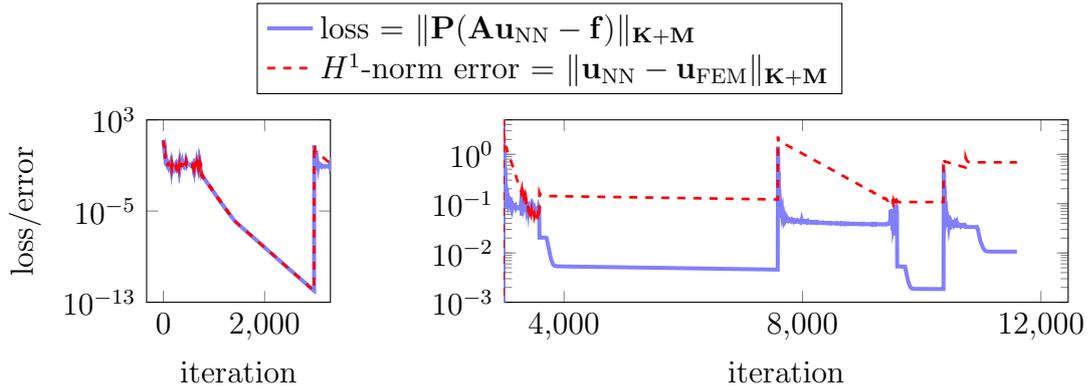
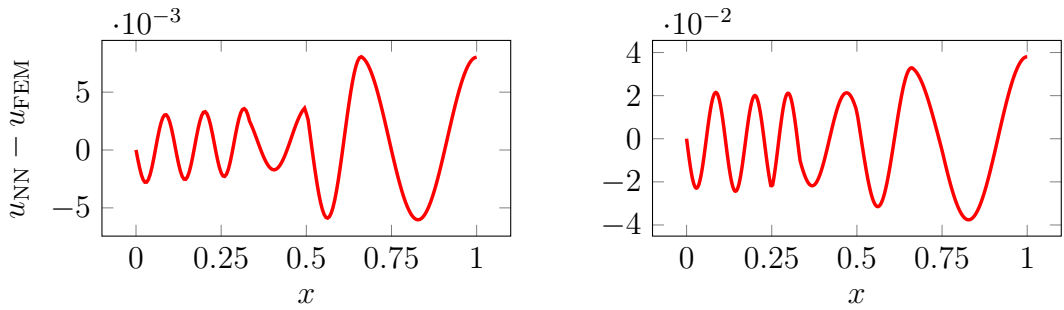

\Cref{fig:helmholtz_piecewiseconstant} shows the NN predictions and \Cref{fig:loss_helmholtz_piecewiseconstant_3refinements_endtoend} shows the loss function evolution during training. The NN nicely converges when employing the inverse as a preconditioner, but in the subsequent steps, we observe a stagnation of the loss, as was previously the case. In addition, there is a decreasing influence of the loss on the $H^1$-norm error, which remains constant throughout the last training step, while the loss decreases by two orders of magnitude. \Cref{fig:helmholtz_piecewiseconstant_errorfunction_193nodes_endtoend} and \Cref{fig:helmholtz_piecewiseconstant_errorfunction_385nodes_endtoend} show the errors at the end of steps three and four, respectively.

\newpage
\subsection{Parametric boundary value problems}
\label[section]{section3.6.4}

In this section, we naturally extend the \ac{DeepFEM} in its parametric variant: we train the \ac{NN} to learn the \ac{FEM} solutions from one mesh to another for a family of PDE coefficients with fixed boundary conditions. We consider experiments varying the $\alpha$ coefficient with a constant parametric behavior for the following \ac{BVP}:

\begin{equation} \label{eq:parmetric bvp}
\begin{cases}
-u'' + \alpha u = 0, &\\
u(0)=0, u'(1)=2\pi.
\end{cases}
\end{equation}  We solve it by performing three uniform refinements. For the training, we select databases of randomly selected samples of $\alpha$ coefficients. We train the NN on the entire database without batch partitioning. Our loss computes the norm over the preconditioned residual that depends on each sample. Specifically, we optimize with respect to an averaged sum of norms of the residual---recall \Cref{loss_residual}.

\subsubsection{Reaction-diffusion parametric equation: \boldmath{$0<\alpha<200$}}
\label[section]{section3.6.4.1}

The exact solution is $u^*(x)=C(e^{\sqrt{\alpha}x}-e^{-\sqrt{\alpha}x})$ with $C=\frac{2\pi e^{\sqrt{\alpha}}}{\sqrt{\alpha}(e^{2\sqrt{\alpha}}+1)}$ and $u^*(x) \to 2\pi x$ as $\alpha\to 0$. At each step, we establish one-layer depth and $20$-neuron width training block architecture for the NN. We start the \ac{DeepFEM} with a uniform eight-element mesh and finish with a $64$-element mesh. At each step, we employ block-Jacobi preconditioners with blocks of size eight. Since the parametric problem is symmetric and positive-definite, we consider the norm induced by the preconditioner for the optimization.

\begin{figure}[htbp]
\centering
\begin{subfigure}[b]{\textwidth}
\centering
 \begin{tikzpicture}
\pgfplotsset{A/.append style={
	    hide axis,
	    xmin=-0.1,
	    xmax=1.1,
	    ymin=-0.1,
	    ymax=1.1,
	    legend style={fill=none, anchor=north, /tikz/every even column/.append style={column sep=0.2cm}},
	    legend columns = 5
            }
}
\begin{axis}[A]
\addlegendimage{color=green!30!white, line width=1., mark=*, mark options={scale=0.8, solid}}
\addlegendentry{$u_{\text{FEM}}(0)$}
\addlegendimage{color=red!30!white, line width=1., mark=*, mark options={scale=0.8, solid}}
\addlegendentry{$u_{\text{FEM}}(3)$}
\addlegendimage{color=blue!30!white, line width=1., mark=*, mark options={scale=0.8, solid}}
\addlegendentry{$u_{\text{FEM}}(15)$}
\addlegendimage{color=orange!30!white, line width=1., mark=*, mark options={scale=0.8, solid}}
\addlegendentry{$u_{\text{FEM}}(50)$}
\addlegendimage{color=gray!30!white, line width=1., mark=*, mark options={scale=0.8, solid}}
\addlegendentry{$u_{\text{FEM}}(200)$}
\addlegendimage{color=green!80!black, style=dashed, line width=0.7 ,mark=asterisk, mark options = {scale=1, solid}}
\addlegendentry{$u_{\text{NN}}(0)$};
\addlegendimage{color=red!80!black, style=dashed, line width=0.7 ,mark=asterisk, mark options = {scale=1, solid}}
\addlegendentry{$u_{\text{NN}}(3)$};
\addlegendimage{color=blue!80!black, style=dashed, line width=0.7 ,mark=asterisk, mark options = {scale=1, solid}}
\addlegendentry{$u_{\text{NN}}(15)$};
\addlegendimage{color=orange!80!black, style=dashed, line width=0.7 ,mark=asterisk, mark options = {scale=1, solid}}
\addlegendentry{$u_{\text{NN}}(50)$};
\addlegendimage{color=gray!80!black, style=dashed, line width=0.7 ,mark=asterisk, mark options = {scale=1, solid}}
\addlegendentry{$u_{\text{NN}}(200)$};
\end{axis}
\end{tikzpicture}
\end{subfigure}
\vskip 0.5em
\begin{subfigure}[b]{0.49\textwidth}
\centering
\begin{tikzpicture}
\newcommand{\fileDL}{Chapters/3.Chapter/figures/numerical_results/DATA/results_reactiondiffusion_0-200_3refinements_endtoend/start/test_prediction_9nodes.csv}
\newcommand{\fileFEM}{Chapters/3.Chapter/figures/numerical_results/DATA/results_reactiondiffusion_0-200_3refinements_endtoend/fem/test_FEM_9nodes.csv}
\pgfplotsset{A/.append style={
	    xmin=-0.1,
	    xmax=1.1,
	    ymin=-0.1,
	    ymax=2.2,
	    ylabel = {$u(x)$},
	    height=0.59*\textwidth,
	    width=\textwidth,
	    xtick={0,0.25,...,1},
            }
}
\begin{axis}[A]
\addplot[color=green!30!white, line width=1., mark=*, mark options={scale=.8, solid}] table[x expr=\thisrow{x},y=u_FEM0]{\fileFEM};
\addplot[color=red!30!white, line width=1., mark=*, mark options={scale=.8, solid}] table[x expr=\thisrow{x},y=u_FEM1]{\fileFEM};
\addplot[color=blue!30!white, line width=1., mark=*, mark options={scale=.8, solid}] table[x expr=\thisrow{x},y=u_FEM2]{\fileFEM};
\addplot[color=orange!30!white, line width=1., mark=*, mark options={scale=.8, solid}] table[x expr=\thisrow{x},y=u_FEM3]{\fileFEM};
\addplot[color=gray!30!white, line width=1., mark=*, mark options={scale=.8, solid}] table[x expr=\thisrow{x},y=u_FEM4]{\fileFEM};
\addplot[color=green!80!black, style=dashed, line width=0.7 ,mark=asterisk, mark options = {scale=1, solid}] table[x expr=\thisrow{x},y=u_pred0]{\fileDL};
\addplot[color=red!80!black, style=dashed, line width=0.7 ,mark=asterisk, mark options = {scale=1, solid}] table[x expr=\thisrow{x},y=u_pred1]{\fileDL};
\addplot[color=blue!80!black, style=dashed, line width=0.7 ,mark=asterisk, mark options = {scale=1, solid}] table[x expr=\thisrow{x},y=u_pred2]{\fileDL};
\addplot[color=orange!80!black, style=dashed, line width=0.7 ,mark=asterisk, mark options = {scale=1, solid}] table[x expr=\thisrow{x},y=u_pred3]{\fileDL};
\addplot[color=gray!80!black, style=dashed, line width=0.7 ,mark=asterisk, mark options = {scale=1, solid}] table[x expr=\thisrow{x},y=u_pred4]{\fileDL};
\end{axis}
\end{tikzpicture}
\caption{Step $s=1$ (before training).}
\label{fig:step1_start_reactiondiffusion_parametric_exponential}
\end{subfigure} \hfill
\begin{subfigure}[b]{0.49\textwidth}
\centering
\begin{tikzpicture}
\newcommand{\fileDL}{Chapters/3.Chapter/figures/numerical_results/DATA/results_reactiondiffusion_0-200_3refinements_endtoend/end/test_prediction_9nodes.csv}
\newcommand{\fileFEM}{Chapters/3.Chapter/figures/numerical_results/DATA/results_reactiondiffusion_0-200_3refinements_endtoend/fem/test_FEM_9nodes.csv}
\pgfplotsset{A/.append style={
	    xmin=-0.1,
	    xmax=1.1,
	    ymin=-0.1,
	    ymax=2.2,
	    height=0.59*\textwidth,
	    width=\textwidth,
		xtick={0,0.25,...,1},
	 	legend style={fill=none, at={(1.3,1.6)},anchor=north, /tikz/every even column/.append style={column sep=0.5cm}},
	 	legend columns = 3
            }
}
\begin{axis}[A]
\addplot[color=green!30!white, line width=1., mark=*, mark options={scale=.8, solid}] table[x expr=\thisrow{x},y=u_FEM0]{\fileFEM};
\addplot[color=red!30!white, line width=1., mark=*, mark options={scale=.8, solid}] table[x expr=\thisrow{x},y=u_FEM1]{\fileFEM};
\addplot[color=blue!30!white, line width=1., mark=*, mark options={scale=.8, solid}] table[x expr=\thisrow{x},y=u_FEM2]{\fileFEM};
\addplot[color=orange!30!white, line width=1., mark=*, mark options={scale=.8, solid}] table[x expr=\thisrow{x},y=u_FEM3]{\fileFEM};
\addplot[color=gray!30!white, line width=1., mark=*, mark options={scale=.8, solid}] table[x expr=\thisrow{x},y=u_FEM4]{\fileFEM};
\addplot[color=green!80!black, style=dashed, line width=0.7 ,mark=asterisk, mark options = {scale=1, solid}] table[x expr=\thisrow{x},y=u_pred0]{\fileDL};
\addplot[color=red!80!black, style=dashed, line width=0.7 ,mark=asterisk, mark options = {scale=1, solid}] table[x expr=\thisrow{x},y=u_pred1]{\fileDL};
\addplot[color=blue!80!black, style=dashed, line width=0.7 ,mark=asterisk, mark options = {scale=1, solid}] table[x expr=\thisrow{x},y=u_pred2]{\fileDL};
\addplot[color=orange!80!black, style=dashed, line width=0.7 ,mark=asterisk, mark options = {scale=1, solid}] table[x expr=\thisrow{x},y=u_pred3]{\fileDL};
\addplot[color=gray!80!black, style=dashed, line width=0.7 ,mark=asterisk, mark options = {scale=1, solid}] table[x expr=\thisrow{x},y=u_pred4]{\fileDL};
\end{axis}
\end{tikzpicture}
\caption{Step $s=1$ (after training).}
\label{fig:step1_end_reactiondiffusion_parametric_exponential}
\end{subfigure}
\vskip 0.5em
\begin{subfigure}[b]{0.49\textwidth}
\centering
\begin{tikzpicture}
\newcommand{\fileDL}{Chapters/3.Chapter/figures/numerical_results/DATA/results_reactiondiffusion_0-200_3refinements_endtoend/start/test_prediction_17nodes.csv}
\newcommand{\fileFEM}{Chapters/3.Chapter/figures/numerical_results/DATA/results_reactiondiffusion_0-200_3refinements_endtoend/fem/test_FEM_17nodes.csv}
\pgfplotsset{A/.append style={
	    xmin=-0.1,
	    xmax=1.1,
	    ymin=-0.1,
	    ymax=2.2,
	    ylabel={$u(x)$},
	    height=0.59*\textwidth,
	    width=\textwidth,
		xtick={0,0.25,...,1},
            }
}
\begin{axis}[A]
\addplot[color=green!30!white, line width=1., mark=*, mark options={scale=.8, solid}] table[x expr=\thisrow{x},y=u_FEM0]{\fileFEM};
\addplot[color=red!30!white, line width=1., mark=*, mark options={scale=.8, solid}] table[x expr=\thisrow{x},y=u_FEM1]{\fileFEM};
\addplot[color=blue!30!white, line width=1., mark=*, mark options={scale=.8, solid}] table[x expr=\thisrow{x},y=u_FEM2]{\fileFEM};
\addplot[color=orange!30!white, line width=1., mark=*, mark options={scale=.8, solid}] table[x expr=\thisrow{x},y=u_FEM3]{\fileFEM};
\addplot[color=gray!30!white, line width=1., mark=*, mark options={scale=.8, solid}] table[x expr=\thisrow{x},y=u_FEM4]{\fileFEM};
\addplot[color=green!80!black, style=dashed, line width=0.7 ,mark=asterisk, mark options = {scale=1, solid}] table[x expr=\thisrow{x},y=u_pred0]{\fileDL};
\addplot[color=red!80!black, style=dashed, line width=0.7 ,mark=asterisk, mark options = {scale=1, solid}] table[x expr=\thisrow{x},y=u_pred1]{\fileDL};
\addplot[color=blue!80!black, style=dashed, line width=0.7 ,mark=asterisk, mark options = {scale=1, solid}] table[x expr=\thisrow{x},y=u_pred2]{\fileDL};
\addplot[color=orange!80!black, style=dashed, line width=0.7 ,mark=asterisk, mark options = {scale=1, solid}] table[x expr=\thisrow{x},y=u_pred3]{\fileDL};
\addplot[color=gray!80!black, style=dashed, line width=0.7 ,mark=asterisk, mark options = {scale=1, solid}] table[x expr=\thisrow{x},y=u_pred4]{\fileDL};
\end{axis}
\end{tikzpicture}
\caption{Step $s=2$ (before training).}
\label{fig:step2_start_reactiondiffusion_parametric_exponential}
\end{subfigure} \hfill
\begin{subfigure}[b]{0.49\textwidth}
\centering
\begin{tikzpicture}
\newcommand{\fileDL}{Chapters/3.Chapter/figures/numerical_results/DATA/results_reactiondiffusion_0-200_3refinements_endtoend/end/test_prediction_17nodes.csv}
\newcommand{\fileFEM}{Chapters/3.Chapter/figures/numerical_results/DATA/results_reactiondiffusion_0-200_3refinements_endtoend/fem/test_FEM_17nodes.csv}
\pgfplotsset{A/.append style={
	    xmin=-0.1,
	    xmax=1.1,
	    ymin=-0.1,
	    ymax=2.2,
	    height=0.59*\textwidth,
	    width=\textwidth,
		xtick={0,0.25,...,1},
            }
}
\begin{axis}[A]
\addplot[color=green!30!white, line width=1., mark=*, mark options={scale=.8, solid}] table[x expr=\thisrow{x},y=u_FEM0]{\fileFEM};
\addplot[color=red!30!white, line width=1., mark=*, mark options={scale=.8, solid}] table[x expr=\thisrow{x},y=u_FEM1]{\fileFEM};
\addplot[color=blue!30!white, line width=1., mark=*, mark options={scale=.8, solid}] table[x expr=\thisrow{x},y=u_FEM2]{\fileFEM};
\addplot[color=orange!30!white, line width=1., mark=*, mark options={scale=.8, solid}] table[x expr=\thisrow{x},y=u_FEM3]{\fileFEM};
\addplot[color=gray!30!white, line width=1., mark=*, mark options={scale=.8, solid}] table[x expr=\thisrow{x},y=u_FEM4]{\fileFEM};
\addplot[color=green!80!black, style=dashed, line width=0.7 ,mark=asterisk, mark options = {scale=1, solid}] table[x expr=\thisrow{x},y=u_pred0]{\fileDL};
\addplot[color=red!80!black, style=dashed, line width=0.7 ,mark=asterisk, mark options = {scale=1, solid}] table[x expr=\thisrow{x},y=u_pred1]{\fileDL};
\addplot[color=blue!80!black, style=dashed, line width=0.7 ,mark=asterisk, mark options = {scale=1, solid}] table[x expr=\thisrow{x},y=u_pred2]{\fileDL};
\addplot[color=orange!80!black, style=dashed, line width=0.7 ,mark=asterisk, mark options = {scale=1, solid}] table[x expr=\thisrow{x},y=u_pred3]{\fileDL};
\addplot[color=gray!80!black, style=dashed, line width=0.7 ,mark=asterisk, mark options = {scale=1, solid}] table[x expr=\thisrow{x},y=u_pred4]{\fileDL};
\end{axis}
\end{tikzpicture}
\caption{Step $s=2$ (after training).}
\label{fig:step2_end_reactiondiffusion_parametric_exponential}
\end{subfigure}
\vskip 0.5em
\begin{subfigure}[b]{0.49\textwidth}
\centering
\begin{tikzpicture}
\newcommand{\fileDL}{Chapters/3.Chapter/figures/numerical_results/DATA/results_reactiondiffusion_0-200_3refinements_endtoend/start/test_prediction_33nodes.csv}
\newcommand{\fileFEM}{Chapters/3.Chapter/figures/numerical_results/DATA/results_reactiondiffusion_0-200_3refinements_endtoend/fem/test_FEM_33nodes.csv}
\pgfplotsset{A/.append style={
	    xmin=-0.1,
	    xmax=1.1,
	    ymin=-0.1,
	    ymax=2.2,
	    ylabel={$u(x)$},
	    height=0.59*\textwidth,
	    width=\textwidth,
		xtick={0,0.25,...,1},
            }
}
\begin{axis}[A]
\addplot[color=green!30!white, line width=1., mark=*, mark options={scale=.8, solid}] table[x expr=\thisrow{x},y=u_FEM0]{\fileFEM};
\addplot[color=red!30!white, line width=1., mark=*, mark options={scale=.8, solid}] table[x expr=\thisrow{x},y=u_FEM1]{\fileFEM};
\addplot[color=blue!30!white, line width=1., mark=*, mark options={scale=.8, solid}] table[x expr=\thisrow{x},y=u_FEM2]{\fileFEM};
\addplot[color=orange!30!white, line width=1., mark=*, mark options={scale=.8, solid}] table[x expr=\thisrow{x},y=u_FEM3]{\fileFEM};
\addplot[color=gray!30!white, line width=1., mark=*, mark options={scale=.8, solid}] table[x expr=\thisrow{x},y=u_FEM4]{\fileFEM};
\addplot[color=green!80!black, style=dashed, line width=0.7 ,mark=asterisk, mark options = {scale=1, solid}] table[x expr=\thisrow{x},y=u_pred0]{\fileDL};
\addplot[color=red!80!black, style=dashed, line width=0.7 ,mark=asterisk, mark options = {scale=1, solid}] table[x expr=\thisrow{x},y=u_pred1]{\fileDL};
\addplot[color=blue!80!black, style=dashed, line width=0.7 ,mark=asterisk, mark options = {scale=1, solid}] table[x expr=\thisrow{x},y=u_pred2]{\fileDL};
\addplot[color=orange!80!black, style=dashed, line width=0.7 ,mark=asterisk, mark options = {scale=1, solid}] table[x expr=\thisrow{x},y=u_pred3]{\fileDL};
\addplot[color=gray!80!black, style=dashed, line width=0.7 ,mark=asterisk, mark options = {scale=1, solid}] table[x expr=\thisrow{x},y=u_pred4]{\fileDL};
\end{axis}
\end{tikzpicture}
\caption{Step $s=3$ (before training).}
\label{fig:step3_start_reactiondiffusion_parametric_exponential}
\end{subfigure} \hfill
\begin{subfigure}[b]{0.49\textwidth}
\centering
\begin{tikzpicture}
\newcommand{\fileDL}{Chapters/3.Chapter/figures/numerical_results/DATA/results_reactiondiffusion_0-200_3refinements_endtoend/end/test_prediction_33nodes.csv}
\newcommand{\fileFEM}{Chapters/3.Chapter/figures/numerical_results/DATA/results_reactiondiffusion_0-200_3refinements_endtoend/fem/test_FEM_33nodes.csv}
\pgfplotsset{A/.append style={
	    xmin=-0.1,
	    xmax=1.1,
	    ymin=-0.1,
	    ymax=2.2,
	    height=0.59*\textwidth,
	    width=\textwidth,
		xtick={0,0.25,...,1},
            }
}
\begin{axis}[A]
\addplot[color=green!30!white, line width=1., mark=*, mark options={scale=.8, solid}] table[x expr=\thisrow{x},y=u_FEM0]{\fileFEM};
\addplot[color=red!30!white, line width=1., mark=*, mark options={scale=.8, solid}] table[x expr=\thisrow{x},y=u_FEM1]{\fileFEM};
\addplot[color=blue!30!white, line width=1., mark=*, mark options={scale=.8, solid}] table[x expr=\thisrow{x},y=u_FEM2]{\fileFEM};
\addplot[color=orange!30!white, line width=1., mark=*, mark options={scale=.8, solid}] table[x expr=\thisrow{x},y=u_FEM3]{\fileFEM};
\addplot[color=gray!30!white, line width=1., mark=*, mark options={scale=.8, solid}] table[x expr=\thisrow{x},y=u_FEM4]{\fileFEM};
\addplot[color=green!80!black, style=dashed, line width=0.7 ,mark=asterisk, mark options = {scale=1, solid}] table[x expr=\thisrow{x},y=u_pred0]{\fileDL};
\addplot[color=red!80!black, style=dashed, line width=0.7 ,mark=asterisk, mark options = {scale=1, solid}] table[x expr=\thisrow{x},y=u_pred1]{\fileDL};
\addplot[color=blue!80!black, style=dashed, line width=0.7 ,mark=asterisk, mark options = {scale=1, solid}] table[x expr=\thisrow{x},y=u_pred2]{\fileDL};
\addplot[color=orange!80!black, style=dashed, line width=0.7 ,mark=asterisk, mark options = {scale=1, solid}] table[x expr=\thisrow{x},y=u_pred3]{\fileDL};
\addplot[color=gray!80!black, style=dashed, line width=0.7 ,mark=asterisk, mark options = {scale=1, solid}] table[x expr=\thisrow{x},y=u_pred4]{\fileDL};
\end{axis}
\end{tikzpicture}
\caption{Step $s=3$ (after training).}
\label{fig:step3_end_reactiondiffusion_parametric_exponential}
\end{subfigure}
\vskip 0.5em
\begin{subfigure}[b]{0.49\textwidth}
\centering
\begin{tikzpicture}
\newcommand{\fileDL}{Chapters/3.Chapter/figures/numerical_results/DATA/results_reactiondiffusion_0-200_3refinements_endtoend/start/test_prediction_65nodes.csv}
\newcommand{\fileFEM}{Chapters/3.Chapter/figures/numerical_results/DATA/results_reactiondiffusion_0-200_3refinements_endtoend/fem/test_FEM_65nodes.csv}
\pgfplotsset{A/.append style={
	    xmin=-0.1,
	    xmax=1.1,
	    xlabel={$x$},
	    ymin=-0.1,
	    ymax=2.2,
	    ylabel={$u(x)$},
	    height=0.59*\textwidth,
	    width=\textwidth,
		xtick={0,0.25,...,1},
            }
}
\begin{axis}[A]
\addplot[color=green!30!white, line width=1., mark=*, mark options={scale=.8, solid}] table[x expr=\thisrow{x},y=u_FEM0]{\fileFEM};
\addplot[color=red!30!white, line width=1., mark=*, mark options={scale=.8, solid}] table[x expr=\thisrow{x},y=u_FEM1]{\fileFEM};
\addplot[color=blue!30!white, line width=1., mark=*, mark options={scale=.8, solid}] table[x expr=\thisrow{x},y=u_FEM2]{\fileFEM};
\addplot[color=orange!30!white, line width=1., mark=*, mark options={scale=.8, solid}] table[x expr=\thisrow{x},y=u_FEM3]{\fileFEM};
\addplot[color=gray!30!white, line width=1., mark=*, mark options={scale=.8, solid}] table[x expr=\thisrow{x},y=u_FEM4]{\fileFEM};
\addplot[color=green!80!black, style=dashed, line width=0.7 ,mark=asterisk, mark options = {scale=1, solid}] table[x expr=\thisrow{x},y=u_pred0]{\fileDL};
\addplot[color=red!80!black, style=dashed, line width=0.7 ,mark=asterisk, mark options = {scale=1, solid}] table[x expr=\thisrow{x},y=u_pred1]{\fileDL};
\addplot[color=blue!80!black, style=dashed, line width=0.7 ,mark=asterisk, mark options = {scale=1, solid}] table[x expr=\thisrow{x},y=u_pred2]{\fileDL};
\addplot[color=orange!80!black, style=dashed, line width=0.7 ,mark=asterisk, mark options = {scale=1, solid}] table[x expr=\thisrow{x},y=u_pred3]{\fileDL};
\addplot[color=gray!80!black, style=dashed, line width=0.7 ,mark=asterisk, mark options = {scale=1, solid}] table[x expr=\thisrow{x},y=u_pred4]{\fileDL};
\end{axis}
\end{tikzpicture}
\caption{Step $s=4$ (before training).}
\label{fig:step4_start_reactiondiffusion_parametric_exponential}
\end{subfigure} \hfill
\begin{subfigure}[b]{0.49\textwidth}
\centering
\begin{tikzpicture}
\newcommand{\fileDL}{Chapters/3.Chapter/figures/numerical_results/DATA/results_reactiondiffusion_0-200_3refinements_endtoend/end/test_prediction_65nodes.csv}
\newcommand{\fileFEM}{Chapters/3.Chapter/figures/numerical_results/DATA/results_reactiondiffusion_0-200_3refinements_endtoend/fem/test_FEM_65nodes.csv}
\pgfplotsset{A/.append style={
	    xmin=-0.1,
	    xmax=1.1,
	    xlabel={$x$},
	    ymin=-0.1,
	    ymax=2.2,
	    height=0.59*\textwidth,
	    width=\textwidth,
		xtick={0,0.25,...,1},
            }
}
\begin{axis}[A]
\addplot[color=green!30!white, line width=1., mark=*, mark options={scale=.8, solid}] table[x expr=\thisrow{x},y=u_FEM0]{\fileFEM};
\addplot[color=red!30!white, line width=1., mark=*, mark options={scale=.8, solid}] table[x expr=\thisrow{x},y=u_FEM1]{\fileFEM};
\addplot[color=blue!30!white, line width=1., mark=*, mark options={scale=.8, solid}] table[x expr=\thisrow{x},y=u_FEM2]{\fileFEM};
\addplot[color=orange!30!white, line width=1., mark=*, mark options={scale=.8, solid}] table[x expr=\thisrow{x},y=u_FEM3]{\fileFEM};
\addplot[color=gray!30!white, line width=1., mark=*, mark options={scale=.8, solid}] table[x expr=\thisrow{x},y=u_FEM4]{\fileFEM};
\addplot[color=green!80!black, style=dashed, line width=0.7 ,mark=asterisk, mark options = {scale=1, solid}] table[x expr=\thisrow{x},y=u_pred0]{\fileDL};
\addplot[color=red!80!black, style=dashed, line width=0.7 ,mark=asterisk, mark options = {scale=1, solid}] table[x expr=\thisrow{x},y=u_pred1]{\fileDL};
\addplot[color=blue!80!black, style=dashed, line width=0.7 ,mark=asterisk, mark options = {scale=1, solid}] table[x expr=\thisrow{x},y=u_pred2]{\fileDL};
\addplot[color=orange!80!black, style=dashed, line width=0.7 ,mark=asterisk, mark options = {scale=1, solid}] table[x expr=\thisrow{x},y=u_pred3]{\fileDL};
\addplot[color=gray!80!black, style=dashed, line width=0.7 ,mark=asterisk, mark options = {scale=1, solid}] table[x expr=\thisrow{x},y=u_pred4]{\fileDL};
\end{axis}
\end{tikzpicture}
\caption{Step $s=4$ (after training).}
\label{fig:step4_end_reactiondiffusion_parametric_exponential}
\end{subfigure}
\caption{Four steps of the \acs{DeepFEM} for problem \eqref{eq:parmetric bvp} with $0<\alpha<200$ employing the energy norm in the loss function. $u_{\text{FEM}}(\alpha)$ and $u_{\text{NN}}(\alpha)$ denote the FEM solution and \acs{DeepFEM} prediction for the $\alpha$ parameter coefficient, respectively.}
\label{fig:reactiondiffusion_parametric_exponential}
\end{figure}
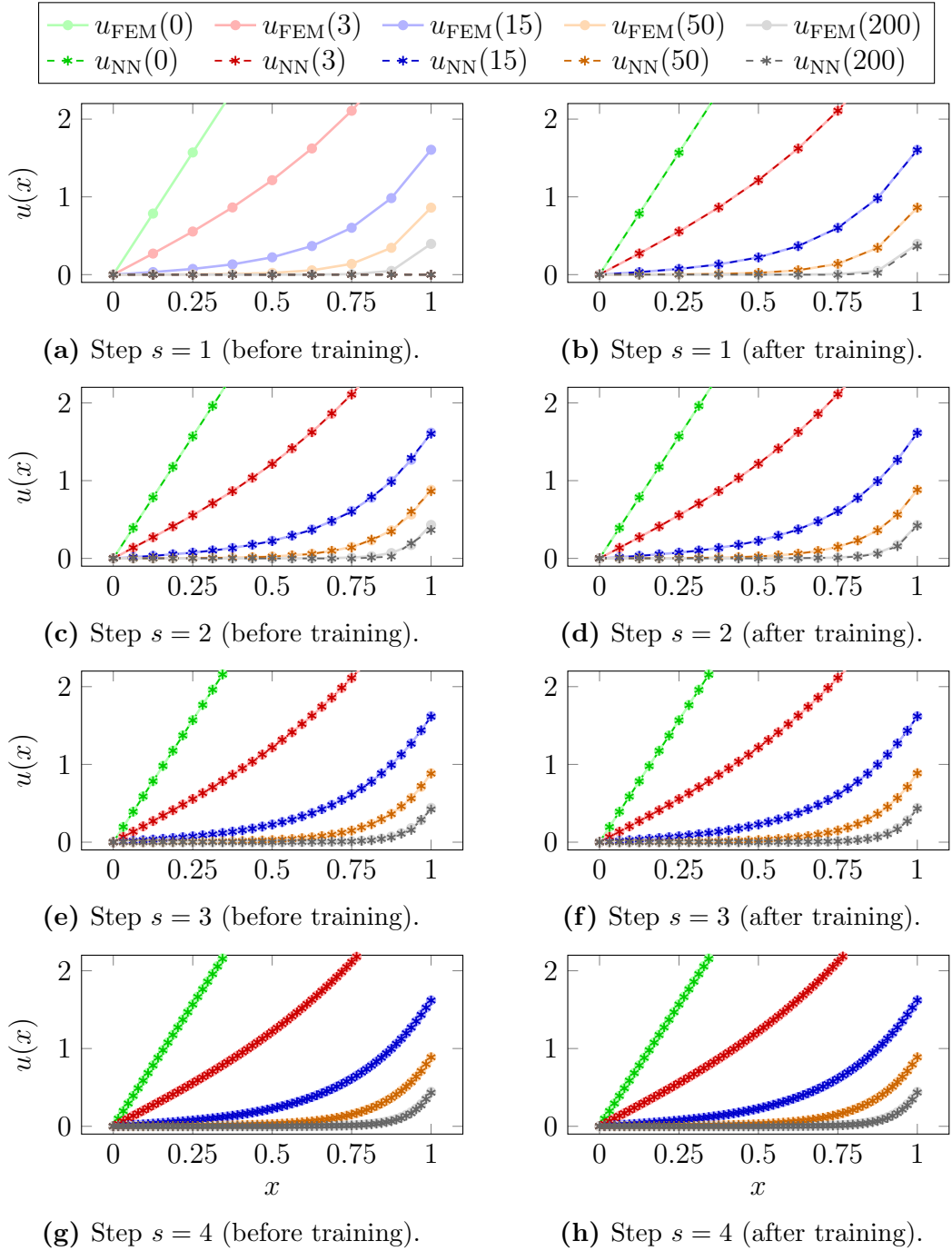

We select a database of $100$ samples randomly and logarithmically distributed to train the NN. To visualize the NN performance after the training, we consider $\{0,3,15,50,200\}$ as test data, which does not take part during training. \Cref{fig:reactiondiffusion_parametric_exponential} shows the parametric mesh-by-mesh adaptivity of the NN over these test data. We observe a good behavior of the NN for all values of $\alpha$. \Cref{table:reaction-diffusion_losserror} displays the losses and energy-norm errors evaluated on the test data points at the end of the training.

\begin{table}[htb]
\centering
\begin{tabular}{cccccc}
\hline
$\alpha$ & $0$ & $3$ & $15$ & $50$ & $200$ \\ \hline
$\lVert \mathbf{A}(\alpha) \mathbf{u}_\text{NN}(\alpha)-\mathbf{f}\rVert_{\mathbf{P}(\alpha)}$ & $0.0015$ & $0.013$ & $0.0046$ & $0.012$ & $0.082$ \\
$\lVert \mathbf{u}_\text{NN}(\alpha)-\mathbf{u}_\text{FEM}(\alpha)\rVert_{\mathbf{A}(\alpha)}$ & $0.0160$ & $0.014$ & $0.0061$ & $0.014$ & $0.095$ \\ \hline
\end{tabular}
\caption{Samplewise residual and error values in the energy norm of the test data at the end of the fourth step for problem \eqref{eq:parmetric bvp} with $0 < \alpha < 200$ when employing eight-size blocks in the preconditioners.}
\label{table:reaction-diffusion_losserror}
\end{table}

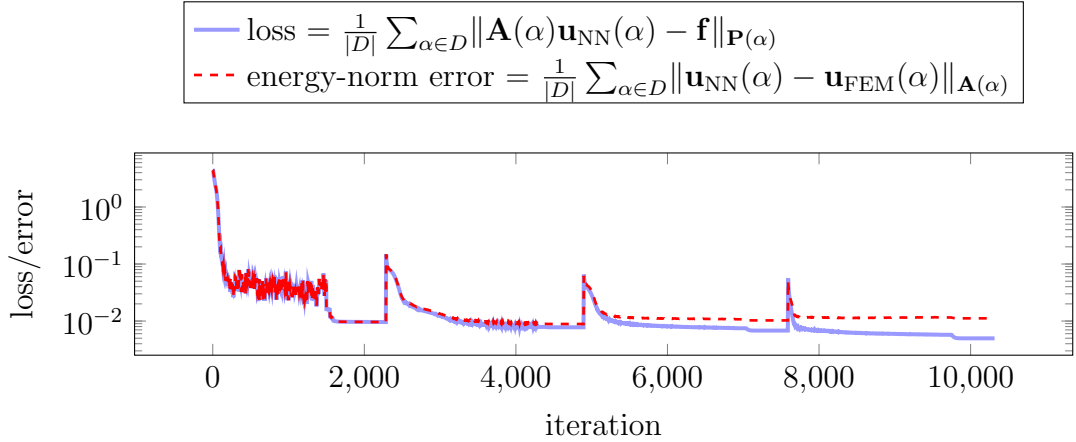
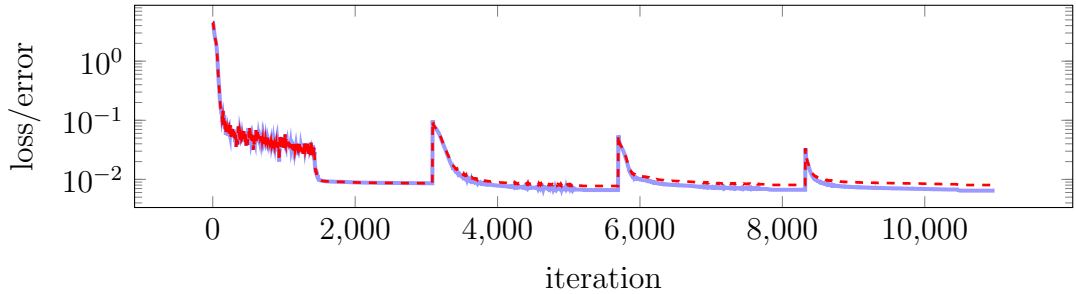
\begin{figure}[htbp]
\begin{subfigure}{0.98\textwidth}
\centering
\begin{tikzpicture}
\pgfplotsset{loss/.append style={ 
     xlabel = {iteration}, 
     xtick={0,2000,...,10000},
     ylabel = {loss/error},
    scaled x ticks=false,
     height=0.3*\textwidth,    
     width=0.98*\textwidth,    
     legend style={fill=none, at={(0.5,1.75)}, anchor=north, minimum height = 0.6cm, /tikz/every even column/.append style={column sep=0.5cm}},
	 legend cell align={left}    
     }
}
\begin{semilogyaxis}[loss]
\addplot[color=blue!40, line width=1.5 ,] table[x expr=\thisrow{iteration},y=loss_history]{Chapters/3.Chapter/figures/numerical_results/DATA/results_reactiondiffusion_0-200_3refinements_endtoend/loss_history_reduced2.csv};
\addlegendentry{loss = $\frac{1}{|D|}\sum_{\alpha\in D} \lVert \mathbf{A}(\alpha) \mathbf{u}_\text{NN}(\alpha)-\mathbf{f}\rVert_{\mathbf{P}(\alpha)}$};
\addplot[color=red!, line width=1. , style=dashed] table[x expr=\thisrow{iteration},y=error_history]{Chapters/3.Chapter/figures/numerical_results/DATA/results_reactiondiffusion_0-200_3refinements_endtoend/error_history_reduced2.csv};
\addlegendentry{energy-norm error = $\frac{1}{|D|}\sum_{\alpha\in D} \lVert \mathbf{u}_\text{NN}(\alpha)-\mathbf{u}_\text{FEM}(\alpha)\rVert_{\mathbf{A}(\alpha)}$};
\end{semilogyaxis}
\end{tikzpicture}
\caption{Employing preconditioners with blocks of size eight in all steps.}
\label{fig:loss_reactiondiffusion_endtoend}
\end{subfigure}
\vskip 0.5em
\begin{subfigure}{0.98\textwidth}
\centering
\begin{tikzpicture}
\pgfplotsset{loss/.append style={ 
     xlabel = {iteration},
     xtick={0,2000,...,10000},
     ylabel = {loss/error},
    scaled x ticks=false,
     height=0.3*\textwidth,    
     width=0.98*\textwidth,    
     legend style={fill=none, at={(0.5,1.65)}, anchor=north, minimum height = 0.6cm, /tikz/every even column/.append style={column sep=0.5cm}},
	 legend cell align={left}    
     }
}
\begin{semilogyaxis}[loss]
\addplot[color=blue!40, line width=1.5 ,] table[x expr=\thisrow{iteration},y=loss_history]{Chapters/3.Chapter/figures/numerical_results/DATA/results_reactiondiffusion_0-200_3refinements_moresmoothed_endtoend/loss_history_reduced2.csv};
\addplot[color=red!, line width=1. , style=dashed] table[x expr=\thisrow{iteration},y=error_history]{Chapters/3.Chapter/figures/numerical_results/DATA/results_reactiondiffusion_0-200_3refinements_moresmoothed_endtoend/error_history_reduced2.csv};

\end{semilogyaxis}
\end{tikzpicture}
\caption{Employing preconditioners with blocks of sizes $8$, $8$, $16$, and $32$ at steps one, two, three, and four, respectively.}
\label{fig:loss_reactiondiffusion_moresmoothed_endtoend}
\end{subfigure}
\caption{Loss function evolution along with the energy-norm error during the four-step training at problem \eqref{eq:parmetric bvp} with $0\leq\alpha\leq 200$.}
\label{fig:loss_helmholtz_parametric}
\end{figure}

Analyzing the evolution of the training in \Cref{fig:loss_reactiondiffusion_endtoend}, we observe that the loss decrease is more significant than the error in the energy norm decrease, as expected. By increasing the size of the blocks in the preconditioners step by step (with sizes $8$, $8$, $16$, and $32$ at steps one, two, three, and four, respectively), the energy-norm error does decrease in consistency with the loss (see \Cref{fig:loss_reactiondiffusion_moresmoothed_endtoend}). \Cref{table:reaction-diffusion_losserror_moresmoothed} displays the losses and energy-norm errors evaluated on the test data at the end of this enhanced training with preconditioning.

\begin{table}[htb]
\centering
\begin{tabular}{cccccc}
\hline
$\alpha$ & $0$ & $3$ & $15$ & $50$ & $200$ \\ \hline
$\lVert \mathbf{A}(\alpha) \mathbf{u}_\text{NN}(\alpha)-\mathbf{f}\rVert_{\mathbf{P}(\alpha)}$ & $0.0022$ & $0.0059$ & $0.0043$ & $0.0086$ & $0.011$ \\
$\lVert \mathbf{u}_\text{NN}(\alpha)-\mathbf{u}_\text{FEM}(\alpha)\rVert_{\mathbf{A}(\alpha)}$ & $0.0053$ & $0.0075$ & $0.0058$ & $0.0101$ & $0.013$ \\ \hline
\end{tabular}
\caption{Samplewise residual and error values in the energy norm of the test data at the end of the fourth step for problem \eqref{eq:parmetric bvp} with $0 < \alpha < 200$ when employing size-increasing blocks in the preconditioners.}
\label{table:reaction-diffusion_losserror_moresmoothed}
\end{table}

\subsubsection{Helmholtz's parametric equation: \boldmath{$-50 < \alpha < -30$}}
\label[section]{section3.6.4.1}

The exact solution is $u^*(x)=C \sin(\sqrt{\alpha}x)$ with $C=\frac{2\pi}{\sqrt{\alpha}\cos(\sqrt{\alpha})}$. We consider $100$ samples randomly and uniformly distributed as the training data and 
\begin{equation}
\{-50, -45, -40, -35, -30\}
\end{equation} as the test data. We consider the $H^1$-norm for the loss and maintain the same trainable block architectures as above. If we train the NN employing the same increasing block-size criterion for the preconditioners as in \Cref{section3.6.4.1}, we observe that the loss does not decrease the $H^1$-norm error in this occasion (\Cref{fig:loss_helmholtz_parametric_-50--30_smoothed2020}). If we employ the inverses as the preconditioners (\Cref{fig:loss_helmholtz_parametric_-50--30_inverse2020}), the NN converges sufficiently to show adequate results. \Cref{fig:helmholtz_parametric} shows the step-by-step predictions for the test data when the training is performed according to \Cref{fig:loss_helmholtz_parametric_-50--30_inverse2020}.

\begin{figure}[htbp]
\begin{subfigure}{0.98\textwidth}
\centering
\begin{tikzpicture}
\pgfplotsset{loss/.append style={ 
     xmax=12000,      
     xlabel = {iteration},
     xtick={0,2000,...,10000},
     ymin=1e-2,     
     ymax=1e1, 
     ylabel = {loss/error},
    scaled x ticks=false,
     height=0.3*\textwidth,    
     width=\textwidth,    
     legend style={fill=none, at={(0.5,1.75)}, anchor=north, minimum height = 0.6cm, /tikz/every even column/.append style={column sep=0.5cm}},
     legend cell align={left}    
     }
}\begin{semilogyaxis}[loss]
\addplot[color=blue!40, line width=1.5 ,] table[x expr=\thisrow{iteration},y=loss_history]{Chapters/3.Chapter/figures/numerical_results/DATA/results_helmholtz_-50--30_smoothed2020_3refinements_endtoend/loss_history_reduced1.csv};
\addlegendentry{loss = $\frac{1}{|D|}\sum_{\alpha\in D} \lVert \mathbf{P}(\alpha)\{\mathbf{A}(\alpha) \mathbf{u}_\text{NN}(\alpha)-\mathbf{f}\}\rVert_{\mathbf{K}+\mathbf{M}}$};
\addplot[color=red!, line width=1. , style=dashed] table[x expr=\thisrow{iteration},y=error_history]{Chapters/3.Chapter/figures/numerical_results/DATA/results_helmholtz_-50--30_smoothed2020_3refinements_endtoend/error_history_reduced1.csv};
\addlegendentry{energy-norm error = $\frac{1}{|D|}\sum_{\alpha\in D} \lVert \mathbf{u}_\text{NN}(\alpha)-\mathbf{u}_\text{FEM}(\alpha)\rVert_{\mathbf{A}(\alpha)}$};
\end{semilogyaxis}
\end{tikzpicture}
\caption{Employing block-Jacobi preconditioners of sizes equal to $8$, $8$, $16$ and $32$ at the first, second, third, and fourth steps, respectively.}
\label{fig:loss_helmholtz_parametric_-50--30_smoothed2020}
\end{subfigure}
\vskip 0.5em
\begin{subfigure}{0.98\textwidth}
\centering
\begin{tikzpicture}
\pgfplotsset{loss/.append style={ 
     xmax=12000,      
     xlabel = {iteration},
     xtick={0,2000,...,10000},
     ymin=1e-2,     
     ymax=1e1, 
     ylabel = {loss/error},
    scaled x ticks=false,
     height=0.3*\textwidth,    
     width=\textwidth,    
     }
}
\begin{semilogyaxis}[loss]
\addplot[color=blue!40, line width=1.5 ,] table[x expr=\thisrow{iteration},y=loss_history]{Chapters/3.Chapter/figures/numerical_results/DATA/results_helmholtz_-50--30_inverse2020_3refinements_endtoend/loss_history_reduced1.csv};
\addplot[color=red!, line width=1. , style=dashed] table[x expr=\thisrow{iteration},y=error_history]{Chapters/3.Chapter/figures/numerical_results/DATA/results_helmholtz_-50--30_inverse2020_3refinements_endtoend/error_history_reduced1.csv};

\end{semilogyaxis}
\end{tikzpicture}
\caption{Employing inverses as preconditioners.}
\label{fig:loss_helmholtz_parametric_-50--30_inverse2020}
\end{subfigure}
\caption{Loss function evolution along $H^1$-norm error during the four-step training at problem \eqref{eq:parmetric bvp} with $-50\leq\alpha\leq -30$.}
\label{}
\end{figure}
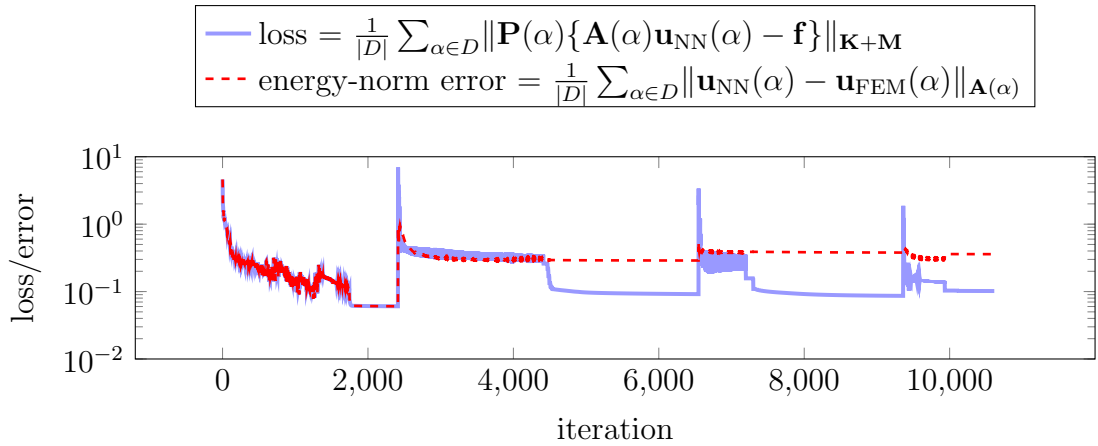
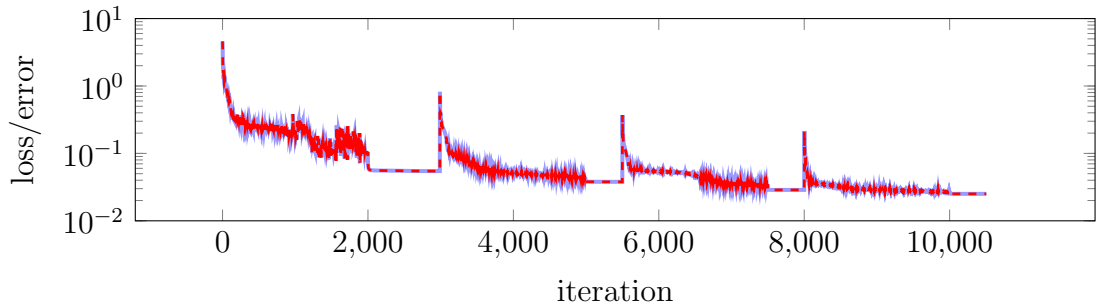

\begin{figure}[htbp]
\centering
\begin{subfigure}[b]{\textwidth}
\centering
 \begin{tikzpicture}
\pgfplotsset{A/.append style={
		hide axis,
	    xmin=-0.1,
	    xmax=1.1,
	    ymin=-2.2,
	    ymax=2.2,
	 	legend style={fill=none, anchor=north, /tikz/every even column/.append style={column sep=0.1cm}},
	 	legend columns = 5
            }
}
\begin{axis}[A]
\addlegendimage{color=green!30!white, line width=1., mark=*, mark options={scale=0.8, solid}}
\addlegendentry{$u_{\text{FEM}}(-50)$}
\addlegendimage{color=red!30!white, line width=1., mark=*, mark options={scale=0.8, solid}}
\addlegendentry{$u_{\text{FEM}}(-45)$}
\addlegendimage{color=blue!30!white, line width=1., mark=*, mark options={scale=0.8, solid}}
\addlegendentry{$u_{\text{FEM}}(-40)$}
\addlegendimage{color=orange!30!white, line width=1., mark=*, mark options={scale=0.8, solid}}
\addlegendentry{$u_{\text{FEM}}(-35)$}
\addlegendimage{color=gray!30!white, line width=1., mark=*, mark options={scale=0.8, solid}}
\addlegendentry{$u_{\text{FEM}}(-30)$}
\addlegendimage{color=green!80!black, style=dashed, line width=0.7 ,mark=asterisk, mark options = {scale=1, solid}}
\addlegendentry{$u_{\text{NN}}(-50)$};
\addlegendimage{color=red!80!black, style=dashed, line width=0.7 ,mark=asterisk, mark options = {scale=1, solid}}
\addlegendentry{$u_{\text{NN}}(-45)$};
\addlegendimage{color=blue!80!black, style=dashed, line width=0.7 ,mark=asterisk, mark options = {scale=1, solid}}
\addlegendentry{$u_{\text{NN}}(-40)$};
\addlegendimage{color=orange!80!black, style=dashed, line width=0.7 ,mark=asterisk, mark options = {scale=1, solid}}
\addlegendentry{$u_{\text{NN}}(-35)$};
\addlegendimage{color=gray!80!black, style=dashed, line width=0.7 ,mark=asterisk, mark options = {scale=1, solid}}
\addlegendentry{$u_{\text{NN}}(-30)$};
\end{axis}
\end{tikzpicture}
\end{subfigure}
\vskip 0.5em
\begin{subfigure}[b]{0.49\textwidth}
\centering
\begin{tikzpicture}
\newcommand{\fileDL}{Chapters/3.Chapter/figures/numerical_results/DATA/results_helmholtz_-50--30_inverse2020_3refinements_endtoend/start/test_prediction_9nodes.csv}
\newcommand{\fileFEM}{Chapters/3.Chapter/figures/numerical_results/DATA/results_helmholtz_-50--30_inverse2020_3refinements_endtoend/fem/test_FEM_9nodes.csv}
\pgfplotsset{A/.append style={
	    xmin=-0.1,
	    xmax=1.1,
	    ymin=-2.2,
	    ymax=2.2,
	    ylabel = {$u(x)$},
	    height=0.59*\textwidth,
	    width=\textwidth,
		xtick={0,0.25,...,1},
            }
}
\begin{axis}[A]
\addplot[color=green!30!white, line width=1., mark=*, mark options={scale=.8, solid}] table[x expr=\thisrow{x},y=u_FEM0]{\fileFEM};
\addplot[color=red!30!white, line width=1., mark=*, mark options={scale=.8, solid}] table[x expr=\thisrow{x},y=u_FEM1]{\fileFEM};
\addplot[color=blue!30!white, line width=1., mark=*, mark options={scale=.8, solid}] table[x expr=\thisrow{x},y=u_FEM2]{\fileFEM};
\addplot[color=orange!30!white, line width=1., mark=*, mark options={scale=.8, solid}] table[x expr=\thisrow{x},y=u_FEM3]{\fileFEM};
\addplot[color=gray!30!white, line width=1., mark=*, mark options={scale=.8, solid}] table[x expr=\thisrow{x},y=u_FEM4]{\fileFEM};
\addplot[color=green!80!black, style=dashed, line width=0.7 ,mark=asterisk, mark options = {scale=1, solid}] table[x expr=\thisrow{x},y=u_pred0]{\fileDL};
\addplot[color=red!80!black, style=dashed, line width=0.7 ,mark=asterisk, mark options = {scale=1, solid}] table[x expr=\thisrow{x},y=u_pred1]{\fileDL};
\addplot[color=blue!80!black, style=dashed, line width=0.7 ,mark=asterisk, mark options = {scale=1, solid}] table[x expr=\thisrow{x},y=u_pred2]{\fileDL};
\addplot[color=orange!80!black, style=dashed, line width=0.7 ,mark=asterisk, mark options = {scale=1, solid}] table[x expr=\thisrow{x},y=u_pred3]{\fileDL};
\addplot[color=gray!80!black, style=dashed, line width=0.7 ,mark=asterisk, mark options = {scale=1, solid}] table[x expr=\thisrow{x},y=u_pred4]{\fileDL};
\end{axis}
\end{tikzpicture}
\caption{Step $s=1$ (before training).}
\label{fig:step1_start_helmholtz_parametric}
\end{subfigure} \hfill
\begin{subfigure}[b]{0.49\textwidth}
\centering
\begin{tikzpicture}
\newcommand{\fileDL}{Chapters/3.Chapter/figures/numerical_results/DATA/results_helmholtz_-50--30_inverse100100_3refinements_endtoend/end/test_prediction_9nodes.csv}
\newcommand{\fileFEM}{Chapters/3.Chapter/figures/numerical_results/DATA/results_helmholtz_-50--30_inverse100100_3refinements_endtoend/fem/test_FEM_9nodes.csv}
\pgfplotsset{A/.append style={
	    xmin=-0.1,
	    xmax=1.1,
	    ymin=-2.2,
	    ymax=2.2,
	    height=0.59*\textwidth,
	    width=\textwidth,
		xtick={0,0.25,...,1},
	 	legend style={fill=none, at={(1.3,1.6)},anchor=north, /tikz/every even column/.append style={column sep=0.5cm}},
	 	legend columns = 3
            }
}
\begin{axis}[A]
\addplot[color=green!30!white, line width=1., mark=*, mark options={scale=.8, solid}] table[x expr=\thisrow{x},y=u_FEM0]{\fileFEM};
\addplot[color=red!30!white, line width=1., mark=*, mark options={scale=.8, solid}] table[x expr=\thisrow{x},y=u_FEM1]{\fileFEM};
\addplot[color=blue!30!white, line width=1., mark=*, mark options={scale=.8, solid}] table[x expr=\thisrow{x},y=u_FEM2]{\fileFEM};
\addplot[color=orange!30!white, line width=1., mark=*, mark options={scale=.8, solid}] table[x expr=\thisrow{x},y=u_FEM3]{\fileFEM};
\addplot[color=gray!30!white, line width=1., mark=*, mark options={scale=.8, solid}] table[x expr=\thisrow{x},y=u_FEM4]{\fileFEM};
\addplot[color=green!80!black, style=dashed, line width=0.7 ,mark=asterisk, mark options = {scale=1, solid}] table[x expr=\thisrow{x},y=u_pred0]{\fileDL};
\addplot[color=red!80!black, style=dashed, line width=0.7 ,mark=asterisk, mark options = {scale=1, solid}] table[x expr=\thisrow{x},y=u_pred1]{\fileDL};
\addplot[color=blue!80!black, style=dashed, line width=0.7 ,mark=asterisk, mark options = {scale=1, solid}] table[x expr=\thisrow{x},y=u_pred2]{\fileDL};
\addplot[color=orange!80!black, style=dashed, line width=0.7 ,mark=asterisk, mark options = {scale=1, solid}] table[x expr=\thisrow{x},y=u_pred3]{\fileDL};
\addplot[color=gray!80!black, style=dashed, line width=0.7 ,mark=asterisk, mark options = {scale=1, solid}] table[x expr=\thisrow{x},y=u_pred4]{\fileDL};
\end{axis}
\end{tikzpicture}
\caption{Step $s=1$ (after training).}
\label{fig:step1_end_helmholtz_parametric}
\end{subfigure}
\vskip 0.5em    
\begin{subfigure}[b]{0.49\textwidth}
\centering
\begin{tikzpicture}
\newcommand{\fileDL}{Chapters/3.Chapter/figures/numerical_results/DATA/results_helmholtz_-50--30_inverse100100_3refinements_endtoend/start/test_prediction_17nodes.csv}
\newcommand{\fileFEM}{Chapters/3.Chapter/figures/numerical_results/DATA/results_helmholtz_-50--30_inverse100100_3refinements_endtoend/fem/test_FEM_17nodes.csv}
\pgfplotsset{A/.append style={
	    xmin=-0.1,
	    xmax=1.1,
	    ymin=-2.2,
	    ymax=2.2,
	    ylabel={$u(x)$},
	    height=0.59*\textwidth,
	    width=\textwidth,
		xtick={0,0.25,...,1},
            }
}
\begin{axis}[A]
\addplot[color=green!30!white, line width=1., mark=*, mark options={scale=.8, solid}] table[x expr=\thisrow{x},y=u_FEM0]{\fileFEM};
\addplot[color=red!30!white, line width=1., mark=*, mark options={scale=.8, solid}] table[x expr=\thisrow{x},y=u_FEM1]{\fileFEM};
\addplot[color=blue!30!white, line width=1., mark=*, mark options={scale=.8, solid}] table[x expr=\thisrow{x},y=u_FEM2]{\fileFEM};
\addplot[color=orange!30!white, line width=1., mark=*, mark options={scale=.8, solid}] table[x expr=\thisrow{x},y=u_FEM3]{\fileFEM};
\addplot[color=gray!30!white, line width=1., mark=*, mark options={scale=.8, solid}] table[x expr=\thisrow{x},y=u_FEM4]{\fileFEM};
\addplot[color=green!80!black, style=dashed, line width=0.7 ,mark=asterisk, mark options = {scale=1, solid}] table[x expr=\thisrow{x},y=u_pred0]{\fileDL};
\addplot[color=red!80!black, style=dashed, line width=0.7 ,mark=asterisk, mark options = {scale=1, solid}] table[x expr=\thisrow{x},y=u_pred1]{\fileDL};
\addplot[color=blue!80!black, style=dashed, line width=0.7 ,mark=asterisk, mark options = {scale=1, solid}] table[x expr=\thisrow{x},y=u_pred2]{\fileDL};
\addplot[color=orange!80!black, style=dashed, line width=0.7 ,mark=asterisk, mark options = {scale=1, solid}] table[x expr=\thisrow{x},y=u_pred3]{\fileDL};
\addplot[color=gray!80!black, style=dashed, line width=0.7 ,mark=asterisk, mark options = {scale=1, solid}] table[x expr=\thisrow{x},y=u_pred4]{\fileDL};
\end{axis}
\end{tikzpicture}
\caption{Step $s=2$ (before training).}
\label{fig:step2_start_helmholtz_parametric}
\end{subfigure} \hfill
\begin{subfigure}[b]{0.49\textwidth}
\centering
\begin{tikzpicture}
\newcommand{\fileDL}{Chapters/3.Chapter/figures/numerical_results/DATA/results_helmholtz_-50--30_inverse100100_3refinements_endtoend/end/test_prediction_17nodes.csv}
\newcommand{\fileFEM}{Chapters/3.Chapter/figures/numerical_results/DATA/results_helmholtz_-50--30_inverse100100_3refinements_endtoend/fem/test_FEM_17nodes.csv}
\pgfplotsset{A/.append style={
	    xmin=-0.1,
	    xmax=1.1,
	    ymin=-2.2,
	    ymax=2.2,
	    height=0.59*\textwidth,
	    width=\textwidth,
		xtick={0,0.25,...,1},
            }
}
\begin{axis}[A]
\addplot[color=green!30!white, line width=1., mark=*, mark options={scale=.8, solid}] table[x expr=\thisrow{x},y=u_FEM0]{\fileFEM};
\addplot[color=red!30!white, line width=1., mark=*, mark options={scale=.8, solid}] table[x expr=\thisrow{x},y=u_FEM1]{\fileFEM};
\addplot[color=blue!30!white, line width=1., mark=*, mark options={scale=.8, solid}] table[x expr=\thisrow{x},y=u_FEM2]{\fileFEM};
\addplot[color=orange!30!white, line width=1., mark=*, mark options={scale=.8, solid}] table[x expr=\thisrow{x},y=u_FEM3]{\fileFEM};
\addplot[color=gray!30!white, line width=1., mark=*, mark options={scale=.8, solid}] table[x expr=\thisrow{x},y=u_FEM4]{\fileFEM};
\addplot[color=green!80!black, style=dashed, line width=0.7 ,mark=asterisk, mark options = {scale=1, solid}] table[x expr=\thisrow{x},y=u_pred0]{\fileDL};
\addplot[color=red!80!black, style=dashed, line width=0.7 ,mark=asterisk, mark options = {scale=1, solid}] table[x expr=\thisrow{x},y=u_pred1]{\fileDL};
\addplot[color=blue!80!black, style=dashed, line width=0.7 ,mark=asterisk, mark options = {scale=1, solid}] table[x expr=\thisrow{x},y=u_pred2]{\fileDL};
\addplot[color=orange!80!black, style=dashed, line width=0.7 ,mark=asterisk, mark options = {scale=1, solid}] table[x expr=\thisrow{x},y=u_pred3]{\fileDL};
\addplot[color=gray!80!black, style=dashed, line width=0.7 ,mark=asterisk, mark options = {scale=1, solid}] table[x expr=\thisrow{x},y=u_pred4]{\fileDL};
\end{axis}
\end{tikzpicture}
\caption{Step $s=2$ (after training).}
\label{fig:step2_end_helmholtz_parametric}
\end{subfigure}
\vskip 0.5em
\begin{subfigure}[b]{0.49\textwidth}
\centering
\begin{tikzpicture}
\newcommand{\fileDL}{Chapters/3.Chapter/figures/numerical_results/DATA/results_helmholtz_-50--30_inverse100100_3refinements_endtoend/start/test_prediction_33nodes.csv}
\newcommand{\fileFEM}{Chapters/3.Chapter/figures/numerical_results/DATA/results_helmholtz_-50--30_inverse100100_3refinements_endtoend/fem/test_FEM_33nodes.csv}
\pgfplotsset{A/.append style={
	    xmin=-0.1,
	    xmax=1.1,
	    ymin=-2.2,
	    ymax=2.2,
	    ylabel={$u(x)$},
	    height=0.59*\textwidth,
	    width=\textwidth,
		xtick={0,0.25,...,1},
            }
}
\begin{axis}[A]
\addplot[color=green!30!white, line width=1., mark=*, mark options={scale=.8, solid}] table[x expr=\thisrow{x},y=u_FEM0]{\fileFEM};
\addplot[color=red!30!white, line width=1., mark=*, mark options={scale=.8, solid}] table[x expr=\thisrow{x},y=u_FEM1]{\fileFEM};
\addplot[color=blue!30!white, line width=1., mark=*, mark options={scale=.8, solid}] table[x expr=\thisrow{x},y=u_FEM2]{\fileFEM};
\addplot[color=orange!30!white, line width=1., mark=*, mark options={scale=.8, solid}] table[x expr=\thisrow{x},y=u_FEM3]{\fileFEM};
\addplot[color=gray!30!white, line width=1., mark=*, mark options={scale=.8, solid}] table[x expr=\thisrow{x},y=u_FEM4]{\fileFEM};
\addplot[color=green!80!black, style=dashed, line width=0.7 ,mark=asterisk, mark options = {scale=1, solid}] table[x expr=\thisrow{x},y=u_pred0]{\fileDL};
\addplot[color=red!80!black, style=dashed, line width=0.7 ,mark=asterisk, mark options = {scale=1, solid}] table[x expr=\thisrow{x},y=u_pred1]{\fileDL};
\addplot[color=blue!80!black, style=dashed, line width=0.7 ,mark=asterisk, mark options = {scale=1, solid}] table[x expr=\thisrow{x},y=u_pred2]{\fileDL};
\addplot[color=orange!80!black, style=dashed, line width=0.7 ,mark=asterisk, mark options = {scale=1, solid}] table[x expr=\thisrow{x},y=u_pred3]{\fileDL};
\addplot[color=gray!80!black, style=dashed, line width=0.7 ,mark=asterisk, mark options = {scale=1, solid}] table[x expr=\thisrow{x},y=u_pred4]{\fileDL};
\end{axis}
\end{tikzpicture}
\caption{Step $s=3$ (before training).}
\label{fig:step3_start_helmholtz_parametric}
\end{subfigure} \hfill
\begin{subfigure}[b]{0.49\textwidth}
\centering
\begin{tikzpicture}
\newcommand{\fileDL}{Chapters/3.Chapter/figures/numerical_results/DATA/results_helmholtz_-50--30_inverse100100_3refinements_endtoend/end/test_prediction_33nodes.csv}
\newcommand{\fileFEM}{Chapters/3.Chapter/figures/numerical_results/DATA/results_helmholtz_-50--30_inverse100100_3refinements_endtoend/fem/test_FEM_33nodes.csv}
\pgfplotsset{A/.append style={
	    xmin=-0.1,
	    xmax=1.1,
	    ymin=-2.2,
	    ymax=2.2,
	    height=0.59*\textwidth,
	    width=\textwidth,
		xtick={0,0.25,...,1},
            }
}
\begin{axis}[A]
\addplot[color=green!30!white, line width=1., mark=*, mark options={scale=.8, solid}] table[x expr=\thisrow{x},y=u_FEM0]{\fileFEM};
\addplot[color=red!30!white, line width=1., mark=*, mark options={scale=.8, solid}] table[x expr=\thisrow{x},y=u_FEM1]{\fileFEM};
\addplot[color=blue!30!white, line width=1., mark=*, mark options={scale=.8, solid}] table[x expr=\thisrow{x},y=u_FEM2]{\fileFEM};
\addplot[color=orange!30!white, line width=1., mark=*, mark options={scale=.8, solid}] table[x expr=\thisrow{x},y=u_FEM3]{\fileFEM};
\addplot[color=gray!30!white, line width=1., mark=*, mark options={scale=.8, solid}] table[x expr=\thisrow{x},y=u_FEM4]{\fileFEM};
\addplot[color=green!80!black, style=dashed, line width=0.7 ,mark=asterisk, mark options = {scale=1, solid}] table[x expr=\thisrow{x},y=u_pred0]{\fileDL};
\addplot[color=red!80!black, style=dashed, line width=0.7 ,mark=asterisk, mark options = {scale=1, solid}] table[x expr=\thisrow{x},y=u_pred1]{\fileDL};
\addplot[color=blue!80!black, style=dashed, line width=0.7 ,mark=asterisk, mark options = {scale=1, solid}] table[x expr=\thisrow{x},y=u_pred2]{\fileDL};
\addplot[color=orange!80!black, style=dashed, line width=0.7 ,mark=asterisk, mark options = {scale=1, solid}] table[x expr=\thisrow{x},y=u_pred3]{\fileDL};
\addplot[color=gray!80!black, style=dashed, line width=0.7 ,mark=asterisk, mark options = {scale=1, solid}] table[x expr=\thisrow{x},y=u_pred4]{\fileDL};
\end{axis}
\end{tikzpicture}
\caption{Step $s=3$ (after training).}
\label{fig:step3_end_helmholtz_parametric}
\end{subfigure}
\vskip 0.5em
\begin{subfigure}[b]{0.49\textwidth}
\centering
\begin{tikzpicture}
\newcommand{\fileDL}{Chapters/3.Chapter/figures/numerical_results/DATA/results_helmholtz_-50--30_inverse100100_3refinements_endtoend/start/test_prediction_65nodes.csv}
\newcommand{\fileFEM}{Chapters/3.Chapter/figures/numerical_results/DATA/results_helmholtz_-50--30_inverse100100_3refinements_endtoend/fem/test_FEM_65nodes.csv}
\pgfplotsset{A/.append style={
	    xmin=-0.1,
	    xmax=1.1,
	    xlabel={$x$},
	    ymin=-2.2,
	    ymax=2.2,
	    ylabel={$u(x)$},
	    height=0.59*\textwidth,
	    width=\textwidth,
		xtick={0,0.25,...,1},
            }
}
\begin{axis}[A]
\addplot[color=green!30!white, line width=1., mark=*, mark options={scale=.8, solid}] table[x expr=\thisrow{x},y=u_FEM0]{\fileFEM};
\addplot[color=red!30!white, line width=1., mark=*, mark options={scale=.8, solid}] table[x expr=\thisrow{x},y=u_FEM1]{\fileFEM};
\addplot[color=blue!30!white, line width=1., mark=*, mark options={scale=.8, solid}] table[x expr=\thisrow{x},y=u_FEM2]{\fileFEM};
\addplot[color=orange!30!white, line width=1., mark=*, mark options={scale=.8, solid}] table[x expr=\thisrow{x},y=u_FEM3]{\fileFEM};
\addplot[color=gray!30!white, line width=1., mark=*, mark options={scale=.8, solid}] table[x expr=\thisrow{x},y=u_FEM4]{\fileFEM};
\addplot[color=green!80!black, style=dashed, line width=0.7 ,mark=asterisk, mark options = {scale=1, solid}] table[x expr=\thisrow{x},y=u_pred0]{\fileDL};
\addplot[color=red!80!black, style=dashed, line width=0.7 ,mark=asterisk, mark options = {scale=1, solid}] table[x expr=\thisrow{x},y=u_pred1]{\fileDL};
\addplot[color=blue!80!black, style=dashed, line width=0.7 ,mark=asterisk, mark options = {scale=1, solid}] table[x expr=\thisrow{x},y=u_pred2]{\fileDL};
\addplot[color=orange!80!black, style=dashed, line width=0.7 ,mark=asterisk, mark options = {scale=1, solid}] table[x expr=\thisrow{x},y=u_pred3]{\fileDL};
\addplot[color=gray!80!black, style=dashed, line width=0.7 ,mark=asterisk, mark options = {scale=1, solid}] table[x expr=\thisrow{x},y=u_pred4]{\fileDL};
\end{axis}
\end{tikzpicture}
\caption{Step $s=4$ (before training).}
\label{fig:step4_start_helmholtz_parametric}
\end{subfigure} \hfill
\begin{subfigure}[b]{0.49\textwidth}
\centering
\begin{tikzpicture}
\newcommand{\fileDL}{Chapters/3.Chapter/figures/numerical_results/DATA/results_helmholtz_-50--30_inverse100100_3refinements_endtoend/end/test_prediction_65nodes.csv}
\newcommand{\fileFEM}{Chapters/3.Chapter/figures/numerical_results/DATA/results_helmholtz_-50--30_inverse100100_3refinements_endtoend/fem/test_FEM_65nodes.csv}
\pgfplotsset{A/.append style={
	    xmin=-0.1,
	    xmax=1.1,
	    xlabel={$x$},
	    ymin=-2.2,
	    ymax=2.2,
	    height=0.59*\textwidth,
	    width=\textwidth,
		xtick={0,0.25,...,1},
            }
}
\begin{axis}[A]
\addplot[color=green!30!white, line width=1., mark=*, mark options={scale=.8, solid}] table[x expr=\thisrow{x},y=u_FEM0]{\fileFEM};
\addplot[color=red!30!white, line width=1., mark=*, mark options={scale=.8, solid}] table[x expr=\thisrow{x},y=u_FEM1]{\fileFEM};
\addplot[color=blue!30!white, line width=1., mark=*, mark options={scale=.8, solid}] table[x expr=\thisrow{x},y=u_FEM2]{\fileFEM};
\addplot[color=orange!30!white, line width=1., mark=*, mark options={scale=.8, solid}] table[x expr=\thisrow{x},y=u_FEM3]{\fileFEM};
\addplot[color=gray!30!white, line width=1., mark=*, mark options={scale=.8, solid}] table[x expr=\thisrow{x},y=u_FEM4]{\fileFEM};
\addplot[color=green!80!black, style=dashed, line width=0.7 ,mark=asterisk, mark options = {scale=1, solid}] table[x expr=\thisrow{x},y=u_pred0]{\fileDL};
\addplot[color=red!80!black, style=dashed, line width=0.7 ,mark=asterisk, mark options = {scale=1, solid}] table[x expr=\thisrow{x},y=u_pred1]{\fileDL};
\addplot[color=blue!80!black, style=dashed, line width=0.7 ,mark=asterisk, mark options = {scale=1, solid}] table[x expr=\thisrow{x},y=u_pred2]{\fileDL};
\addplot[color=orange!80!black, style=dashed, line width=0.7 ,mark=asterisk, mark options = {scale=1, solid}] table[x expr=\thisrow{x},y=u_pred3]{\fileDL};
\addplot[color=gray!80!black, style=dashed, line width=0.7 ,mark=asterisk, mark options = {scale=1, solid}] table[x expr=\thisrow{x},y=u_pred4]{\fileDL};
\end{axis}
\end{tikzpicture}
\caption{Step $s=4$ (after training).}
\label{fig:step4_end_helmholtz_parametric}
\end{subfigure}
\caption{Four steps of the \acs{DeepFEM} for problem \eqref{eq:parmetric bvp} with $-50<\alpha<-30$ employing the $H^1$-norm error in the loss function. $u_\text{FEM}(\alpha)$ and $u_\text{NN}(\alpha)$ denote the FEM solution and DeepFEM prediction for the $\alpha$ parameter coefficient, respectively.}
\label{fig:helmholtz_parametric}
\end{figure}

Utilizing the inverses is equivalent to applying the loss:
\begin{equation} \label{error alternative}
\mathcal{L}(\theta; D) = \frac{1}{|D|} \sum_{\alpha\in D} \lVert \mathbf{u}_{\text{NN}}(\alpha)  - \mathbf{u}_{\text{FEM}}(\alpha) \rVert_{\mathbf{K}+\mathbf{M}},
\end{equation} where $\mathbf{u}_{\text{NN}}(\alpha)$ stands for the \acs{DeepFEM} prediction and $\mathbf{u}_{\text{FEM}}(\alpha)$ for the FEM solution vectors for the $\alpha$ parameter. Precalculating all these vectors and employing \eqref{error alternative} as the loss during training could be significantly less expensive than computing/dealing with the inverse matrices. Nevertheless, we have computed inverses to be consistent with the presented arguing, and we refrained from reporting times since an efficient solver should avoid inverting a matrix explicitly \cite{abur1988parallel, pardo2004integration, moshfegh2017direct}.

\newcommand{\R}{\mathbb{R}}
\newcommand{\V}{\mathbb{V}}
\newcommand{\M}{\mathbb{M}}
\newcommand{\n}{\text{NN}}
\renewcommand{\r}{\text{RNE}}
\newcommand{\J}{\mathcal{J}}
\newcommand{\D}{\mathcal{D}}
\newcommand{\I}{\mathcal{I}}
\newcommand{\F}{\mathcal{F}}
\newcommand{\0}{\mathbf{0}}
\newcommand{\W}{\mathbf{W}}
\newcommand{\bb}{\mathbf{b}}

\chapter[The Deep Double Ritz Method]{The Deep Double Ritz Method}\label{chapter4}

\begin{quote}
\textbf{Summary.} Residual minimization is a widely used technique for solving Partial Differential Equations in variational form. It minimizes the dual norm of the residual, which naturally yields a saddle-point (min–max) problem over the trial and test spaces. In the context of Neural Networks, we can address this min–max approach by employing one Neural Network to seek the trial minimum while another Neural Network seeks the test maximizers. However, the resulting method is numerically unstable as we approach the solution. To overcome this, we reformulate the residual minimization as an equivalent minimization of a Ritz functional fed by optimal test functions computed from another Ritz functional minimization. We call the resulting scheme the \emph{Deep Double Ritz Method}, which combines two Neural Networks for approximating trial and optimal test functions along a nested double Ritz minimization strategy. Numerical results on different diffusion and convection problems support the robustness of our method up to the approximation properties of the considered Neural Networks and the training capacity of the optimizers. \emph{Refer to \cite{uriarte2023deep} for the published version}.
\end{quote}

\section{Introduction}

Within the variational framework introduced in \Cref{section1.5.2}, residual minimization reads as a saddle-point (min-max) problem as follows:
\begin{equation}\label{intro1}
	\min_{u\in \mathbb{U}} \max_{v\in\mathbb{V}\setminus\{0\}} \frac{\langle Bu-l, v\rangle_{\mathbb{V}'\times \mathbb{V}}}{\Vert v\Vert_\mathbb{V}},
\end{equation} where $\mathbb{U}$ and $\mathbb{V}$ are the trial and test spaces, $B:\mathbb{U}\longrightarrow\mathbb{V}'$ is the differential operator governing the considered \ac{BVP} in variational form, $\mathbb{V}'$ is the dual space of $\mathbb{V}$, and $l\in\mathbb{V}'$ is the right-hand side. In \cite{zang2020weak, bao2020numerical}, the authors proposed addressing this optimization scheme employing Generative Adversarial Networks (GANs) \cite{goodfellow2014generative, goodfellow2020generative} by approximating $u$ and $v$ with two NNs. Unfortunately, this approach presents a severe numerical limitation: the Lipschitz continuity constant of the test maximizers with respect to the trial functions might become arbitrarily large when approaching the exact solution. Moreover, the corresponding test maximizer is highly non-unique in the limit. In consequence, we end up with an inherent lack of numerical stability that is easily confirmed by numerical experiments with simple model problems.

To overcome the above limitations, we reformulate residual minimization as a minimization of a Ritz functional fed by optimal test functions \cite{demkowicz2011class,  demkowicz2014overview}. Since optimal test functions are generally unavailable, we compute them for each trial function using another Ritz method. Thus, the resulting scheme is a nested double-loop Ritz minimization method: the outer loop seeks the trial solution, while the inner loop seeks the optimal test function for each trial function.  We call the resulting scheme the \emph{Double Ritz Method}.

In some occasions, the trial-to-test operator that maps each trial function with the corresponding optimal test function is available. For example, when the problem is symmetric and positive-definite, and we consider the norm induced by the bilinear form for the trial and test spaces,  the trial-to-test operator is the identity; or when selecting the strong variational formulation, the trial-to-test operator is the one given by the \ac{PDE} operator. In these cases, the Double Ritz Method reduces to a single-loop Ritz minimization (this will be further discussed in \Cref{Generalized Ritz method}).  Thus, the Double Ritz Method is a general method for solving \acp{PDE} in different variational forms, which in some particular cases simplifies into a single-loop Ritz minimization method.

Thanks to \acp{NN}, we find a simple and advantageous computational framework to approximate and connect the trial and test functions between the inner- and outer-loop minimizations in the Double Ritz Method, a task that is challenging to tackle with traditional numerical methods. Herein, we propose using one \ac{NN} to represent the trial functions and another \ac{NN} to represent the local actions of the trial-to-test operator. Thus, the composition of both \acp{NN} represents the (optimal) test functions,  and we preserve the trial dependence of the test functions during the entire process.  We call the resulting \ac{NN}-based method the \emph{\acf{D2RM}}.

While the \ac{D2RM} replicates existing residual minimization methods in the context of \acp{NN}, we fall short of providing a detailed mathematical convergence analysis due to the manifold structure of \acp{NN} that departs from the traditional vector-space-based mathematical approach (recall \Cref{chapter2}). Related to this, we encounter the usual drawbacks of lack of convexity between the objective/loss function with respect to the trainable parameters,  which prevents us from making a proper diagnosis of the optimizer during training. 

The remainder of this chapter is organized as follows. \Cref{section:Mathematical framework} formalizes the variational setting introduced in \Cref{section1.5.2} and derives the Double Ritz Method at the continuum level. \Cref{section:Adversarial Neural Networks} introduces the \ac{D2RM} within the \ac{NN} framework.  \Cref{section:Implementation} provides implementation details of the addressed methods and \Cref{section:Numerical results} develops on numerical experimentation. 

\section{From residual to Ritz minimization}
\label{section:Mathematical framework}

We introduce the residual minimization framework, followed by a saddle-point reformulation and an alternative Double Ritz scheme at the continuum level.  Subsequently, we describe three particular cases for which the Double Ritz method simplifies into a single Ritz method.
\subsection{Residual minimization}

Let
\begin{subequations}\label{PG}
\begin{equation}
\displaystyle{\left\{
\begin{tabular}{l}
Find $u^* \in \mathbb{U}$ such that\\
$b(u^*,v)=l(v), \; \forall v\in \V$,
\end{tabular}
\right.}
\end{equation} where $\mathbb{U}$ and $\V$ are real Hilbert trial and test spaces, respectively, $b:\mathbb{U}\times \V\longrightarrow\R$ is a bilinear form, and $l:\V\longrightarrow\R$ is a continuous linear functional. Equivalently, in operator form:
\begin{equation}
\displaystyle{\left\{
\begin{tabular}{l}
Find $u^*\in \mathbb{U}$ such that\\
$Bu^*=l$,
\end{tabular}
\right.} 
\end{equation}
\end{subequations}
where $B:\mathbb{U}\longrightarrow \V'$ is the operator defined by $\langle Bu, v\rangle_{\mathbb{V}'\times\mathbb{V}} := b(u,v)$, $\V'$ denotes the topological dual of $\V$, and $l\in \V'$.  The equivalent residual minimization formulation reads as
\begin{equation}\label{DualRes}
    u^* = \arg\min_{u\in \mathbb{U}} \Vert Bu-l\Vert_{\V'},
\end{equation} where $Bu-l\in\V'$ is the \emph{residual} for each trial function $u\in\mathbb{U}$.

To guarantee well-posedness of \eqref{PG}--\eqref{DualRes}, we assume the hypotheses of the Babu{\v{s}}ka–Lax–Milgram Theorem \cite{babuvska1971error}, which according to our presentation 
translate into that $B$ is an isomorphism that is bounded from above and below, i.e., there exist some positive constants $\gamma \leq M$ such that it satisfies
\begin{equation}
\gamma \Vert u\Vert_\mathbb{U} \leq \Vert Bu\Vert_{\mathbb{V}'} \leq M \Vert u\Vert_\mathbb{U}, \qquad u\in\mathbb{U}.
\end{equation}  Then, the error in $\mathbb{U}$ is equivalent to the residual in $\V'$ in the following sense:
\begin{equation}\label{residual_error_relation}
\frac{1}{M} \Vert Bu-l \Vert_{\mathbb{V}'} \leq \Vert u-u^*\Vert_{\mathbb{U}} \leq \frac{1}{\gamma} \Vert Bu-l \Vert_{\mathbb{V}'}, \qquad u\in\mathbb{U}.
\end{equation} Below, we examine two alternatives to evaluate and minimize the residual in its dual norm.

\subsection{Saddle-point problem (min-max optimization)}\label{section:Approach1}
The norm in $\mathbb{V}'$ is defined in terms of the norm in $\V$ as 
\begin{equation}\label{dual}
\Vert t\Vert_{\mathbb{V}'}:=\sup_{v\in\mathbb{V}\setminus\{0\}} \frac{\langle t, v\rangle_\mathbb{V}}{\Vert v\Vert_{\mathbb{V}}} = \sup_{v\in\mathbb{V}\setminus\{0\}} \left\langle t, \frac{v}{\Vert v\Vert_\mathbb{V}}\right\rangle_{\mathbb{V}'\times\mathbb{V}} ,  \qquad t\in\mathbb{V}'.
\end{equation} Hence, combining \eqref{DualRes} and \eqref{dual} yields
\begin{subequations}
\begin{equation}
    u^* = \arg\min_{u\in \mathbb{U}} \max_{v\in\mathbb{V}\setminus\{0\}} \mathcal{F}_\text{min}^\text{max}(u,v)\label{infsup_residual}
\end{equation} with 
\begin{equation}
\mathcal{F}_\text{min}^\text{max}(u,v):=\left\langle Bu-l, \frac{v}{\Vert v\Vert_\mathbb{V}}\right\rangle_{\mathbb{V}'\times\mathbb{V}} = b\left(u,\frac{v}{\Vert v\Vert_\mathbb{V}}\right) - l\left(\frac{v}{\Vert v\Vert_\mathbb{V}}\right).
\end{equation}
\end{subequations} When $\mathcal{F}_\text{min}^\text{max}$ is fed with the exact solution $u^*$, the resulting operator becomes the null functional, i.e., $\mathcal{F}_\text{min}^\text{max}(u^*, \cdot)=0$, which implies a highly non-uniqueness of the test maximizer in the limit---indeed, any element in $\mathbb{V}\setminus\{0\}$ is a test maximizer for $u^*$. Outside this singular case, the operator that maps each trial function $u\in\mathbb{U}$ to its unitary test maximizer is well-defined but not Lipschitz continuous when approaching the exact solution $u^*$, leading to an unstable numerical method. We formalize this in the following two items: 
\begin{itemize}
	\item Let $v_\text{max}:\mathbb{U}\setminus\{u^*\}\longrightarrow\V$ be the mapping that for each trial function returns the test maximizer of the actions of the residual $Bu-l\in\V'$ over the unitary sphere, i.e.,
\begin{equation}
v_\text{max}(u) := \arg\max_{\Vert v\Vert_\V = 1} \langle Bu-l, v\rangle_{\mathbb{V}'\times\mathbb{V}}.\end{equation} Then, $v_\text{max}(u)\in\mathbb{V}$ is unique for each $u\in\mathbb{U}\setminus\{u^*\}$. 
\begin{proof} Let $u\in\mathbb{U}\setminus\{u^*\}$. By the Riesz Representation Theorem, there exists a unique $r_u\in\V$ such that
\begin{equation}\label{appendix1}
(r_u,v)_\mathbb{V}=\langle Bu-l, v\rangle_{\mathbb{V}'\times\mathbb{V}},\qquad \forall v\in\V,
\end{equation} and
\begin{equation}\label{test_maximizer_norm}
\Vert r_u\Vert_\mathbb{V} = \Vert Bu-l\Vert_{\mathbb{V}'} = \max_{\Vert v\Vert_\mathbb{V} = 1} \langle Bu-l, v\rangle_{\mathbb{V}'\times\mathbb{V}}.
\end{equation} Let $v_\text{max}=v_\text{max}(u)$ be a test maximizer of \eqref{test_maximizer_norm} in the unitary sphere of $\mathbb{V}$, i.e.,
\begin{equation}
\langle Bu-l, v_\text{max}\rangle_{\mathbb{V}'\times\mathbb{V}}=\max_{\Vert v\Vert_{\mathbb{V}} = 1} \langle Bu-l, v\rangle_{\mathbb{V}'\times\mathbb{V}}.
\end{equation} By the Cauchy-Schwarz inequality:
\begin{subequations}
\begin{align}
0 &< \Vert r_u\Vert_\mathbb{V} = \Vert Bu-l\Vert_{\mathbb{V}'} = \langle Bu-l, v_\text{max}\rangle_{\mathbb{V}'\times\mathbb{V}}=(r_u, v_\text{max})_\mathbb{V}\\
&\leq \Vert r_u\Vert_{\mathbb{V}} \; \Vert v_\text{max}\Vert_{\mathbb{V}},
\end{align} 
\end{subequations} where the equality holds if and only if $v_\text{max} = \lambda r_u$ for some $\lambda> 0$.  Then,  $v_\text{max} = \frac{r_u}{\Vert r_u\Vert_\V}$ is unique.
\end{proof}

\item $v_\text{max}$ is not Lipschitz continuous around any (reduced) neighborhood of $u^*$, i.e.,  there does not exist a constant $0<C<\infty$ such that 
\begin{equation}\label{appendix continuity}
\Vert v_\text{max}(u_1)-v_\text{max}(u_2)\Vert_\V \leq C\Vert u_1-u_2\Vert_\mathbb{U}, \qquad\forall u_1,u_2\in\mathbb{U}\setminus\{u^*\}.
\end{equation}
\begin{proof}
Assume by contradiction that there exists $0<C<\infty$ such that \eqref{appendix continuity} holds, and let $u_2=2u^*-u_1$ with $u_1\neq u^*$. Then,
\begin{equation}
2 = \left\Vert v_\text{max}(u_1) - v_\text{max}(u_2)\right\Vert_{\V} \leq C \Vert u_1-u_2\Vert_\mathbb{U} = 2 C \Vert u_1-u^*\Vert_\mathbb{U}.
\end{equation} Letting  $u_1\to u^*$, we obtain $C\to\infty$.
\end{proof}
\end{itemize}

\subsection{The Double Ritz Method with optimal test functions}\label{section:Double Ritz Method with optimal test functions}
The Riesz Representation Theorem allows us to work isometrically in the test space instead of in its dual. In particular, for the residual, we have
\begin{equation}
\Vert Bu-l\Vert_{\mathbb{V}'} = \Vert \mathfrak{R}^{-1}_\mathbb{V}(Bu-l)\Vert_{\mathbb{V}},\qquad u\in\mathbb{U},
\end{equation} where $\mathfrak{R}_\mathbb{V}: \V\ni v\longmapsto (v,\cdot)_\mathbb{V}\in\mathbb{V}'$ denotes the Riesz operator.  This relation suggests considering the \emph{trial-to-test} operator $T:\mathbb{U}\longrightarrow\mathbb{V}$ defined by $T:=\mathfrak{R}^{-1}_\mathbb{V} B$ \cite{demkowicz2014overview} since it relates the error in $\mathbb{U}$ with the Riesz representative of the residual in $\mathbb{V}$, 
\begin{equation} \label{error_residual_representation}
T(u-u^*) = \mathfrak{R}^{-1}_\mathbb{V} (Bu-l) \in \mathbb{V},\qquad u\in\mathbb{U}.
\end{equation} The images of trial functions through $T$ are known as \emph{optimal test functions} \cite{demkowicz2011class}, and they allow us to rewrite \eqref{PG} in terms of the following symmetric and positive-definite variational problem:
\begin{equation}\label{Tproblem}
\displaystyle{\left\{
\begin{tabular}{l}
Find $u^* \in \mathbb{U}$ such that\\
$(T u^*, T u)_\mathbb{V} = l(T u), \; \forall u\in \mathbb{U}$.
\end{tabular}
\right.} 
\end{equation}
\begin{proof} On the one hand, $b(u,v)=\langle Bu,v\rangle_{\mathbb{V}'\times\mathbb{V}} = (Tu,v)_\mathbb{V}$ for all $u\in\mathbb{U}$ and all $v\in\mathbb{V}$. On the other hand, $T$ is an isomorphism because so is $B$. Hence, testing with all $v\in\mathbb{V}$ is equivalent to testing with $Tu\in\mathbb{V}$ for all $u\in\mathbb{U}$.\end{proof} 

Because the bilinear form $(T\cdot,T\cdot):\mathbb{U}\times\mathbb{U}\longrightarrow\mathbb{R}$ of \eqref{Tproblem} is symmetric and positive definite, we can reformulate \eqref{Tproblem} in terms of a quadratic-functional minimization combining the ``usual'' Ritz functional and the trial-to-test operator as follows:
\begin{subequations}\label{Ritz_u}
\begin{equation} 
    u^* = \arg\min_{u\in \mathbb{U}} \mathcal{F}_T(u),
\end{equation} with
\begin{equation}
\mathcal{F}_T(u):=\frac{1}{2}\Vert T u\Vert_\mathbb{V}^2 -l(T u).
\end{equation}
\end{subequations}
One might interpret that $\mathcal{F}_T$ acts as a generalization of the ``usual'' Ritz functional\footnote{The Ritz method \cite{ritz1909neue} is typically reserved for symmetric and positive-definite problems due to the natural ``energy'' concept arising from those kinds of problems. In this way, the Ritz functional is typically presented as $\mathcal{F}(\cdot)=\frac{1}{2}b(\cdot, \cdot)-l(\cdot)$, where the inner product of the underlying Hilbert space is selected as the bilinear form, i.e., $(\cdot,\cdot)_\mathbb{V}=b(\cdot,\cdot)$.} $\mathcal{F}(\cdot)=\frac{1}{2}\Vert\cdot\Vert_\mathbb{V}^2 -l (\cdot)$ into problems that are not necessarily symmetric or positive-definite. The relation between residual minimization \eqref{DualRes} and our proposed generalized Ritz minimization \eqref{Ritz_u} is
\begin{align}\label{residual_min_Ritzu}
\mathcal{F}_T(u)-\mathcal{F}_T(u^*) = \frac{1}{2}\Vert Bu-l\Vert_{\mathbb{V}'}^2.
\end{align}

\begin{proof}
According to \eqref{Tproblem}, we can write $\mathcal{F}_T(u)=\frac{1}{2}(Tu,Tu)_\mathbb{V}-(Tu^*,Tu)_{\mathbb{V}}$ for all $u\in\mathbb{U}$. Straightforward calculations lead to the desired identity.
\end{proof}

Now,  the challenge is to compute $Tu\in\V$ when iterating along $u\in \mathbb{U}$. Because finding $Tu\in\V$ is equivalent to solving the symmetric and positive-definite variational problem
\begin{equation}\label{VarVh}
\displaystyle{\left\{
\begin{tabular}{l}
Find $Tu \in \V$ such that\\
$(T u,  v)_\V = b(u,v), \; \forall v\in \V$,
\end{tabular}
\right.} 
\end{equation} following the same arguing as before, we can reformulate it as a corresponding Ritz minimization for each given/fixed $u\in\mathbb{U}$ as follows:
\begin{subequations}
\begin{equation} \label{Ritz_v} 
   T u = \arg\min_{v\in \mathbb{V}} \mathcal{F}_u^{\text{opt}}(v),
 \end{equation} with
 \begin{equation}
  \mathcal{F}_u^{\text{opt}}(v):=\frac{1}{2}\Vert v\Vert_\mathbb{V}^2 - b(u,v).
 \end{equation}
 \end{subequations} However,  characterizing $Tu$ by means of \eqref{Ritz_v} remains impractical when iterating along $u\in\mathbb{U}$ to minimize \eqref{Ritz_u} in the absence of knowledge of $T$\footnote{Intuitively, the minimizer in \eqref{Ritz_v} is delivered with an \emph{implicit} dependence on $u\in\mathbb{U}$. Consequently, trying to iterate along $\mathbb{U}$ via \eqref{Ritz_u} becomes challenging due to the absence of explicit dependencies to deal with.}. To overcome this,  instead of iterating along $\mathbb{V}$ to seek the optimal test function for a given $u\in\mathbb{U}$, we consider seeking along a ``convenient'' family $\mathbb{M}$ of operators from $\mathbb{U}$ to $\mathbb{V}$\footnote{At a continuum level, we should ensure that $\mathbb{M}$ is able to produce optimal test functions for any ``iterable'' trial function, e.g., by assuming that there exists a family of operators $L_u$ in $\mathbb{M}$ such that $Tu\in\text{Im}(L_u)$ for all ``iterable'' $u\in\mathbb{U}$ (note that this does not imply $T\in\mathbb{M}$).}. Then, we reformulate \eqref{Ritz_v} as seeking a candidate $\tau_u\in\mathbb{M}$ that when acting on $u\in\mathbb{U}$ delivers $T u\in\mathbb{V}$, namely, given $u\in\mathbb{U}$,
\begin{equation} \label{Ritz_T}
    \tau_{u} = \arg\min_{\tau\in\M} \mathcal{F}_u^{\text{opt}}(\tau(u)).
\end{equation} This provides a framework where a minimizer $\tau_u$ of \eqref{Ritz_T} acts as $T$ \emph{only} at the given $u\in\mathbb{U}$,  i.e., $\tau_u(u)=Tu$ but $\tau_u(w)$ might differ from $T w$ whenever $u\neq w$.  Moreover, depending on the construction of $\mathbb{M}$,  $\tau_u$ is possibly non-unique.

Hence,  by performing a nested minimization of \eqref{Ritz_T} within \eqref{Ritz_u},  we solve the problem at hand without explicitly dealing with the full operator $T$.  We call this the \emph{Double Ritz Method}. \Cref{alg:outer-inner loop} depicts its nested-loop optimization strategy that iterates separately either with elements in $\mathbb{U}$ or in $\mathbb{M}$.

\begin{algorithm}
\caption{Double Ritz Method}\label[algorithm]{alg:outer-inner loop}
Initialize $u\in \mathbb{U}$ and $\tau\in\M$\;
\While{\upshape not converged}{
    Find $\displaystyle\tau_{u} \in \arg\min_{\tau\in \M} \mathcal{F}_u^{\text{opt}}(\tau(u))$\;
    $u =$ following candidate in $\mathbb{U}$ minimizing $\mathcal{F}_{\tau_u}(u)$\;
}
\Return $u$
\end{algorithm}

\subsection{Generalized Ritz Methods}\label{Generalized Ritz method}

We show three scenarios where the Double Ritz Method simplifies into a single Ritz minimization.

\begin{enumerate}
\item[(a)] \textbf{(Traditional) Ritz Method.} If the bilinear form $b(\cdot,\cdot)$ is symmetric and positive definite, we have $\V=\mathbb{U}$. Considering the inner product as the bilinear form yields $T u = u$ for all $u\in \mathbb{U}$,  i.e., $T$ is the identity operator.

\item[(b)] \textbf{Strong formulation.} Let $A:\D(A)\longrightarrow L^2(\Omega)$ denote a PDE operator with domain $\D(A)$.  Considering $b(u,v)=(Au,v)_{L^2(\Omega)}$, $\mathbb{U}=\D(A)$ equipped with the graph norm,  and $\V=L^2(\Omega)$, we obtain $Tu = Au$ for all $u\in \mathbb{U}$, i.e., $T$ is the PDE operator.

\item[(c)] \textbf{Ultraweak formulation.} Let $A':\D(A')\longrightarrow L^2(\Omega)$ denote the adjoint operator of the PDE operator $A:\D(A)\longrightarrow L^2(\Omega)$, i.e., $A'$ satisfies $(Au,v)_{L^2(\Omega)}=(u,A' v)_{L^2(\Omega)}$ for all $u\in\mathcal{D}(A)$ and $v\in\mathcal{D}(A')$. Considering $b(u,v)=(u,A' v)_{L^2(\Omega)}$, $\mathbb{U}=L^2(\Omega)$,  and $\V=\D(A')$ equipped with the graph norm, we obtain $A' Tu = u$ for all $u\in \mathbb{U}$.

Following the ideas of \cite{brunken2019parametrized, demkowicz2020dpg}, we can first solve for the optimal test function of the trial solution:
\begin{equation}\label{adjoint_minimization}
Tu^*=\arg\min_{v\in \V} \F'(v), \qquad \F'(v):=\frac{1}{2}\Vert A' v\Vert^2_{L^2(\Omega)}-l(v), \end{equation}
and then apply $A'$ to the minimizer to recover the trial solution: $u^*=A'Tu^*$.  We call \emph{Adjoint Ritz Method} to this Ritz minimization with post-processing based on the use of the $A'$ operator.
\end{enumerate}

To illustrate the above three cases, we show a simple \ac{PDE} problem with different variational formulations. 

\begin{example}[Poisson's Equation]\label{Poisson's Equation}
Let $f\in L^2(0,1)$ and consider the following 1D pure diffusion equation with homogeneous Dirichlet boundary conditions over $\Omega=(0,1)$:
\begin{equation}\label{1DLap}
\begin{cases}
-u''=f, &\\
u(0)=u(1)=0. &
\end{cases}
\end{equation} We multiply the PDE by a test function and integrate over $\Omega$. Depending on the number of times we integrate by parts to derive the variational formulation, we obtain the following scenarios:
\begin{enumerate}
\item \textbf{Strong formulation.} Without integration by parts. Then, $\mathbb{U}=H^2(0,1)\cap H^1_0(0,1)$, $\mathbb{V}=L^2(0,1)$, $b(u,v)=\int_0^1 -u''v$, and $Tu=-u''$ for all $u\in \mathbb{U}$. This is the case {\upshape (b)} above.
\item \textbf{Weak formulation.} We integrate by parts once, passing one derivative from the trial to the test function. Then, $\mathbb{U}=H_0^1(0,1)=\V$, $b(u,v)=\int_0^1 u'v'$, $(v_1, v_2)_\mathbb{V}=b(v_1,v_2)$, and $Tu=u$ for all $u\in \mathbb{U}$.  This is the case {\upshape (a)} above.
\item \textbf{Ultraweak formulation.}  We integrate by parts twice, passing the two derivatives of the trial to the test function. Then, $\mathbb{U}=L^2(0,1)$, $\mathbb{V}=H^2(0,1)\cap H^1_0(0,1)$,  $(v_1,v_2)_\mathbb{V}=(-v_1'',-v_2'')_{L^2(\Omega)}$,  $b(u,v)=\int_0^1-uv''$,  and $-(Tu)''=u$ for all $u\in \mathbb{U}$. This is the case {\upshape (c)} above.
\end{enumerate}
\end{example}

All three variational formulations are valid, although they exhibit different convergence behaviors that go beyond the scope of this work.

\section{Approximation with Neural Networks}
\label{section:Adversarial Neural Networks}

In practice,  instead of seeking in $\mathbb{U}$, $\mathbb{V}$,  and $\mathbb{M}$, we seek along corresponding computationally accessible subsets generated by \acp{NN}.  This section elaborates on the NN setting for the computability of the proposed methods following the presentation and notation of \Cref{chapter2}.

\subsection{Parameterization}

Let $u_\text{NN}$, $v_\text{NN}$ and $\tau_\text{NN}$ be \ac{FFNN} architectures with corresponding sets of learnable parameters $\theta_u$, $\theta_v$ and $\theta_\tau$ and sets of realizations $\mathbb{U}_\text{NN}$, $\mathbb{V}_\text{NN}$, and $\mathbb{M}_\text{NN}$, respectively.

To ensure the containment of the set of realizations within the trial space (i.e., $\mathbb{U}_\text{NN}\subset \mathbb{U}$), we select smooth activation functions (e.g., $\tanh$) and strongly impose Dirichlet boundary conditions by considering $u_\text{NN}^\xi:\Omega\longrightarrow\mathbb{R}$ defined by
\begin{equation}\label{neural_network_architecture_dirichlet}
u_\text{NN}^\xi(x):=\xi(x) u_\text{NN}(x),\qquad x\in\Omega,
\end{equation} where $u_\text{NN}:\Omega\longrightarrow\mathbb{R}$ is the previously considered boundary-free \ac{FFNN}  architecture and $\xi:\Omega\longrightarrow\mathbb{R}$ is a non-trainable smooth cut-off function that vanishes only on the Dirichlet boundary. For simplicity, we drop the superscript $\xi$ notation, meaning that $u_\text{NN}$ will directly refer to $u_\text{NN}^\xi$ from now on. We apply similar criteria to $v_\text{NN}$ and $\tau_\text{NN}$.

\Cref{table_methods_NNs} summarizes some of the methods introduced in \Cref{section:Mathematical framework} within the framework of NNs.

\begin{table}[htb]
\centering
\resizebox{\textwidth}{!}{%
\begin{tabular}{@{}|r|l|l|@{}}
\toprule
\multicolumn{1}{|c|}{\textbf{Method}} & \multicolumn{1}{c|}{\textbf{Objective function(s)}}                                                                                                                                                                                                                    & \multicolumn{1}{c|}{\textbf{Optimization}}                                                                                                                                                                                  \\ \midrule
Weak Adversarial Networks             & $\displaystyle \mathcal{F}_\text{min}^\text{max}(u_\text{NN},v_\text{NN})=b\left(u_\text{NN},\frac{v_\text{NN}}{\Vert v_\text{NN}\Vert_\mathbb{V}}\right)- l\left(\frac{v_\text{NN}}{\Vert v_\text{NN}\Vert_\mathbb{V}}\right)$                                                                                                                                                                     & $\displaystyle u^*_\text{NN} = \arg\min_{u_\text{NN}\in\mathbb{U}_\text{NN}} \max_{v_\text{NN}\in\mathbb{V}_\text{NN}}  \mathcal{F}_\text{min}^\text{max}(u_\text{NN},v_\text{NN})$                                                                                                                         \\ \midrule
Generalized Deep Ritz Method          & \begin{tabular}[c]{@{}l@{}}$\displaystyle \;\;\;\;\;\;\;\;\;\;\;\mathcal{F}_T(u_\text{NN})=\tfrac{1}{2} \lVert Tu_\text{NN}\rVert_V^2 - l(Tu_\text{NN})$\vspace{0.2cm}  \\ \;\;\;\;\;$T$ is the ``available'' trial-to-test operator                                                                                                                                                                                  \end{tabular} & $\displaystyle u^*_\text{NN} = \arg\min_{u_\text{NN}\in\mathbb{U}_\text{NN}} \F_T(u_\text{NN})$                                                                                                                                                                 \\ \midrule
Deep Double Ritz Method               & \begin{tabular}[c]{@{}l@{}}$\displaystyle \;\;\;\;\;\;\;\;\mathcal{F}_{\tau_\text{NN}}(u_\text{NN})=\tfrac{1}{2} \lVert \tau_\text{NN}(u_\text{NN})\rVert_V^2 - l(\tau_\text{NN}(u_\text{NN}))$\vspace{0.2cm}  \\ $\displaystyle \mathcal{F}_{u_\text{NN}}^{\text{opt}}(\tau_\text{NN}(u_\text{NN})) =\tfrac{1}{2} \lVert \tau_\text{NN}(u_\text{NN})\rVert_V^2 - b(u_\text{NN}, \tau_\text{NN}(u_\text{NN}))$\end{tabular} & \begin{tabular}[c]{@{}l@{}}$\displaystyle u^*_\text{NN} = \arg\min_{u_\text{NN}\in\mathbb{U}_\text{NN}} \mathcal{F}_{\tau_{u_\text{NN}}}(u_\text{NN})$ \\ $\displaystyle \tau_{u_\text{NN}} = \arg\min_{\tau_\text{NN}\in\M_\text{NN}} \mathcal{F}_{u_\text{NN}}^{\text{opt}}(\tau_\text{NN}(u_\text{NN}))$\end{tabular} \\ \bottomrule
\end{tabular}%
}
\caption{The saddle-point approach,  the Generalized Ritz Method, and the Double Ritz Method in the context of NNs: \emph{\acf{WANs}}, \emph{the \acf{GDRM}}, and \emph{the \acf{D2RM}}. Here, the ``arg'' and ``min/max'' terms should be understood loosely according to the possible lack of existence of minimizers/maximizers discussed in \Cref{chapter2}.}
\label{table_methods_NNs}
\end{table}

\subsection{Training}\label{section:Quadrature rules}

We approximate the corresponding integrals via quadrature rules, thus producing the loss functions (recall \Cref{chapter2}). We replace the calligraphic $\mathcal{F}$ with a calligraphic $\mathcal{L}$, maintaining the remaining symbology, to denote the loss functions, namely,\begin{subequations}\label{equation:integral_functional vs losses}\begin{alignat}{4}
\mathcal{F}_\text{min}^\text{max}(u_\text{NN}, v_\text{NN})&\approx \mathcal{L}_\text{min}^\text{max}(\theta_u,\theta_v; \{x_j\}_{j=1}^N;), \label{equation:integral_functional vs losses a}\\
\mathcal{F}_T(u_\text{NN})&\approx \mathcal{L}_T(\theta_u; \{x_j\}_{j=1}^N),\\
\mathcal{F}^\text{opt}_{u_\text{NN}}(\tau_\text{NN}(u_\text{NN}))&\approx \mathcal{L}^\text{opt}_{u_\text{NN}}(\theta_\tau; \{x_j\}_{j=1}^N), \label{equation:integral_functional vs losses d}
\end{alignat} 
\end{subequations} where $\{x_j\}_{j=1}^N$ denotes the set of integration points employed during integral approximation (recall \Cref{section2.2}).

To maintain the advantages of stochastic integration, reduce the runtime and sample size---from thousands to hundreds of points, and control the integration error around singularities, we consider a \emph{randomized (composite) intermediate-point quadrature rule} that generates random integration points following a beta $\beta(a,b)$ probability distribution \cite{gupta2004handbook, weinzierl2000introduction}. We conveniently tune the hyperparameters $a$ and $b$ according to our problem specifications. \Cref{alg:intermediate-point quadrature rule} describes the randomized intermediate-point quadrature rule in $\Omega=(0,1)$.

\begin{algorithm}
\caption{Randomized intermediate-point quadrature in $\Omega=(0,1)$}\label[algorithm]{alg:intermediate-point quadrature rule}
Generate $x_i\in (0,1)$ for $1\leq i\leq N$\;
Sort $\{x_i:1\leq i\leq N\}$ so that $x_{i-1}<x_i$ for all $1\leq i\leq N$\;
Evaluate $I$ in $\{x_i:1\leq i\leq N\}$\;
Define $m_0:=0$,  $m_i:=(x_{i-1}+x_i)/2$ for $1\leq i\leq N-1$, and $m_{N}:=1$\;
Define $\omega_i:=m_i-m_{i-1}$ for $1\leq i\leq N$\;
\Return $\sum_{i=1}^N \omega_i\cdot I(x_i)$, which approximates $\int_0^1 I(x) dx$\;
\end{algorithm}

Below, we detail the different optimization strategies when using NNs for the minimization schemes described in \Cref{section:Mathematical framework}:
\begin{itemize}
\item \textbf{Weak Adversarial Networks (WANs).} We use two stochastic gradient-based optimizers \cite{ruder2016overview} to adjust the learnable parameters during the min-max optimization. For simplicity, we illustrate the functioning with the classical \acf{SGD} optimizer, 
\begin{subequations}
\begin{align}
    \theta_u &= \theta_u - \lambda_u \frac{\partial\mathcal{L}_\text{min}^\text{max}}{\partial\theta_u}(\theta_u, \theta_v; \{x_j\}_{j=1}^N), \label{gradient-descent}
\end{align} and its maximization analogue, the \acf{SGA},
\begin{align}    
    \theta_v &= \theta_v + \lambda_v \frac{\partial\mathcal{L}_\text{min}^\text{max}}{\partial\theta_v}(\theta_u, \theta_v; \{x_j\}_{j=1}^N). \label{gradient-ascent}
\end{align}
\end{subequations} Above, $\lambda_u,\lambda_v>0$ denote the learning rates, and we dropped the iteration (sub)indexes for simplicity---recall \eqref{Chapter1_SGD}.  We perform multiple iterations on the ascent for each iteration on the descent (see \Cref{alg:nested_minmax}). We call \emph{outer} and \emph{inner loops} to the minimization and maximization processes,  respectively,  because of the nested optimization structure.

\begin{algorithm}[htb]
\caption{Training of \acf{WANs}}\label[algorithm]{alg:nested_minmax}
Initialize $\theta_u\in\Theta_u$ and $\theta_v\in\Theta_v$\;
\tcc{Outer-loop}
\While{\upshape not converged}{
Randomly sample $\{x_j\}_{j=1}^N\subset\Omega$\;
$\theta_u = \theta_u - \lambda_u\displaystyle\frac{\partial\mathcal{L}_\text{min}^\text{max}}{\partial\theta_u}(\theta_u, \theta_v; \{x_j\}_{j=1}^N)$\;
 \tcc{Inner-loop}
\While{\upshape not converged}{
Randomly sample $\{x_j\}_{j=1}^N\subset\Omega$\;
$\theta_v = \theta_v + \lambda_v \displaystyle\frac{\partial\mathcal{L}_\text{min}^\text{max}}{\partial\theta_v}(\theta_u, \theta_v; \{x_j\}_{j=1}^N)$\;
  }
}
\Return $\theta_u$
\end{algorithm}

\item \textbf{Deep Double Ritz Method (D$\boldmath{^2}$RM).} The training is similar to \Cref{alg:nested_minmax}, but modifying the loss when jumping from the outer to the inner loop, and performing only gradient-descents for minimizations (see \Cref{alg:nested_minmin}). 

\begin{algorithm}[htb]
\caption{Training of the \acf{D2RM}}\label[algorithm]{alg:nested_minmin}
Initialize $\theta_u\in\Theta_u$ and $\theta_\tau\in\Theta_\tau$\;
\tcc{Outer-loop}
\While{\upshape not converged}{
Randomly sample $\{x_j\}_{j=1}^N\subset\Omega$\;
$\theta_u = \theta_u - \lambda_u \displaystyle\frac{\partial\mathcal{L}_{\tau_\text{NN}}}{\partial\theta_u}(\theta_u,\theta_\tau; \{x_j\}_{j=1}^N)$\;
\tcc{Inner-loop}
\While{\upshape not converged}{
Randomly sample $\{x_j\}_{j=1}^N\subset\Omega$\;
$\theta_\tau = \theta_\tau - \lambda_\tau \displaystyle\frac{\partial\mathcal{L}_{u_\text{NN}}^\text{opt}}{\partial\theta_\tau}(\theta_u,\theta_\tau; \{x_j\}_{j=1}^N)$\;
  }
}
\Return $\theta_u$
\end{algorithm}

\item \textbf{Generalized Deep Ritz Method (GDRM).} The optimization consists of a single-loop minimization (see \Cref{alg:single_min}) when the trial-to-test operator $T$ is ``computably available''.  In particular, when the bilinear form is symmetric and positive definite (recall item (a) in \Cref{Generalized Ritz method}), this method is the well-known \acf{DRM} \cite{e2018deep}.

\begin{algorithm}[htb]
\caption{Training of the \acf{GDRM}}\label[algorithm]{alg:single_min}
Initialize $\theta_u\in\Theta_u$\;
\While{\upshape not converged}{
Randomly sample $\{x_j\}_{j=1}^N\subset\Omega$\;
$\theta_u = \theta_u - \lambda_u \displaystyle\frac{\partial\mathcal{L}_T}{\partial\theta_u}(\theta_u; \{x_j\}_{j=1}^N)$\;
}
\Return $\theta_u$
\end{algorithm}

\item \textbf{Adjoint Deep Ritz Method (DRM$)'$.} The optimization consists of a single-loop minimization (as in \Cref{alg:single_min}) with post-processing and is only valid in ultraweak formulations (recall item (c) in \Cref{Generalized Ritz method}).

\end{itemize}

We consider the Adam optimizer \cite{kingma2014adam} to carry out our experiments. In the case of nested optimizations, we select fully independent Adam optimizers for the inner and outer loops. We establish an accumulated maximum number of iterations for both nested loops, and we fix four inner-loop iterations for each outer-loop iteration\footnote{We made this decision motivated by the continuity of the elements involved during the process: a slight modification in $\theta_u$ translates into a slight variation in $u_\text{NN}$ (continuity of the realization mapping---recall \Cref{section2.1}) that produces a slight variation in $Tu_\text{NN}$ (continuity of $T$). Therefore, we expect (have the hope) that $\tau_\text{NN}(u_\text{NN})$ can provide a good approximation of $Tu_{\text{NN}}$ after a small number of iterations in $\theta_\tau$.} (as considered in \cite{zang2020weak,  bao2020numerical} for \ac{WANs}\footnote{The lack of Lipschitz continuity in the min-max approach suggests that this is a poor strategy for optimizing \ac{WANs}.}) unless otherwise specified.
 
\newpage
\section{Implementation}
\label{section:Implementation}

We use the \acf{TF2} library \cite{tensorflow2015-whitepaper, abadi2016tensorflow} within Python to implement our neural network architectures and loss functions, and to manage the random creation and flow of data. Specifically, we accommodate all of our implementations to Keras (\texttt{tf.keras}).

\subsection{Samples generation, input batch flow, encapsulation of models, and optimization}

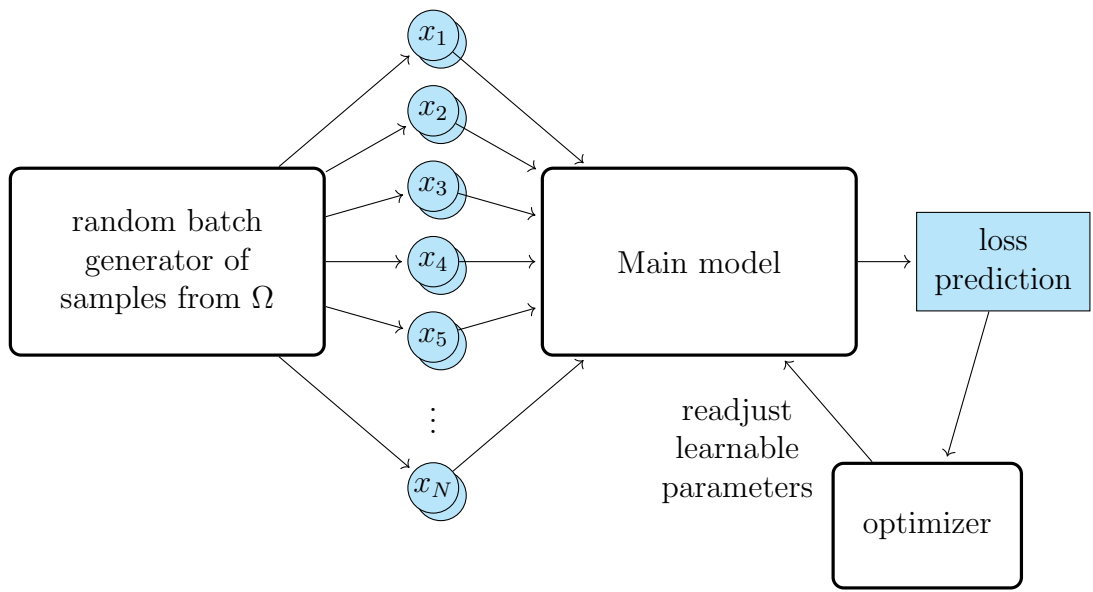
\begin{figure}[htb]
\centering
\begin{tikzpicture}[
    shorten >=2pt,->,
    draw=black,
    node distance=\layersep,
    every pin edge/.style={<-,shorten <=1pt},
    neuron/.style={circle, draw=black, fill=black,minimum size=19pt,inner sep=0pt},
    element/.style={rectangle, draw=black, fill=brown!25!, minimum size=30pt,inner sep=6pt},
    input neuron/.style={neuron, fill=cyan!25},
    input network/.style={neuron, fill=green!25!},
    white neuron/.style={neuron, fill=white, minimum size=5pt},
    output neuron/.style={neuron, fill=red!35, inner sep=2pt},
    hidden neuron/.style={neuron, fill=purple!50},
    annot/.style={text width=4em, text centered},
    punkt/.style={
           rectangle,
           rounded corners,
           draw=black, very thick,
           text width=4em,
           minimum height=2em,
           text centered},
    pil/.style={
           ->,
           thick,
           shorten <=2pt,
           shorten >=2pt,}
]
    \def\layersep{1.5cm}
    
   \node[punkt, minimum width = 10em, minimum height = 6em, text width=8em] (batch) at (-4,-3) {random batch generator of samples from $\Omega$};
   
    \node[input neuron] (I-1aux) at (-0.4,-0.1) {};
    \node[input neuron] (I-1) at (-0.5,0) {$x_1$};
    \node[input neuron] (I-2aux) at (-0.4,-1.1) {};
    \node[input neuron] (I-2) at (-0.5,-1) {$x_2$};
    \node[input neuron] (I-3aux) at (-0.4,-2.1) {};
    \node[input neuron] (I-3) at (-0.5,-2) {$x_3$};
    \node[input neuron] (I-4aux) at (-0.4,-3.1) {};
    \node[input neuron] (I-4) at (-0.5,-3) {$x_4$};
    \node[input neuron] (I-5aux) at (-0.4,-4.1) {};
    \node[input neuron] (I-5) at (-0.5,-4) {$x_5$};
    \node[annot, text width=0.5em] (I-6) at (-0.5,-5) {$\vdots$};
    \node[input neuron] (I-7aux) at (-0.4,-6.1) {};
    \node[input neuron] (I-7) at (-0.5,-6) {$x_N$};
    
    \path (batch) edge (I-1);
    \path (batch) edge (I-2);
    \path (batch) edge (I-3);
    \path (batch) edge (I-4);
    \path (batch) edge (I-5);
    \path (batch) edge (I-7);
    
   \node[punkt, minimum width = 10em, minimum height = 6em, text width=6em] (main) at (3,-3) {Main model};
   
    \path (I-1) edge (main);
    \path (I-2) edge (main);
    \path (I-3) edge (main);
    \path (I-4) edge (main);
    \path (I-5) edge (main);
    \path (I-7) edge (main);

    \node[element, fill=cyan!25, text width=4.5em, text centered] (loss_out) at (7, -3) {loss\\ prediction};
    \path (main) edge (loss_out);
    
   \node[punkt, minimum width = 6em, minimum height = 4em, text width=5em] (optimizer) at (6,-6.5) {optimizer};

    \path (loss_out) edge (optimizer);
   	\path (optimizer) edge (main) ;
   	\node[annot, text width=6em] at (3.5,-5.5) {readjust learnable parameters};

\end{tikzpicture}
\caption{Implementation sketch of the general flux for the proposed methods at each training iteration.}
\label{figure:implementation_sketch_optimization}
\end{figure}

Our inputs to networks are samples over the domain $\Omega$.  After feeding our networks with the batch of inputs, we combine the batch of outputs in a single loss prediction. We encapsulate the networks and losses inside a main model whose input and output are the batch of samples and the loss prediction, respectively. From the loss prediction, we optimize the learnable parameters of the networks (i.e., we fit the learnable parameters to the data) with the Keras built-in Adam optimizer.  \Cref{figure:implementation_sketch_optimization} illustrates the described process at each training iteration.  

For the maximization involved in \ac{WANs}, we reverse the sign of the gradients so that when feeding the (default) gradient-descent-based optimizer, it instead performs a gradient ascent. 

\subsection{Design of models by methods}

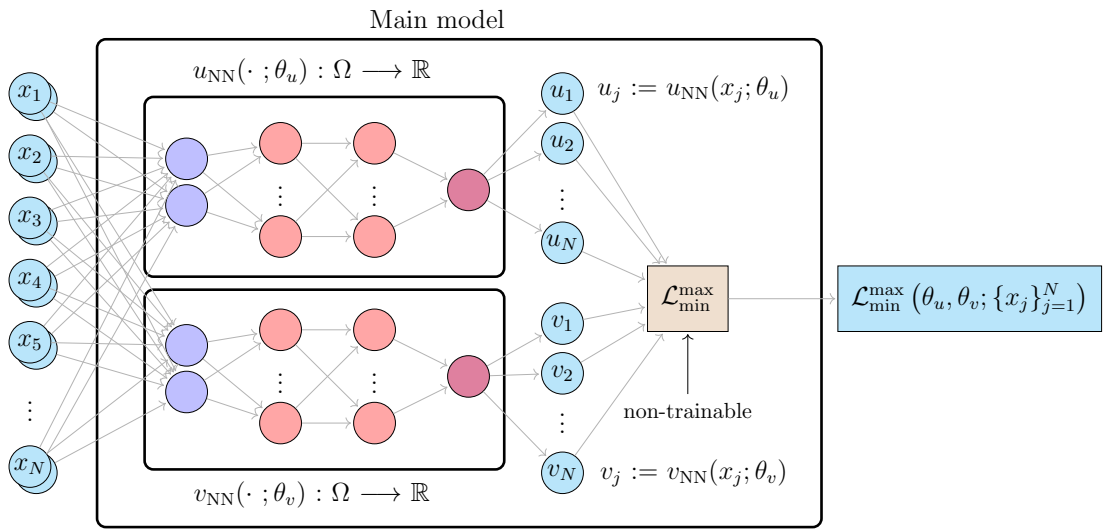
\begin{figure}[htbp]
\centering
\resizebox{\textwidth}{!}{%
\begin{tikzpicture}[    
    shorten >=1pt,->,
    draw=gray!60,
    every pin edge/.style={<-,shorten <=1pt},
    neuron/.style={circle, draw=black, fill=black,minimum size=19pt,inner sep=0pt},
    element/.style={rectangle, draw=black, fill=brown!25!, minimum size=30pt,inner sep=6pt},
    input neuron/.style={neuron, fill=cyan!25},
    input network/.style={neuron, fill=blue!25!},
    white neuron/.style={neuron, fill=white, minimum size=5pt},
    output neuron/.style={neuron, fill=purple!50, inner sep=2pt},
    hidden neuron/.style={neuron, fill=red!35},
    annot/.style={text width=4em, text centered},
    punkt/.style={
           rectangle,
           rounded corners,
           draw=black, very thick,
           text width=4em,
           minimum height=2em,
           text centered},
    pil/.style={
           ->,
           thick,
           shorten <=2pt,
           shorten >=2pt,}
]
    \node[punkt, minimum width = 28em, minimum height = 19em, right] (block) at (0.55,-2.75) {};
    \node[] at (6,1.5) {Main model};

    \node[input neuron] (I-1aux) at (-0.4,0.2) {};
    \node[input neuron] (I-1) at (-0.5,0.3) {$x_1$};
    \node[input neuron] (I-2aux) at (-0.4,-0.8) {};
    \node[input neuron] (I-2) at (-0.5,-0.7) {$x_2$};
    \node[input neuron] (I-3aux) at (-0.4,-1.8) {};
    \node[input neuron] (I-3) at (-0.5,-1.7) {$x_3$};
    \node[input neuron] (I-4aux) at (-0.4,-2.8) {};
    \node[input neuron] (I-4) at (-0.5,-2.7) {$x_4$};
    \node[input neuron] (I-5aux) at (-0.4,-3.8) {};
    \node[input neuron] (I-5) at (-0.5,-3.7) {$x_5$};
    \node[annot, text width=0.5em] (I-6) at (-0.5,-4.7) {$\vdots$};
    \node[input neuron] (I-7aux) at (-0.4,-5.8) {};
    \node[input neuron] (I-7) at (-0.5,-5.7) {$x_N$};
    
    \node[punkt, minimum width = 13.9em, minimum height = 7em] (block) at (4.2,-1.2) {};
    \node[input network] (H1-1) at (2,-0.75) {};
    \node[input network] (H1-2) at (2,-1.5) {};
    \node[hidden neuron] (H2-1) at (3.5,-0.5) {};
    \node[] (H2-2) at (3.5,-1.25) {$\vdots$};
    \node[hidden neuron] (H2-3) at (3.5,-2) {};
    \node[hidden neuron] (H3-1) at (5,-0.5) {};
    \node[] (H3-2) at (5,-1.25) {$\vdots$};
    \node[hidden neuron] (H3-3) at (5,-2) {};
    \node[output neuron] (O-1) at (6.5,-1.25) {};
    
    \node[annot, text width=10em] (label) at (4, 0.65) {$u_\text{NN}(\cdot\; ;\theta_u): \Omega\longrightarrow\mathbb{R}$};
    \path (I-1) edge (H1-1);
    \path (I-2) edge (H1-1);
    \path (I-3) edge (H1-1);
    \path (I-4) edge (H1-1);
    \path (I-5) edge (H1-1);
    \path (I-7) edge (H1-1);
    \path (I-1aux) edge (H1-2);
    \path (I-2aux) edge (H1-2);
    \path (I-3aux) edge (H1-2);
    \path (I-4aux) edge (H1-2);
    \path (I-5aux) edge (H1-2);
    \path (I-7aux) edge (H1-2);
    \path (H1-1) edge (H2-1);
    \path (H1-1) edge (H2-3);
    \path (H1-2) edge (H2-1);
    \path (H1-2) edge (H2-3);
    \path (H2-1) edge (H3-1);
    \path (H2-1) edge (H3-3);
    \path (H2-3) edge (H3-1);
    \path (H2-3) edge (H3-3);
    \path (H3-1) edge (O-1);
    \path (H3-3) edge (O-1);
    
    \node[input neuron] (Out-1) at (8,0.3) {$u_1$};
    \node[input neuron] (Out-2) at (8,-0.5) {$u_2$};
    \node[] (Out-3) at (8,-1.3) {$\vdots$};
    \node[input neuron] (Out-4) at (8,-2.1) {$u_N$};
    
    \path (O-1) edge (Out-1);
    \path (O-1) edge (Out-2);
    \path (O-1) edge (Out-4);
   
    \node[annot,text width=8em,right] (label) at (8.3,0.37) {$u_j := u_\text{NN}(x_j; \theta_u)$};
    \node[element] (loss) at (10,-3) {$\mathcal{L}^\text{max}_\text{min}$};
    
    \node[punkt, minimum width = 13.9em, minimum height = 7em] (block) at (4.2,-4.3) {};
    \node[input network] (H12-1) at (2,-3.75) {};
    \node[input network] (H12-2) at (2,-4.5) {};
    \node[hidden neuron] (H22-1) at (3.5,-3.5) {};
    \node[] (H22-2) at (3.5,-4.25) {$\vdots$};
    \node[hidden neuron] (H22-3) at (3.5,-5) {};
    \node[hidden neuron] (H32-1) at (5,-3.5) {};
    \node[] (H32-2) at (5,-4.25) {$\vdots$};
    \node[hidden neuron] (H32-3) at (5,-5) {};
    \node[output neuron] (O2-1) at (6.5,-4.25) {};
    
    \path (I-1) edge (H12-1);
    \path (I-2) edge (H12-1);
    \path (I-3) edge (H12-1);
    \path (I-4) edge (H12-1);
    \path (I-5) edge (H12-1);
    \path (I-7) edge (H12-1);
    \path (I-1aux) edge (H12-2);
    \path (I-2aux) edge (H12-2);
    \path (I-3aux) edge (H12-2);
    \path (I-4aux) edge (H12-2);
    \path (I-5aux) edge (H12-2);
    \path (I-7aux) edge (H12-2);
    \path (H12-1) edge (H22-1);
    \path (H12-1) edge (H22-3);
    \path (H12-2) edge (H22-1);
    \path (H12-2) edge (H22-3);
    \path (H22-1) edge (H32-1);
    \path (H22-1) edge (H32-3);
    \path (H22-3) edge (H32-1);
    \path (H22-3) edge (H32-3);
    \path (H32-1) edge (O2-1);
    \path (H32-3) edge (O2-1);
    
    \node[input neuron] (Out2-1) at (8,-3.4) {$v_1$};
    \node[input neuron] (Out2-2) at (8,-4.2) {$v_2$};
    \node[] (Out2-6) at (8,-4.9) {$\vdots$};
    \node[input neuron] (Out2-4) at (8,-5.8) {$v_N$};
    \node[annot,text width=8em,right] (label) at (8.3,-5.8) {$v_j := v_\text{NN}(x_j; \theta_v)$};
    \node[annot, text width=10em] (label) at (4, -6.15) {$v_\text{NN}(\cdot\; ;\theta_v): \Omega\longrightarrow\mathbb{R}$};
    
    \path (O2-1) edge (Out2-1);
    \path (O2-1) edge (Out2-2);
    \path (O2-1) edge (Out2-4);
    
    \node[annot, text width=12em] (label_loss) at (10,-4.8) {\footnotesize non-trainable};
    \path (label_loss) edge[black] (loss);
    
    \path (Out-1) edge (loss);
    \path (Out-2) edge (loss);
    \path (Out-4) edge (loss);
    \path (Out2-1) edge (loss);
    \path (Out2-2) edge (loss);
    \path (Out2-4) edge (loss);
    
    \node[element, fill=cyan!25] (loss_out) at (14.5, -3) {$\mathcal{L}^\text{max}_\text{min}\left(\theta_u, \theta_v; \{x_j\}_{j=1}^N \right)$};
    \path (loss) edge (loss_out);
\end{tikzpicture}
}%
\caption{Main model architecture for \acs{WANs}. It consists of two independent NNs, $u_\text{NN}$ and $v_\text{NN}$,  combined via the loss function $\mathcal{L}_\text{min}^\text{max}$.}
\label{figure:implementation_sketch_MIN_MAX}
\end{figure}

We implement networks, losses, and operators as models and layers in \ac{TF2} from the redefinitions of the corresponding Keras base classes (\texttt{tf.keras.Model} and \texttt{tf.keras.Layer}). In \ac{WANs}, we implement $u_\text{NN}$ and $v_\text{NN}$ as two independent models that are combined via a common non-trainable layer for the loss (see \Cref{figure:implementation_sketch_MIN_MAX}). In the \ac{DRM}, we implement $u_\text{NN}$ as a model that subsequently connects with two non-trainable layers for the trial-to-test operator and the loss function, respectively (see \Cref{figure:implementation_sketch_MIN}).  In the \ac{D2RM}, we implement $u_\text{NN}$ and $\tau_\text{NN}$ as two sequential models whose output feeds into two separate losses (see \Cref{figure:implementation_sketch_MIN_MIN}).  The loss functions are implemented as latent outputs of the main model, which, together with a TF2-suitable boolean variable, activate and deactivate alternatively. Despite the significant compilation time that the two-branch model (for the \ac{D2RM}) takes compared with one-branch models (for \ac{WANs} and the \ac{DRM})\footnote{Around one minute for the \ac{D2RM} vs. a couple of seconds for \ac{WANs} and the \ac{DRM}.},  its fitting execution in graph mode is as fast as that of one-branch models.

\begin{figure}[htb]
\centering
\resizebox{\textwidth}{!}{
\begin{tikzpicture}[
	shorten >=1pt,->,
    draw=gray!60,
    every pin edge/.style={<-,shorten <=1pt},
    neuron/.style={circle, draw=black, fill=black,minimum size=19pt,inner sep=0pt},
    element/.style={rectangle, draw=black, fill=brown!25!, minimum size=30pt,inner sep=6pt},
    input neuron/.style={neuron, fill=cyan!25},
    input network/.style={neuron, fill=blue!25},
    white neuron/.style={neuron, fill=white, minimum size=5pt},
    output neuron/.style={neuron, fill=purple!50, inner sep=2pt},
    hidden neuron/.style={neuron, fill=red!35},
    annot/.style={text width=4em, text centered},
    punkt/.style={
           rectangle,
           rounded corners,
           draw=black, very thick,
           text width=4em,
           minimum height=2em,
           text centered},
    pil/.style={
           ->,
           thick,
           shorten <=2pt,
           shorten >=2pt,}
]
    
    \node[punkt, minimum width = 38em, minimum height = 19em, right] (block) at (0.55,-3) {};
    \node[] at (8,1.5) {Main model};

    \node[input neuron] (I-1aux) at (-0.4,-0.1) {};
    \node[input neuron] (I-1) at (-0.5,0) {$x_1$};
    \node[input neuron] (I-2aux) at (-0.4,-1.1) {};
    \node[input neuron] (I-2) at (-0.5,-1) {$x_2$};
    \node[input neuron] (I-3aux) at (-0.4,-2.1) {};
    \node[input neuron] (I-3) at (-0.5,-2) {$x_3$};
    \node[input neuron] (I-4aux) at (-0.4,-3.1) {};
    \node[input neuron] (I-4) at (-0.5,-3) {$x_4$};
    \node[input neuron] (I-5aux) at (-0.4,-4.1) {};
    \node[input neuron] (I-5) at (-0.5,-4) {$x_5$};
    \node[annot, text width=0.5em] (I-6) at (-0.5,-5) {$\vdots$};
    \node[input neuron] (I-7aux) at (-0.4,-6.1) {};
    \node[input neuron] (I-7) at (-0.5,-6) {$x_N$};
    \node[punkt, minimum width = 15em, minimum height = 10em] (block) at (4,-3) {};
    \node[input network] (H1-1) at (2,-2) {};
    \node[input network] (H1-2) at (2,-4) {};
    \node[hidden neuron] (H2-1) at (3.5,-1.5) {};
    \node[hidden neuron] (H2-2) at (3.5,-2.5) {};
    \node[] (H2-3) at (3.5,-3.5) {\vspace{-12em} $\vdots$};
    \node[hidden neuron] (H2-4) at (3.5,-4.5) {};
    \node[hidden neuron] (H3-1) at (5,-1.5) {};
    \node[hidden neuron] (H3-2) at (5,-2.5) {};
    \node[] (H3-3) at (5,-3.5) {\vspace{-12em} $\vdots$};
    \node[hidden neuron] (H3-4) at (5,-4.5) {};
    \node[output neuron] (O-1) at (6.5,-3) {};
    
        \node[annot, text width=10em] (label) at (4.25,-0.6) {$u_\text{NN}(\cdot\; ; \theta_u):\Omega\longrightarrow\mathbb{R}$};

    \path (I-1) edge (H1-1);
    \path (I-2) edge (H1-1);
    \path (I-3) edge (H1-1);
    \path (I-4) edge (H1-1);
    \path (I-5) edge (H1-1);
    \path (I-7) edge (H1-1);
    \path (I-1aux) edge (H1-2);
    \path (I-2aux) edge (H1-2);
    \path (I-3aux) edge (H1-2);
    \path (I-4aux) edge (H1-2);
    \path (I-5aux) edge (H1-2);
    \path (I-7aux) edge (H1-2);
    \path (H1-1) edge (H2-1);
    \path (H1-1) edge (H2-2);
    \path (H1-1) edge (H2-4);    
    \path (H1-2) edge (H2-1);
    \path (H1-2) edge (H2-2);
    \path (H1-2) edge (H2-4);     
    \path (H2-1) edge (H3-1);
    \path (H2-1) edge (H3-2);
    \path (H2-1) edge (H3-4);
    \path (H2-2) edge (H3-1);
    \path (H2-2) edge (H3-2);
    \path (H2-2) edge (H3-4);
    \path (H2-4) edge (H3-1);
    \path (H2-4) edge (H3-2);
    \path (H2-4) edge (H3-4);
    \path (H3-1) edge (O-1);
    \path (H3-2) edge (O-1);
    \path (H3-4) edge (O-1);
    
    \node[input neuron] (Out-1) at (8,0) {$u_1$};
    \node[input neuron] (Out-2) at (8,-1) {$u_2$};
    \node[input neuron] (Out-3) at (8,-2) {$u_3$};
    \node[input neuron] (Out-4) at (8,-3) {$u_4$};
    \node[input neuron] (Out-5) at (8,-4) {$u_5$};
    \node[annot, text width=0.5em] (Out-6) at (8,-5) {$\vdots$};
    \node[input neuron] (Out-7) at (8,-6) {$u_N$};
    
    \path (O-1) edge (Out-1);
    \path (O-1) edge (Out-2);
    \path (O-1) edge (Out-3);
    \path (O-1) edge (Out-4);
    \path (O-1) edge (Out-5);
    \path (O-1) edge (Out-7);
    
    \node[annot,text width=8em,right] (label) at (8.5,0) {$u_j := u_\text{NN}(x_j; \theta_u)$};
    \node[element] (T_h) at (10,-3) {$T$};
    
    \node[annot, text width=12em] (label_T_h) at (10,-4.8) {\footnotesize non-trainable};
    \path (label_T_h) edge[black] (T_h);
    
    \path (Out-1) edge (T_h);
    \path (Out-2) edge (T_h);
    \path (Out-3) edge (T_h);
    \path (Out-4) edge (T_h);
    \path (Out-5) edge (T_h);
    \path (Out-7) edge (T_h);
    
    \node[input neuron] (T-1) at (12.5,0) {$T_1$};
    \node[input neuron] (T-2) at (12.5,-1) {$T_2$};
    \node[input neuron] (T-3) at (12.5,-2) {$T_3$};
    \node[input neuron] (T-4) at (12.5,-3) {$T_4$};
    \node[input neuron] (T-5) at (12.5,-4) {$T_5$};
    \node[annot, text width=0.5em] (T-6) at (12.5,-5) {$\vdots$};
    \node[input neuron] (T-7) at (12.5,-6) {$T_N$};
    
    \path (T_h) edge (T-1);
    \path (T_h) edge (T-2);
    \path (T_h) edge (T-3);
    \path (T_h) edge (T-4);
    \path (T_h) edge (T-5);
    \path (T_h) edge (T-7);
    
   \node[annot,text width=10em, right] (label) at (12.25,0) {$T_j := T(u_j)$};
   
    \node[element] (loss) at (15,-3) {$\mathcal{L}$};

    \path (T-1) edge (loss);
    \path (T-2) edge (loss);
    \path (T-3) edge (loss);
    \path (T-4) edge (loss);
    \path (T-5) edge (loss);
    \path (T-7) edge (loss);
    
    \node[element, fill=cyan!25] (loss_out) at (18.5, -3) {$\mathcal{L}_T\left(\theta_u; \{x_j\}_{j=1}^N \right)$};
    \path (loss) edge (loss_out);
    
    \node[annot, text width=6em] (label_losses) at (15,-4.8) {\footnotesize non-trainable};
    \path (label_losses) edge[black] (loss);

\end{tikzpicture}
}%
\caption{Main model architecture for the \acs{DRM}. It consists of a NN, $u_\text{NN}$, composed with the trial-to-test operator $T$ and the loss function $\mathcal{L}$.}
\label{figure:implementation_sketch_MIN}
\end{figure}
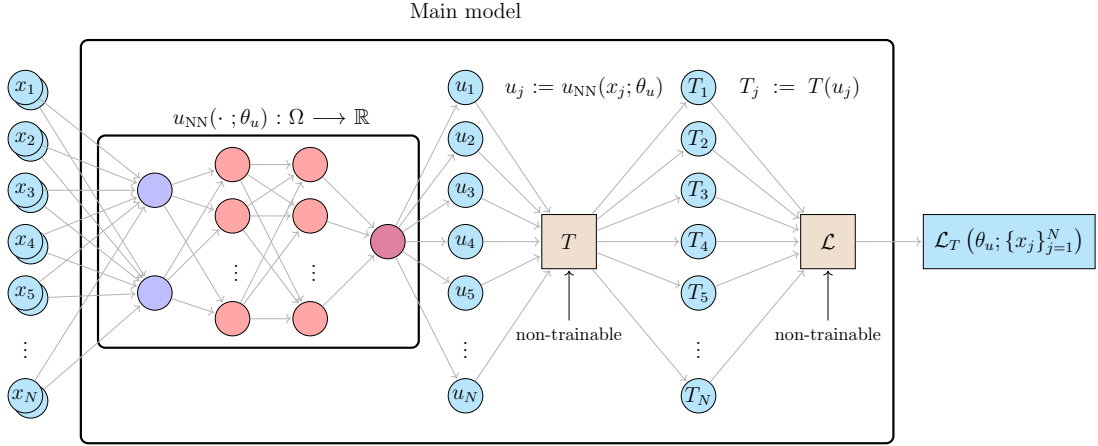
\begin{figure}[htb]
\centering
\resizebox{\textwidth}{!}{%
\begin{tikzpicture}[
	shorten >=1pt,->,
    draw=gray!60,
    node distance=\layersep,
    every pin edge/.style={<-,shorten <=1pt},
    neuron/.style={circle, draw=black, fill=black,minimum size=19pt,inner sep=0pt},
    element/.style={rectangle, draw=black, fill=brown!25!, minimum size=30pt,inner sep=6pt},
    input neuron/.style={neuron, fill=cyan!25},
    input network/.style={neuron, fill=blue!25!},
    white neuron/.style={neuron, fill=white, minimum size=5pt},
    output neuron/.style={neuron, fill=purple!50, inner sep=2pt},
    hidden neuron/.style={neuron, fill=red!35},
    annot/.style={text width=4em, text centered},
    punkt/.style={
           rectangle,
           rounded corners,
           draw=black, very thick,
           text width=4em,
           minimum height=2em,
           text centered},
    pil/.style={
           ->,
           thick,
           shorten <=2pt,
           shorten >=2pt,}
]
    \def\layersep{1.5cm}
    
    \node[punkt, minimum width = 50.25em, minimum height = 18.5em, right] (block) at (0.25,-3) {};
    \node[annot, text width=8em] at (10,1.5) {Main model};

    \node[input neuron] (I-1aux) at (-0.4,-0.1) {};
    \node[input neuron] (I-1) at (-0.5,0) {$x_1$};
    \node[input neuron] (I-2aux) at (-0.4,-1.1) {};
    \node[input neuron] (I-2) at (-0.5,-1) {$x_2$};
    \node[input neuron] (I-3aux) at (-0.4,-2.1) {};
    \node[input neuron] (I-3) at (-0.5,-2) {$x_3$};
    \node[input neuron] (I-4aux) at (-0.4,-3.1) {};
    \node[input neuron] (I-4) at (-0.5,-3) {$x_4$};
    \node[input neuron] (I-5aux) at (-0.4,-4.1) {};
    \node[input neuron] (I-5) at (-0.5,-4) {$x_5$};
    \node[annot, text width=0.5em] (I-6) at (-0.5,-5) {$\vdots$};
    \node[input neuron] (I-7aux) at (-0.4,-6.1) {};
    \node[input neuron] (I-7) at (-0.5,-6) {$x_N$};
    \node[punkt, minimum width = 16em, minimum height = 10em] (block) at (4,-3) {};
    \node[input network] (H1-1) at (2,-2) {};
    \node[input network] (H1-2) at (2,-4) {};
    \node[hidden neuron] (H2-1) at (3.5,-1.5) {};
    \node[hidden neuron] (H2-2) at (3.5,-2.5) {};
    \node[] (H2-3) at (3.5,-3.5) {\vspace{-12em} $\vdots$};
    \node[hidden neuron] (H2-4) at (3.5,-4.5) {};
    \node[hidden neuron] (H3-1) at (5,-1.5) {};
    \node[hidden neuron] (H3-2) at (5,-2.5) {};
    \node[] (H3-3) at (5,-3.5) {\vspace{-12em} $\vdots$};
    \node[hidden neuron] (H3-4) at (5,-4.5) {};
    \node[output neuron] (O-1) at (6.5,-3) {};

    \node[annot, text width=10em] (label) at (4.25,-0.6) {$u_\text{NN}(\cdot\; ; \theta_u):\Omega\longrightarrow\mathbb{R}$};
    
    \path (I-1) edge (H1-1);
    \path (I-2) edge (H1-1);
    \path (I-3) edge (H1-1);
    \path (I-4) edge (H1-1);
    \path (I-5) edge (H1-1);
    \path (I-7) edge (H1-1);
    \path (I-1aux) edge (H1-2);
    \path (I-2aux) edge (H1-2);
    \path (I-3aux) edge (H1-2);
    \path (I-4aux) edge (H1-2);
    \path (I-5aux) edge (H1-2);
    \path (I-7aux) edge (H1-2);
    \path (H1-1) edge (H2-1);
    \path (H1-1) edge (H2-2);
    \path (H1-1) edge (H2-4);    
    \path (H1-2) edge (H2-1);
    \path (H1-2) edge (H2-2);
    \path (H1-2) edge (H2-4);     
    \path (H2-1) edge (H3-1);
    \path (H2-1) edge (H3-2);
    \path (H2-1) edge (H3-4);
    \path (H2-2) edge (H3-1);
    \path (H2-2) edge (H3-2);
    \path (H2-2) edge (H3-4);
    \path (H2-4) edge (H3-1);
    \path (H2-4) edge (H3-2);
    \path (H2-4) edge (H3-4);
    \path (H3-1) edge (O-1);
    \path (H3-2) edge (O-1);
    \path (H3-4) edge (O-1);
    
    \node[input neuron] (Out-1) at (8,0) {$u_1$};
    \node[input neuron] (Out-2) at (8,-1) {$u_2$};
    \node[input neuron] (Out-3) at (8,-2) {$u_3$};
    \node[input neuron] (Out-4) at (8,-3) {$u_4$};
    \node[input neuron] (Out-5) at (8,-4) {$u_5$};
    \node[annot, text width=0.5em] (Out-6) at (8,-5) {$\vdots$};
    \node[input neuron] (Out-7) at (8,-6) {$u_N$};
    
    \path (O-1) edge (Out-1);
    \path (O-1) edge (Out-2);
    \path (O-1) edge (Out-3);
    \path (O-1) edge (Out-4);
    \path (O-1) edge (Out-5);
    \path (O-1) edge (Out-7);
    
    \node[annot,text width=10em,right] (label) at (8.25,-6) {$u_j := u_\text{NN}(x_j; \theta_u)$};
    \node[punkt, minimum width = 16em, minimum height = 10em] (block) at (12,-3) {};
    \node[input network] (H4-1) at (10,-3) {};
    \node[hidden neuron] (H5-1) at (11.5,-1.5) {};
    \node[hidden neuron] (H5-2) at (11.5,-2.5) {};
    \node[] (H5-3) at (11.5,-3.5) {\vspace{-12em} $\vdots$};
    \node[hidden neuron] (H5-4) at (11.5,-4.5) {};
    \node[hidden neuron] (H6-1) at (13,-1.5) {};
    \node[hidden neuron] (H6-2) at (13,-2.5) {};
    \node[] (H6-3) at (13,-3.5) {\vspace{-12em} $\vdots$};
    \node[hidden neuron] (H6-4) at (13,-4.5) {};
    \node[output neuron] (O1-1) at (14.5,-3) {};

    \node[annot, text width=14em] (label) at (12,-0.6) {$\tau_\text{NN}(\cdot\; ; \theta_\tau):u_\text{NN}(\Omega; \theta_u)\longrightarrow\mathbb{R}$};
    
    \path (Out-1) edge (H4-1);
    \path (Out-2) edge (H4-1);
    \path (Out-3) edge (H4-1);
    \path (Out-4) edge (H4-1);
    \path (Out-5) edge (H4-1);
    \path (Out-7) edge (H4-1);
    \path (H4-1) edge (H5-1);
    \path (H4-1) edge (H5-2);
    \path (H4-1) edge (H5-4);     
    \path (H5-1) edge (H6-1);
    \path (H5-1) edge (H6-2);
    \path (H5-1) edge (H6-4);
    \path (H5-2) edge (H6-1);
    \path (H5-2) edge (H6-2);
    \path (H5-2) edge (H6-4);
    \path (H5-4) edge (H6-1);
    \path (H5-4) edge (H6-2);
    \path (H5-4) edge (H6-4);
    \path (H6-1) edge (O1-1);
    \path (H6-2) edge (O1-1);
    \path (H6-4) edge (O1-1);
    
    \node[input neuron] (T-1) at (16,0) {$\tau_1$};
    \node[input neuron] (T-2) at (16,-1) {$\tau_2$};
    \node[input neuron] (T-3) at (16,-2) {$\tau_3$};
    \node[input neuron] (T-4) at (16,-3) {$\tau_4$};
    \node[input neuron] (T-5) at (16,-4) {$\tau_5$};
    \node[annot, text width=0.5em] (T-6) at (16,-5) {$\vdots$};
    \node[input neuron] (T-7) at (16,-6) {$\tau_N$};
    
    \path (O1-1) edge (T-1);
    \path (O1-1) edge (T-2);
    \path (O1-1) edge (T-3);
    \path (O1-1) edge (T-4);
    \path (O1-1) edge (T-5);
    \path (O1-1) edge (T-7);
    
   \node[annot,text width=9em,right] (label) at (16.25,-6) {$\tau_j := \tau_\text{NN}(u_j; \theta_\tau)$};
   
    \node[element] (loss_u) at (18,-1.5) {$\mathcal{L}_{\tau_\text{NN}}$};
    \node[element] (loss_T) at (18,-4.5) {$\mathcal{L}_{u_\text{NN}}^\text{opt}$};

    \path (T-1) edge (loss_u);
    \path (T-2) edge (loss_u);
    \path (T-3) edge (loss_u);
    \path (T-4) edge (loss_u);
    \path (T-5) edge (loss_u);
    \path (T-7) edge (loss_u);
    
    \path (T-1) edge (loss_T);
    \path (T-2) edge (loss_T);
    \path (T-3) edge (loss_T);
    \path (T-4) edge (loss_T);
    \path (T-5) edge (loss_T);
    \path (T-7) edge (loss_T);
    
    \node[element, fill=cyan!25] (loss_u_out) at (23.5, -1.5) {$\mathcal{L}_{\tau_\text{NN}}\left(\theta_u; \{x_j\}_{j=1}^N \right)$};
    \node[element, fill=cyan!25] (loss_T_out) at (23.5, -4.5) {$\mathcal{L}_{u_\text{NN}}^{\text{opt}}\left(\theta_\tau; \{x_j\}_{j=1}^N \right)$};
    
    \node[annot, text width=5em] (label_losses) at (18.8,-3) {\footnotesize non-trainable};
    \path (label_losses) edge[black] (loss_u);
    \path (label_losses) edge[black] (loss_T);
    
    \path (loss_u) edge node[above=0mm] {\footnotesize if $u_\text{NN}$} node[below=0mm] {\footnotesize trainable} (loss_u_out);
    \path (loss_T) edge node[above=0mm] {\footnotesize if $\tau_\text{NN}$} node[below=0mm] {\footnotesize trainable} (loss_T_out);
    
\end{tikzpicture}
}%
\caption{Main model architecture for the \acs{D2RM}. It consists of two NNs, $u_\text{NN}$ and $\tau_\text{NN}$,  equipped with the loss functions $\mathcal{L}_{\tau_\text{NN}}$ and $\mathcal{L}_{u_\text{NN}}^{\text{opt}}$.}
\label{figure:implementation_sketch_MIN_MIN}
\end{figure}
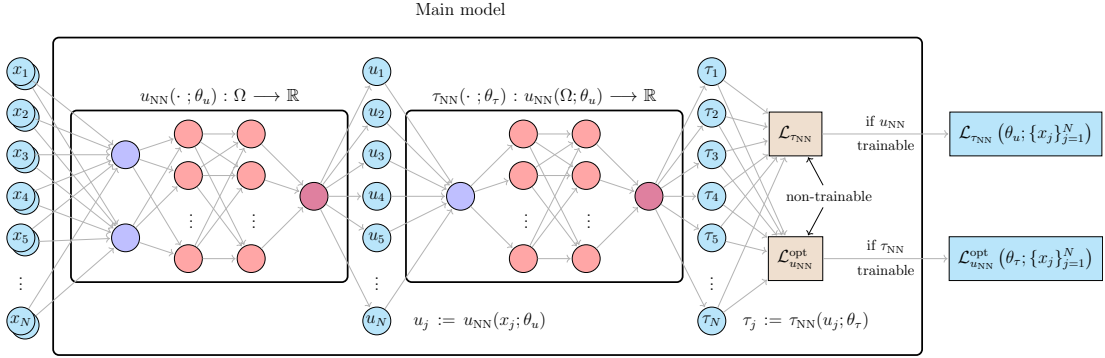

\subsection{Graph-mode execution dynamics}

We employ \emph{callbacks} to avoid interrupting the graph execution mode carried out by the Keras fitting instruction (\texttt{.fit}). Callbacks act during fitting and enable accessing certain elements of (main) models and modifying them. We utilize callbacks for loss monitoring (for \ac{WANs},  the \ac{DRM},  and the \ac{D2RM}), to activate or deactivate the trainability of the networks and switch between optimizers (for \ac{WANs} and the \ac{D2RM}), or to interchange losses (for the \ac{D2RM}) during training when iterating either over the outer- or the inner-loop.
\newpage

\section{Numerical results}
\label{section:Numerical results}

We show numerical experiments to compare the methods introduced above.  \Cref{section51} considers a simple model problem and compares \ac{WANs}, the \ac{DRM}, and the \ac{D2RM}.  \Cref{section:Comparison of performances by approaches} makes a more profound comparison considering a parametric model problem.  \Cref{section:Diffusion equation with a Dirac delta source} and \Cref{section:Advection equation in ultraweak formulation} consider sources that lead to singular problems in pure diffusion and convection equations, respectively. Finally, \Cref{section:Convection equation in 2D} considers a pure convection equation in 2D.
 
\subsection{Initial comparison of \acs{WANs},  the \acs{DRM}, and the \acs{D2RM} on a simple problem}\label{section51}

We select model problem \eqref{1DLap} in weak form with source $f=-2$, so the exact solution is $u^*=x(x-1)$. Hence, $\mathbb{U}=H^1_0(0,1)=\V$ and $T$ is the identity operator.  We solve this problem by employing \ac{WANs},  the \ac{DRM}, and the \ac{D2RM}.

We select a two-layer fully-connected \ac{NN} with 20 neurons on each layer and $\tanh$ activation functions for the architectures of $u_\text{NN}$, $v_\text{NN}$,  and $\tau_\text{NN}$.  We perform $200$ iterations for $u_\text{NN}$ in the three methods: \ac{WANs},  the \ac{DRM}, and the \ac{D2RM}.  Since in \ac{WANs} and the \ac{D2RM} we established four iterations to approximate the test maximizers, we end up with a total of $1\mathord{,}000$ training iterations ($200$ for the trial function, and $800$ for the test functions). Moreover, we select batches of size $200$ for the training and a uniform distribution for the sample generation. 

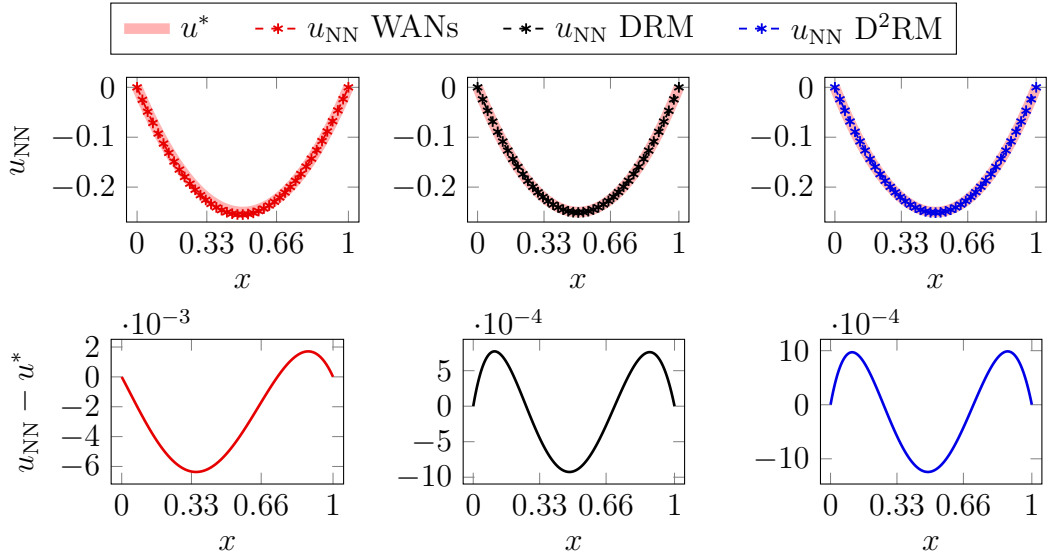
\begin{figure}[htbp]
\centering
\begin{subfigure}{\textwidth}
\centering
\begin{tikzpicture}
\begin{axis}[hide axis,
	    xmin=-0.1,
	    xmax=1.1,
	    ymin=-0.1,
	    ymax=1.1,
	 	legend style={anchor=north, /tikz/every even column/.append style={column sep=0.5cm}},
	 	legend columns = -1
            ]

    	\addlegendimage{color=red!30!white, line width=4.};
	\addlegendentry{$u^*$};	
	\addlegendimage{color=red!90!black, style=dashed, line width=.7 ,mark=asterisk, mark options = {scale=1, solid}};
	\addlegendentry{$u_{\text{NN}}$ WANs};
	\addlegendimage{color=black, style=dashed, line width=.7 ,mark=asterisk, mark options = {scale=1, solid}};
	\addlegendentry{$u_{\text{NN}}$ \ac{DRM}};
	\addlegendimage{color=blue!90!black, style=dashed, line width=.7 ,mark=asterisk, mark options = {scale=1, solid}};
	\addlegendentry{$u_{\text{NN}}$ \ac{D2RM}};
\end{axis}
\end{tikzpicture}
\end{subfigure}\vskip 0.3em%
\begin{subfigure}[t]{0.32\textwidth}
\centering
\begin{tikzpicture}
\begin{axis}[
		xmin=-0.05,
	    xmax=1.05,
	    xlabel = {$x$},
	    ymin=-0.27,
	    ymax=0.02,
	    ylabel = {$u_\text{NN}$},
	    height=3.5cm,
	    width=\textwidth,
		xtick={0,0.33,0.66,1},
	    yticklabel style={
        /pgf/number format/fixed,
        /pgf/number format/precision=2},
        scaled y ticks=false,
	    legend style={at={(0.5,0.95)}, anchor=north, /tikz/every even column/.append style={column sep=0.5cm}},
	 	]
  
\addplot[color=red!30!white, line width=4] table[x expr=\thisrow{x},y=u_analytic_plot]{Chapters/4.Chapter/FIGURES/numerical_results/section52/WAN_x(x-1)/u_analytic_plot.csv};
\addplot[color=red!90!black, style=dashed, line width=.7 ,mark=asterisk, mark options = {scale=1, solid},] table[x expr=\thisrow{x},y=u_net]{Chapters/4.Chapter/FIGURES/numerical_results/section52/WAN_x(x-1)/u_net_plot.csv};

\end{axis}
\end{tikzpicture}
\end{subfigure}\hskip 0.5em%
\begin{subfigure}[t]{0.31\textwidth}
\centering
\begin{tikzpicture}
\begin{axis}[
		xmin=-0.05,
	    xmax=1.05,
	    xlabel = {$x$},
	    ymin=-0.27,
	    ymax=0.02,
	    height=3.5cm,
	    width=\textwidth,
		xtick={0,0.33,0.66,1},
	    yticklabel style={
        /pgf/number format/fixed,
        /pgf/number format/precision=2},
        scaled y ticks=false,
	    legend style={at={(0.5,0.95)}, anchor=north, /tikz/every even column/.append style={column sep=0.5cm}},
	 	]
  
\addplot[color=red!30!white, line width=4,each nth point=2] table[x expr=\thisrow{x},y=u_analytic]{Chapters/4.Chapter/FIGURES/numerical_results/section52/WAN_x(x-1)/u_analytic.csv};
\addplot[color=black, style=dashed, line width=.7 ,mark=asterisk, mark options = {scale=1, solid}] table[x expr=\thisrow{x},y=u_net]{Chapters/4.Chapter/FIGURES/numerical_results/section51/DeepRitz_mid_1000_x(x-1)/u_net_plot.csv};

\end{axis}
\end{tikzpicture}
\end{subfigure}\hskip 0.5em%
\begin{subfigure}[t]{0.31\textwidth}
\centering
\begin{tikzpicture}
\begin{axis}[
		xmin=-0.05,
	    xmax=1.05,
	    xlabel = {$x$},
	    ymin=-0.27,
	    ymax=0.02,
	    height=3.5cm,
	    width=\textwidth,
		xtick={0,0.33,0.66,1},
	    yticklabel style={
        /pgf/number format/fixed,
        /pgf/number format/precision=2},
        scaled y ticks=false,
	    legend style={at={(0.5,0.95)}, anchor=north, /tikz/every even column/.append style={column sep=0.5cm}},
	 	]
  
\addplot[color=red!30!white, line width=4,each nth point=2] table[x expr=\thisrow{x},y=u_analytic]{Chapters/4.Chapter/FIGURES/numerical_results/section52/WAN_x(x-1)/u_analytic.csv};
\addplot[color=blue!90!black, style=dashed, line width=.7 ,mark=asterisk, mark options = {scale=1, solid},] table[x expr=\thisrow{x},y=u_net]{Chapters/4.Chapter/FIGURES/numerical_results/section52/DeepDRitz_x(x-1)/u_net_plot.csv};

\end{axis}
\end{tikzpicture}
\end{subfigure}\vskip 0.25em%
\begin{subfigure}[t]{0.32\textwidth}
\centering
\begin{tikzpicture}
\begin{axis}[
	    xmin=-0.05,
	    xmax=1.05,
	    xlabel = {$x$},
	    ylabel = {$u_\text{NN}-u^*$},
	    height=3.5cm,
	    width=\textwidth,
	    xtick={0,0.33,0.66,1},
	    scaled y ticks=true,
	 	]
  
\addplot[color=red!90!black, line width=1,each nth point=2] table[x expr=\thisrow{x},y=e_u]{Chapters/4.Chapter/FIGURES/numerical_results/section52/WAN_x(x-1)/e_u.csv};

\end{axis}
\end{tikzpicture}
\end{subfigure}\hskip 0.5em%
\begin{subfigure}[t]{0.31\textwidth}
\centering
\begin{tikzpicture}
\begin{axis}[
	    xmin=-0.05,
	    xmax=1.05,
	    xlabel = {$x$},
	    height=3.5cm,
	    width=\textwidth,
	    xtick={0,0.33,0.66,1},
	    scaled y ticks=base 10:4,
	 	]
  
\addplot[color=black, line width=1,each nth point=2] table[x expr=\thisrow{x},y=e_u]{Chapters/4.Chapter/FIGURES/numerical_results/section51/DeepRitz_mid_1000_x(x-1)/e_u.csv};

\end{axis}
\end{tikzpicture}
\end{subfigure}\hskip 0.5em%
\begin{subfigure}[t]{0.31\textwidth}
\centering
\begin{tikzpicture}
\begin{axis}[
		xmin=-0.05,
	    xmax=1.05,
	    xlabel = {$x$},
	    height=3.5cm,
	    width=\textwidth,
		xtick={0,0.33,0.66,1},
	    yticklabel style={
        /pgf/number format/fixed,
        /pgf/number format/precision=1},
        scaled y ticks=base 10:4,
	 	]
  
\addplot[color=blue!90!black, line width=1,each nth point=2] table[x expr=\thisrow{x},y=e_u]{Chapters/4.Chapter/FIGURES/numerical_results/section52/DeepDRitz_x(x-1)/e_u.csv};
\end{axis}
\end{tikzpicture}
\end{subfigure}
\caption{Trial network predictions and errors in \acs{WANs}, the \acs{DRM}, and the \acs{D2RM} at the end of the training in model problem \eqref{1DLap} with exact solution $u^*=x(x-1)$.}
\label{figure:WANs_DRM_DDRM_x(x-1)}
\end{figure}

\Cref{figure:WANs_DRM_DDRM_x(x-1)} shows the $u_\text{NN}$ network predictions and errors of the three methods at the end of the training.  We observe that \ac{WANs} produce a larger error than the \ac{DRM} and the \ac{D2RM}.

\begin{figure}[htbp]
\centering
\begin{subfigure}[b]{\textwidth}
\centering
\begin{tikzpicture}
\begin{axis}[hide axis,
	    xmin=-0.1,
	    xmax=1.1,
	    ymin=-0.1,
	    ymax=1.1,
	 	legend style={anchor=north, /tikz/every even column/.append style={column sep=0.5cm}},
	 	legend columns = -1
            ]

	\addlegendimage{line width=1, color=blue!50!white, opacity=0.8};
	\addlegendentry{$\mathcal{L}_\text{min}^\text{max}(u_\text{NN}, v_\text{NN})$};
	\addlegendimage{color=red, style=dashed, line width=0.8};
	\addlegendentry{$0$};
\end{axis}
\end{tikzpicture}
\end{subfigure}\vskip 0.5em%
\begin{subfigure}[b]{\textwidth}
\centering
\begin{tikzpicture}[]         
 \begin{axis}[ 
     xlabel = {iteration},
     scaled x ticks=false,
     ylabel = {loss},
     height=4.25cm,    
     width=0.98*\textwidth,    
     xtick={0,100,200,...,1000},
     legend style={at={(0.5,1.5)}, anchor=north, /tikz/every even column/.append style={column sep=0.5cm}},
     legend columns = -1,
     yticklabel style={
        /pgf/number format/fixed,
        /pgf/number format/precision=4},
     select coords between index/.style 2 args={
    x filter/.code={
        \ifnum\coordindex<#1\def\pgfmathresult{}\fi
        \ifnum\coordindex>#2\def\pgfmathresult{}\fi
    }
},
     ]        
\addplot[line width=1, color=blue!50!white, opacity=0.8, select coords between index={0}{1000}] table[x expr=\thisrow{iteration},y=loss_uv]{Chapters/4.Chapter/FIGURES/numerical_results/section52/WAN_x(x-1)/loss.csv};
\addplot [domain = 0:1000, line width=0.8, dashed, red]{0.};                 
\end{axis}
\end{tikzpicture}
\end{subfigure}\vskip 0.25em
\begin{subfigure}[b]{0.53\textwidth}
\centering
\begin{tikzpicture}[]         
 \begin{axis}[ 
     xlabel = {iteration}, 
     xtick={0,20,...,100},
     ylabel = {loss},
     height=4.25cm,    
     width=\textwidth,    
     yticklabel style={
        /pgf/number format/fixed,
        /pgf/number format/precision=4},
     select coords between index/.style 2 args={
    x filter/.code={
        \ifnum\coordindex<#1\def\pgfmathresult{}\fi
        \ifnum\coordindex>#2\def\pgfmathresult{}\fi
    }
},
     ]        
\addplot[line width=1, color=blue!50!white, opacity=0.8, select coords between index={0}{100}] table[x expr=\thisrow{iteration},y=loss_uv]{Chapters/4.Chapter/FIGURES/numerical_results/section52/WAN_x(x-1)/loss.csv};
\end{axis}
\end{tikzpicture}
\end{subfigure}\hfill
\begin{subfigure}[b]{0.45\textwidth}
\centering
\begin{tikzpicture}[]         
 \begin{axis}[ 
     xlabel = {iteration},
     scaled x ticks=false,
     height=4.25cm,    
     width=\textwidth,    
     xtick={500,750,...,1000},
     yticklabel style={
        /pgf/number format/fixed,
        /pgf/number format/precision=4},
     select coords between index/.style 2 args={
    x filter/.code={
        \ifnum\coordindex<#1\def\pgfmathresult{}\fi
        \ifnum\coordindex>#2\def\pgfmathresult{}\fi
    }
},
     ]        
\addplot[line width=1, color=blue!50!white, opacity=0.8, select coords between index={500}{1000}] table[x expr=\thisrow{iteration},y=loss_uv]{Chapters/4.Chapter/FIGURES/numerical_results/section52/WAN_x(x-1)/loss.csv};
\addplot [domain = 500:1000, line width=0.8, dashed, red]{0.};           
\end{axis}
\end{tikzpicture}
\end{subfigure}
\caption{Loss evolution during the \acs{WANs} training for model problem \eqref{1DLap} with exact solution $u^*=x(x-1)$.}
\label{figure:WAN_x(x-1)_loss}
\end{figure}
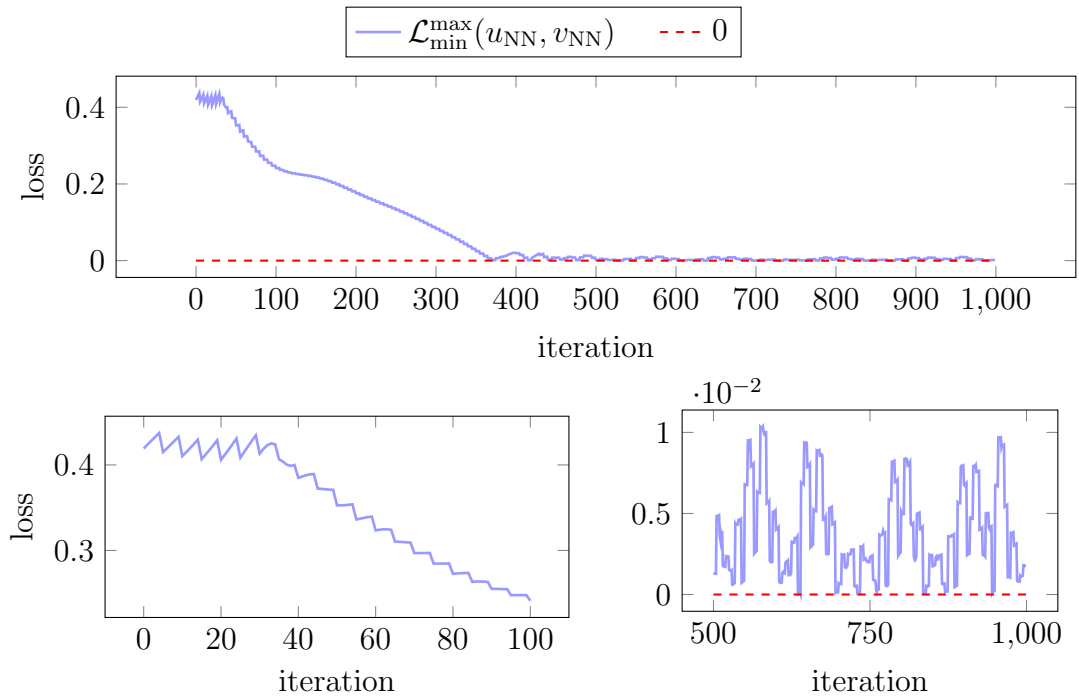

\Cref{figure:WAN_x(x-1)_loss} shows the loss evolution for \ac{WANs}.  Every five iterations, the loss decreases (minimization in $u_\text{NN}$) and increases in the remaining iterations (maximization in $v_\text{NN}$ with fixed $u_\text{NN}$).  From iteration $500$ onwards, the loss stops improving and oscillates above the optimal value. 

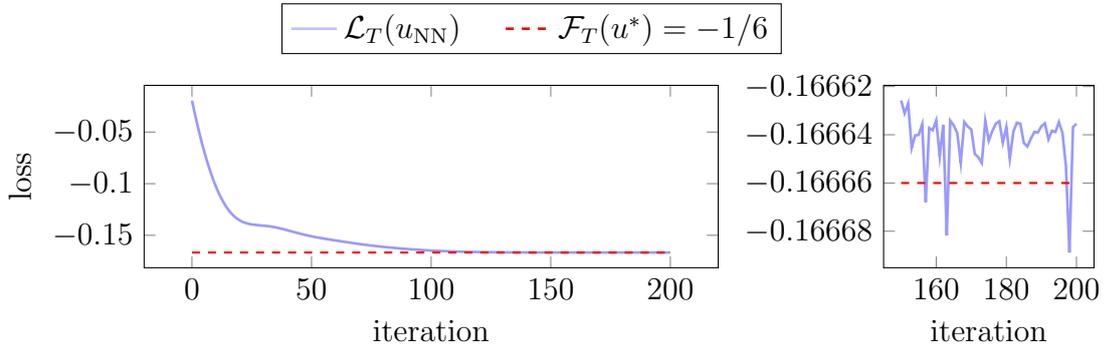
\begin{figure}[htbp]
\centering
\begin{subfigure}{\textwidth}
\centering
\begin{tikzpicture}
\begin{axis}[hide axis,
	    xmin=-0.1,
	    xmax=1.1,
	    ymin=-0.1,
	    ymax=1.1,
	 	legend style={anchor=north, /tikz/every even column/.append style={column sep=0.5cm}},
	 	legend columns = -1
            ]

	\addlegendimage{color=blue!50!white, opacity=.5, line width=1};
	\addlegendentry{$\mathcal{L}_T(u_\text{NN})$};
	\addlegendimage{color=red, style=dashed, line width=1.};
	\addlegendentry{$\mathcal{F}_T(u^*)=-1/6$};
	
\end{axis}
\end{tikzpicture}
\end{subfigure}\vskip 0.25em%
\centering
\begin{subfigure}[t]{0.63\textwidth}
\centering
\begin{tikzpicture}[]         
 \begin{axis}[
     xlabel = {iteration},
     scaled x ticks=false,
     ylabel = {loss},
     height=4cm,    
     width=\textwidth,    
     xtick={0,50,100,...,200},
     legend style={at={(0.5,1.5)}, anchor=north, /tikz/every even column/.append style={column sep=0.5cm}},
     legend columns = -1,
     yticklabel style={
        /pgf/number format/fixed,
        /pgf/number format/precision=4},
     select coords between index/.style 2 args={
    x filter/.code={
        \ifnum\coordindex<#1\def\pgfmathresult{}\fi
        \ifnum\coordindex>#2\def\pgfmathresult{}\fi
    }
},
     ]        

\addplot[line width=1, color=blue!50!white, opacity=0.8, select coords between index={0}{200}] table[x expr=\thisrow{iteration},y=loss]{Chapters/4.Chapter/FIGURES/numerical_results/section51/DeepRitz_mid_1000_x(x-1)/loss.csv};
\addplot [domain = 0:200, line width=0.8, dashed, red]{-1/6};
                 
\end{axis}
\end{tikzpicture}
\label{figure:DRM_loss_global}
\end{subfigure} \hskip 1em%
\begin{subfigure}[t]{0.30\textwidth}
\centering
\begin{tikzpicture}[]         
 \begin{axis}[
     xlabel = {iteration},
     scaled x ticks=false,
     height=4cm,    
     width=\textwidth,    
     scaled y ticks=true,
     legend style={at={(0.5,1.5)}, anchor=north, /tikz/every even column/.append style={column sep=0.5cm}},
     legend columns = -1,
     yticklabel style={
        /pgf/number format/fixed,
        /pgf/number format/precision=5},
     select coords between index/.style 2 args={
    x filter/.code={
        \ifnum\coordindex<#1\def\pgfmathresult{}\fi
        \ifnum\coordindex>#2\def\pgfmathresult{}\fi
    }
},
     ]        

\addplot[line width=1, color=blue!50!white, opacity=0.8, select coords between index={150}{200}] table[x expr=\thisrow{iteration},y=loss]{Chapters/4.Chapter/FIGURES/numerical_results/section51/DeepRitz_mid_1000_x(x-1)/loss.csv};
\addplot [domain = 150:200, line width=0.8, dashed, red]{-1/6};
\end{axis}
\end{tikzpicture}
\label{figure:DRM_loss_end}
\end{subfigure}
\caption{Loss evolution during the DRM training in model problem \eqref{1DLap} with exact solution $u^*=x(x-1)$.}
\label{figure:DRM_x(x-1)_loss}
\end{figure}

\Cref{figure:DRM_x(x-1)_loss} shows the evolution of the loss for the \ac{DRM}.  Here, we have a single minimization, so the observed noisy behavior of the loss towards the end of the training is attributed to the optimizer performance. Due to the lower complexity of the training, we achieve better convergence performance than with \ac{WANs}.

\Cref{figure:DDRM_x(x-1)_loss} shows the loss evolutions for the \ac{D2RM}.  Here, we have a nested min-min optimization.  At each iteration, we evaluate both $\mathcal{L}_{\tau_\text{NN}}$ and $\mathcal{L}_{u_\text{NN}}^\text{opt}$, even if we are only optimizing with respect to one of them. We superimpose both losses, each one with its own scale (the left vertical axis corresponds to $\mathcal{L}_{\tau_\text{NN}}$ and the right vertical axis to $\mathcal{L}_{u_\text{NN}}^\text{opt}$).  Both losses exhibit a decreasing staircase shape with downward-sloping steps.  Jumps occur when the optimization is performed with respect to $\mathcal{L}_{\tau_\text{NN}}$, which suggests that $\mathcal{L}_{u_\text{NN}}^\text{opt}$ takes longer to converge as it depends on the convergence of $\mathcal{L}_{\tau_\text{NN}}$ (inner- vs. outer-loop). 

\begin{figure}[htbp]
\centering
\begin{subfigure}{.98\textwidth}
\centering
\hskip 5em \begin{tikzpicture}
\begin{axis}[
		hide axis,
	    xmin=-0.1,
	    xmax=1.1,
	    ymin=-0.1,
	    ymax=1.1,
	    legend style={ anchor=north, /tikz/every even column/.append style={column sep=0.5cm}},
	    legend columns = -1,]

	\addlegendimage{line width=1, color=blue!50!white, opacity=0.8};
	\addlegendentry{$\mathcal{L}_{\tau_\text{NN}}(u_\text{NN})$};
	\addlegendimage{line width=1, color=purple!50!white, opacity=0.8};
	\addlegendentry{$\mathcal{L}_{u_\text{NN}}^\text{opt}(\tau_\text{NN}(u_\text{NN}))$};
	\addlegendimage{color=red, style=dashed, line width=0.8};
	\addlegendentry{$\mathcal{F}_{u^*}^\text{opt}(Tu^*)=\mathcal{F}_T(u^*)=-1/6$};
\end{axis}
\end{tikzpicture}
\end{subfigure}\vskip 0.35em
\begin{subfigure}[t]{\textwidth}
\centering
\begin{tikzpicture}[]         
 \begin{axis}[
     xlabel = {iteration},
     scaled x ticks=false,
     ylabel = {loss},
     height=4cm,    
     width=0.98*\textwidth,    
     xtick={0,100,200,...,1000},
     legend style={at={(0.5,1.5)}, anchor=north, /tikz/every even column/.append style={column sep=0.5cm}},
     legend columns = -1,
     yticklabel style={
        /pgf/number format/fixed,
        /pgf/number format/precision=4},
     select coords between index/.style 2 args={
    x filter/.code={
        \ifnum\coordindex<#1\def\pgfmathresult{}\fi
        \ifnum\coordindex>#2\def\pgfmathresult{}\fi
    }
},
     ]        

\addplot[line width=1, color=blue!50!white, opacity=0.8, select coords between index={0}{1000}, each nth point=2] table[x expr=\thisrow{iteration},y=loss_u]{Chapters/4.Chapter/FIGURES/numerical_results/section52/DeepDRitz_x(x-1)/loss.csv};
\addplot[line width=1, color=purple!50!white, opacity=0.8, select coords between index={0}{1000}, each nth point=2] table[x expr=\thisrow{iteration},y=loss_T]{Chapters/4.Chapter/FIGURES/numerical_results/section52/DeepDRitz_x(x-1)/loss.csv};
\addplot [domain = 0:1000, line width=0.8, dashed, red]{-1/6};
\end{axis}
\end{tikzpicture}
\end{subfigure}\vskip 0.25em
\begin{subfigure}[t]{0.5\textwidth}
\centering
\begin{tikzpicture}

\begin{axis}[
name=plot1,
xlabel=iteration,
xtick={0,10,20,...,50},
ylabel=loss,
axis y line*=left,
height=4cm,    
width=\textwidth,    
yticklabel style={blue},
select coords between index/.style 2 args={x filter/.code={
        \ifnum\coordindex<#1\def\pgfmathresult{}\fi
        \ifnum\coordindex>#2\def\pgfmathresult{}\fi}},
]
\draw [dashed, gray!50!white] (0,-1) -- (0,1);
\draw [dashed, gray!50!white] (5,-1) -- (5,1);
\draw [dashed, gray!50!white] (10,-1) -- (10,1);
\draw [dashed, gray!50!white] (15,-1) -- (15,1);
\draw [dashed, gray!50!white] (20,-1) -- (20,1);
\draw [dashed, gray!50!white] (25,-1) -- (25,1);
\draw [dashed, gray!50!white] (30,-1) -- (30,1);
\draw [dashed, gray!50!white] (35,-1) -- (35,1);
\draw [dashed, gray!50!white] (40,-1) -- (40,1);
\draw [dashed, gray!50!white] (45,-1) -- (45,1);
\draw [dashed, gray!50!white] (50,-1) -- (50,1);
\node[text width=5em, text centered] (min l_u) at (8,-0.05) {minimizing $\mathcal{L}_{\tau_\text{NN}}$};
\path[->] (min l_u) edge (0,-0.004);
\path[->] (min l_u) edge (5,-0.006);
\path[->] (min l_u) edge (10,-0.008);
\path[->] (min l_u) edge (15,-0.01);
\node[text width=5em, text centered] (min l_T) at (40,-0.012) {minimizing $\phantom{fvghtg}\mathcal{L}_{u_\text{NN}}^\text{opt}$};
\path[<->] (30.5,-0.012) edge (34.5,-0.0155);
\path[<->] (35.5,-0.02) edge (39.5,-0.0245);
\path[<->] (40.5,-0.032) edge (44.5,-0.0375);
\addplot[line width=1, color=blue!50!white, opacity=0.8, select coords between index={0}{50}] table[x expr=\thisrow{iteration},y=loss_u]{Chapters/4.Chapter/FIGURES/numerical_results/section52/DeepDRitz_x(x-1)/loss.csv};
\end{axis}

\begin{axis}[
name=plot2,
axis y line*=right,
axis x line=none,
height=4cm,    
width=\textwidth,   
yticklabel style={red},
select coords between index/.style 2 args={x filter/.code={
        \ifnum\coordindex<#1\def\pgfmathresult{}\fi
        \ifnum\coordindex>#2\def\pgfmathresult{}\fi}},
]

\addplot[line width=1, color=purple!50!white, opacity=0.8, select coords between index={0}{50}] table[x expr=\thisrow{iteration},y=loss_T]{Chapters/4.Chapter/FIGURES/numerical_results/section52/DeepDRitz_x(x-1)/loss.csv};
\end{axis}
\end{tikzpicture}
\end{subfigure}\hfill
\begin{subfigure}[t]{0.44\textwidth}
\centering
\begin{tikzpicture}[]         
 \begin{axis}[
     xlabel = {iteration},
     scaled x ticks=false, 
     height=4cm,    
     width=\textwidth,     
     xtick={500,750,...,1000},
     legend style={at={(0.5,1.5)}, anchor=north, /tikz/every even column/.append style={column sep=0.5cm}},
     legend columns = -1,
     yticklabel style={
        /pgf/number format/fixed,
        /pgf/number format/precision=5},
     select coords between index/.style 2 args={
    x filter/.code={
        \ifnum\coordindex<#1\def\pgfmathresult{}\fi
        \ifnum\coordindex>#2\def\pgfmathresult{}\fi
    }
},
     ]        

\addplot[line width=1, color=blue!50!white, opacity=0.8, select coords between index={500}{1000}] table[x expr=\thisrow{iteration},y=loss_u]{Chapters/4.Chapter/FIGURES/numerical_results/section52/DeepDRitz_x(x-1)/loss.csv};
\addplot[line width=1, color=purple!50!white, opacity=0.8, select coords between index={500}{1000}] table[x expr=\thisrow{iteration},y=loss_T]{Chapters/4.Chapter/FIGURES/numerical_results/section52/DeepDRitz_x(x-1)/loss.csv};
\addplot [domain = 500:1000, line width=0.8, dashed, red]{-1/6};
\end{axis}
\end{tikzpicture}
\end{subfigure}
\caption{Loss evolution of the \acs{D2RM} training for model problem \eqref{1DLap} with exact solution $u^*=x(x-1)$.}
\label{figure:DDRM_x(x-1)_loss}
\end{figure}

The relative errors of the trial network predictions\footnote{To approximate $\frac{\Vert u_{\text{NN}}-u^*\Vert_\mathbb{U}}{\Vert u^*\Vert_\mathbb{U}}\times 100$, we perform a composite intermediate-point rule with $10^4$ integration nodes for the numerator, and analytically calculate $\Vert u^*\Vert_\mathbb{U} = \sqrt{3}/3$ for the denominator.} at the end of the training are $3.12\%$, $0.99\%$,  and $1.31\%$ in \ac{WANs}, the \ac{DRM}, and the \ac{D2RM}, respectively.

\subsection{Comparison of WANs,  the DRM, and the D$^2$RM on singular problems}\label{section:Comparison of performances by approaches}

Now,  we focus on the evolution of the relative error for the previous weak formulation of \eqref{1DLap}, but selecting the source so that the solution varies according to a parameter:
\begin{equation}
u_\alpha^* = x^{\alpha}(x-1)\in H^1_0(0,1), \qquad \alpha>1/2.
\end{equation}

\subsubsection{Without singularities: \boldmath{$\alpha\geq 1$}}
We experiment individually for $\alpha \in\{2, 5,10\}$ with $5\mathord{,}000$ training iterations for $u_\text{NN}$ in \ac{WANs}, the \ac{DRM}, and the \ac{D2RM}. Note that this corresponds to a total of $25\mathord{,}000$ iterations for \ac{WANs} and the \ac{D2RM} when taking into account the iterations dedicated to the test maximizers. 

\Cref{table:experiments2} displays the relative errors along different stages of the training, and \Cref{figure:xalpha(x-1)_predictions} shows the trial network predictions and error functions at the end of training.

\begin{table}[htbp]
\centering
\begin{tabular}{|c|c||c|c|c|c|c|}
\toprule
\multicolumn{2}{|c||}{\textbf{Training progress}} & $\mathbf{4}\bf{\%}$ & $\mathbf{20}\bf{\%}$ & $\mathbf{40}\bf{\%}$ & $\mathbf{60}\bf{\%}$ & $\mathbf{100}\bf{\%}$\\ \bottomrule
\toprule
\textbf{Method} & $\boldsymbol{\alpha}$ & \multicolumn{5}{c|}{$\frac{\Vert u_\text{NN}-u^*\Vert_\mathbb{U}}{\Vert u^*\Vert_\mathbb{U}}\times 100$}\\ \midrule
\multicolumn{1}{|l||}{\multirow{3}{*}{WANs}}     & $2$ &  $2.49\%$ & $3.43\%$  & $5.92\%$ & $7.40\%$ & $10.07\%$\\ \cline{2-7} 
\multicolumn{1}{|l||}{}                       & $5$ & $58.69\%$ & $41.78\%$ & $63.03\%$ & $74.06\%$ & $40.33\%$ \\ \cline{2-7} 
\multicolumn{1}{|l||}{}                       & $10$ & $93.15\%$ & $93.95\%$ & $89.04\%$ & $68.85\%$ & $367.61\%$ \\ \midrule
\multicolumn{1}{|l||}{\multirow{3}{*}{\ac{DRM}}}  & $2$ & $3.40\%$ & $2.47\%$ & $1.17\%$ & $0.30\%$ & $0.23\%$ \\ \cline{2-7} 
\multicolumn{1}{|l||}{}                       & $5$ &  $50.50\%$&$8.53\%$&$2.83\%$&$2.27\%$&$1.59\%$\\ \cline{2-7} 
\multicolumn{1}{|l||}{}                       & $10$ &  $58.55\%$&$9.51\%$&$3.39\%$&$2.83\%$&$1.69\%$\\ \midrule
\multicolumn{1}{|l||}{\multirow{3}{*}{\ac{D2RM}}} & $2$ &  $3.27\%$&$1.84\%$&$0.31\%$&$0.27\%$&$0.54\%$\\ \cline{2-7} 
\multicolumn{1}{|l||}{}                       & $5$ & $59.13\%$&$12.93\%$&$2.92\%$&$2.52\%$&$1.56\%$\\ \cline{2-7} 
\multicolumn{1}{|l||}{}                       & $10$ & $82.31\%$ & $20.18\%$ & $6.24\%$& $2.68\%$& $2.60\%$\\ \bottomrule
\toprule
\textbf{Method} & $\boldsymbol{\alpha}$ & \multicolumn{5}{c|}{$\frac{\Vert \tau_\text{NN}(u_\text{NN}) - Tu^*\Vert_\V}{\Vert Tu^*\Vert_\V}\times 100$}\\ \midrule
\multicolumn{1}{|l||}{\multirow{3}{*}{\ac{D2RM}}} & $2$ &  $3.36\%$&$1.94\%$&$0.48\%$&$0.45\%$&$0.58\%$\\ \cline{2-7} 
\multicolumn{1}{|l||}{}                       & $5$ &  $59.13\%$&$12.93\%$&$2.90\%$&$2.11\%$&$1.55\%$\\ \cline{2-7} 
\multicolumn{1}{|l||}{}                       & $10$ &  $83.51\%$&$20.18\%$&$6.23\%$&$2.68\%$&$2.61\%$\\ \bottomrule
\end{tabular}%
\caption{Relative errors of $u_\text{NN}$ (in \acs{WANs}, the \acs{DRM}, and the \acs{D2RM}) and $\tau_\text{NN} (u_\text{NN})$ (in the \acs{D2RM}) along different stages of the training progress in problem \eqref{1DLap} with exact solution $u_\alpha^*=x^{\alpha}(x-1)$ and $\alpha\in\{2,5,10\}$. }
\label{table:experiments2}
\end{table}
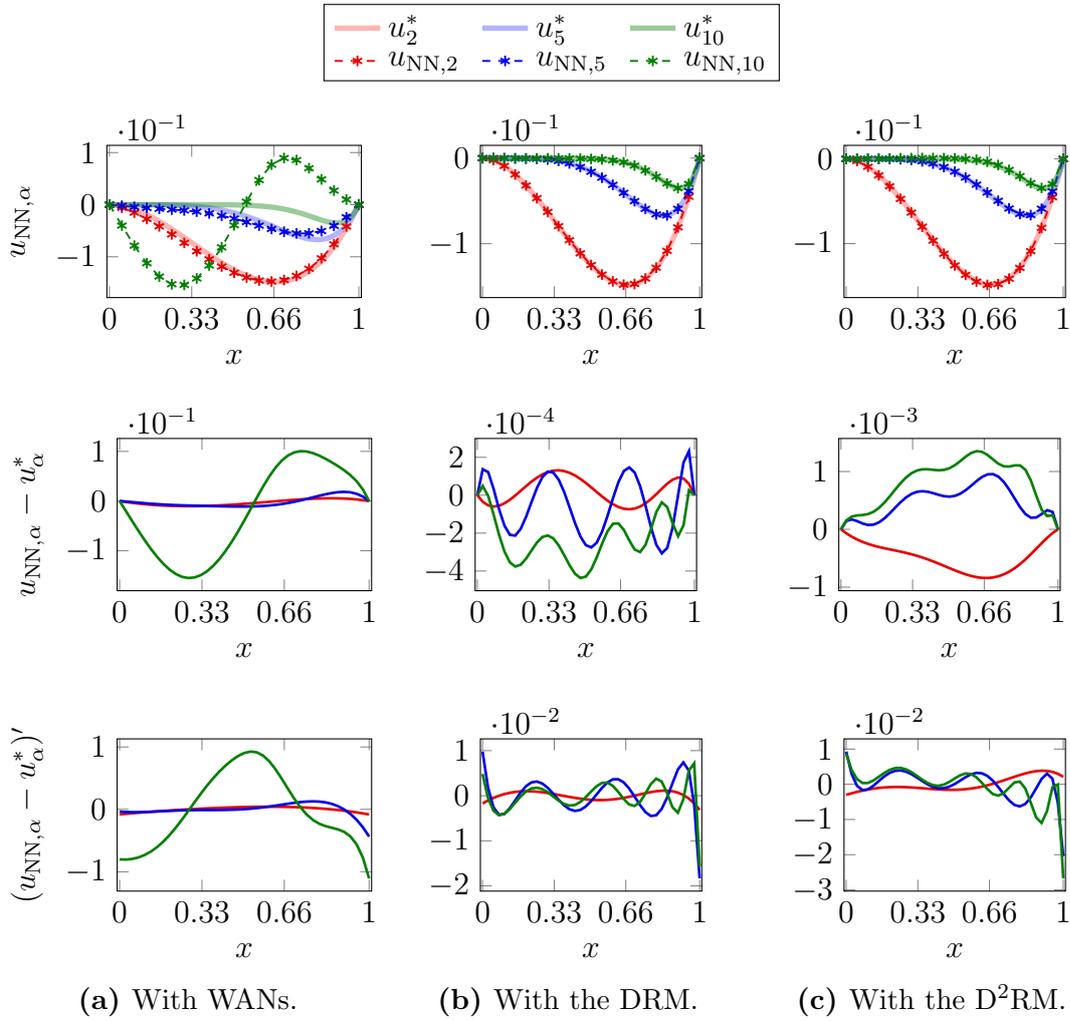
\begin{figure}[htbp]
\centering
\begin{subfigure}[t]{\textwidth}
\centering
 \begin{tikzpicture}
\pgfplotsset{A/.append style={
		hide axis,
	    xmin=-0.1,
	    xmax=1.1,
	    ymin=-0.1,
	    ymax=1.1,
	 	legend style={anchor=north, /tikz/every even column/.append style={column sep=0.2cm}},
	 	legend columns = 3,
	 	legend cell align={left},
            }
}
\begin{axis}[A]
\addlegendimage{line width=2, color=red!30!white}
\addlegendentry{$u_{2}^*$}
\addlegendimage{line width=2, color=blue!30!white}
\addlegendentry{$u_{5}^*$};
\addlegendimage{line width=2, color=green!50!black, opacity=0.4}
\addlegendentry{$u_{10}^*$};

\addlegendimage{dashed, color=red!90!black,  line width=0.7, mark=asterisk, mark options = {scale=1, solid}}
\addlegendentry{$u_{\text{NN},2}$}
\addlegendimage{dashed, color=blue!90!black,  line width=0.7, mark=asterisk, mark options = {scale=1, solid}}
\addlegendentry{$u_{\text{NN},5}$};
\addlegendimage{dashed, color=green!50!black,  line width=0.7, mark=asterisk, mark options = {scale=1, solid}}
\addlegendentry{$u_{\text{NN},10}$};
\end{axis}
\end{tikzpicture}
\end{subfigure}\vskip 1em%
\begin{subfigure}[t]{0.34\textwidth}
\centering
\begin{tikzpicture}         
 \begin{axis}[xmin=-0.01,         
     xmax=1.01,      
     xlabel = {$x$},  
     ylabel = {$u_{\text{NN},\alpha}$},
     xtick={0,0.33,0.66,1},
     height=3.6cm,    
     width=\textwidth,
     scaled y ticks=base 10:1,    
     yticklabel style={
        /pgf/number format/fixed,
        /pgf/number format/precision=2}, ]        

\addplot [domain = 0:1, samples=1000, line width=2, color=red!30!white]{x^2*(x-1)};
\addplot[dashed, color=red!90!black,  line width=0.7, mark=asterisk, mark options = {scale=1, solid}, each nth point=2] table[x expr=\thisrow{x},y=u_net]{Chapters/4.Chapter/FIGURES/numerical_results/section53/WAN_x2(x-1)_25000/u_net_plot.csv};
\addplot [domain = 0:1, samples=1000, line width=2, color=blue!30!white]{x^5*(x-1)};
\addplot[dashed, color=blue!90!black,  line width=0.7, mark=asterisk, mark options = {scale=1, solid}, each nth point=2] table[x expr=\thisrow{x},y=u_net]{Chapters/4.Chapter/FIGURES/numerical_results/section53/WAN_x5(x-1)_25000/u_net_plot.csv};
\addplot [domain = 0:1, samples=1000, line width=2,green!50!black, opacity=0.4]{x^10*(x-1)};
\addplot[dashed, color=green!50!black,  line width=0.7, mark=asterisk, mark options = {scale=1, solid}, each nth point=2] table[x expr=\thisrow{x},y=u_net]{Chapters/4.Chapter/FIGURES/numerical_results/section53/WAN_x10(x-1)_25000/u_net_plot.csv};
\end{axis}
\end{tikzpicture}
\end{subfigure}\hfill%
\begin{subfigure}[t]{0.31\textwidth}
\centering
\begin{tikzpicture}         
 \begin{axis}[xmin=-0.01,         
     xmax=1.01,      
     xlabel = {$x$},
     ytick={0,-0.1},
     xtick={0,0.33,0.66,1},
     height=3.6cm,    
     width=\textwidth,
     scaled y ticks=base 10:1,   
     yticklabel style={
        /pgf/number format/fixed,
        /pgf/number format/precision=1}, ]        
\addplot [domain = 0:1, samples=1000, line width=2, red!30!white]{x^2*(x-1)};
\addplot[dashed, color=red!90!black,  line width=.7, mark=asterisk, mark options = {scale=1, solid}, each nth point=2] table[x expr=\thisrow{x},y=u_net]{Chapters/4.Chapter/FIGURES/numerical_results/section53/DeepRitz_x2(x-1)_5000/u_net_plot.csv};
\addplot [domain = 0:1, samples=1000, line width=2, blue!30!white]{x^5*(x-1)};
\addplot[dashed, color=blue!90!black,  line width=.7, mark=asterisk, mark options = {scale=1, solid}, each nth point=2] table[x expr=\thisrow{x},y=u_net]{Chapters/4.Chapter/FIGURES/numerical_results/section53/DeepRitz_x5(x-1)_5000/u_net_plot.csv};
\addplot [domain = 0:1, samples=1000, line width=2,green!50!black, opacity=0.4]{x^10*(x-1)};
\addplot[dashed, color=green!50!black,  line width=.7, mark=asterisk, mark options = {scale=1, solid}, each nth point=2] table[x expr=\thisrow{x},y=u_net]{Chapters/4.Chapter/FIGURES/numerical_results/section53/DeepRitz_x10(x-1)_5000/u_net_plot.csv};
\end{axis}
\end{tikzpicture}
\end{subfigure}\hfill%
\begin{subfigure}[t]{0.31\textwidth}
\centering
\begin{tikzpicture}         
 \begin{axis}[xmin=-0.01,         
     xmax=1.01,      
     xlabel = {$x$},
     xtick={0,0.33,0.66,1},
     ytick={0,-0.1},
     height=3.6cm,    
     width=\textwidth,
     scaled y ticks=base 10:1,   
     yticklabel style={
        /pgf/number format/fixed,
        /pgf/number format/precision=1}, ]        
\addplot [domain = 0:1, samples=1000, line width=2, red!30!white]{x^2*(x-1)};
\addplot[dashed, color=red!90!black,  line width=.7, mark=asterisk, mark options = {scale=1, solid}, each nth point=2] table[x expr=\thisrow{x},y=u_net]{Chapters/4.Chapter/FIGURES/numerical_results/section53/DeepDRitz_x2(x-1)_25000/u_net_plot.csv};
\addplot [domain = 0:1, samples=1000, line width=2, blue!30!white]{x^5*(x-1)};
\addplot[dashed, color=blue!90!black,  line width=.7, mark=asterisk, mark options = {scale=1, solid}, each nth point=2] table[x expr=\thisrow{x},y=u_net]{Chapters/4.Chapter/FIGURES/numerical_results/section53/DeepDRitz_x5(x-1)_25000/u_net_plot.csv};
\addplot [domain = 0:1, samples=1000, line width=2,green!50!black, opacity=0.4]{x^10*(x-1)};
\addplot[dashed,color=green!50!black,  line width=.7, mark=asterisk, mark options = {scale=1, solid}, each nth point=2] table[x expr=\thisrow{x},y=u_net]{Chapters/4.Chapter/FIGURES/numerical_results/section53/DeepDRitz_x10(x-1)_25000/u_net_plot.csv};
\end{axis}
\end{tikzpicture}
\end{subfigure}
\vskip 1em
\begin{subfigure}[t]{0.34\textwidth}
\centering
\begin{tikzpicture}         
 \begin{axis}[xmin=-0.01,         
     xmax=1.01,      
     xlabel = {$x$},
     ylabel = {$u_{\text{NN}, \alpha}-u_\alpha^*\phantom{)'}$},
     xtick={0,0.33,0.66,1},
     height=3.6cm,    
     width=\textwidth,    
     scaled y ticks=base 10:1,
     yticklabel style={
        /pgf/number format/fixed,
        /pgf/number format/precision=1}, ]        

\addplot[color=red!90!black,  line width=1] table[x expr=\thisrow{x},y=e_u]{Chapters/4.Chapter/FIGURES/numerical_results/section53/WAN_x2(x-1)_25000/e_u_plot.csv};
\addplot[color=blue!90!black,  line width=1] table[x expr=\thisrow{x},y=e_u]{Chapters/4.Chapter/FIGURES/numerical_results/section53/WAN_x5(x-1)_25000/e_u_plot.csv};
\addplot[color=green!50!black,  line width=1] table[x expr=\thisrow{x},y=e_u]{Chapters/4.Chapter/FIGURES/numerical_results/section53/WAN_x10(x-1)_25000/e_u_plot.csv};
\end{axis}
\end{tikzpicture}
\label{figure:WAN_xalpha(x-1)_25000}
\end{subfigure}\hfill%
\begin{subfigure}[t]{0.31\textwidth}
\centering
\begin{tikzpicture}         
 \begin{axis}[xmin=-0.01,         
     xmax=1.01,      
     xlabel = {$x$},
     xtick={0,0.33,0.66,1},
     height=3.6cm,    
     width=\textwidth,    
     yticklabel style={
        /pgf/number format/fixed,
        /pgf/number format/precision=1}, ]        
\addplot[color=red!90!black,  line width=1] table[x expr=\thisrow{x},y=e_u]{Chapters/4.Chapter/FIGURES/numerical_results/section53/DeepRitz_x2(x-1)_5000/e_u_plot.csv};
\addplot[color=blue!90!black,  line width=1] table[x expr=\thisrow{x},y=e_u]{Chapters/4.Chapter/FIGURES/numerical_results/section53/DeepRitz_x5(x-1)_5000/e_u_plot.csv};
\addplot[color=green!50!black,  line width=1] table[x expr=\thisrow{x},y=e_u]{Chapters/4.Chapter/FIGURES/numerical_results/section53/DeepRitz_x10(x-1)_5000/e_u_plot.csv};
\end{axis}
\end{tikzpicture}
\label{figure:DRM_xalpha(x-1)_5000}
\end{subfigure}\hfill%
\begin{subfigure}[t]{0.31\textwidth}
\centering
\begin{tikzpicture}         
 \begin{axis}[xmin=-0.01,         
     xmax=1.01,      
     xlabel = {$x$},
     xtick={0,0.33,0.66,1},  
     height=3.6cm,    
     width=\textwidth,    
     yticklabel style={
        /pgf/number format/fixed,
        /pgf/number format/precision=1}, ]        
\addplot[color=red!90!black,  line width=1] table[x expr=\thisrow{x},y=e_u]{Chapters/4.Chapter/FIGURES/numerical_results/section53/DeepDRitz_x2(x-1)_25000/e_u_plot.csv};
\addplot[color=blue!90!black,  line width=1] table[x expr=\thisrow{x},y=e_u]{Chapters/4.Chapter/FIGURES/numerical_results/section53/DeepDRitz_x5(x-1)_25000/e_u_plot.csv};
\addplot[color=green!50!black,  line width=1] table[x expr=\thisrow{x},y=e_u]{Chapters/4.Chapter/FIGURES/numerical_results/section53/DeepDRitz_x10(x-1)_25000/e_u_plot.csv};
\end{axis}
\end{tikzpicture}
\label{figure:DDRM_xalpha(x-1)_25000}
\end{subfigure}
\begin{subfigure}[t]{0.34\textwidth}
\centering
\begin{tikzpicture}         
 \begin{axis}[xmin=-0.01,         
     xmax=1.01,      
     xlabel = {$x$}, 
     ylabel = {$(u_{\text{NN},\alpha}-u_\alpha^*)'$},
     xtick={0,0.33,0.66,1},
     height=3.6cm,    
     width=\textwidth,    
     scaled y ticks=base 10:0,
     yticklabel style={
        /pgf/number format/fixed,
        /pgf/number format/precision=2}, ]        
\addplot[color=red!90!black,  line width=1] table[x expr=\thisrow{x},y=de_u]{Chapters/4.Chapter/FIGURES/numerical_results/section53/WAN_x2(x-1)_25000/de_u_plot.csv};
\addplot[color=blue!90!black,  line width=1] table[x expr=\thisrow{x},y=de_u]{Chapters/4.Chapter/FIGURES/numerical_results/section53/WAN_x5(x-1)_25000/de_u_plot.csv};
\addplot[color=green!50!black,  line width=1] table[x expr=\thisrow{x},y=de_u]{Chapters/4.Chapter/FIGURES/numerical_results/section53/WAN_x10(x-1)_25000/de_u_plot.csv};
\end{axis}
\end{tikzpicture}
\caption{With \acs{WANs}.}
\end{subfigure}\hfill%
\begin{subfigure}[t]{0.31\textwidth}
\centering
\begin{tikzpicture}         
 \begin{axis}[xmin=-0.01,         
     xmax=1.01,      
     xlabel = {$x$},
     xtick={0,0.33,0.66,1}, 
     height=3.6cm,    
     width=\textwidth,    
     yticklabel style={
        /pgf/number format/fixed,
        /pgf/number format/precision=2}, ]        
\addplot[color=red!90!black,  line width=1] table[x expr=\thisrow{x},y=de_u]{Chapters/4.Chapter/FIGURES/numerical_results/section53/DeepRitz_x2(x-1)_5000/de_u_plot.csv};
\addplot[color=blue!90!black,  line width=1] table[x expr=\thisrow{x},y=de_u]{Chapters/4.Chapter/FIGURES/numerical_results/section53/DeepRitz_x5(x-1)_5000/de_u_plot.csv};
\addplot[color=green!50!black,  line width=1] table[x expr=\thisrow{x},y=de_u]{Chapters/4.Chapter/FIGURES/numerical_results/section53/DeepRitz_x10(x-1)_5000/de_u_plot.csv};
\end{axis}
\end{tikzpicture}
\caption{With the \acs{DRM}.}
\end{subfigure}\hfill%
\begin{subfigure}[t]{0.31\textwidth}
\centering
\begin{tikzpicture}         
 \begin{axis}[xmin=-0.01,         
     xmax=1.01,      
     xlabel = {$x$},
     xtick={0,0.33,0.66,1},     
     height=3.6cm,    
     width=\textwidth,    
     yticklabel style={
        /pgf/number format/fixed,
        /pgf/number format/precision=2}, ]        
\addplot[color=red!90!black,  line width=1] table[x expr=\thisrow{x},y=de_u]{Chapters/4.Chapter/FIGURES/numerical_results/section53/DeepDRitz_x2(x-1)_25000/de_u_plot.csv};
\addplot[color=blue!90!black,  line width=1] table[x expr=\thisrow{x},y=de_u]{Chapters/4.Chapter/FIGURES/numerical_results/section53/DeepDRitz_x5(x-1)_25000/de_u_plot.csv};
\addplot[color=green!50!black,  line width=1] table[x expr=\thisrow{x},y=de_u]{Chapters/4.Chapter/FIGURES/numerical_results/section53/DeepDRitz_x10(x-1)_25000/de_u_plot.csv};
\end{axis}
\end{tikzpicture}
\caption{With the \acs{D2RM}.}
\end{subfigure}
\caption{Trial network predictions and errors for \acs{WANs}, the \acs{DRM}, and the \acs{D2RM} in model problem \eqref{1DLap} with exact solution $u_\alpha^*=x^{\alpha}(x-1)$ for $\alpha\in\{2,5,10\}$. }
\label{figure:xalpha(x-1)_predictions}
\end{figure}

\ac{WANs} show poor results in all the cases,  with a clear non-convergent tendency as the training progresses, possibly justified by the unstable behavior of the method at the continuous level. Thus, we discard the WANs for the remainder of the experiments and focus on the other two methods.

In the \ac{DRM} and the \ac{D2RM}, we observe a decreasing behavior of the relative error during training.  We highlight the nearly identical behavior of the relative norm errors of $u_\text{NN}$ and $\tau_\text{NN}(u_\text{NN})$ in the \ac{D2RM},  which suggests that the \ac{D2RM} behaves as the \ac{DRM}, as desired. 

\newpage

\subsubsection{With singularities: \boldmath{$1/2 < \alpha < 1$}}
Here,  $(u^*_\alpha)'(x)\to-\infty$ as $x\to 0$, which suggests that a large portion of the integration nodes should concentrate on a neighborhood of zero to appropriately capture the singular trend of the derivative at that point.  To this end,  we divide the batch as the union of two equal-sized different $\beta(a,b)$ samples:  one with $a=1=b$, and the other with $a=10^4$ and $b=1$. We individually experiment for $\alpha\in\{0.6, 0.7, 0.8\}$ with $100\mathord{,}000$ iterations for approximating $u_\text{NN}$, which implies $500\mathord{,}000$ for the \ac{D2RM} when taking into account the iterations dedicated to approximate the optimal test functions. 

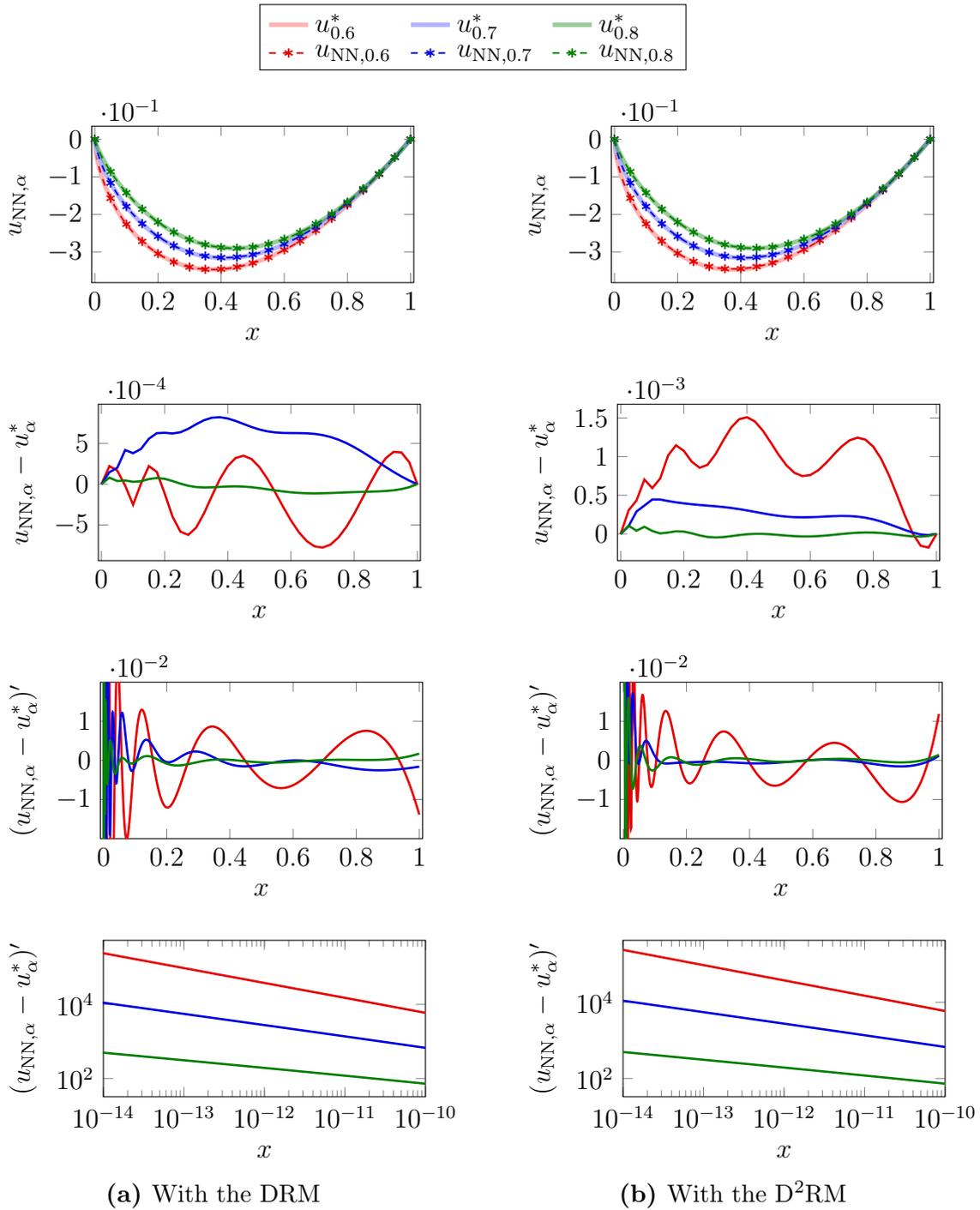
\begin{figure}[htbp]
\centering
\begin{subfigure}[t]{\textwidth}
\centering
 \begin{tikzpicture}
\pgfplotsset{A/.append style={
		hide axis,
	    xmin=-0.1,
	    xmax=1.1,
	    ymin=-0.1,
	    ymax=1.1,
	 	legend style={anchor=north, /tikz/every even column/.append style={column sep=0.2cm}},
	 	legend columns = 3,
	 	legend cell align={left},
            }
}
\begin{axis}[A]
\addlegendimage{line width=2, color=red!30!white}
\addlegendentry{$u_{0.6}^*$}
\addlegendimage{line width=2, blue!30!white}
\addlegendentry{$u_{0.7}^*$};
\addlegendimage{line width=2, green!50!black, opacity=0.4}
\addlegendentry{$u_{0.8}^*$};

\addlegendimage{dashed, color=red!90!black,  line width=.7, mark=asterisk, mark options = {scale=1, solid}}
\addlegendentry{$u_{\text{NN},0.6}$}
\addlegendimage{dashed, color=blue!90!black,  line width=.7, mark=asterisk, mark options = {scale=1, solid}}
\addlegendentry{$u_{\text{NN},0.7}$};
\addlegendimage{dashed, color=green!50!black,  line width=.7, mark=asterisk, mark options = {scale=1, solid}}
\addlegendentry{$u_{\text{NN},0.8}$};
\end{axis}
\end{tikzpicture}
\end{subfigure}\vskip 1em%
\begin{subfigure}[t]{0.45\textwidth}
\centering
\begin{tikzpicture}         
 \begin{axis}[xmin=-0.01,         
     xmax=1.01,      
     xlabel = {$x$},  
     ylabel = {$u_{\text{NN},\alpha}$},
     height=4cm,    
     width=\textwidth,
     scaled y ticks=base 10:1,    
     yticklabel style={
        /pgf/number format/fixed,
        /pgf/number format/precision=2}, ]        

\addplot [domain = 0:1, samples=1000, line width=2, red!30!white]{x^(0.6)*(x-1)};
\addplot[only marks, color=red!90!black,  line width=.7, mark=asterisk, mark options = {scale=1, solid}, each nth point=2] table[x expr=\thisrow{x},y=u_net]{Chapters/4.Chapter/FIGURES/numerical_results/section53/DeepRitz_x0.6(x-1)_100000_a20/u_net_plot.csv};
\addplot[dashed, color=red!90!black,  line width=.7, each nth point=3] table[x expr=\thisrow{x},y=u_net]{Chapters/4.Chapter/FIGURES/numerical_results/section53/DeepRitz_x0.6(x-1)_100000_a20/u_net.csv};
\addplot [domain = 0:1, samples=1000, line width=2, blue!30!white]{x^(0.7)*(x-1)};
\addplot[only marks, color=blue!90!black,  line width=.7, mark=asterisk, mark options = {scale=1, solid}, each nth point=2] table[x expr=\thisrow{x},y=u_net]{Chapters/4.Chapter/FIGURES/numerical_results/section53/DeepRitz_x0.7(x-1)_100000_a20/u_net_plot.csv};
\addplot[dashed, color=blue!90!black,  line width=.7, each nth point=3] table[x expr=\thisrow{x},y=u_net]{Chapters/4.Chapter/FIGURES/numerical_results/section53/DeepRitz_x0.7(x-1)_100000_a20/u_net.csv};
\addplot [domain = 0:1, samples=1000, line width=2,green!50!black, opacity=0.4]{x^(0.8)*(x-1)};
\addplot[only marks, color=green!50!black,  line width=.7, mark=asterisk, mark options = {scale=1, solid}, each nth point=2] table[x expr=\thisrow{x},y=u_net]{Chapters/4.Chapter/FIGURES/numerical_results/section53/DeepRitz_x0.8(x-1)_100000_a20/u_net_plot.csv};
\addplot[dashed, color=green!50!black,  line width=.7, each nth point=3] table[x expr=\thisrow{x},y=u_net]{Chapters/4.Chapter/FIGURES/numerical_results/section53/DeepRitz_x0.8(x-1)_100000_a20/u_net.csv};
\end{axis}
\end{tikzpicture}
\end{subfigure}\hfill%
\begin{subfigure}[t]{0.45\textwidth}
\centering
\begin{tikzpicture}         
 \begin{axis}[xmin=-0.01,         
     xmax=1.01,      
     xlabel = {$x$}, 
     ylabel = {$u_{\text{NN},\alpha}$},
     height=4cm,    
     width=\textwidth,
     scaled y ticks=base 10:1,    
     yticklabel style={
        /pgf/number format/fixed,
        /pgf/number format/precision=2}, ]        
\addplot [domain = 0:1, samples=1000, line width=2, red!30!white]{x^(0.6)*(x-1)};
\addplot[only marks, color=red!90!black,  line width=.7, mark=asterisk, mark options = {scale=1, solid}, each nth point=2] table[x expr=\thisrow{x},y=u_net]{Chapters/4.Chapter/FIGURES/numerical_results/section53/DeepDRitz_x0.6(x-1)_500000_a20/u_net_plot.csv};
\addplot[dashed, color=red!90!black,  line width=.7, each nth point=3] table[x expr=\thisrow{x},y=u_net]{Chapters/4.Chapter/FIGURES/numerical_results/section53/DeepDRitz_x0.6(x-1)_500000_a20/u_net.csv};
\addplot [domain = 0:1, samples=1000, line width=2, blue!30!white]{x^(0.7)*(x-1)};
\addplot[only marks, color=blue!90!black,  line width=.7, mark=asterisk, mark options = {scale=1, solid}, each nth point=2] table[x expr=\thisrow{x},y=u_net]{Chapters/4.Chapter/FIGURES/numerical_results/section53/DeepDRitz_x0.7(x-1)_500000_a20/u_net_plot.csv};
\addplot[dashed, color=blue!90!black,  line width=.7,each nth point=3] table[x expr=\thisrow{x},y=u_net]{Chapters/4.Chapter/FIGURES/numerical_results/section53/DeepDRitz_x0.7(x-1)_500000_a20/u_net.csv};
\addplot [domain = 0:1, samples=1000, line width=2,green!50!black, opacity=0.4]{x^(0.8)*(x-1)};
\addplot[only marks, color=green!50!black,  line width=.7, mark=asterisk, mark options = {scale=1, solid}, each nth point=2] table[x expr=\thisrow{x},y=u_net]{Chapters/4.Chapter/FIGURES/numerical_results/section53/DeepDRitz_x0.8(x-1)_500000_a20/u_net_plot.csv};
\addplot[dashed, color=green!50!black,  line width=.7, each nth point=3] table[x expr=\thisrow{x},y=u_net]{Chapters/4.Chapter/FIGURES/numerical_results/section53/DeepDRitz_x0.8(x-1)_500000_a20/u_net.csv};
\end{axis}
\end{tikzpicture}
\end{subfigure}\vskip 1em%
\begin{subfigure}[t]{0.45\textwidth}
\centering
\begin{tikzpicture}         
 \begin{axis}[xmin=-0.01,
     xmax=1.01,      
     xlabel = {$x$},  
     ylabel = {$u_{\text{NN}, \alpha}-u_\alpha^*$},
     height=4cm,    
     width=\textwidth,    
     yticklabel style={
        /pgf/number format/fixed,
        /pgf/number format/precision=1}, ]        
\addplot[color=red!90!black,  line width=1] table[x expr=\thisrow{x},y=e_u]{Chapters/4.Chapter/FIGURES/numerical_results/section53/DeepRitz_x0.6(x-1)_100000_a20/e_u_plot.csv};
\addplot[color=blue!90!black,  line width=1] table[x expr=\thisrow{x},y=e_u]{Chapters/4.Chapter/FIGURES/numerical_results/section53/DeepRitz_x0.7(x-1)_100000_a20/e_u_plot.csv};
\addplot[color=green!50!black,  line width=1,] table[x expr=\thisrow{x},y=e_u]{Chapters/4.Chapter/FIGURES/numerical_results/section53/DeepRitz_x0.8(x-1)_100000_a20/e_u_plot.csv};
\end{axis}
\end{tikzpicture}
\end{subfigure}\hfill%
\begin{subfigure}[t]{0.45\textwidth}
\centering
\begin{tikzpicture}         
 \begin{axis}[xmin=-0.01,         
     xmax=1.01,      
     xlabel = {$x$},  
     ylabel = {$u_{\text{NN}, \alpha}-u_\alpha^*$},
     height=4cm,    
     width=\textwidth,    
     yticklabel style={
        /pgf/number format/fixed,
        /pgf/number format/precision=1}, ]        
\addplot[color=red!90!black,  line width=1] table[x expr=\thisrow{x},y=e_u]{Chapters/4.Chapter/FIGURES/numerical_results/section53/DeepDRitz_x0.6(x-1)_500000_a20/e_u_plot.csv};
\addplot[color=blue!90!black,  line width=1] table[x expr=\thisrow{x},y=e_u]{Chapters/4.Chapter/FIGURES/numerical_results/section53/DeepDRitz_x0.7(x-1)_500000_a20/e_u_plot.csv};
\addplot[color=green!50!black,  line width=1] table[x expr=\thisrow{x},y=e_u]{Chapters/4.Chapter/FIGURES/numerical_results/section53/DeepDRitz_x0.8(x-1)_500000_a20/e_u_plot.csv};
\end{axis}
\end{tikzpicture}
\end{subfigure}\vskip 1 em%
\begin{subfigure}[t]{0.45\textwidth}
\centering
\begin{tikzpicture}         
 \begin{axis}[
     xmin=-0.01,         
     xmax=1.01,      
     xlabel = {$x$}, 
     ymin=-2e-2,     
     ymax=2e-2,      
     ylabel = {$(u_{\text{NN}, \alpha}-u_\alpha^*)'$},
     height=4cm,    
     width=\textwidth,    
     ytick={-0.01,0,...,0.01},
     yticklabel style={
        /pgf/number format/fixed,
        /pgf/number format/precision=2}, 
    select coords between index/.style 2 args={
    x filter/.code={
        \ifnum\coordindex<#1\def\pgfmathresult{}\fi
        \ifnum\coordindex>#2\def\pgfmathresult{}\fi
    }
},]        
\addplot[color=red!90!black,  line width=1, select coords between index={4000}{10000}, each nth point=2] table[x expr=\thisrow{x},y=de_u]{Chapters/4.Chapter/FIGURES/numerical_results/section53/DeepRitz_x0.6(x-1)_100000_a20/de_u.csv};
\addplot[color=blue!90!black,  line width=1, select coords between index={4000}{10000}, each nth point=2] table[x expr=\thisrow{x},y=de_u]{Chapters/4.Chapter/FIGURES/numerical_results/section53/DeepRitz_x0.7(x-1)_100000_a20/de_u.csv};
\addplot[color=green!50!black,  line width=1, select coords between index={4000}{10000}, each nth point=2] table[x expr=\thisrow{x},y=de_u]{Chapters/4.Chapter/FIGURES/numerical_results/section53/DeepRitz_x0.8(x-1)_100000_a20/de_u.csv};
\end{axis}
\end{tikzpicture}
\end{subfigure}\hfill%
\begin{subfigure}[t]{0.45\textwidth}
\centering
\begin{tikzpicture}         
 \begin{axis}[
     xmin=-0.01,         
     xmax=1.01,      
     xlabel = {$x$}, 
     ymin=-2e-2,     
     ymax=2e-2,      
     ylabel = {$(u_{\text{NN}, \alpha}-u_\alpha^*)'$},
     ytick={-0.01,0,...,0.01},
     height=4cm,    
     width=\textwidth,    
     yticklabel style={
        /pgf/number format/fixed,
        /pgf/number format/precision=2}, 
    select coords between index/.style 2 args={
    x filter/.code={
        \ifnum\coordindex<#1\def\pgfmathresult{}\fi
        \ifnum\coordindex>#2\def\pgfmathresult{}\fi
    }
},]        
\addplot[color=red!90!black,  line width=1, select coords between index={4000}{10000}, each nth point=2] table[x expr=\thisrow{x},y=de_u]{Chapters/4.Chapter/FIGURES/numerical_results/section53/DeepDRitz_x0.6(x-1)_500000_a20/de_u.csv};
\addplot[color=blue!90!black,  line width=1, select coords between index={4000}{10000}, each nth point=2] table[x expr=\thisrow{x},y=de_u]{Chapters/4.Chapter/FIGURES/numerical_results/section53/DeepDRitz_x0.7(x-1)_500000_a20/de_u.csv};
\addplot[color=green!50!black,  line width=1, select coords between index={4000}{10000}, each nth point=2] table[x expr=\thisrow{x},y=de_u]{Chapters/4.Chapter/FIGURES/numerical_results/section53/DeepDRitz_x0.8(x-1)_500000_a20/de_u.csv};
\end{axis}
\end{tikzpicture}
\end{subfigure}\vskip 1 em%
\begin{subfigure}[t]{0.45\textwidth}
\centering
\begin{tikzpicture}         
 \begin{axis}[
     xmode=log,
     ymode=log,
     xmin=1e-14,         
     xmax=1e-10,      
     xlabel = {$x$}, 
     ylabel = {$(u_{\text{NN}, \alpha}-u_\alpha^*)'$},
     height=4cm,    
     width=\textwidth, 
     yticklabel style={
        /pgf/number format/fixed,
        /pgf/number format/precision=2}, ]        
\addplot[color=red!90!black,  line width=1, each nth point=2] table[x expr=\thisrow{x},y=de_u]{Chapters/4.Chapter/FIGURES/numerical_results/section53/DeepRitz_x0.6(x-1)_100000_a20/de_u.csv};
\addplot[color=blue!90!black,  line width=1, each nth point=2] table[x expr=\thisrow{x},y=de_u]{Chapters/4.Chapter/FIGURES/numerical_results/section53/DeepRitz_x0.7(x-1)_100000_a20/de_u.csv};
\addplot[color=green!50!black,  line width=1, each nth point=2] table[x expr=\thisrow{x},y=de_u]{Chapters/4.Chapter/FIGURES/numerical_results/section53/DeepRitz_x0.8(x-1)_100000_a20/de_u.csv};
\end{axis}
\end{tikzpicture}\caption{With the \acs{DRM}}
\end{subfigure}\hfill%
\begin{subfigure}[t]{0.45\textwidth}
\centering
\begin{tikzpicture}         
 \begin{axis}[
     xmode=log,
     ymode=log,
     xmin=1e-14,         
     xmax=1e-10,      
     xlabel = {$x$},  
     ylabel = {$(u_{\text{NN}, \alpha}-u_\alpha^*)'$},
     height=4cm,    
     width=\textwidth,
     yticklabel style={
        /pgf/number format/fixed,
        /pgf/number format/precision=2}, ]        
\addplot[color=red!90!black,  line width=1, each nth point=2] table[x expr=\thisrow{x},y=de_u]{Chapters/4.Chapter/FIGURES/numerical_results/section53/DeepDRitz_x0.6(x-1)_500000_a20/de_u.csv};
\addplot[color=blue!90!black,  line width=1, each nth point=2] table[x expr=\thisrow{x},y=de_u]{Chapters/4.Chapter/FIGURES/numerical_results/section53/DeepDRitz_x0.7(x-1)_500000_a20/de_u.csv};
\addplot[color=green!50!black,  line width=1, each nth point=2] table[x expr=\thisrow{x},y=de_u]{Chapters/4.Chapter/FIGURES/numerical_results/section53/DeepDRitz_x0.8(x-1)_500000_a20/de_u.csv};
\end{axis}
\end{tikzpicture}
\caption{With the \acs{D2RM}}
\end{subfigure}
\caption{Trial network predictions and errors for the \acs{DRM} and the \acs{D2RM} in model problem \eqref{1DLap} with exact solution $u_\alpha^*=x^{\alpha}(x-1)$ for $\alpha\in\{0.6,0.7,0.8\}$.  The last row is a zoomed in version of the third one in a reduced neighborhood of zero.}
\label{figure:xalpha(x-1)_singular_predictions}
\end{figure}

\Cref{table:experiments3} displays the record of the relative error estimates along different stages of the training, and \Cref{figure:xalpha(x-1)_singular_predictions} shows the $u_\text{NN}$ predictions and error functions at the end of the training.

\begin{table}[htbp]
\centering
\begin{tabular}{|c|c||c|c|c|c|c|}
\toprule
\multicolumn{2}{|c||}{\textbf{Training progress}} & $\mathbf{4}\bf{\%}$ & $\mathbf{20}\bf{\%}$ & $\mathbf{40}\bf{\%}$ & $\mathbf{60}\bf{\%}$ & $\mathbf{100}\bf{\%}$\\ \bottomrule
\toprule
\textbf{Method} & $\boldsymbol{\alpha}$ & \multicolumn{5}{c|}{$\frac{\Vert u_\text{NN}-u^*\Vert_\mathbb{U}}{\Vert u^*\Vert_\mathbb{U}}\times 100$} \\ \midrule
\multicolumn{1}{|l||}{\multirow{3}{*}{\ac{DRM}}}  & $0.6$ & $42.23\%$ & $30.29\%$ & $26.91\%$ & $25.34\%$ & $23.84\%$\\ \cline{2-7} 
\multicolumn{1}{|l||}{}                       & $0.7$ &  $18.47\%$ &$9.49\%$ &$7.42\%$ &$6.65\%$ &$5.95\%$\\ \cline{2-7} 
\multicolumn{1}{|l||}{}                       & $0.8$ &  $8.78\%$ &$3.50\%$ &$2.46\%$ &$2.10\%$ &$1.81\%$\\ \midrule
\multicolumn{1}{|l||}{\multirow{3}{*}{\ac{D2RM}}} & $0.6$ &  $41.67\%$ &$30.25\%$ &$27.05\%$ &$25.76\%$ &$24.40\%$\\ \cline{2-7} 
\multicolumn{1}{|l||}{}                       & $0.7$ & $14.12\%$ &$9.99\%$ &$7.91\%$ &$6.98\%$ &$6.15\%$ \\ \cline{2-7} 
\multicolumn{1}{|l||}{}                       & $0.8$ & $6.00\%$ & $3.62\%$ & $2.78\%$ & $2.55\%$ & $2.13\%$\\ \bottomrule
\toprule
\textbf{Method} & $\boldsymbol{\alpha}$ & \multicolumn{5}{c|}{$\frac{\Vert \tau_\text{NN}(u_\text{NN})-Tu^*\Vert_\V}{\Vert Tu^*\Vert_\V}\times 100$} \\ \midrule
\multicolumn{1}{|l||}{\multirow{3}{*}{\ac{D2RM}}} & $0.6$ &  $41.68\%$ &$30.25\%$ &$27.05\%$ &$25.77\%$ &$24.40\%$ \\ \cline{2-7} 
\multicolumn{1}{|l||}{}                       & $0.7$ &  $14.12\%$ &$10.00\%$ &$7.91\%$ &$6.98\%$ &$6.15\%$ \\ \cline{2-7} 
\multicolumn{1}{|l||}{}                       & $0.8$ &  $6.01\%$ &$3.62\%$ &$2.78\%$ &$2.54\%$ &$2.13\%$\\ \bottomrule
\end{tabular}%
\caption{Relative errors of $u_\text{NN}$ (in the \acs{DRM} and the \acs{D2RM}) and $\tau_\text{NN}(u_\text{NN})$ (in the \acs{D2RM}) along different stages of the training progress in problem \eqref{1DLap} with exact solution $u_\alpha^*=x^{\alpha}(x-1)$ and $\alpha\in\{0.6, 0.7, 0.8\}$. }
\label{table:experiments3}
\end{table}

The \ac{DRM} and the \ac{D2RM} decrease the relative errors at similar rates, slowing down notably from the $40\%$ of the training progress onwards.  For $\alpha \in\{ 0.7,  0.8\}$, we achieve final relative errors below $7.5\%$ at the end of the training, with loss predictions around $-0.3632$ and $-0.2562$ where the optimal values are around $-0.3646$ and $-0.2564$, respectively.  For $\alpha = 0.6$, the relative error is above $20\%$,  with a final loss prediction around $-0.6429$ where the optimal value is around $-0.6818$. Intending to concentrate more points around the singularity, we performed several experiments modifying the $a$ parameter of the beta distribution, but without success.

\subsection{Pure diffusion with a Dirac delta source}\label{section:Diffusion equation with a Dirac delta source}
We now select the source $l=4 \delta_{1/2}\in H^{-1}(0,1)\setminus (L^2(0,1))'$, where $\delta_{1/2}$ denotes the Dirac delta distribution at $1/2$, i.e., $\delta_{1/2}(v)=v(1/2)$ for all $v\in H^1_0(0,1)$. The corresponding exact solution is
 \begin{equation}\label{equation:corner_solution}
 u^*=\begin{cases}2x, &\text{if }x<1/2,\\
 2(1-x), &\text{if } x>1/2.
 \end{cases}
 \end{equation} 

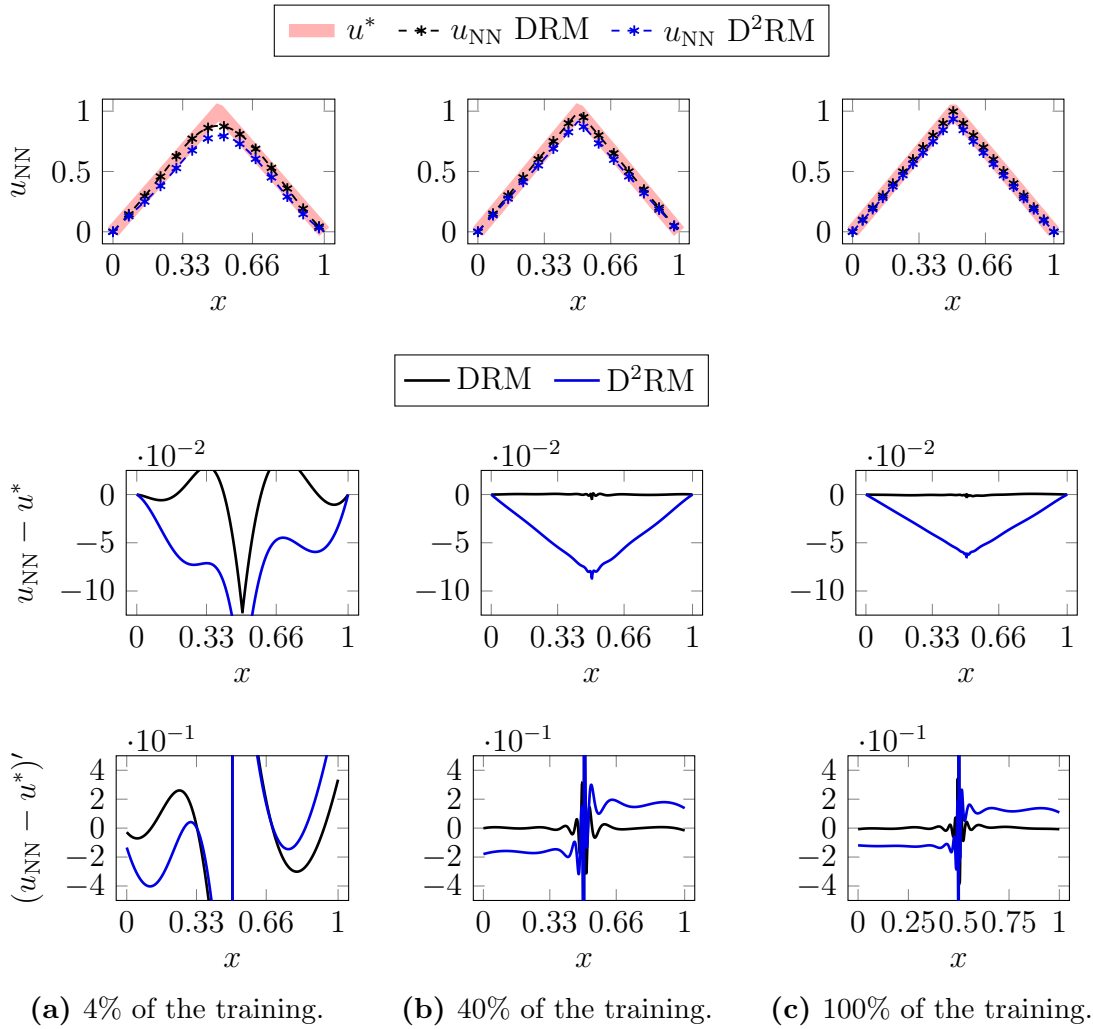
\begin{figure}[htbp]
\centering
\begin{subfigure}[t]{\textwidth}
\centering
 \begin{tikzpicture}
\pgfplotsset{A/.append style={
		hide axis,
	    xmin=-0.1,
	    xmax=1.1,
	    ymin=-0.1,
	    ymax=1.1,
	 	legend style={anchor=north, /tikz/every even column/.append style={column sep=0.2cm}},
	 	legend columns = 3
            }
}
\begin{axis}[A]
\addlegendimage{color=red!30!white, line width=5}
\addlegendentry{$u^*$}
\addlegendimage{style=dashed, color=black, line width=.7, mark=asterisk, mark options = {scale=1, solid},}
\addlegendentry{$u_\text{NN}$ \ac{DRM}}
\addlegendimage{color=blue!90!black, style=dashed,line width=.7, mark=asterisk, mark options = {scale=1, solid},}
\addlegendentry{$u_\text{NN}$ \ac{D2RM}}
\end{axis}
\end{tikzpicture}
\end{subfigure}\vskip 1em
\begin{subfigure}[t]{0.32\textwidth}
\centering
\begin{tikzpicture}
\begin{axis}[xmin=-0.05,
	    xmax=1.05,
	    xlabel = {$x$},
	    ylabel = {$u_\text{NN}$},
	    height=3.5cm,
	    width=\textwidth,
	    xtick={0,0.33,0.66,1},
	    yticklabel style={
        /pgf/number format/fixed,
        /pgf/number format/precision=2},
        scaled y ticks=false,
	    legend style={at={(0.5,0.95)}, anchor=north, /tikz/every even column/.append style={column sep=0.5cm}},
	 	]  
	\addplot[color=red!30!white, line width=5,each nth point=3] table[x expr=\thisrow{x},y=u_analytic]{Chapters/4.Chapter/FIGURES/numerical_results/section54/DeepRitz_corner_800/u_analytic.csv};
	\addplot[color=black!90!black,  style=dashed,line width=.7,each nth point=3] table[x expr=\thisrow{x},y=u_net]{Chapters/4.Chapter/FIGURES/numerical_results/section54/DeepRitz_corner_800/u_net.csv};
	\addplot[only marks, color=black, line width=.7, mark=asterisk, mark options = {scale=1, solid},each nth point=3] table[x expr=\thisrow{x},y=u_net]{Chapters/4.Chapter/FIGURES/numerical_results/section54/DeepRitz_corner_800/u_net_plot.csv};
    \addplot[color=blue!90!black, style=dashed,line width=.7,each nth point=3] table[x expr=\thisrow{x},y=u_net]{Chapters/4.Chapter/FIGURES/numerical_results/section54/DeepDRitz_corner_4000/u_net.csv};
	\addplot[only marks, color=blue!90!black,line width=.7, mark=asterisk, mark options = {scale=1, solid}, each nth point=3] table[x expr=\thisrow{x},y=u_net]{Chapters/4.Chapter/FIGURES/numerical_results/section54/DeepDRitz_corner_4000/u_net_plot.csv};
\end{axis}
\end{tikzpicture}
\label{figure:DeepRitz_corner_prediction}
\end{subfigure}\hfill
\begin{subfigure}[t]{0.31\textwidth}
\centering
\begin{tikzpicture}
\begin{axis}[xmin=-0.05,
	    xmax=1.05,
	    xlabel = {$x$},
	    height=3.5cm,
	    width=\textwidth,
	    xtick={0,0.33,.66,1},
	    yticklabel style={
        /pgf/number format/fixed,
        /pgf/number format/precision=2},
        scaled y ticks=false,
	    legend style={at={(0.5,0.95)}, anchor=north, /tikz/every even column/.append style={column sep=0.5cm}},
	 	]  
	\addplot[color=red!30!white, line width=5,each nth point=3] table[x expr=\thisrow{x},y=u_analytic]{Chapters/4.Chapter/FIGURES/numerical_results/section54/DeepRitz_corner_8000/u_analytic.csv};
	\addplot[color=black!90!black,  style=dashed,line width=.7, each nth point=3] table[x expr=\thisrow{x},y=u_net]{Chapters/4.Chapter/FIGURES/numerical_results/section54/DeepRitz_corner_8000/u_net.csv};
	\addplot[only marks, color=black, line width=.7, mark=asterisk, mark options = {scale=1, solid}, each nth point=3] table[x expr=\thisrow{x},y=u_net]{Chapters/4.Chapter/FIGURES/numerical_results/section54/DeepRitz_corner_8000/u_net_plot.csv};
    \addplot[color=blue!90!black, style=dashed,line width=.7,each nth point=3] table[x expr=\thisrow{x},y=u_net]{Chapters/4.Chapter/FIGURES/numerical_results/section54/DeepDRitz_corner_40000/u_net.csv};
	\addplot[only marks, color=blue!90!black,line width=.7, mark=asterisk, mark options = {scale=1, solid}, each nth point=3] table[x expr=\thisrow{x},y=u_net]{Chapters/4.Chapter/FIGURES/numerical_results/section54/DeepDRitz_corner_40000/u_net_plot.csv};
\end{axis}
\end{tikzpicture}
\label{figure:DeepDRitz_corner_prediction}
\end{subfigure}\hfill
\begin{subfigure}[t]{0.31\textwidth}
\centering
\begin{tikzpicture}
\begin{axis}[xmin=-0.05,
	    xmax=1.05,
	    xlabel = {$x$},
	    height=3.5cm,
	    width=\textwidth,
	    xtick={0,0.33,0.66,1},
	    yticklabel style={
        /pgf/number format/fixed,
        /pgf/number format/precision=2},
        scaled y ticks=false,
	 	]  
	\addplot[color=red!30!white, line width=5,each nth point=2] table[x expr=\thisrow{x},y=u_analytic]{Chapters/4.Chapter/FIGURES/numerical_results/section54/DeepRitz_corner_20000/u_analytic.csv};
	\addplot[color=black!90!black,  style=dashed,line width=.7,each nth point=2] table[x expr=\thisrow{x},y=u_net, each nth point=2]{Chapters/4.Chapter/FIGURES/numerical_results/section54/DeepRitz_corner_20000/u_net.csv};
	\addplot[only marks, color=black, line width=.7,mark=asterisk, mark options = {scale=1, solid},each nth point=2] table[x expr=\thisrow{x},y=u_net]{Chapters/4.Chapter/FIGURES/numerical_results/section54/DeepRitz_corner_20000/u_net_plot.csv};
    \addplot[color=blue!90!black, style=dashed,line width=.7,each nth point=2] table[x expr=\thisrow{x},y=u_net]{Chapters/4.Chapter/FIGURES/numerical_results/section54/DeepDRitz_corner_100000/u_net.csv};
	\addplot[only marks, color=blue!90!black,line width=.7, mark=asterisk, mark options = {scale=1, solid},each nth point=2] table[x expr=\thisrow{x},y=u_net]{Chapters/4.Chapter/FIGURES/numerical_results/section54/DeepDRitz_corner_100000/u_net_plot.csv};
\end{axis}
\end{tikzpicture}
\label{figure:DeepDRitz_corner_prediction}
\end{subfigure}\vskip1em%
\begin{subfigure}[t]{\textwidth}
\centering
 \begin{tikzpicture}
\pgfplotsset{A/.append style={hide axis,
	    xmin=-0.1,
	    xmax=1.1,
	    ymin=-0.1,
	    ymax=1.1,
	    legend style={anchor=north, /tikz/every even column/.append style={column sep=0.2cm}},
	    legend columns = 2
            }
}
\begin{axis}[A]
\addlegendimage{color=black, line width=1}
\addlegendentry{\ac{DRM}};
\addlegendimage{color=blue!90!black, line width=1}
\addlegendentry{\ac{D2RM}};
\end{axis}
\end{tikzpicture}
\end{subfigure}\vskip1em
\begin{subfigure}[t]{0.32\textwidth}
\centering
\begin{tikzpicture}
\begin{axis}[
		xmin=-0.05,
	    xmax=1.05,
	    xlabel = {$x$},
	    ymin=-0.125,
	    ymax=0.025,
	    ylabel = {$u_\text{NN}-u^*$},
	    height=3.5cm,
	    width=\textwidth,
	    xtick={0,0.33,0.66,1},
	    scaled y ticks=base 10:2,   
	 	]
	\addplot[color=black!90!black, line width=1] table[x expr=\thisrow{x},y=e_u, each nth point=2]{Chapters/4.Chapter/FIGURES/numerical_results/section54/DeepRitz_corner_800/e_u.csv};
	\addplot[color=blue!90!black, line width=1] table[x expr=\thisrow{x},y=e_u, each nth point=2]{Chapters/4.Chapter/FIGURES/numerical_results/section54/DeepDRitz_corner_4000/e_u.csv};
\end{axis}
\end{tikzpicture}
\end{subfigure}\hfill%
\begin{subfigure}[t]{0.31\textwidth}
\centering
\begin{tikzpicture}
\begin{axis}[xmin=-0.05,
	    xmax=1.05,
	    xlabel = {$x$},
	    ymin=-0.125,
	    ymax=0.025,
	    height=3.5cm,
	    width=\textwidth,
	    xtick={0,0.33,0.66,1},
	    scaled y ticks=base 10:2,   
	 	]
	\addplot[color=black!90!black, line width=1] table[x expr=\thisrow{x},y=e_u, each nth point=2]{Chapters/4.Chapter/FIGURES/numerical_results/section54/DeepRitz_corner_8000/e_u.csv};
	\addplot[color=blue!90!black, line width=1] table[x expr=\thisrow{x},y=e_u, each nth point=2]{Chapters/4.Chapter/FIGURES/numerical_results/section54/DeepDRitz_corner_40000/e_u.csv};
\end{axis}
\end{tikzpicture}
\end{subfigure}\hfill%
\begin{subfigure}[t]{0.31\textwidth}
\centering
\begin{tikzpicture}
\begin{axis}[xmin=-0.05,
	    xmax=1.05,
	    xlabel = {$x$},
	    ymin=-0.125,
	    ymax=0.025,
	    height=3.5cm,
	    width=\textwidth,
	    xtick={0,0.33,0.66,1},
	    ytick={-0.1,-0.05,0},
	    scaled y ticks=base 10:2,   
	 	]
	\addplot[color=black!90!black, line width=1] table[x expr=\thisrow{x},y=e_u, each nth point=2]{Chapters/4.Chapter/FIGURES/numerical_results/section54/DeepRitz_corner_20000/e_u.csv};
	\addplot[color=blue!90!black, line width=1] table[x expr=\thisrow{x},y=e_u, each nth point=2]{Chapters/4.Chapter/FIGURES/numerical_results/section54/DeepDRitz_corner_100000/e_u.csv};
\end{axis}
\end{tikzpicture}
\end{subfigure}\vskip 1em%
\begin{subfigure}[t]{0.32\textwidth}
\centering
\begin{tikzpicture}
\begin{axis}[xmin=-0.05,
	    xmax=1.05,
	    xlabel = {$x$},
	    ymin=-0.5,
	    ymax=0.5,
	    ylabel = {$(u_\text{NN}-u^*)'$},
	    height=3.5cm,
	    width=\textwidth,
	    xtick={0,0.33,0.66,1},
	    scaled y ticks=base 10:1,   
	 	]
	\addplot[color=black!90!black, line width=1] table[x expr=\thisrow{x},y=de_u, each nth point=2]{Chapters/4.Chapter/FIGURES/numerical_results/section54/DeepRitz_corner_800/de_u.csv};
	\addplot[color=blue!90!black, line width=1] table[x expr=\thisrow{x},y=de_u, each nth point=2]{Chapters/4.Chapter/FIGURES/numerical_results/section54/DeepDRitz_corner_4000/de_u.csv};
\end{axis}
\end{tikzpicture}
\caption{$4\%$ of the training.}
\label{figure:corner_derror1}
\end{subfigure}\hfill%
\begin{subfigure}[t]{0.31\textwidth}
\centering
\begin{tikzpicture}
\begin{axis}[xmin=-0.05,
	    xmax=1.05,
	    xlabel = {$x$},
	    ymin=-0.5,
	    ymax=0.5,
	    height=3.5cm,
	    width=\textwidth,
	    xtick={0,0.33,0.66,1},
	    scaled y ticks=base 10:1,   
	 	]
	\addplot[color=black!90!black, line width=1] table[x expr=\thisrow{x},y=de_u, each nth point=2]{Chapters/4.Chapter/FIGURES/numerical_results/section54/DeepRitz_corner_8000/de_u.csv};
	\addplot[color=blue!90!black, line width=1] table[x expr=\thisrow{x},y=de_u, each nth point=2]{Chapters/4.Chapter/FIGURES/numerical_results/section54/DeepDRitz_corner_40000/de_u.csv};
\end{axis}
\end{tikzpicture}
\caption{$40\%$ of the training.}
\label{figure:corner_derror2}
\end{subfigure}\hfill%
\begin{subfigure}[t]{0.31\textwidth}
\centering
\begin{tikzpicture}
\begin{axis}[xmin=-0.05,
	    xmax=1.05,
	    xlabel = {$x$},
	    ymin=-0.5,
	    ymax=0.5,
	    height=3.5cm,
	    width=\textwidth,
	    xtick={0,0.25,...,1},
	    scaled y ticks=base 10:1,   
	 	]
	\addplot[color=black!90!black, line width=1] table[x expr=\thisrow{x},y=de_u, each nth point=2]{Chapters/4.Chapter/FIGURES/numerical_results/section54/DeepRitz_corner_20000/de_u.csv};
	\addplot[color=blue!90!black, line width=1] table[x expr=\thisrow{x},y=de_u, each nth point=2]{Chapters/4.Chapter/FIGURES/numerical_results/section54/DeepDRitz_corner_100000/de_u.csv};
\end{axis}
\end{tikzpicture}
\caption{$100\%$ of the training.}
\label{figure:corner_derror3}
\end{subfigure}
\caption{Trial error functions in model problem \eqref{1DLap} with with solution \eqref{equation:corner_solution} at different stages of the training progress for the \acs{DRM} and the \acs{D2RM}.}
\label{figure:corner_error_u}
\end{figure}

We train $u_\text{NN}$ for $20\mathord{,}000$ iterations, selecting $a=1=b$ and $a=10=b$ in the $\beta(a,b)$ probability distribution for the two-sampled batch. \Cref{table:experiments4} records the relative errors at three different stages of the training progress. \Cref{figure:corner_error_u} shows the trial network predictions, errors, and derivative of the errors at those checkpoints. As we observe, the \ac{D2RM} shows higher errors than the \ac{DRM}.  

\begin{table}[htbp]
\centering
\begin{tabular}{|c||c|c|c|c|c|}
\toprule
\textbf{Training progress} & $\mathbf{4}\bf{\%}$ & $\mathbf{20}\bf{\%}$ & $\mathbf{40}\bf{\%}$ & $\mathbf{60}\bf{\%}$ & $\mathbf{100}\bf{\%}$\\ \bottomrule
\toprule
\textbf{Method} & \multicolumn{5}{c|}{$\frac{\Vert u_\text{NN}-u^*\Vert_\mathbb{U}}{\Vert u^*\Vert_\mathbb{U}}\times 100$}\\ \midrule
\ac{DRM} & $32.32\%$ & $10.10\%$ & $6.22\%$ & $5.08\%$ & $4.34\%$\\ \midrule
\ac{D2RM} & $32.43\%$&$13.95\%$&$10.87\%$&$9.41\%$&$7.95\%$\\ \bottomrule
\toprule
\textbf{Method} & \multicolumn{5}{c|}{$\frac{\Vert \tau_\text{NN}(u_\text{NN})-Tu^*\Vert_\V}{\Vert Tu^*\Vert_\V}\times 100$}\\ \midrule
\ac{D2RM} & $32.68\%$&$14.09\%$&$10.84\%$&$9.53\%$&$8.16\%$\\ \bottomrule
\end{tabular}%
\caption{Relative errors of $u_\text{NN}$ (in the \acs{DRM} and the \acs{D2RM}) and $\tau_\text{NN}(u_\text{NN})$ (in the \acs{D2RM}) along different stages of the training progress in model problem \eqref{1DLap} with exact solution \eqref{equation:corner_solution}.}
\label{table:experiments4}
\end{table}

\subsection{Pure convection with Dirac delta source}\label{section:Advection equation in ultraweak formulation}
We consider the spatial convection equation
\begin{equation}\label{equation:advection}
\begin{cases}
u'=\delta_{1/2}, &\text{in } (0,1),\\
u(0)=0.
\end{cases}
\end{equation} Here, the (ultra)weak formulation is appropriate because $\delta_{1/2}$ does not belong to $(L^2(0,1))'$.  Integration by parts yields:
\begin{equation}
b(u,v):=-\int_0^1 uv',\qquad l(v):=v(1/2),\qquad u\in\mathbb{U}, v\in\V,
\end{equation} with $\mathbb{U}=L^2(0,1)$ and $\mathbb{V}=H^1_{0)}:=\{v\in H^1(0,1):v(1)=0\}$.  The trial-to-test operator is no longer the identity and has the following integral form \cite{gopalakrishnan2013five, munoz2021equivalence}:
\begin{equation}\label{trial-to-test convection}
(Tu)(x)=\int_x^1 u(s) ds, \qquad u\in L^2(0,1).
\end{equation} For the exact solution and its corresponding optimal test function, we have
\begin{subequations}\label{equation:advection_solution}
\begin{align}
u^*&=\begin{cases}
0, &\text{if } 0<x<1/2,\\
1, &\text{if } 1/2<x<1,
\end{cases} \in L^2(0,1),\\
Tu^*&=\begin{cases}
1/2, &\text{if } 0<x<1/2,\\
1-x, &\text{if } 1/2<x<1,
\end{cases} \in H^1_{0)}(0,1).
\end{align}
\end{subequations}

In our context of NNs, considering $T$ as an available operator is challenging, so we discard employing the \ac{GDRM} here, and alternatively employ the (DRM$)'$---recall \eqref{adjoint_minimization}---and the \ac{D2RM}:

\begin{itemize}
\item \textbf{Adjoint Deep Ritz Method.} Following item (c) in \Cref{Generalized Ritz method}, we minimize $\mathcal{F}'$ to find an approximation to the optimal test function of the trial solution, i.e., 
\begin{equation}
Tu^* \approx \arg\min_{v_\text{NN}\in\mathbb{V}_\text{NN}} \mathcal{F}'(v_\text{NN})\footnote{The ``arg'' and ``min'' terms should be understood loosely according to the possible lack of existence discussed in \Cref{chapter2}.}. 
\end{equation} Then, we post-process the obtained minimizer by applying the available adjoint operator $A'=-d/dx$. 

We perform $50\mathord{,}000$ training iterations.  \Cref{table:experiments5} records the evolution of the approximated relative norm errors $\frac{\Vert v_\text{NN}-Tu^*\Vert_\V}{\Vert Tu^*\Vert_\V}$ and $\frac{\Vert A' v_\text{NN}-u^* \Vert_\mathbb{U}}{\Vert u^*\Vert_\mathbb{U}}$ along the training progress. \Cref{figure:Adjoint DRM_break} shows the predictions and errors of both $v_\text{NN}$ and $A' v_\text{NN}$ along different stages of the training. 

\begin{table}[htbp]
\centering
\begin{tabular}{|c|c|c|c|c|c|}
\toprule
\textbf{Training progress} & $\mathbf{4}\bf{\%}$ & $\mathbf{20}\bf{\%}$ & $\mathbf{40}\bf{\%}$ &$\mathbf{60}\bf{\%}$ & $\mathbf{100}\bf{\%}$\\
\toprule
$\frac{\Vert v_\text{NN}-Tu^*\Vert_\V}{\Vert Tu^*\Vert_\V}\times 100$  & $19.03\%$ & $4.58\%$ & $3.55\%$ & $3.20\%$ & $2.83\%$ \\ \midrule
$\frac{\Vert A' v_\text{NN}-u^*\Vert_\mathbb{U}}{\Vert u^*\Vert_\mathbb{U}}\times 100$  & $21.95\%$ & $5.28\%$ & $4.10\%$ & $3.69\%$ & $3.27\%$ \\ \bottomrule
\end{tabular}%
\caption{Relative errors of $v_\text{NN}$ and $A' v_\text{NN}$ along different stages of the training progress in problem \eqref{equation:advection} with exact solution \eqref{equation:advection_solution}. }
\label{table:experiments5}
\end{table}
\begin{figure}[htbp]
\centering
\begin{subfigure}[t]{\textwidth}
\centering
 \begin{tikzpicture}
\pgfplotsset{A/.append style={
		hide axis,
	    xmin=-0.1,
	    xmax=1.1,
	    ymin=-0.1,
	    ymax=1.1,
	 	legend style={anchor=north, /tikz/every even column/.append style={column sep=0.2cm}},
	 	legend columns = 3
            }
}
\begin{axis}[A]
\addlegendimage{color=red!30!white, line width=5}
\addlegendentry{$Tu^*$}
\addlegendimage{dashed, color=red!50!black, line width=.7, mark=asterisk, mark options = {scale=1, solid},}
\addlegendentry{$v_\text{NN}$}
\addlegendimage{color=red!50!black, line width=1}
\addlegendentry{$v_\text{NN}-Tu^*$}
\addlegendimage{color=blue!30!white, line width=5}
\addlegendentry{$u^*$}
\addlegendimage{dashed, color=blue!50!black,line width=.7, mark=asterisk, mark options = {scale=1, solid},}
\addlegendentry{$A' v_\text{NN}$}
\addlegendimage{color=blue!50!black, line width=1}
\addlegendentry{$A' v_\text{NN}-u^*$}
\end{axis}
\end{tikzpicture}
\end{subfigure}\vskip 1em%
\begin{subfigure}[t]{0.34\textwidth}
\centering
\begin{tikzpicture}
\begin{axis}[xmin=-0.05,
	    xmax=1.05,
	    xlabel = {$x$},
	    ylabel = {$v_\text{NN}$},
	    height=3.6cm,
	    width=\textwidth,
	    xtick={0,0.33,0.66,1},
	    yticklabel style={
        /pgf/number format/fixed,
        /pgf/number format/precision=2},
        scaled y ticks=false,
	    legend style={at={(0.5,0.95)}, anchor=north, /tikz/every even column/.append style={column sep=0.5cm}},
	 	]
  
	\addplot[color=red!30!white, line width=5,each nth point=2] table[x expr=\thisrow{x},y=v_analytic]{Chapters/4.Chapter/FIGURES/numerical_results/section55/AdjointDeepRitz_break_2000/v_analytic.csv};
	\addplot[dashed,color=red!50!black,  line width=.7,each nth point=2] table[x expr=\thisrow{x},y=v_net]{Chapters/4.Chapter/FIGURES/numerical_results/section55/AdjointDeepRitz_break_2000/v_net.csv};
	\addplot[only marks, color=red!50!black, line width=1, mark=asterisk, mark options = {scale=1, solid},each nth point=2] table[x expr=\thisrow{x},y=v_net]{Chapters/4.Chapter/FIGURES/numerical_results/section55/AdjointDeepRitz_break_2000/v_net_plot.csv};
\end{axis}
\end{tikzpicture}
\end{subfigure}\hfill%
\begin{subfigure}[t]{0.32\textwidth}
\centering
\begin{tikzpicture}
\begin{axis}[xmin=-0.05,
	    xmax=1.05,
	    xlabel = {$x$},
	    height=3.6cm,
	    width=\textwidth,
	    xtick={0,0.33,0.66,1},
	    yticklabel style={
        /pgf/number format/fixed,
        /pgf/number format/precision=2},
        scaled y ticks=false,
	    legend style={at={(0.5,0.95)}, anchor=north, /tikz/every even column/.append style={column sep=0.5cm}},
	 	]
  
	\addplot[color=red!30!white, line width=5,each nth point=2] table[x expr=\thisrow{x},y=v_analytic]{Chapters/4.Chapter/FIGURES/numerical_results/section55/AdjointDeepRitz_break_20000/v_analytic.csv};
	\addplot[dashed,color=red!50!black,line width=.7,each nth point=2] table[x expr=\thisrow{x},y=v_net]{Chapters/4.Chapter/FIGURES/numerical_results/section55/AdjointDeepRitz_break_20000/v_net.csv};
	\addplot[only marks, color=red!50!black, line width=1, mark=asterisk, mark options = {scale=1, solid},each nth point=2] table[x expr=\thisrow{x},y=v_net]{Chapters/4.Chapter/FIGURES/numerical_results/section55/AdjointDeepRitz_break_20000/v_net_plot.csv};

\end{axis}
\end{tikzpicture}
\end{subfigure}\hfill%
\begin{subfigure}[t]{0.32\textwidth}
\centering
\begin{tikzpicture}
\begin{axis}[xmin=-0.05,
	    xmax=1.05,
	    xlabel = {$x$},
	    height=3.6cm,
	    width=\textwidth,
	    xtick={0,0.33,0.66,1},
	    yticklabel style={
        /pgf/number format/fixed,
        /pgf/number format/precision=2},
        scaled y ticks=false,
	    legend style={at={(0.5,0.95)}, anchor=north, /tikz/every even column/.append style={column sep=0.5cm}},
	 	]
  
	\addplot[color=red!30!white, line width=5,each nth point=2] table[x expr=\thisrow{x},y=v_analytic]{Chapters/4.Chapter/FIGURES/numerical_results/section55/AdjointDeepRitz_break_50000/v_analytic.csv};
	\addplot[dashed,color=red!50!black,line width=.7,each nth point=2] table[x expr=\thisrow{x},y=v_net]{Chapters/4.Chapter/FIGURES/numerical_results/section55/AdjointDeepRitz_break_50000/v_net.csv};
	\addplot[only marks, color=red!50!black, line width=1, mark=asterisk, mark options = {scale=1, solid},each nth point=2] table[x expr=\thisrow{x},y=v_net]{Chapters/4.Chapter/FIGURES/numerical_results/section55/AdjointDeepRitz_break_50000/v_net_plot.csv};
\end{axis}
\end{tikzpicture}
\end{subfigure}\vskip 1em%
\begin{subfigure}[t]{0.34\textwidth}
\centering
\begin{tikzpicture}
\begin{axis}[xmin=-0.05,
	    xmax=1.05,
	    xlabel = {$x$},
	    ylabel = {$v_\text{NN}-Tu^*$},
	    height=3.6cm,
	    width=\textwidth,
	    xtick={0,0.33,0.66,1},
	    yticklabel style={
        /pgf/number format/fixed,
        /pgf/number format/precision=1},
        scaled y ticks=true,
	    legend style={at={(0.5,0.95)}, anchor=north, /tikz/every even column/.append style={column sep=0.5cm}},
	 	]
  
	\addplot[color=red!50!black,  line width=1,each nth point=2] table[x expr=\thisrow{x},y=e_v]{Chapters/4.Chapter/FIGURES/numerical_results/section55/AdjointDeepRitz_break_2000/e_v.csv};
\end{axis}
\end{tikzpicture}
\end{subfigure}\hfill
\begin{subfigure}[t]{0.32\textwidth}
\centering
\begin{tikzpicture}
\begin{axis}[xmin=-0.05,
	    xmax=1.05,
	    xlabel = {$x$},
	    height=3.6cm,
	    width=\textwidth,
	    xtick={0,0.33,0.66,1},
	    yticklabel style={
        /pgf/number format/fixed,
        /pgf/number format/precision=1},
        scaled y ticks=true,
	    legend style={at={(0.5,0.95)}, anchor=north, /tikz/every even column/.append style={column sep=0.5cm}},
	 	]
  
	\addplot[color=red!50!black,  line width=1,each nth point=2] table[x expr=\thisrow{x},y=e_v]{Chapters/4.Chapter/FIGURES/numerical_results/section55/AdjointDeepRitz_break_20000/e_v.csv};
\end{axis}
\end{tikzpicture}
\end{subfigure}\hfill%
\begin{subfigure}[t]{0.32\textwidth}
\centering
\begin{tikzpicture}
\begin{axis}[xmin=-0.05,
	    xmax=1.05,
	    xlabel = {$x$},
	    height=3.6cm,
	    width=\textwidth,
	    xtick={0,0.33,0.66,1},
	    yticklabel style={
        /pgf/number format/fixed,
        /pgf/number format/precision=1},
        scaled y ticks=true,
	    legend style={at={(0.5,0.95)}, anchor=north, /tikz/every even column/.append style={column sep=0.5cm}},
	 	]
  
	\addplot[color=red!50!black,  line width=1,each nth point=2] table[x expr=\thisrow{x},y=e_v]{Chapters/4.Chapter/FIGURES/numerical_results/section55/AdjointDeepRitz_break_50000/e_v.csv};
\end{axis}
\end{tikzpicture}
\end{subfigure}\vskip 1em%
\begin{subfigure}[t]{0.34\textwidth}
\centering
\begin{tikzpicture}
\begin{axis}[xmin=-0.05,
	    xmax=1.05,
	    xlabel = {$x$},
	    ylabel = {$A' v_\text{NN}$},
	    height=3.6cm,
	    width=\textwidth,
	    xtick={0,0.33,0.66,1},
	    yticklabel style={
        /pgf/number format/fixed,
        /pgf/number format/precision=2},
        scaled y ticks=false,
	    legend style={at={(0.5,0.95)}, anchor=north, /tikz/every even column/.append style={column sep=0.5cm}},
	 	]
  
	\addplot[color=blue!30!white, line width=5,each nth point=2,  domain=0:1/2] {0.};
	\addplot[color=blue!30!white, line width=5,each nth point=2, domain=1/2:1] {1.};
	\addplot[dashed,color=blue!50!black,  line width=.7,each nth point=2] table[x expr=\thisrow{x},y=u_net]{Chapters/4.Chapter/FIGURES/numerical_results/section55/AdjointDeepRitz_break_2000/u_net.csv};
	\addplot[only marks, color=blue!50!black, line width=1, mark=asterisk, mark options = {scale=1, solid},each nth point=2] table[x expr=\thisrow{x},y=u_net]{Chapters/4.Chapter/FIGURES/numerical_results/section55/AdjointDeepRitz_break_2000/u_net_plot.csv};
\end{axis}
\end{tikzpicture}
\end{subfigure}\hfill%
\begin{subfigure}[t]{0.32\textwidth}
\centering
\begin{tikzpicture}
\begin{axis}[xmin=-0.05,
	    xmax=1.05,
	    xlabel = {$x$},
	    height=3.6cm,
	    width=\textwidth,
	    xtick={0,0.33,0.66,1},
	    yticklabel style={
        /pgf/number format/fixed,
        /pgf/number format/precision=2},
        scaled y ticks=false,
	    legend style={at={(0.5,0.95)}, anchor=north, /tikz/every even column/.append style={column sep=0.5cm}},
	 	]
  
	\addplot[color=blue!30!white, line width=5,each nth point=2,  domain=0:1/2] {0.};
	\addplot[color=blue!30!white, line width=5,each nth point=2, domain=1/2:1] {1.};
	\addplot[dashed,color=blue!50!black,  line width=.7,each nth point=2] table[x expr=\thisrow{x},y=u_net]{Chapters/4.Chapter/FIGURES/numerical_results/section55/AdjointDeepRitz_break_20000/u_net.csv};
	\addplot[only marks, color=blue!50!black, line width=1, mark=asterisk, mark options = {scale=1, solid},each nth point=2] table[x expr=\thisrow{x},y=u_net]{Chapters/4.Chapter/FIGURES/numerical_results/section55/AdjointDeepRitz_break_20000/u_net_plot.csv};
\end{axis}
\end{tikzpicture}
\end{subfigure}\hfill%
\begin{subfigure}[t]{0.32\textwidth}
\centering
\begin{tikzpicture}
\begin{axis}[xmin=-0.05,
	    xmax=1.05,
	    xlabel = {$x$},
	    height=3.6cm,
	    width=\textwidth,
	    xtick={0,0.33,0.66,1},
	    yticklabel style={
        /pgf/number format/fixed,
        /pgf/number format/precision=2},
        scaled y ticks=false,
	    legend style={at={(0.5,0.95)}, anchor=north, /tikz/every even column/.append style={column sep=0.5cm}},
	 	]
  
	\addplot[color=blue!30!white, line width=5,each nth point=2,  domain=0:1/2] {0.};
	\addplot[color=blue!30!white, line width=5,each nth point=2, domain=1/2:1] {1.};
	\addplot[dashed,color=blue!50!black,  line width=.7,each nth point=2] table[x expr=\thisrow{x},y=u_net]{Chapters/4.Chapter/FIGURES/numerical_results/section55/AdjointDeepRitz_break_50000/u_net.csv};
	\addplot[only marks, color=blue!50!black, line width=1, mark=asterisk, mark options = {scale=1, solid},each nth point=2] table[x expr=\thisrow{x},y=u_net]{Chapters/4.Chapter/FIGURES/numerical_results/section55/AdjointDeepRitz_break_50000/u_net_plot.csv};
\end{axis}
\end{tikzpicture}
\end{subfigure}\vskip 1em%
\begin{subfigure}[t]{0.34\textwidth}
\centering
\begin{tikzpicture}
\begin{axis}[xmin=-0.05,
	    xmax=1.05,
	    xlabel = {$x$},
	    ylabel = {$A' v_\text{NN}-u^*$},
	    height=3.6cm,
	    width=\textwidth,
	    xtick={0,0.33,0.66,1},
	    yticklabel style={
        /pgf/number format/fixed,
        /pgf/number format/precision=1},
        scaled y ticks=true,
	    legend style={at={(0.5,0.95)}, anchor=north, /tikz/every even column/.append style={column sep=0.5cm}},
	 	]
  
	\addplot[color=blue!50!black,  line width=1,each nth point=2] table[x expr=\thisrow{x},y=e_u]{Chapters/4.Chapter/FIGURES/numerical_results/section55/AdjointDeepRitz_break_2000/e_u.csv};
\end{axis}
\end{tikzpicture}
\caption{$4\%$ of the training.}
\end{subfigure}\hfill%
\begin{subfigure}[t]{0.32\textwidth}
\centering
\begin{tikzpicture}
\begin{axis}[xmin=-0.05,
	    xmax=1.05,
	    xlabel = {$x$},
	    height=3.6cm,
	    width=\textwidth,
	    xtick={0,0.33,0.66,1},
	    yticklabel style={
        /pgf/number format/fixed,
        /pgf/number format/precision=1},
        scaled y ticks=true,
	    legend style={at={(0.5,0.95)}, anchor=north, /tikz/every even column/.append style={column sep=0.5cm}},
	 	]
  
	\addplot[color=blue!50!black,  line width=1,each nth point=2] table[x expr=\thisrow{x},y=e_u]{Chapters/4.Chapter/FIGURES/numerical_results/section55/AdjointDeepRitz_break_20000/e_u.csv};
\end{axis}
\end{tikzpicture}
\caption{$40\%$ of the training.}
\end{subfigure}\hfill%
\begin{subfigure}[t]{0.32\textwidth}
\centering
\begin{tikzpicture}
\begin{axis}[xmin=-0.05,
	    xmax=1.05,
	    xlabel = {$x$},
	    height=3.6cm,
	    width=\textwidth,
	    xtick={0,0.33,0.66,1},
	    yticklabel style={
        /pgf/number format/fixed,
        /pgf/number format/precision=1},
        scaled y ticks=true,
	    legend style={at={(0.5,0.95)}, anchor=north, /tikz/every even column/.append style={column sep=0.5cm}},
	 	]
  
	\addplot[color=blue!50!black,  line width=1,each nth point=2] table[x expr=\thisrow{x},y=e_u]{Chapters/4.Chapter/FIGURES/numerical_results/section55/AdjointDeepRitz_break_50000/e_u.csv};
\end{axis}
\end{tikzpicture}
\caption{$100\%$ of the training.}
\end{subfigure}
\caption{Test network predicitions $v_\text{NN}$,  post-processings $A' v_\text{NN}$,  and error functions for the (DRM$)'$ in model problem \eqref{equation:advection} with exact solution \eqref{equation:advection_solution} along different stages of the training progress.}
\label{figure:Adjoint DRM_break}
\end{figure}
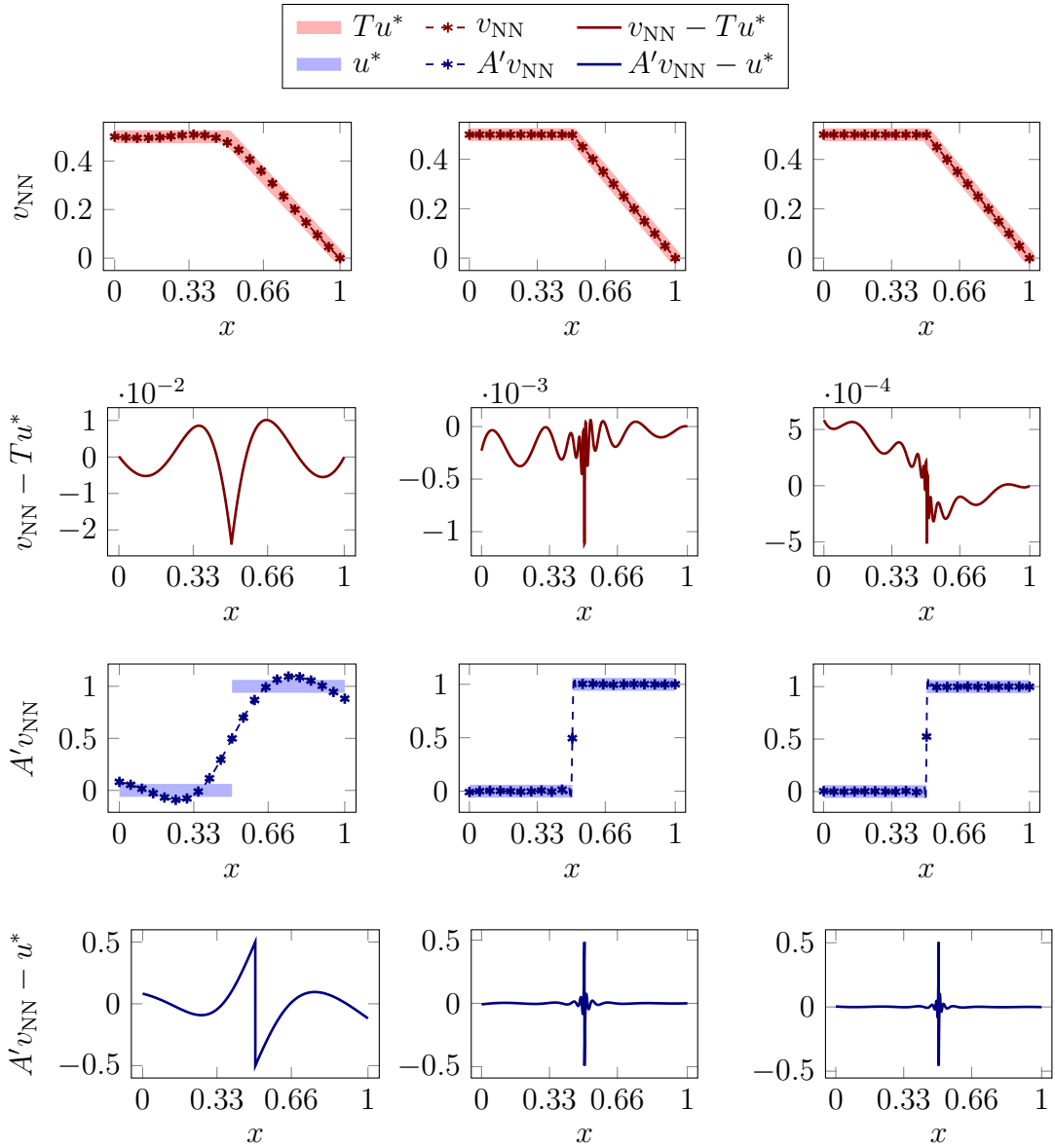

\item  \textbf{Deep Double Ritz Method.} We impose the outflow-boundary condition only to $\tau_\text{NN}$, perform $50\mathord{,}000$ iterations for $u_\text{NN}$ and increase the number of inner-loop iterations from four to nine (i.e., a total of $500\mathord{,}000$ iterations if we account those dedicated to the test functions).  \Cref{table:experiments6} records the relative errors for $u_\text{NN}$ and $\tau_\text{NN} u_\text{NN}$ along different stages of the training. \Cref{figure:DDRM_break} shows the trial and optimal test network predictions and error functions at the end of the training.
\begin{table}[htbp]
\centering
\begin{tabular}{|c|c|c|c|c|c|}
\toprule
\textbf{Training progress} & $\mathbf{4}\bf{\%}$ & $\mathbf{20}\bf{\%}$ & $\mathbf{40}\bf{\%}$ &$\mathbf{60}\bf{\%}$ & $\mathbf{100}\bf{\%}$\\
\toprule
$\frac{\Vert u_\text{NN} - u^*\Vert_\mathbb{U}}{\Vert u^*\Vert_\mathbb{U}}$ & $42.05\%$ & $29.69\%$ & $18.42\%$ & $11.20\%$ & $8.93\%$ \\ \midrule
$\frac{\Vert \tau_\text{NN}(u_\text{NN}) - Tu^*\Vert_\V}{\Vert Tu^*\Vert_\V}$ & $37.27\%$ & $28.25\%$ & $19.02\%$ & $9.55\%$ & $6.13\%$ \\ \bottomrule
\end{tabular}%
\caption{Relative errors of $u_\text{NN}$ and $\tau_\text{NN}(u_\text{NN})$ along different stages of the training progress in problem \eqref{equation:advection} with exact solution \eqref{equation:advection_solution}. }
\label{table:experiments6}
\end{table}

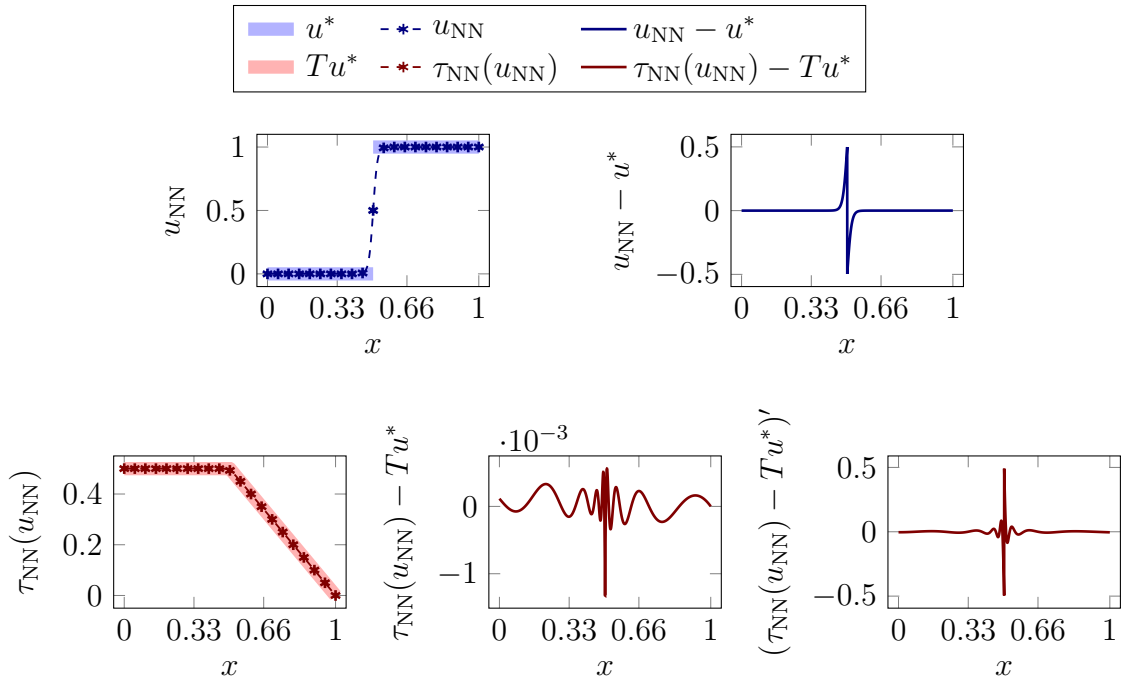
\begin{figure}[htbp]
\centering
\begin{subfigure}[t]{\textwidth}
\centering
 \begin{tikzpicture}
\pgfplotsset{A/.append style={hide axis,
	    xmin=-0.1,
	    xmax=1.1,
	    ymin=-0.1,
	    ymax=1.1,
	    legend style={anchor=north, /tikz/every even column/.append style={column sep=0.2cm}},
	    legend columns = 3
            }
}
\begin{axis}[A]
\addlegendimage{color=blue!30!white, line width=5}
\addlegendentry{$u^*$}
\addlegendimage{dashed,color=blue!50!black,line width=.7, mark=asterisk, mark options = {scale=1, solid},}
\addlegendentry{$u_\text{NN}$}
\addlegendimage{color=blue!50!black, line width=1}
\addlegendentry{$u_\text{NN}-u^*$}
\addlegendimage{color=red!30!white, line width=5}
\addlegendentry{$Tu^*$}
\addlegendimage{dashed,color=red!50!black, line width=.7, mark=asterisk, mark options = {scale=1, solid},}
\addlegendentry{$\tau_\text{NN}(u_\text{NN})$}
\addlegendimage{color=red!50!black, line width=1}
\addlegendentry{$\tau_\text{NN}(u_\text{NN})-Tu^*$}
\end{axis}
\end{tikzpicture}
\end{subfigure}\vskip 1em%
\begin{subfigure}[t]{0.32\textwidth}
\centering
\begin{tikzpicture}
\begin{axis}[
		xmin=-0.05,
	    xmax=1.05,
	    xlabel = {$x$},
	    ylabel = {$u_\text{NN}$},
	    height=3.6cm,
	    width=\textwidth,
	    xtick={0,0.33,0.66,1},
	    yticklabel style={
        /pgf/number format/fixed,
        /pgf/number format/precision=2},
        scaled y ticks=false,
	    legend style={at={(0.5,0.95)}, anchor=north, /tikz/every even column/.append style={column sep=0.5cm}},
	 	]
  
	\addplot[color=blue!30!white, line width=5,each nth point=2,  domain=0:1/2] {0.};
	\addplot[color=blue!30!white, line width=5,each nth point=2, domain=1/2:1] {1.};
	\addplot[dashed, color=blue!50!black,  line width=.7,each nth point=2] table[x expr=\thisrow{x},y=u_net]{Chapters/4.Chapter/FIGURES/numerical_results/section55/c/u_net.csv};
	\addplot[only marks, color=blue!50!black, line width=1, mark=asterisk, mark options = {scale=1, solid},each nth point=2] table[x expr=\thisrow{x},y=u_net]{Chapters/4.Chapter/FIGURES/numerical_results/section55/c/u_net_plot.csv};
\end{axis}
\end{tikzpicture}
\end{subfigure}\hskip 3em%
\begin{subfigure}[t]{0.32\textwidth}
\centering
\begin{tikzpicture}
\begin{axis}[xmin=-0.05,
	    xmax=1.05,
	    xlabel = {$x$},
	    ylabel = {$u_\text{NN}-u^*$},
	    height=3.6cm,
	    width=\textwidth,
	    xtick={0,0.33,0.66,1},
	    yticklabel style={
        /pgf/number format/fixed,
        /pgf/number format/precision=1},
        scaled y ticks=true,
	    legend style={at={(0.5,0.95)}, anchor=north, /tikz/every even column/.append style={column sep=0.5cm}},
	 	]  
	\addplot[color=blue!50!black,  line width=1,each nth point=2] table[x expr=\thisrow{x},y=e_u]{Chapters/4.Chapter/FIGURES/numerical_results/section55/c/e_u.csv};
\end{axis}
\end{tikzpicture}
\end{subfigure}\vskip 1em%
\begin{subfigure}[t]{0.32\textwidth}
\centering
\begin{tikzpicture}
\begin{axis}[xmin=-0.05,
	    xmax=1.05,
	    xlabel = {$x$},
	    ylabel = {$\tau_\text{NN}(u_\text{NN})$},
	    height=3.6cm,
	    width=\textwidth,
	    xtick={0,0.33,0.66,1},
	    yticklabel style={
        /pgf/number format/fixed,
        /pgf/number format/precision=2},
        scaled y ticks=false,
	    legend style={at={(0.5,0.95)}, anchor=north, /tikz/every even column/.append style={column sep=0.5cm}},
	 	]
  
	\addplot[color=red!30!white, line width=5,each nth point=2] table[x expr=\thisrow{x},y=Tu_analytic]{Chapters/4.Chapter/FIGURES/numerical_results/section55/a/Tu_analytic.csv};
	\addplot[dashed, color=red!50!black,line width=.7,each nth point=2] table[x expr=\thisrow{x},y=Tu_net]{Chapters/4.Chapter/FIGURES/numerical_results/section55/c/Tu_net.csv};
	\addplot[only marks, color=red!50!black, line width=1, mark=asterisk, mark options = {scale=1, solid},each nth point=2] table[x expr=\thisrow{x},y=Tu_net]{Chapters/4.Chapter/FIGURES/numerical_results/section55/c/Tu_net_plot.csv};
\end{axis}
\end{tikzpicture}
\end{subfigure}\hfill%
\begin{subfigure}[t]{0.32\textwidth}
\centering
\begin{tikzpicture}
\begin{axis}[xmin=-0.05,
	    xmax=1.05,
	    xlabel = {$x$},
	    ylabel = {$\tau_\text{NN}(u_\text{NN})-Tu^*$},
	    height=3.6cm,
	    width=\textwidth,
	    xtick={0,0.33,0.66,1},
	    yticklabel style={
        /pgf/number format/fixed,
        /pgf/number format/precision=1},
        scaled y ticks=true,
	    legend style={at={(0.5,0.95)}, anchor=north, /tikz/every even column/.append style={column sep=0.5cm}},
	 	]  
	\addplot[color=red!50!black,  line width=1,each nth point=2] table[x expr=\thisrow{x},y=e_v]{Chapters/4.Chapter/FIGURES/numerical_results/section55/AdjointDeepRitz_break_10000/e_v.csv};
\end{axis}
\end{tikzpicture}
\end{subfigure}\hfill%
\begin{subfigure}[t]{0.32\textwidth}
\centering
\begin{tikzpicture}
\begin{axis}[xmin=-0.05,
	    xmax=1.05,
	    xlabel = {$x$},
	    ylabel = {$(\tau_\text{NN}(u_\text{NN}) - Tu^*)'$},
	    height=3.6cm,
	    width=\textwidth,
	    xtick={0,0.33,0.66,1},
	    yticklabel style={
        /pgf/number format/fixed,
        /pgf/number format/precision=1},
        scaled y ticks=true,
	    legend style={at={(0.5,0.95)}, anchor=north, /tikz/every even column/.append style={column sep=0.5cm}},
	 	]
  
	\addplot[color=red!50!black,  line width=1,each nth point=2] table[x expr=\thisrow{x},y=de_v]{Chapters/4.Chapter/FIGURES/numerical_results/section55/AdjointDeepRitz_break_10000/de_v.csv};
\end{axis}
\end{tikzpicture}
\end{subfigure}
\caption{$u_\text{NN}$ and $\tau_\text{NN}(u_\text{NN})$ predictions,  errors, and derivative of the errors for the \acs{D2RM} in problem \eqref{equation:advection} with exact solution \eqref{equation:advection_solution}.}
\label{figure:DDRM_break}
\end{figure}

\end{itemize}

\subsection{Pure convection in 2D}\label{section:Convection equation in 2D}

Let
\begin{equation} \label{equation:advection2D}
\begin{cases}
\displaystyle \frac{\partial u}{\partial x} + \frac{\partial u}{\partial y} =k\pi \sin\left(k\pi (x+y)\right), &\text{in } \Omega=(0,1)\times(0,1),\\
u(x,0)=u(0,y)=0,& 0\leq x,y\leq 1,
\end{cases}
\end{equation} and consider its strong variational formulation
\begin{equation}
b(u,v):=\int_{\Omega} \left(\frac{\partial u}{\partial x} + \frac{\partial u}{\partial y}\right)v, \quad l(v):=k\pi \int_{\Omega}  \sin(k\pi (x+y)) v,
\end{equation} for $u\in\mathbb{U}$ and $v\in\mathbb{V}$ such that
\begin{subequations}
\begin{align}
\mathbb{U}&=\left\lbrace u\in L^2(\Omega): \frac{\partial u}{\partial x} + \frac{\partial u}{\partial y} \in L^2(\Omega) \text{ and }u(x,y)=0 \text{ when } xy=0\right\rbrace,\\
\mathbb{V}&=L^2(\Omega).
\end{align}
\end{subequations} Its exact solution is $u^*=\sin(kx)\sin(ky)$ and the trial-to-test operator is the PDE operator.

For $k=3/2$, we perform $200\mathord{,}000$ iterations in the \ac{D2RM} with a training regime of nine iterations in the inner loop for each iteration in the outer loop. We run the inner loop for $2\mathord{,}000$ iterations before the first outer-loop iteration.  For integration, we consider $50$ nodes on each axis (i.e., $250$ integration points on the entire domain due to the cartesian-product structure). We increase the NN architecture to three layers of $50$-neuron width.

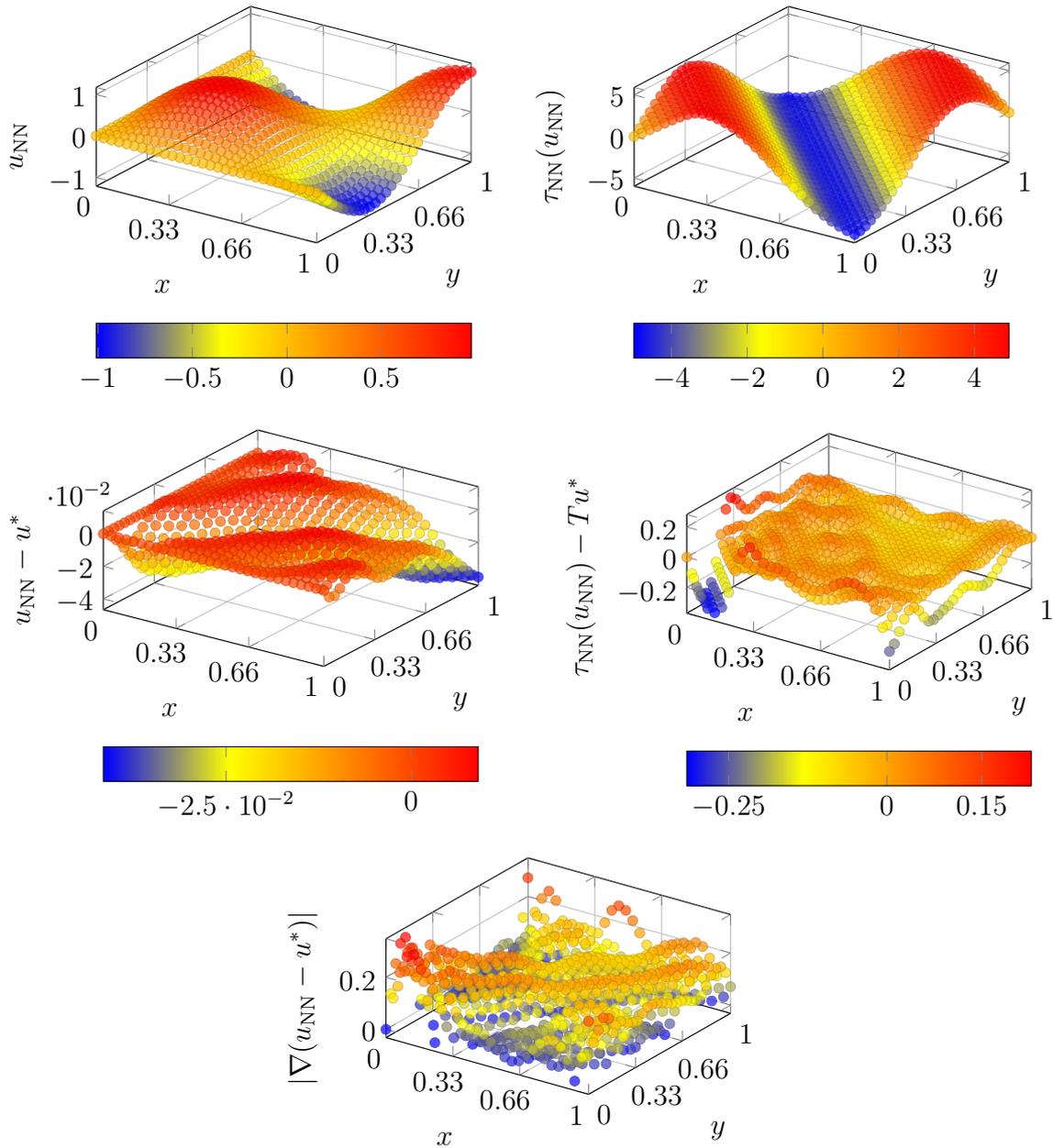
\begin{figure}[htbp]
\centering
\begin{subfigure}[t]{0.48\textwidth}
\centering
\begin{tikzpicture}
\begin{axis}[xrange=0:1,
        yrange=0:1,
        xlabel = {$x$},
        ylabel = {$y$},
        zlabel={$u_\text{NN}$},
        height=5cm,
        width=\textwidth,
        view = {35}{45},
        grid=major,
        colorbar horizontal,
        xtick={0,0.33,0.66,1},
        ytick={0,0.33,0.66,1},
        yticklabel style={
        /pgf/number format/fixed,
        /pgf/number format/precision=3},
        scaled y ticks=true,
             select coords between index/.style 2 args={
    x filter/.code={
        \ifnum\coordindex<#1\def\pgfmathresult{}\fi
        \ifnum\coordindex>#2\def\pgfmathresult{}\fi
    }
},
	    legend style={at={(0.5,0.95)}, anchor=north, /tikz/every even column/.append style={column sep=0.5cm}},
	 	]	 	
        \addplot3[domain=0:1, y domain=0:1, surf,scatter, only marks, mark=*,mesh/ordering=y varies, opacity=0.7] table[x expr=\thisrow{x},y expr=\thisrow{y}, z=u_net]{Chapters/4.Chapter/FIGURES/numerical_results/section56/solutions_2D/u_net_plot.csv};

\end{axis}
\end{tikzpicture}
\end{subfigure}\hfill%
\begin{subfigure}[t]{0.48\textwidth}
\centering
\begin{tikzpicture}
\begin{axis}[xrange=0:1,
	yrange=0:1,
	xlabel = {$x$},
	ylabel = {$y$},
	zlabel={$\tau_\text{NN}(u_\text{NN})$},
        height=5cm,
	width=\textwidth,
	view = {35}{45},
        grid=major,
        colorbar horizontal,
        xtick={0,0.33,0.66,1},
        ytick={0,0.33,0.66,1},
        yticklabel style={
        /pgf/number format/fixed,
        /pgf/number format/precision=3},
        scaled y ticks=true,
             select coords between index/.style 2 args={
    x filter/.code={
        \ifnum\coordindex<#1\def\pgfmathresult{}\fi
        \ifnum\coordindex>#2\def\pgfmathresult{}\fi
    }
},
	    legend style={at={(0.5,0.95)}, anchor=north, /tikz/every even column/.append style={column sep=0.5cm}},
	 	]

        \addplot3[domain=0:1, y domain=0:1, surf,scatter, only marks, mark=*,mesh/ordering=y varies, opacity=0.7] table[x expr=\thisrow{x},y expr=\thisrow{y}, z=Tu_net]{Chapters/4.Chapter/FIGURES/numerical_results/section56/solutions_2D/Tu_net_plot.csv};
\end{axis}
\end{tikzpicture}
\end{subfigure}\vskip 1em%
\begin{subfigure}[t]{0.48\textwidth}
\centering
\begin{tikzpicture}
\begin{axis}[xrange=0:1,
        yrange=0:1,
        xlabel = {$x$},
        ylabel = {$y$},
        zlabel = {$u_\text{NN}-u^*$},
        height=5cm,
        width=\textwidth,
        view = {35}{45},
        grid=major,
        colorbar horizontal,
        colorbar style={scaled x ticks=false},
        colorbar style={xtick={-0.025, 0, 0.0175}},
        xtick={0,0.33,0.66,1},
        ytick={0,0.33,0.66,1},
	    yticklabel style={
        /pgf/number format/fixed,
        /pgf/number format/precision=4},
        scaled y ticks=true,
             select coords between index/.style 2 args={
    x filter/.code={
        \ifnum\coordindex<#1\def\pgfmathresult{}\fi
        \ifnum\coordindex>#2\def\pgfmathresult{}\fi
    }
},
	    legend style={at={(0.5,0.95)}, anchor=north, /tikz/every even column/.append style={column sep=0.5cm}},
	 	]

        \addplot3[domain=0:1, y domain=0:1, surf,scatter, only marks, mark=*,mesh/ordering=y varies, opacity=0.7] table[x expr=\thisrow{x},y expr=\thisrow{y}, z=e_u]{Chapters/4.Chapter/FIGURES/numerical_results/section56/solutions_2D/e_u_plot.csv};
\end{axis}
\end{tikzpicture}
\end{subfigure}\hfill%
\begin{subfigure}[t]{0.45\textwidth}
\centering
\begin{tikzpicture}
\begin{axis}[xrange=0:1,
        yrange=0:1,
        xlabel = {$x$},
        ylabel = {$y$},
        zlabel = {$\tau_\text{NN}(u_\text{NN})-Tu^*$},
        height=5cm,
        width=\textwidth,
        view = {35}{45},
        grid=major,
        colorbar horizontal,
        colorbar style={xtick={-0.25, 0, 0.15}},
        xtick={0,0.33,0.66,1},
        ytick={0,0.33,0.66,1},
	    yticklabel style={
        /pgf/number format/fixed,
        /pgf/number format/precision=3},
        scaled y ticks=true,
             select coords between index/.style 2 args={
    x filter/.code={
        \ifnum\coordindex<#1\def\pgfmathresult{}\fi
        \ifnum\coordindex>#2\def\pgfmathresult{}\fi
    }
},
	    legend style={at={(0.5,0.95)}, anchor=north, /tikz/every even column/.append style={column sep=0.5cm}},
	 	]

        \addplot3[domain=0:1, y domain=0:1, surf,scatter, only marks, mark=*,mesh/ordering=y varies, opacity=0.7] table[x expr=\thisrow{x},y expr=\thisrow{y}, z=e_Tu]{Chapters/4.Chapter/FIGURES/numerical_results/section56/solutions_2D/e_Tu_plot.csv};
\end{axis}
\end{tikzpicture}
\end{subfigure}\vskip 1em%
\begin{subfigure}[t]{0.45\textwidth}
\centering
\begin{tikzpicture}
\begin{axis}[xrange=0:1,
        yrange=0:1,
        xlabel = {$x$},
        ylabel = {$y$},
        zlabel = {$\vert\nabla(u_\text{NN}-u^*)\vert$},
        height=5cm,
        width=\textwidth,
        view = {35}{45},
        grid=major,
        xtick={0,0.33,0.66,1},
        ytick={0,0.33,0.66,1},
	    yticklabel style={
        /pgf/number format/fixed,
        /pgf/number format/precision=3},
        scaled y ticks=true,
             select coords between index/.style 2 args={
    x filter/.code={
        \ifnum\coordindex<#1\def\pgfmathresult{}\fi
        \ifnum\coordindex>#2\def\pgfmathresult{}\fi
    }
},
	    legend style={at={(0.5,0.95)}, anchor=north, /tikz/every even column/.append style={column sep=0.5cm}},
	 	]

        \addplot3[domain=0:1, y domain=0:1, surf,scatter, only marks, mark=*,mesh/ordering=y varies, opacity=0.7] table[x expr=\thisrow{x},y expr=\thisrow{y}, z=de_u]{Chapters/4.Chapter/FIGURES/numerical_results/section56/solutions_2D/de_u_plot.csv};
\end{axis}
\end{tikzpicture}
\end{subfigure}
\caption{Trial and optimal test predictions, errors, and derivative of the trial error for the \acs{D2RM} in model problem \eqref{equation:advection} with exact solution \eqref{equation:advection_solution}.}
\label{figure:DDRM_2D}
\end{figure}

\Cref{figure:DDRM_2D} shows the trial and optimal test predictions with corresponding error functions at the end of the training. The resulting relative errors are $3.62\%$ and $1.70\%$ for $u_\n$ and $\tau_\n(u_\n)$, respectively.

\chapter{Memory-based Monte Carlo integration} \label{chapter5}

\begin{quote}
\textbf{Summary.} Monte Carlo integration is a widely used quadrature rule to solve Partial Differential Equations using Neural Networks due to its ability to limit overfitting and high-dimensional scalability.  However, this stochastic method produces noisy losses and gradients during training, which hinders a proper convergence diagnosis. Typically, this is overcome using an immense (disproportionate) amount of integration points, which deteriorates the training performance. This work proposes a memory-based Monte Carlo integration method that produces accurate integral approximations without requiring the high computational costs of processing large samples during training.  \emph{Refer to \cite{uriarte2023memory} for the published version}.
\end{quote}

\section{Introduction}

Using a deterministic quadrature rule with fixed integration points may allow misbehavior of the NN away from the integration points (an overfitting issue) \cite{rivera2022quadrature},  producing significant integration errors and, consequently, poor solutions.  To overcome this,  \emph{\acf{MC} integration} is a popular and suitable choice of quadrature rule due to the mesh-free and stochastic sampling of the integration points during training \cite{leobacher2014introduction, grohs2022deep,dick2013high,chen2021quasi}. However,  MC integration error is of order $\mathcal{O}(1/\sqrt{N})$, where $N$ is the number of integration points \cite{newman1999monte}. Thus,  in practice, this may require tens or hundreds of thousands of integration points to obtain an acceptable error per integral approximation---even for one-dimensional problems,  which deteriorates the training speed.

In this work, we propose a memory-based approach that approximates definite integrals involving NNs by taking advantage of the information gained in previous iterations.  As long as the expected value of these integrals does not change significantly, this technique reduces the expected integration error and leads to better approximations.  Moreover, since gradients are also described in terms of definite integrals, we apply this approach to the gradient computations,  which reinterprets the well-known momentum method \cite{polyak1964some} when we appropriately modify the hyperparameters of the optimizer.

The remainder of this chapter is organized as follows. \Cref{section5.2} reviews the approximation and optimization frameworks on \acp{NN}, \Cref{section5.3} proposes the memory-based integration and optimization strategies, and \Cref{section5.4} relates these proposals with the momentum method.

\section{Review of the approximation, discretization, parameterization, and optimization setups}\label[section]{section5.2}

We summarize many of the concepts already discussed in \Cref{chapter2} to (re)define notation and present the method in question more fluently for the reader.

Let us consider a well-defined minimization problem of the form
\begin{equation}\label{minimization_problem}
u^*=\arg\min_{u\in\mathbb{U}} \mathcal{F}(u),
\end{equation} where $\mathbb{U}$ denotes the search space of functions with domain $\Omega$, $\mathcal{F}:\mathbb{U}\longrightarrow\mathbb{R}$ is the objective function governing our minimization problem, and $u^*$ is the exact solution.

Let $u_\theta:\Omega\longrightarrow\mathbb{R}$ denote a neural network architecture parameterized by the set of trainable parameters $\theta$ with domain $\Theta$. Then, a \ac{NN} approximation of problem \eqref{minimization_problem} consists in replacing the continuum-level search space $\mathbb{U}$ with the parameterized search space $\mathbb{U}_\Theta =  \{u_\theta: \theta\in\Theta\}$.

To carry out the minimization, we resort to a first-order gradient-descent scheme that is described as follows:
\begin{align}
\theta_{t+1} & = \theta_t -\lambda \frac{\partial\mathcal{F}}{\partial\theta}(v_{\theta_t}),\label{gradients}
\end{align} where $\lambda>0$ is the learning rate  and $\theta_t$ denotes the trainable parameters at the $t^\text{th}$ iteration. 

For $\mathcal{F}$ in the form of a definite integral, 
\begin{equation}
    \mathcal{F}(v_\theta)=\int_\Omega I(v_\theta)(x) \ dx,
\end{equation} we approximate it by a quadrature rule, producing a \emph{loss function} $\mathcal{L}$.  Considering MC integration as the quadrature rule for the loss, we have
\begin{equation}\label{Monte_Carlo_rule}
\mathcal{F}(v_\theta) \approx \mathcal{L}(v_\theta):= \frac{\text{Vol}(\Omega)}{N}\sum_{i=1}^N I(v_\theta)(x_i),
\end{equation} where $\{x_i\}_{i=1}^N$ is a stochastic set of integration points sampled from a random uniform distribution in $\Omega$.

Similarly, for the gradients, we have
\begin{equation}\label{practical_gradients}
    \frac{\partial\mathcal{F}}{\partial\theta}(v_\theta) \approx g(v_\theta):=\frac{\partial\mathcal{L}}{\partial\theta}(v_\theta) = \frac{\text{Vol}(\Omega)}{N}\sum_{i=1}^N \frac{\partial I(v_\theta)}{\partial\theta}(x_i),
\end{equation} which yields a discretized version of \eqref{gradients}, the SGD optimizer \cite{robbins1951stochastic},
\begin{equation}\label{numeric_gradients}
    \theta_{t+1} := \theta_t - \lambda g(v_{\theta_{t}}).
\end{equation}

From now on, we write $\mathcal{F}(\theta_t)$, $\frac{\partial\mathcal{F}}{\partial\theta}(\theta_t)$, $\mathcal{L}(\theta_t)$, and $g(\theta_t)$ as simplified versions of $\mathcal{F}(v_{\theta_t})$, $\frac{\partial\mathcal{F}}{\partial\theta}(v_{\theta_t})$, $\mathcal{L}(v_{\theta_t})$, and $g(v_{\theta_t})$, respectively\footnote{Note that we have dropped the \emph{realization map} terminology---recall \eqref{realization_map}---when expressing the dependence of the objective function on the parameters. That is, $\mathcal{F}(\theta)$ directly refers to $(\mathcal{F}\circ\Phi_\text{NN})(\theta)$ according to the notation developed in \Cref{chapter2}.}.

\section{Memory-based integration and optimization}\label[section]{section5.3}

If we train the network according to \eqref{numeric_gradients}, we obtain a noisy and oscillatory behavior of the loss and the gradient. This occurs because of the introduced MC integration error at each training iteration.  \Cref{monte_carlo_first_training} (blue curve) illustrates the noisy behavior of MC integration in a network with a single trainable parameter that permits exact calculation of $\mathcal{F}$ (black curve).  

This section is divided into three parts.  \Cref{section5.2.1} introduces the considered model problem for experimentation,  \Cref{section5.2.2} describes the proposed memory-based MC integration rule, and \Cref{section5.2.3} extends such memory-based scheme to the gradients for training.

\subsection{Model problem}\label[section]{section5.2.1}
Let us consider the model problem $-u''=4\delta_{1/2}$ in $\Omega=(0,1)$ with homogeneous Dirichlet conditions on $\partial\Omega=\{0,1\}$.  Then, its weak form reformulation possesses the following bilinear and linear forms:
\begin{equation}
b(u,v)=\int_0^1 u'v', \qquad l(v)=v(1/2),\qquad u\in \mathbb{U}=H^1_0(0,1)=\mathbb{V}\ni v.
\end{equation} Its unique exact solution is $u^*(x)=2x$ in $0\leq x\leq 1/2$ and $u^*(x)=2(1-x)$ in $1/2\leq x\leq 1$.  Moreover,  it admits a Ritz-type minimization reformulation as follows:
\begin{equation}
u^*=\arg\min_{u\in \mathbb{U}} \frac{1}{2} \int_0^1 [u'(x)]^2 dx - u(1/2).
\end{equation}

Noticing that $2\tanh(\theta(x-1/2))$ has arbitrary approximation capacity for $(u^*)'$---recall \Cref{section2.4}, integrating and imposing the boundary conditions, we obtain the (atypical) NN architecture with arbitrary approximation capacity for $u^*$ given by
\begin{equation}
u_\theta(x)=\frac{2}{\theta} \Big\lbrace \log\big[\cosh(\theta/2)\big] - \log\big[\cosh(\theta(1/2-x))\big]\Big\rbrace.
\end{equation} Consequently,  straightforward calculations yield:
\begin{subequations}
\begin{align}
\mathcal{F}(\theta)&=2-\frac{4}{\theta} \tanh(\theta/2) - \frac{8}{\theta}\log\big[\cosh(\theta/2)\big],\\
\frac{\partial\mathcal{F}}{\partial\theta}(\theta) &= \frac{-4(\theta-1)\tanh(\theta/2) - 2\theta\text{sech}^2(\theta/2)+8\log[\cosh(\theta/2)]}{\theta^2},
\end{align}
\end{subequations} for $\theta\neq 0$. Notice that approximation is achieved as $\theta\to\infty$ with corresponding $\mathcal{F}(\theta)\to-2$ and $\frac{\partial\mathcal{F}}{\partial\theta}(\theta)\to 0$.

\subsection{Integration}\label[subsection]{section5.2.2}

In order to decrease the integration error during training, we replace \eqref{Monte_Carlo_rule} with the following recurrence process: 
\begin{subequations}\label{memory_integration}
\begin{equation}
    \mathcal{F}(\theta_t)\approx \mathcal{L}_t := \begin{cases} \mathcal{L}(\theta_0), \qquad &t=0, \\ \alpha_t \mathcal{L}(\theta_t) + (1-\alpha_t) \mathcal{L}_{t-1},\qquad &t\geq 1,\end{cases}
\end{equation} where $\mathcal{L}(\theta_t)$ is the MC estimate at the $t^\text{th}$ iteration,  and $\{\alpha_t\}_{t\geq 0}$ is a selected sequence of coefficients $0<\alpha_t\leq 1$ such that $\alpha_0=1$.  In expanded form,
\begin{equation}\label{memory_expanded}
    \mathcal{L}_t = \sum_{l= 0}^t \alpha_l \left(\prod_{s=1}^{t-l}(1-\alpha_{l+s})\right) \mathcal{L}(\theta_l),
\end{equation}
\end{subequations} which shows that the approximation $\mathcal{L}_t$ of $\mathcal{F}(\theta_t)$ is indeed a linear combination of the current and all previous MC integration estimates. If $\alpha_t=1$ for all $t\geq 0$, we recover the usual MC integration case without memory.

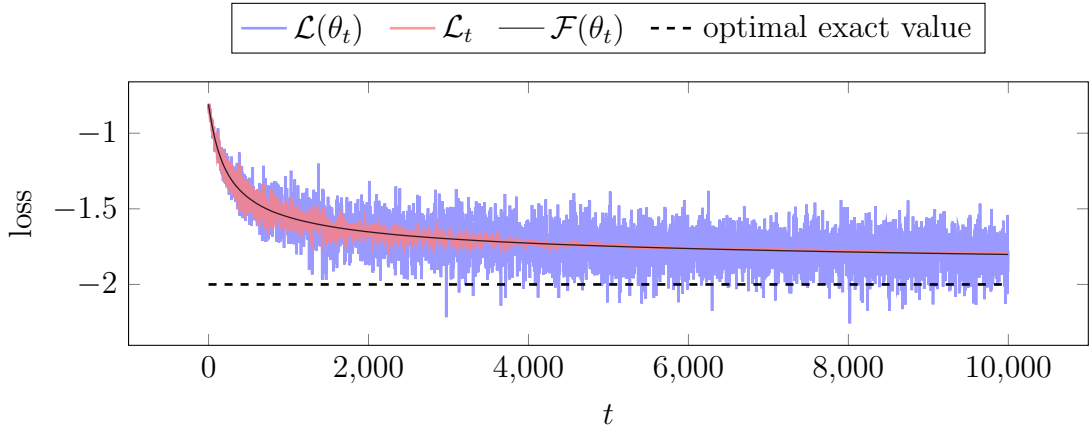
\begin{figure}[htbp]
\centering
\begin{tikzpicture}[]         
 \begin{axis}[
     xlabel = {$t$},
     scaled x ticks=false,
     ylabel = {loss},
     height=0.35*\textwidth,    
     width=0.98\textwidth,    
     xtick={0,2000,4000,...,10000},
     legend style={at={(0.5,1.3)}, anchor=north, /tikz/every even column/.append style={column sep=0.25cm}},
     legend columns = -1,
     yticklabel style={
        /pgf/number format/fixed,
        /pgf/number format/precision=4},
     select coords between index/.style 2 args={
    x filter/.code={
        \ifnum\coordindex<#1\def\pgfmathresult{}\fi
        \ifnum\coordindex>#2\def\pgfmathresult{}\fi
    }
},
     ]        

\addplot[line width=1, color=blue!50!white, opacity=0.8,each nth point=2] table[x expr=\thisrow{iteration},y=numeric]{Chapters/5.Chapter/figures/data/section32/loss.csv};
\addlegendentry{$\mathcal{L}(\theta_t)$};
\addplot[line width=1, color=red!50!white, opacity=0.8,each nth point=2] table[x expr=\thisrow{iteration},y=memory]{Chapters/5.Chapter/figures/data/section32/loss.csv};
\addlegendentry{$\mathcal{L}_t$};
\addplot[line width=.5, color=black, opacity=0.8,each nth point=2] table[x expr=\thisrow{iteration},y=analytic]{Chapters/5.Chapter/figures/data/section32/loss.csv};
\addlegendentry{$\mathcal{F}(\theta_t)$};
\addplot[domain=0:10000,color=black, dashed,line width=1]{-2};
\addlegendentry{optimal exact value};

\end{axis}
\end{tikzpicture}
\caption{Training of the single-trainable-parameter NN presented in \Cref{section5.2.1} whose architecture permits exact calculation of $\mathcal{F}$. The training is performed according to \eqref{numeric_gradients}, and thus, $\mathcal{L}(\theta_t)$ is its associated loss function. $\mathcal{L}_t$ and $\mathcal{F}(\theta_t)$ are computed for monitoring.}
\label{monte_carlo_first_training}
\end{figure}

\Cref{monte_carlo_first_training} (red curve) shows the memory-based loss $\mathcal{L}_t$ evolution along training according to ordinary SGD optimization \eqref{numeric_gradients} and selecting $\alpha_t = e^{-0.001t} + 0.001$.  Initially, we integrate with large errors ($\alpha_t$ is practically one, and therefore,  there is hardly any memory in $\mathcal{L}_t$). However, as we progress in training, we increasingly endow memory to $\mathcal{L}_t$ ($\alpha_t$ becomes small),  and as a consequence,  its integration error decreases.  $\mathcal{L}_t$ produces more accurate approximations than $\mathcal{L}(\theta_t)$,  allowing better convergence monitoring to,  for example,  establish proper stopping criteria during training.

\subsection{Optimization}\label[subsection]{section5.2.3}

While the proposed scheme improves the approximation of $\mathcal{F}$,  we have that a corresponding SGD scheme using \eqref{memory_integration} is equivalent to classical SGD with a learning rate that is multiplied by $\alpha_t$.

We can naturally endow the idea of memory-based integration to the gradients,  as $g(\theta_t)$ is also obtained via MC integration---recall \eqref{practical_gradients},
\begin{subequations}\label{memory-based_optimization}
\begin{equation}\label{memory_gradients}
    \frac{\partial\mathcal{F}}{\partial\theta}(\theta_t)\approx g_t := \begin{cases} g(\theta_0), \qquad &t=0, \\ \gamma_t g(\theta_t) + (1-\gamma_t) g_{t-1},\qquad &t\geq 1,\end{cases} 
\end{equation} where $\{\gamma_t\}_{t\geq 0}$ is a selected sequence of coefficients such that $0<\gamma_t\leq 1$ and $\gamma_0=1$. Then, we obtain a memory-based SGD optimizer employing the $g_t$ term instead of $g(\theta_t)$---recall \eqref{numeric_gradients},
\begin{equation}\label{memory_optimization}
\theta_{t+1} := \theta_{t} - \lambda g_t.
\end{equation} In the expanded form, we have
\begin{equation}\label{SGD_gradient_accumulator}
g_t = \sum_{l=0}^t \gamma_l \left(\prod_{s=1}^{t-l}(1-\gamma_{l+s})\right) g(\theta_l).
\end{equation}
\end{subequations}

\begin{figure}[htbp]
\centering
\begin{tikzpicture}[]         
 \begin{axis}[
     xlabel = {$t$},
     scaled x ticks=false,
     ylabel = {gradient},
     height=0.35*\textwidth,    
     width=0.98\textwidth,    
     xtick={0,2000,4000,...,10000},
     legend style={at={(0.5,1.3)}, anchor=north, /tikz/every even column/.append style={column sep=0.25cm}},
     legend columns = -1,
     yticklabel style={
        /pgf/number format/fixed,
        /pgf/number format/precision=4},
     select coords between index/.style 2 args={
    x filter/.code={
        \ifnum\coordindex<#1\def\pgfmathresult{}\fi
        \ifnum\coordindex>#2\def\pgfmathresult{}\fi
    }
},
     ]        

\addplot[line width=1, color=blue!50!white, opacity=0.8,each nth point=2] table[x expr=\thisrow{iteration},y=numeric]{Chapters/5.Chapter/figures/data/section32/gradient.csv};
\addlegendentry{$g(\theta_t)$}
\addplot[line width=1, color=red!50!white, opacity=0.8,each nth point=2] table[x expr=\thisrow{iteration},y=memory]{Chapters/5.Chapter/figures/data/section32/gradient.csv};
\addlegendentry{$g_t$}
\addplot[line width=.5, color=black, opacity=0.8,each nth point=2] table[x expr=\thisrow{iteration},y=analytic]{Chapters/5.Chapter/figures/data/section32/gradient.csv};
\addlegendentry{$\frac{\partial\mathcal{F}}{\partial\theta}(\theta_t)$}
\addplot[domain=0:10000,color=black, dashed,line width=1]{0};
\addlegendentry{optimal exact value};
                 
\end{axis}
\end{tikzpicture}
\caption{Gradient evolution of the single-parameter model problem in \Cref{monte_carlo_first_training} during training. The optimization is performed according to \eqref{numeric_gradients} using $g(\theta_t)$,  while $g_t$ and $\frac{\partial\mathcal{F}}{\partial\theta}(\theta_t)$ are computed for monitoring. }
\label{monte_carlo_second_training}
\end{figure}
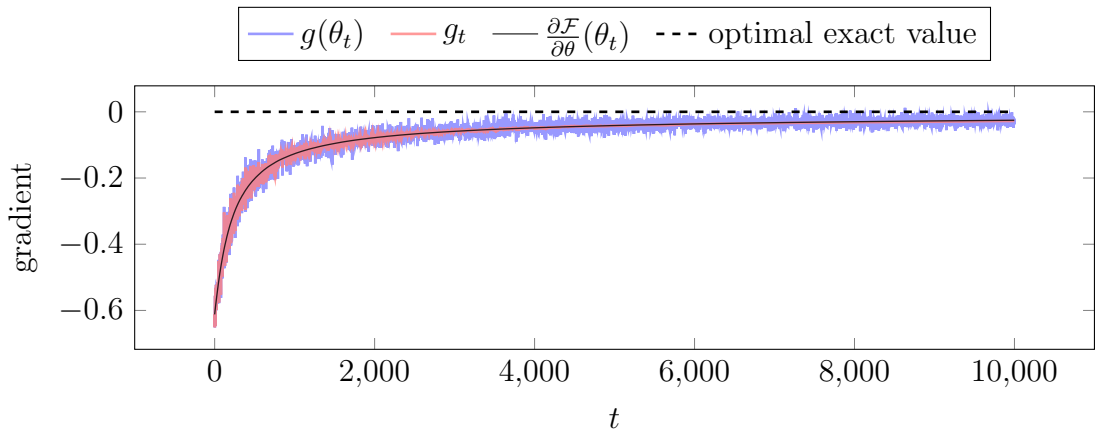

\Cref{monte_carlo_second_training} shows the gradient evolution during training of the previous single-trainable-parameter model problem.   We select $\alpha_t=\gamma_t$ for all $t$,  with the same exponential decay as before.  $g_t$ produces more accurate approximations of the exact gradients than $g(\theta_t)$, minimizing the noise. 

Proper tuning of the coefficients $\alpha_t$ and $\gamma_t$ is critical to maximize integration performance.  Coefficients should be high (low memory) when the involved integrals vary rapidly (e.g.,  at the beginning of training). Conversely, when the approximated solution is near equilibrium and the relevant integrals vary slowly, the coefficients should be low (high memory). 
\section{Relation with the momentum method}\label[section]{section5.4}

The \ac{SGD} optimizer with momentum (\acs{SGDM}) \cite{polyak1964some,sutskever2013momentum} is commonly introduced as the following two-step recursive method:
\begin{subequations}\label{momentum_definition}
\begin{align}
    v_{t+1} &:= \beta v_{t} - g(\theta_t),\\
    \theta_{t+1} &:= \theta_t + \lambda v_{t+1},
\end{align} where $v_t$ is the \emph{momentum accumulator} initialized by $v_0=0$, and $0\leq \beta< 1$ is the momentum coefficient.  If $\beta=0$, we recover the classical \ac{SGD} optimizer \eqref{numeric_gradients}.  Rewriting \eqref{momentum_definition} in terms of the scheme $\theta_{t+1}=\theta_t-\lambda g_t$,  we obtain
\begin{equation}\label{momentum_optimization}
    g_t = g(\theta_t) + \beta g_{t-1} = \sum_{l=0}^t \beta^{t-l} g(\theta_l).
\end{equation}
\end{subequations}

A more sophisticated version of the \acs{SGDM} modifies the momentum coefficient during training (see, e.g., \cite{chen2022demon}),  namely,  defined as in \eqref{momentum_definition} but replacing $\beta$ with $\beta_t\in [0,1)$ for some conveniently selected sequence $\{\beta_t\}_{t\geq 1}$. Then, the $g_t$ term results
\begin{equation}\label{momentum_optimization_beta_variable}
    g_t = g(\theta_t) + \beta_t g_{t-1} = \sum_{l=0}^t \left(\prod_{s=1}^{t-l} \beta_{l+s}\right) g(\theta_l).
\end{equation} Selecting proper hyper-parameters $\beta_t$ during training is challenging.  However,  by readjusting the learning rate and momentum coefficient in the \acs{SGDM} optimizer according to $\gamma_t$ in \eqref{memory-based_optimization} for $t\geq 1$, 
\begin{subequations}
\label{change_variable}
\begin{align}
    \lambda_t &:= \lambda_{t-1}\frac{\gamma_t}{\gamma_{t-1}},\qquad \lambda_0:=\lambda,\\
    \beta_t &:= \gamma_{t-1}\frac{1-\gamma_t}{\gamma_t},
\end{align}
\end{subequations} we recover our memory-based proposal \eqref{memory-based_optimization}.

Both optimizations \eqref{memory-based_optimization} and \eqref{momentum_definition}--\eqref{momentum_optimization_beta_variable} stochastically accumulate gradients to readjust the trainable parameters.  However,  while \eqref{momentum_definition}--\eqref{momentum_optimization_beta_variable} considers a geometrically weighted average of past gradients, our proposal \eqref{memory-based_optimization} re-scales the current and prior gradients so $g_t$ intends to imitate $\frac{\partial\mathcal{F}}{\partial\theta}(\theta_t)$---recall \Cref{monte_carlo_second_training}. 

\eqref{change_variable} reinterprets the \acs{SGDM} as an exact-gradient performer by re-scaling the learning rate.  In contrast,  our memory-based proposal provides the exact-gradient interpretation, leaving the learning rate free.  Consequently,  the learning rate is an independent hyperparameter of the gradient-based optimizer and not an auxiliary element to interpret gradients during training.  Moreover, our optimizer is designed to work in parallel with the memory-based loss \eqref{memory_integration} that approximates the (typically unavailable) objective function during training.

%
%


\chapter{Conclusions and Future Work}\label{chapter6}

One of the main benefits of \acfp{ANN} is their ability to approximate nearly any continuous function to a desired precision. This is achieved through an underlying non-linear parameterization, typically turning the approximation task into a highly non-convex optimization problem. 



In this dissertation, we have used ANNs to numerically solve PDEs. Although the ANN-based framework provides greater approximation capacity compared to traditional methods, it also produces significant drawbacks that hamper the performance during experimentation. Among the most remarkable, we highlight the following two: (a) lack of efficient integration methods for ANNs, especially in the presence of singular solutions where \acf{MC} integration fails; and (b) poor optimizer performance, which often leads to fall in a local minimum whose distance to a global minimum is uncertain.

Below, we collect the most salient aspects of the novel proposals exposed throughout this dissertation and our visions for future work that aim at overcoming the aforementioned limitations.

\section{The Deep Finite Element Method}

We developed a deep-learning-based method to solve PDEs that mimics a FEM setup when applying mesh refinements. Thanks to the FEM-based formulation, the integration is exact: the approximated solution is a linear combination of pre-established basis functions and the ANN only optimizes the involved optimal coefficients. Moreover, thanks to the proposed dynamic architecture, we are able to start from a low-dimensional discretization scheme and step-by-step scale to higher-dimensional spaces (when applying refinements), taking advantage of the solutions found in the previous subspaces. 

The method is capable of solving both parametric and non-parametric PDEs. While solving non-parametric problems with the \ac{DeepFEM} is computationally inefficient (because, essentially, we are unable to beat a traditional linear solver), the execution times when solving parametric problems are comparable to non-parametric ones (i.e., the number of employed iterations is similar), which is where the true power of the method is shown. We developed and implemented the \ac{DeepFEM} in one spatial dimension employing piecewise-linear approximations and uniform refinements. Extending it to more complicated 2D or 3D geometries, with adaptive meshes, and higher-order polynomials is straightforward but challenging in a \acf{TF2} implementation. We leave this as future work.

Since our \ac{DeepFEM} architecture is partially explainable, we identified equivalences between different types of training (end-to-end vs. layer-by-layer), the degrees of freedom arising in the FEM (coefficients of the involved linear combination), and the considered preconditioners and norms. However, the lack of convexity of the loss with respect to the trainable parameters is critical and prevents us from obtaining a robust method.

\section{The Deep Double Ritz Method}

We studied the problem of residual minimization for solving \acp{PDE} using \acp{NN}. Using the dual norm definition as a maximum over the test space, we first obtained a min-max saddle-point problem. This method turned out to be numerically unstable due to a lack of Lipschitz continuity of the test maximizers. To overcome this, we rewrote a general variational problem as a minimization of a Ritz functional employing optimal test functions.  To carry out this minimization while computing the optimal test functions over general problems, we proposed a \acf{D2RM} that combines two nested Ritz optimization loops using ANNs. This novel method constructs local approximations of the trial-to-test operator to express the optimal test functions as dependent on the trial functions.  By doing this, ANNs allowed us to easily combine the two nested Ritz problems in a way that is difficult to treat with traditional numerical methods.

We tested the \ac{D2RM} in several smooth and singular problems. Numerical results illustrate the advantages and limitations of the proposed method. Among the advantages, we encounter the good approximation capabilities of ANNs, the generality of the proposed method, and its application to different linear PDEs, variational formulations, and spatial dimensions. As main limitations, we faced the aforementioned two common difficulties shared by most ANN-based PDE solvers, namely, the lack of efficient integration and optimization. Other than that, numerical results are promising and supported by a solid mathematical framework at the continuum level.

In future work, we will investigate the behavior of optimal test functions to propose enhanced stopping criteria in the inner loop. In particular, we believe that a global approximation of the trial-to-test operator using ANNs would be helpful for this purpose. In parallel, we will investigate reformulations in terms of residuals in replacement of the optimal functions. We believe that its integral formulation involving an integrand that converges to zero is key for stabilizing the challenge of numerical integration when using ANNs.

\section{Memory-based Monte Carlo integration}

With the aim of improving integration accuracy when using ANNs (without incurring prohibitive computational costs), we proposed a memory-based iterative method that conveniently accumulates MC integral estimations in previous iterations when the convergence reaches an equilibrium phase. Appropriately conveying this memory behavior to the gradient computation, we found a reinterpretation of the momentum method in gradient descent optimization.  However, this integration error improvement that makes the loss function resemble the objective function (which similarly occurs with the gradient of the loss towards the gradient of the objective function) is insufficient for enhancing convergence, as we are limited to the inherent lack of convexity of the ANN parameterization. In other words, we are able to overcome the integration challenge but not the optimization one.

\section{Ongoing research}

Our following research steps involve exploiting the interpretation of the ANN as a family of finite-dimensional vector spaces (recall \Cref{appendix2.a}). In this way, we are able to formulate a suitable goal-oriented adaptivity scheme using ANNs and to find an efficient optimization framework for Variational Physics-Informed Neural Networks (VPINNs). These ongoing projects are tentatively titled ``\emph{GOANNs: Goal-Oriented Adaptivity Using Neural Networks}'' and ``\emph{Optimizing (Variational) Physics-Informed Neural Networks Using Least Squares}'', respectively.

\UseRawInputEncoding
\chapter{Achievements} \label{chapter7}

This dissertation contributes to the solution of \acfp{PDE} using \acfp{ANN}, providing three main novel contributions as a result of the successful accomplishment of a series of activities developed throughout the PhD program.

Below, we collect the most relevant research contributions and stays carried out since enrolling in the Mathematics and Statistics PhD program of the University of the Basque Country (UPV/EHU) in October 2019.

\section{Publications}

Peer-reviewed published works.

\begin{description}

\item[2023] {Carlos Uriarte,  Jamie M. Taylor, David Pardo, Oscar A. Rodr\'iguez, Patrick Vega.  \textit{Memory-Based Monte Carlo Integration for Solving Partial Differential Equations Using Neural Networks}. In: Miky\v{s}ka, J., de Mulatier, C., Paszynski, M., Krzhizhanovskaya, V.V., Dongarra, J.J., Sloot, P.M. (eds) Computational Science -- ICCS 2023. Lecture Notes in Computer Science, Vol. 14074, 2023. Springer, Cham. \\\url{https://doi.org/10.1007/978-3-031-36021-3_51}.}

\item[2023]{Carlos Uriarte, David Pardo, Ignacio Muga, Judit Mu\~{n}oz-Matute.
\textit{A Deep Double Ritz Method (D$^2$RM) for solving Partial Differential Equations using Neural Networks}. Computer Methods in Applied Mechanics and Engineering (Q1, Top 1\%), Vol. 405, 115892, 2023.\\\url{https://doi.org/10.1016/j.cma.2023.115892}.}

\item[2022]{Carlos Uriarte,  David Pardo,  \'Angel J. Omella.  \textit{A Finite Element based Deep Learning solver for parametric PDEs}. Computer Methods in Applied Mechanics and Engineering (Q1, Top 1\%), Vol. 391, 114562, 2022.\\\url{https://doi.org/10.1016/j.cma.2021.114562}.}

\end{description}

\section{Conferences}

Contributions to national and international conferences.  Underlined appears the speaker of the presentation.

\begin{description}

\item[2023]{\underline{Carlos Uriarte},  David Pardo, Ignacio Muga, Judit Mu\~{n}oz-Matute.  \textit{The Deep Double Ritz Method: a Deep Learning Residual Minimization Method for solving Partial Differential Equations}. 2nd IACM Mechanistic Machine Learning and Digital Engineering for Computational Science and Technology Conference. University of Texas at El Paso, USA.  September 24-27, 2023.}

\item[2023]{Carlos Uriarte,  Jamie M. Taylor, David Pardo, Oscar A. Rodr\'iguez, \underline{Patrick} \underline{Vega}.  \textit{Memory-based Monte Carlo integration for solving Partial Differential Equations using Neural Networks}.  XXXI COMCA -- Congreso de Matem\'atica Capricornio, Antofagasta, Chile.  August 2-4, 2023.}

\item[2023]{\underline{Carlos Uriarte},  Jamie M. Taylor, David Pardo, Oscar A. Rodr\'iguez,  Patrick Vega.  \textit{Memory-based Monte Carlo integration for solving Partial Differential Equations using Neural Networks}..  ICCS 2023 -- 23rd International Conference on Computational Science. Prague, Czech Republic.  July 3-5, 2023.}

\item[2023]{\underline{Carlos Uriarte},  David Pardo, Ignacio Muga.  \textit{Goal-Oriented Deep Ritz and Least-Squares methods}.  ICCS 2023 -- 23rd International Conference on Computational Science. Prague, Czech Republic.  July 3-5, 2023.}

\item[2022]{\underline{Carlos Uriarte}, David Pardo, and \'{A}ngel J. Omella. \textit{A Finite Element based Deep Learning solver for parametric PDEs}. SIAM MDS22 -- Conference on Mathematics of Data Science. San Diego, California, USA.  September 26-30, 2023.}

\item[2022]{\underline{Carlos Uriarte}, David Pardo, Ignacio Muga, and Judit Mu\~{n}oz-Matute.  \textit{Solving Partial Differential Equations using Adversarial Neural Networks.} CMN 2022 -- Congress on Numerical Methods in Engineering. Las Palmas de Gran Canaria, Spain. September 12-14,2022.}

\item[2022]{\underline{Carlos Uriarte}, David Pardo, Ignacio Muga, and Judit Mu\~{n}oz-Matute. \textit{Adversarial Neural Networks for solving variationally formulated Partial Differential Equations}.  WCCM/APCOM 2022 -- 15th World Congress on Computational Mechanics and 8th Asian Pacific Congress on Computational Mechanics.  Yokohama, Japan.  July 31 -- August 5, 2022.}

\item[2022]{\underline{Carlos Uriarte}, David Pardo, Ignacio Muga, and Judit Mu\~{n}oz-Matute. \textit{An Adversarial Networks approach for solving Partial Differential Equations}.  ICCS 2022 -- 22nd International Conference on Computational Science. London, UK.  June 21-23, 2022.}

\item[2022]{\underline{Carlos Uriarte}, David Pardo, Ignacio Muga, and Judit Mu\~{n}oz-Matute.  \textit{A Generative Adversarial Networks approach for solving Partial Differential Equations}.  ECCOMAS 2022 -- 8th European Congress on Computational Methods in Applied Sciences and Engineering.  Oslo, Norway.  June 5-9, 2022.}

\item[2021]{\underline{Carlos Uriarte},  \'Angel J. Omella,  David Pardo. \textit{A Finite Element based Deep Learning solver
for parametric PDEs}.  ICCS 2021 -- 21st International Conference on Computational Science. Krakow, Poland.  June 16-18, 2021.}
\end{description}

\section{Courses,  Seminars \& Workshops}\label{Achievements:seminars}

Contributions to courses, seminars, and workshops. Underlined appears the main contributor. Lack of underlining means that there is no specific principal contributor.

\begin{description}
\item[2023]{\underline{Carlos Uriarte}, David Pardo.  \textit{(Goal-Oriented) Deep Residual Minimization Methods}.  Seminar at the Faculty of Computer Science, Electronics, and Telecommunications of the AGH University of Science and Technology.  Krakow, Poland.  July 11,  2023.}

\item[2023]{Tom\'as Teijeiro,  \'Angel J. Omella,  Jamie M. Taylor,  Carlos Uriarte,  David Pardo.  \textit{Parametric PDEs using Deep Learning}.  Organizer of the Course and Working Group.  Asturias, Spain.  May 21-27, 2023.}

\item[2023]{\underline{Carlos Uriarte}, David Pardo, Ignacio Muga,  Judit Mu\~{n}oz-Matute. \textit{A Deep Double Ritz Method for solving Partial Differential Equations using Neural Networks}.  Oral presentation at the Workshop Numerical Methods in Geophysics: Present, Future, and Applications held in Valpara\'iso, Chile.  January 12-13, 2023.}

\item[2023]{\'Angel J. Omella, Carlos Uriarte, and David Pardo. \textit{Coding Deep Neural Networks for PDEs}. Professor at the Course organized by the Pontificia Universidad Cat\'olica de Valpara\'iso and held in Olmu\'e, Chile. January 15-20,  2023.}

\item[2022]{\underline{Carlos Uriarte}, David Pardo, Ignacio Muga,  Judit Mu\~{n}oz-Matute.  \textit{A Deep Double Ritz Method for solving Partial Differential Equations}.  Oral presentation at the XC Encuentro Anual de la Sociedad Matem\'atica de Chile in Punta de Tralca, Chile. December 8-10, 2022.}

\item[2022]{\underline{Carlos Uriarte}, David Pardo, Ignacio Muga,  Judit Mu\~{n}oz-Matute. \textit{A Deep Double Ritz Method for solving Partial Differential Equations}.  Seminar at the Pontificia Universidad Cat\'olica de Valpara\'iso, Chile. November 11, 2022.}

\item[2022]{\underline{David Pardo}, Magdalena Strugaru, Jamie M. Taylor,  \'Angel J. Omella,  Jon A. Rivera,  Carlos Uriarte, Ignacio Muga,  Judit Mu\~{n}oz-Matute.  \textit{Deep Learning for Simulation and Inversion Problems}. Plenary talk at the Oden Institute for Computational Engineering and Sciences of the University of Texas at Austin, Texas, USA. October 21,  2022.}

\item[2022]{\underline{David Pardo},  \'Angel J. Omella,  Jamie M. Taylor, Carlos Uriarte, Jon A. Rivera,  Magdalena Strugaru. \textit{Deep Learning for Simulation and Inversion Problems}. Plenary talk at the 9{\textordmasculine} Congreso Metropolitano de Modelado y Simulaci\'on Num\'erica held in Mexico. May 4-6, 2022.}

\item[2021]{David Pardo,  \'Angel J. Omella, Jon Ander River, Carlos Uriarte,  Ana Fern\'adez-Navamuel. \textit{Solving Forward and Inverse Problems with Deep Learning}. Organizer of the Course and Working Group held in Cantabria, Spain.  May 10-18, 2021.}

\item[2021]{\underline{Carlos Uriarte}, David Pardo,  \'Angel J. Omella. \textit{A Finite Element based Deep Learning solver for parametric PDEs}.  Seminar at the Pontificia Universidad Cat\'olica de Valpara\'i­so, Chile. March 26, 2021.}

\item[2021]{\underline{Carlos Uriarte}, David Pardo,  \'Angel J. Omella. \textit{A Finite Element based Deep Learning solver for parametric PDEs}. Lecture at the Faculty of Computer Science, Electronics, and Telecommunications of the AGH University of Science and Technology in Krakow, Poland. January 21, 2021.}

\item[2020]{David Pardo,  \'Angel J. Omella,  Carlos Uriarte. \textit{Deep Learning for Solving Inverse Problems using TF2.0}. Organizer of the Course and Working Group held in Asturias, Spain. June 28 -- July 4,  2020.}
\end{description}

\newpage

\section{Research Stays}

Research stays abroad longer than 30 days.

\begin{description}

\item[2023]{Instituto de Matem\'aticas de la Pontificia Universidad Cat\'olica de Vapara\'iso, Chile.  Supervisor: Prof.  Ignacio Muga (92 days).}

\item[2022]{Oden Institute for Computational Engineering and Sciences, University of Texas at Austin,  USA.  Supervisor: Prof.  Leszek F.  Demkowicz (37 days).}

\end{description}

\section{Disseminating Activities}

Disseminating talks.

\begin{description}

\item[2023]{Elisabete Alberdi, Carlos Uriarte. \textit{Matematika eguneroko bizitzan. Ekuazio diferentzialak mundua azaltzeko}.  Disseminating talk at Uhagon Kulturgunea in Markina-Xemein, Bizkaia, Spain. May 12, 2023.}

\item[2022]{Carlos Uriarte. \textit{Ekuazio diferentzialak mundua azaltzeko eta sare-neuronal artifizialak horiek ebazteko}. Disseminating talk at the 19th Edition of the Mathematics in Everyday Life (Matem\'aticas en la vida cotidiana) lecture series at Bidebarrieta Kulturgunea, Bilbao, Spain. May 19, 2022.}
 
\end{description}

%
%

\cleardoublepage
\phantomsection
\addcontentsline{toc}{part}{Bibliography}
\bibliographystyle{abbrv}
\bibliography{bibliography}

\begin{thebibliography}{100}

\bibitem{aarts2001neural}
L.~P. Aarts and P.~Van Der~Veer.
\newblock {Neural network method for solving partial differential equations}.
\newblock {\em Neural Processing Letters}, 14:261--271, 2001.

\bibitem{abadi2016tensorflow}
M.~Abadi, A.~Agarwal, P.~Barham, E.~Brevdo, Z.~Chen, C.~Citro, G.~S. Corrado,
  A.~Davis, J.~Dean, M.~Devin, et~al.
\newblock {Tensorflow: Large-Scale Machine Learning on Heterogeneous
  Distributed Systems}.
\newblock {\em arXiv preprint arXiv:1603.04467}, 2016.

\bibitem{tensorflow2015-whitepaper}
M.~Abadi, A.~Agarwal, P.~Barham, E.~Brevdo, Z.~Chen, C.~Citro, G.~S. Corrado,
  A.~Davis, J.~Dean, M.~Devin, S.~Ghemawat, I.~Goodfellow, A.~Harp, G.~Irving,
  M.~Isard, Y.~Jia, R.~Jozefowicz, L.~Kaiser, M.~Kudlur, J.~Levenberg,
  D.~Man\'{e}, R.~Monga, S.~Moore, D.~Murray, C.~Olah, M.~Schuster, J.~Shlens,
  B.~Steiner, I.~Sutskever, K.~Talwar, P.~Tucker, V.~Vanhoucke, V.~Vasudevan,
  F.~Vi\'{e}gas, O.~Vinyals, P.~Warden, M.~Wattenberg, M.~Wicke, Y.~Yu, and
  X.~Zheng.
\newblock {TensorFlow: Large-Scale Machine Learning on Heterogeneous Systems},
  2015.
\newblock Software available from tensorflow.org.

\bibitem{abiodun2018state}
O.~I. Abiodun, A.~Jantan, A.~E. Omolara, K.~V. Dada, N.~A. Mohamed, and
  H.~Arshad.
\newblock {State-of-the-art in artificial neural network applications: A
  survey}.
\newblock {\em Heliyon}, 4(11), 2018.

\bibitem{abur1988parallel}
A.~Abur.
\newblock {A parallel scheme for the forward/backward substitutions in solving
  sparse linear equations}.
\newblock {\em IEEE transactions on power systems}, 3(4):1471--1478, 1988.

\bibitem{ainsworth1997posteriori}
M.~Ainsworth and J.~T. Oden.
\newblock {A posteriori error estimation in finite element analysis}.
\newblock {\em Computer methods in applied mechanics and engineering},
  142(1-2):1--88, 1997.

\bibitem{al2022extensions}
A.~Al-Aradi, A.~Correia, G.~Jardim, D.~de~Freitas~Naiff, and Y.~Saporito.
\newblock {Extensions of the deep Galerkin method}.
\newblock {\em Applied Mathematics and Computation}, 430:127287, 2022.

\bibitem{alnaes2015fenics}
M.~Aln{\ae}s, J.~Blechta, J.~Hake, A.~Johansson, B.~Kehlet, A.~Logg,
  C.~Richardson, J.~Ring, M.~E. Rognes, and G.~N. Wells.
\newblock {The FEniCS project version 1.5}.
\newblock {\em Archive of numerical software}, 3(100), 2015.

\bibitem{alvarez2015hp}
J.~Alvarez-Aramberri.
\newblock {\em {hp-Adaptive Simulation and Inversion of Magnetotelluric
  Measurements}}.
\newblock PhD thesis, University of the Basque Country (UPV/EHU) and University
  of Pau (UPPA), 2015.

\bibitem{alvarez-aramberri2013inversion}
J.~Alvarez-Aramberri, D.~Pardo, and H.~Barucq.
\newblock {Inversion of Magnetotelluric Measurements Using Multigoal Oriented
  hp-adaptivity}.
\newblock {\em Procedia Computer Science}, 18:1564--1573, 2013.

\bibitem{alyaev2021modeling}
S.~Alyaev, M.~Shahriari, D.~Pardo, {\'{A}}.~J. Omella, D.~S. Larsen, N.~Jahani,
  and E.~Suter.
\newblock {Modeling extra-deep electromagnetic logs using a deep neural
  network}.
\newblock {\em GEOPHYSICS}, 86(3):E269--E281, 2021.

\bibitem{amari1967theory}
S.~Amari.
\newblock {A Theory of Adaptive Pattern Classifiers}.
\newblock {\em IEEE Transactions on Electronic Computers}, EC-16(3):299--307,
  1967.

\bibitem{aminataei2008numerical}
A.~Aminataei and M.~M. Mazarei.
\newblock {Numerical solution of Poisson’s equation using radial basis
  function networks on the polar coordinate}.
\newblock {\em Computers \& Mathematics with Applications}, 56(11):2887--2895,
  2008.

\bibitem{anderson1995computational}
J.~D. Anderson and J.~Wendt.
\newblock {\em {Computational Fluid Dynamics}}, volume 206.
\newblock Springer, 1995.

\bibitem{arndt2005approximating}
O.~Arndt, T.~Barth, B.~Freisleben, and M.~Grauer.
\newblock {Approximating a finite element model by neural network prediction
  for facility optimization in groundwater engineering}.
\newblock {\em European journal of operational research}, 166(3):769--781,
  2005.

\bibitem{aster2019parameter}
R.~C. Aster, B.~Borchers, and C.~H. Thurber.
\newblock {\em {Parameter Estimation and Inverse Problems}}.
\newblock Elsevier, 2019.

\bibitem{babuvska1971error}
I.~Babu{\v{s}}ka.
\newblock {Error-bounds for finite element method}.
\newblock {\em Numerische Mathematik}, 16(4):322--333, 1971.

\bibitem{babuvvska1978error}
I.~Babuv{\v{s}}ka and W.~C. Rheinboldt.
\newblock {Error estimates for adaptive finite element computations}.
\newblock {\em SIAM Journal on Numerical Analysis}, 15(4):736--754, 1978.

\bibitem{bach2017breaking}
F.~Bach.
\newblock {Breaking the curse of dimensionality with convex neural networks}.
\newblock {\em The Journal of Machine Learning Research}, 18(1):629--681, 2017.

\bibitem{bao2020numerical}
G.~Bao, X.~Ye, Y.~Zang, and H.~Zhou.
\newblock {Numerical solution of inverse problems by weak adversarial
  networks}.
\newblock {\em Inverse Problems}, 36(11):115003, 2020.

\bibitem{barredoarrieta2020explainable}
A.~{Barredo Arrieta}, N.~D{\'{i}}az-Rodr{\'{i}}guez, J.~{Del Ser}, A.~Bennetot,
  S.~Tabik, A.~Barbado, S.~Garcia, S.~Gil-Lopez, D.~Molina, R.~Benjamins,
  R.~Chatila, and F.~Herrera.
\newblock {Explainable Artificial Intelligence (XAI): Concepts, taxonomies,
  opportunities and challenges toward responsible AI}.
\newblock {\em Information Fusion}, 58:82--115, 2020.

\bibitem{barron1993universal}
A.~R. Barron.
\newblock {Universal approximation bounds for superpositions of a sigmoidal
  function}.
\newblock {\em IEEE Transactions on Information theory}, 39(3):930--945, 1993.

\bibitem{batchelor1967introduction}
G.~K. Batchelor.
\newblock {\em {An Introduction to Fluid Dynamics}}.
\newblock Cambridge University Press, 1967.

\bibitem{bauer2019deep}
B.~Bauer and M.~Kohler.
\newblock {On Deep Learning as a remedy for the curse of dimensionality in
  nonparamtric regression}.
\newblock {\em The Annals of Statistics}, 47(4):2261--2285, 2019.

\bibitem{baur1983complexity}
W.~Baur and V.~Strassen.
\newblock {The complexity of partial derivatives}.
\newblock {\em Theoretical computer science}, 22(3):317--330, 1983.

\bibitem{baxter1996financial}
M.~Baxter and A.~Rennie.
\newblock {\em {Financial Calculus: An Introduction to Derivative Pricing}}.
\newblock Cambridge University Press, 1996.

\bibitem{baydin2018automatic}
A.~G. Baydin, B.~A. Pearlmutter, A.~A. Radul, and J.~M. Siskind.
\newblock {Automatic differentiation in machine learning: A survey}.
\newblock {\em Journal of Marchine Learning Research}, 18:1--43, 2018.

\bibitem{beidokhti2009solving}
R.~S. Beidokhti and A.~Malek.
\newblock {Solving initial-boundary value problems for systems of partial
  differential equations using neural networks and optimization techniques}.
\newblock {\em Journal of the Franklin Institute}, 346(9):898--913, 2009.

\bibitem{bellman1957dynamic}
R.~E. Bellman.
\newblock {Dynamic programming}.
\newblock {\em New Jersey Google Scholar}, pages 24--73, 1957.

\bibitem{bellman1961adaptive}
R.~E. Bellman.
\newblock {\em {Adaptive Control Processes: A Guided Tour}}.
\newblock Princeton University Press, 1961.

\bibitem{beltzer2003neural}
A.~I. Beltzer and T.~Sato.
\newblock {Neural classification of finite elements}.
\newblock {\em Computers \& structures}, 81(24-25):2331--2335, 2003.

\bibitem{berger1977nonlinearity}
M.~S. Berger.
\newblock {\em {Nonlinearity and functional analysis: lectures on nonlinear
  problems in mathematical analysis}}, volume~74.
\newblock Academic press, 1977.

\bibitem{bertero2020introduction}
M.~Bertero and P.~Boccacci.
\newblock {\em {Introduction to Inverse Problems in Imaging}}.
\newblock CRC Press, 2020.

\bibitem{bhattacharya2021model}
K.~Bhattacharya, B.~Hosseini, N.~B. Kovachki, and A.~M. Stuart.
\newblock {Model reduction and neural networks for parametric PDEs}.
\newblock {\em The SMAI journal of computational mathematics}, 7:121--157,
  2021.

\bibitem{blechschmidt2021three}
J.~Blechschmidt and O.~G. Ernst.
\newblock {Three ways to solve partial differential equations with neural
  networks—A review}.
\newblock {\em GAMM-Mitteilungen}, 44(2):e202100006, 2021.

\bibitem{boyd2001chebyshev}
J.~P. Boyd.
\newblock {\em {Chebyshev and Fourier spectral methods}}.
\newblock Courier Corporation, 2001.

\bibitem{bramble2019multigrid}
J.~H. Bramble.
\newblock {\em {Multigrid methods}}.
\newblock Chapman and Hall/CRC, jan 2019.

\bibitem{brauer2011mathematical}
F.~Brauer and C.~Castillo-Chavez.
\newblock {\em {Mathematical Models in Population Biology and Epidemiology}}.
\newblock Texts in Applied Mathematics. Springer New York, 2011.

\bibitem{brenner2008mathematical}
S.~C. Brenner.
\newblock {\em {The mathematical theory of finite element methods}}.
\newblock Springer, 2008.

\bibitem{brevis2022neural}
I.~Brevis, I.~Muga, and K.~G. van~der Zee.
\newblock {Neural control of discrete weak formulations: Galerkin, least
  squares \& minimal-residual methods with quasi-optimal weights}.
\newblock {\em Computer Methods in Applied Mechanics and Engineering},
  402:115716, 2022.

\bibitem{brunken2019parametrized}
J.~Brunken, K.~Smetana, and K.~Urban.
\newblock {(Parametrized) First Order Transport Equations: Realization of
  Optimally Stable Petrov--Galerkin Methods}.
\newblock {\em SIAM Journal on Scientific Computing}, 41(1):A592--A621, 2019.

\bibitem{cai2021physics2}
S.~Cai, Z.~Mao, Z.~Wang, M.~Yin, and G.~E. Karniadakis.
\newblock {Physics-informed neural networks (PINNs) for fluid mechanics: A
  review}.
\newblock {\em Acta Mechanica Sinica}, 37(12):1727--1738, 2021.

\bibitem{cai2021physics1}
S.~Cai, Z.~Wang, S.~Wang, P.~Perdikaris, and G.~E. Karniadakis.
\newblock {Physics-informed neural networks for heat transfer problems}.
\newblock {\em Journal of Heat Transfer}, 143(6):060801, 2021.

\bibitem{cai2021least}
Z.~Cai, J.~Chen, and M.~Liu.
\newblock {Least-squares ReLU neural network (LSNN) method for linear
  advection-reaction equation}.
\newblock {\em Journal of Computational Physics}, 443:110514, 2021.

\bibitem{cai2020deep}
Z.~Cai, J.~Chen, M.~Liu, and X.~Liu.
\newblock {Deep least-squares methods: An unsupervised learning-based numerical
  method for solving elliptic PDEs}.
\newblock {\em Journal of Computational Physics}, 420:109707, 2020.

\bibitem{canuto2007spectral}
C.~Canuto, M.~Y. Hussaini, A.~Quarteroni, and T.~A. Zang.
\newblock {\em {Spectral methods: evolution to complex geometries and
  applications to fluid dynamics}}.
\newblock Springer Science \& Business Media, 2007.

\bibitem{chen2011numerical}
H.~Chen, L.~Kong, and W.-J. Leng.
\newblock {Numerical solution of PDEs via integrated radial basis function
  networks with adaptive training algorithm}.
\newblock {\em Applied Soft Computing}, 11(1):855--860, 2011.

\bibitem{chen2021quasi}
J.~Chen, R.~Du, P.~Li, and L.~Lyu.
\newblock Quasi-{M}onte {C}arlo sampling for solving partial differential
  equations by deep neural networks.
\newblock {\em Numerical Mathematics. Theory, Methods and Applications},
  14(2):377--404, 2021.

\bibitem{chen2020comparison}
J.~Chen, R.~Du, and K.~Wu.
\newblock {A comparison study of deep Galerkin method and deep Ritz method for
  elliptic problems with different boundary conditions}.
\newblock {\em arXiv preprint arXiv:2005.04554}, 2020.

\bibitem{chen2022demon}
J.~Chen, C.~Wolfe, Z.~Li, and A.~Kyrillidis.
\newblock {Demon: Improved neural network training with momentum decay}.
\newblock In {\em ICASSP 2022-2022 IEEE International Conference on Acoustics,
  Speech and Signal Processing (ICASSP)}, pages 3958--3962. IEEE, 2022.

\bibitem{claeskens2008model}
G.~Claeskens, N.~L. Hjort, et~al.
\newblock {\em {Model selection and model averaging}}, volume 330.
\newblock Cambridge University Press Cambridge, 2008.

\bibitem{collobert2008unified}
R.~Collobert and J.~Weston.
\newblock {A unified architecture for natural language processing: Deep neural
  networks with multitask learning}.
\newblock In {\em Proceedings of the 25th international conference on Machine
  learning}, pages 160--167, 2008.

\bibitem{cuomo2022scientific}
S.~Cuomo, V.~S. Di~Cola, F.~Giampaolo, G.~Rozza, M.~Raissi, and F.~Piccialli.
\newblock {Scientific machine learning through physics-informed neural
  networks: Where we are and what’s next}.
\newblock {\em Journal of Scientific Computing}, 92(3):88, 2022.

\bibitem{cybenko1989approximation}
G.~Cybenko.
\newblock {Approximation by superpositions of a sigmoidal function}.
\newblock {\em Mathematics of control, signals and systems}, 2(4):303--314,
  1989.

\bibitem{demkowicz2011class}
L.~Demkowicz and J.~Gopalakrishnan.
\newblock {A class of discontinuous Petrov--Galerkin methods. II. Optimal test
  functions}.
\newblock {\em Numerical Methods for Partial Differential Equations},
  27(1):70--105, 2011.

\bibitem{demkowicz2020dpg}
L.~Demkowicz, J.~Gopalakrishnan, and B.~Keith.
\newblock {The DPG-star method}.
\newblock {\em Computers \& Mathematics with Applications}, 79(11):3092--3116,
  2020.

\bibitem{demkowicz2000rham}
L.~Demkowicz, P.~Monk, L.~Vardapetyan, and W.~Rachowicz.
\newblock {deRham Diagram for $hp$ Finite Element Spaces}.
\newblock {\em Computers \& Mathematics with Applications}, 39(7-8):29--38,
  2000.

\bibitem{demkowicz2014overview}
L.~F. Demkowicz and J.~Gopalakrishnan.
\newblock {An overview of the discontinuous Petrov Galerkin method}.
\newblock {\em Recent Developments in Discontinuous Galerkin Finite Element
  Methods for Partial Differential Equations: 2012 John H Barrett Memorial
  Lectures}, pages 149--180, 2014.

\bibitem{deng2009imagenet}
J.~Deng, W.~Dong, R.~Socher, L.-J. Li, K.~Li, and L.~Fei-Fei.
\newblock {Imagenet: A large-scale hierarchical image database}.
\newblock In {\em 2009 IEEE conference on computer vision and pattern
  recognition}, pages 248--255. Ieee, 2009.

\bibitem{deng2003pillar}
J.~Deng, Z.~Yue, L.~Tham, and H.~Zhu.
\newblock {Pillar design by combining finite element methods, neural networks
  and reliability: a case study of the Feng Huangshan copper mine, China}.
\newblock {\em International Journal of Rock Mechanics and Mining Sciences},
  40(4):585--599, 2003.

\bibitem{dennis1977quasi}
J.~E. Dennis, Jr and J.~J. Mor{\'e}.
\newblock {Quasi-Newton methods, motivation and theory}.
\newblock {\em SIAM review}, 19(1):46--89, 1977.

\bibitem{dick2013high}
J.~Dick, F.~Y. Kuo, and I.~H. Sloan.
\newblock High-dimensional integration: the quasi-{M}onte {C}arlo way.
\newblock {\em Acta Numerica}, 22:133--288, 2013.

\bibitem{dissanayake1994neural}
M.~W. M.~G. Dissanayake and N.~Phan-Thien.
\newblock {Neural-network-based approximations for solving Partial Differential
  Equations}.
\newblock {\em Communications in Numerical Methods in Engineering},
  10(3):195--201, 1994.

\bibitem{dozat2016incorporating}
T.~Dozat.
\newblock {Incorporating Nesterov Momentum into Adam}, 2016.

\bibitem{duan2022convergence}
C.~Duan, Y.~Jiao, Y.~Lai, D.~Li, J.~Z. Yang, et~al.
\newblock {Convergence Rate Analysis for Deep Ritz Method}.
\newblock {\em Communications in Computational Physics}, 31(4):1020--1048,
  2022.

\bibitem{duchi2011adaptive}
J.~Duchi, E.~Hazan, and Y.~Singer.
\newblock {Adaptive subgradient methods for online learning and stochastic
  optimization}.
\newblock {\em Journal of machine learning research}, 12(7), 2011.

\bibitem{dupond2019thorough}
S.~Dupond.
\newblock {A thorough review on the current advance of neural network
  structures}.
\newblock {\em Annual Reviews in Control}, 14(14):200--230, 2019.

\bibitem{e2017deep}
W.~E and B.~Yu.
\newblock {The Deep Ritz method: A deep learning-based numerical algorithm for
  solving variational problems}, 2017.

\bibitem{e2018deep}
W.~E and B.~Yu.
\newblock {The deep Ritz method: a deep learning-based numerical algorithm for
  solving variational problems}.
\newblock {\em Communications in Mathematics and Statistics}, 6(1):1--12, 2018.

\bibitem{farlow1993partial}
S.~J. Farlow.
\newblock {\em {Partial Differential Equations for Scientists and Engineers}}.
\newblock Courier Corporation, 1993.

\bibitem{fernandez1994application}
A.~A. Fernandez.
\newblock {The application of feedforward artificial neural networks to
  function approximation and the solution of differential equations}.
\newblock Master's thesis, Rice University, 1994.

\bibitem{feynman1986feynman}
R.~P. Feynman.
\newblock {\em {The Feynman Lectures on Physics Vol l}}.
\newblock Narosa, 1986.

\bibitem{gershenfeld1999nature}
N.~A. Gershenfeld.
\newblock {\em {The nature of mathematical modeling}}.
\newblock {Cambridge University Press}, 1999.

\bibitem{girault2012finite}
V.~Girault and P.-A. Raviart.
\newblock {\em {Finite element methods for Navier-Stokes equations: theory and
  algorithms}}, volume~5.
\newblock Springer Science \& Business Media, 2012.

\bibitem{golbabai2009solving}
A.~Golbabai, M.~Mammadov, and S.~Seifollahi.
\newblock {Solving a system of nonlinear integral equations by an RBF network}.
\newblock {\em Computers \& Mathematics with Applications}, 57(10):1651--1658,
  2009.

\bibitem{golbabai2007radial}
A.~Golbabai and S.~Seifollahi.
\newblock {Radial basis function networks in the numerical solution of linear
  integro-differential equations}.
\newblock {\em Applied Mathematics and Computation}, 188(1):427--432, 2007.

\bibitem{goodfellow2016deep}
I.~Goodfellow, Y.~Bengio, and A.~Courville.
\newblock {\em {Deep learning}}.
\newblock MIT press, 2016.

\bibitem{goodfellow2014generative}
I.~Goodfellow, J.~Pouget-Abadie, M.~Mirza, B.~Xu, D.~Warde-Farley, S.~Ozair,
  A.~Courville, and Y.~Bengio.
\newblock {Generative Adversarial Nets}.
\newblock {\em Advances in neural information processing systems}, 27, 2014.

\bibitem{goodfellow2020generative}
I.~Goodfellow, J.~Pouget-Abadie, M.~Mirza, B.~Xu, D.~Warde-Farley, S.~Ozair,
  A.~Courville, and Y.~Bengio.
\newblock {Generative Adversarial Networks}.
\newblock {\em Communications of the ACM}, 63(11):139--144, 2020.

\bibitem{gopalakrishnan2013five}
J.~Gopalakrishnan.
\newblock {Five lectures on DPG methods}.
\newblock {\em arXiv preprint arXiv:1306.0557}, 2013.

\bibitem{gottlieb1977numerical}
D.~Gottlieb and S.~A. Orszag.
\newblock {Numerical analysis of Spectral Methods}.
\newblock {\em Philadelphia, Pennsylvania}, 19103, 1977.

\bibitem{griewank2008evaluating}
A.~Griewank and A.~Walther.
\newblock {\em {Evaluating derivatives: principles and techniques of
  algorithmic differentiation}}.
\newblock SIAM, 2008.

\bibitem{gripenberg2003approximation}
G.~Gripenberg.
\newblock {Approximation by neural networks with a bounded number of nodes at
  each level}.
\newblock {\em Journal of approximation theory}, 122(2):260--266, 2003.

\bibitem{grohs2022deep}
P.~Grohs, A.~Jentzen, and D.~Salimova.
\newblock {Deep neural network approximations for solutions of {PDE}s based on
  {M}onte {C}arlo algorithms}.
\newblock {\em Partial Differential Equations and Applications}, 3(4):Paper No.
  45, 41, 2022.

\bibitem{guliyev2018approximation}
N.~J. Guliyev and V.~E. Ismailov.
\newblock {On the approximation by single hidden layer feedforward neural
  networks with fixed weights}.
\newblock {\em Neural Networks}, 98:296--304, 2018.

\bibitem{gupta2004handbook}
A.~K. Gupta and S.~Nadarajah.
\newblock {\em {Handbook of beta distribution and its applications}}.
\newblock CRC press, 2004.

\bibitem{haberman2012applied}
R.~Haberman.
\newblock {\em {Applied Partial Differential Equations with Fourier series and
  Boundary Value Problems}}.
\newblock Pearson Higher Ed, 2012.

\bibitem{hahn2012heat}
D.~W. Hahn and M.~N. Ozisik.
\newblock {\em {Heat conduction}}.
\newblock John Wiley \& Sons, 2012.

\bibitem{hanin2019universal}
B.~Hanin.
\newblock {Universal function approximation by deep neural nets with bounded
  width and ReLU activations}.
\newblock {\em Mathematics}, 7(10):992, 2019.

\bibitem{hanin2017approximating}
B.~Hanin and M.~Sellke.
\newblock {Approximating continuous functions by ReLU nets of minimal width}.
\newblock {\em arXiv preprint arXiv:1710.11278}, 2017.

\bibitem{harris2020array}
C.~R. Harris, K.~J. Millman, S.~J. Van Der~Walt, R.~Gommers, P.~Virtanen,
  D.~Cournapeau, E.~Wieser, J.~Taylor, S.~Berg, N.~J. Smith, et~al.
\newblock {Array programming with NumPy}.
\newblock {\em Nature}, 585(7825):357--362, 2020.

\bibitem{hawkins2004problem}
D.~M. Hawkins.
\newblock {The problem of overfitting}.
\newblock {\em Journal of chemical information and computer sciences},
  44(1):1--12, 2004.

\bibitem{he2016deep}
K.~He, X.~Zhang, S.~Ren, and J.~Sun.
\newblock {Deep residual learning for image recognition}.
\newblock In {\em Proceedings of the IEEE conference on computer vision and
  pattern recognition}, pages 770--778, 2016.

\bibitem{he2000multilayer}
S.~He, K.~Reif, and R.~Unbehauen.
\newblock {Multilayer neural networks for solving a class of partial
  differential equations}.
\newblock {\em Neural networks}, 13(3):385--396, 2000.

\bibitem{hinton2012neural}
G.~Hinton, N.~Srivastava, and K.~Swersky.
\newblock {Neural networks for machine learning—Lecture 6e—RMSprop: Divide
  the gradient by a running average of its recent magnitude}, 2012.
\newblock (Last accesed: August 2023).

\bibitem{holmes1994partial}
E.~E. Holmes, M.~A. Lewis, J.~Banks, and R.~Veit.
\newblock {Partial Differential Equations in Ecology: Spatial Interactions and
  Population Dynamics}.
\newblock {\em Ecology}, 75(1):17--29, 1994.

\bibitem{hornik1991approximation}
K.~Hornik.
\newblock {Approximation capabilities of multilayer feedforward networks}.
\newblock {\em Neural networks}, 4(2):251--257, 1991.

\bibitem{hornik1989multilayer}
K.~Hornik, M.~Stinchcombe, and H.~White.
\newblock {Multilayer feedforward networks are universal approximators}.
\newblock {\em Neural networks}, 2(5):359--366, 1989.

\bibitem{huang2022applications}
B.~Huang and J.~Wang.
\newblock {Applications of physics-informed neural networks in power
  systems---A review}.
\newblock {\em IEEE Transactions on Power Systems}, 38(1):572--588, 2022.

\bibitem{hughes1987finite}
T.~Hughes.
\newblock {The Finite Element Analysis}, 1987.

\bibitem{hull2015options}
J.~Hull.
\newblock {\em {Options, Futures, and Other Derivatives}}.
\newblock Pearson, 2015.

\bibitem{incropera1996fundamentals}
F.~P. Incropera, D.~P. DeWitt, T.~L. Bergman, A.~S. Lavine, et~al.
\newblock {\em {Fundamentals of heat and mass transfer}}, volume~6.
\newblock Wiley New York, 1996.

\bibitem{isakov2006inverse}
V.~Isakov.
\newblock {\em {Inverse problems for Partial Differential Equations}}, volume
  127.
\newblock Springer, 2006.

\bibitem{iturraran-viveros2014artificial}
U.~Iturrar{\'{a}}n-Viveros and J.~O. Parra.
\newblock {Artificial Neural Networks applied to estimate permeability,
  porosity and intrinsic attenuation using seismic attributes and well-log
  data}.
\newblock {\em Journal of Applied Geophysics}, 107:45--54, 2014.

\bibitem{ivakhnenko1967cybernetics}
A.~Ivakhnenko and V.~Lapa.
\newblock {\em {Cybernetics and Forecasting Techniques}}.
\newblock Modern analytic and computational methods in science and mathematics.
  American Elsevier Publishing Company, 1967.

\bibitem{jackson1999classical}
J.~D. Jackson.
\newblock {Classical electrodynamics}, 1999.

\bibitem{ameya2020conservative}
A.~D. Jagtap, E.~Kharazmi, and G.~E. Karniadakis.
\newblock {Conservative physics-informed neural networks on discrete domains
  for conservation laws: Applications to forward and inverse problems}.
\newblock {\em Computer Methods in Applied Mechanics and Engineering},
  365:113028, 2020.

\bibitem{jianyu2002numerical}
L.~Jianyu, L.~Siwei, Q.~Yingjian, and H.~Yaping.
\newblock {Numerical solution of differential equations by radial basis
  function neural networks}.
\newblock In {\em Proceedings of the 2002 International Joint Conference on
  Neural Networks. IJCNN'02 (Cat. No. 02CH37290)}, volume~1, pages 773--777.
  IEEE, 2002.

\bibitem{jianyu2003numerical}
L.~Jianyu, L.~Siwei, Q.~Yingjian, and H.~Yaping.
\newblock {Numerical solution of elliptic partial differential equation using
  radial basis function neural networks}.
\newblock {\em Neural Networks}, 16(5-6):729--734, 2003.

\bibitem{jilani2009adaptive}
H.~Jilani, A.~Bahreininejad, and M.~Ahmadi.
\newblock {Adaptive finite element mesh triangulation using self-organizing
  neural networks}.
\newblock {\em Advances in Engineering Software}, 40(11):1097--1103, 2009.

\bibitem{jin2021nsfnets}
X.~Jin, S.~Cai, H.~Li, and G.~E. Karniadakis.
\newblock {NSFnets (Navier-Stokes flow nets): Physics-informed neural networks
  for the incompressible Navier-Stokes equations}.
\newblock {\em Journal of Computational Physics}, 426:109951, 2021.

\bibitem{kachurovskii1960monotone}
R.~Kachurovskii.
\newblock {Monotone operators and convex functionals}.
\newblock {\em Uspekhi Matematicheskikh Nauk}, 15(4):213--215, 1960.

\bibitem{kansa2004volumetric}
E.~Kansa, H.~Power, G.~Fasshauer, and L.~Ling.
\newblock {A volumetric integral radial basis function method for
  time-dependent partial differential equations. I. Formulation}.
\newblock {\em Engineering Analysis with Boundary Elements}, 28(10):1191--1206,
  2004.

\bibitem{karniadakis2021physics}
G.~E. Karniadakis, I.~G. Kevrekidis, L.~Lu, P.~Perdikaris, S.~Wang, and
  L.~Yang.
\newblock {Physics-informed machine learning}.
\newblock {\em Nature Reviews Physics}, 3(6):422--440, 2021.

\bibitem{kharazmi2019variational}
E.~Kharazmi, Z.~Zhang, and G.~E. Karniadakis.
\newblock {VPINNs: Variational Physics-Informed Neural Networks for solving
  Partial Differential Equations}.
\newblock {\em arXiv preprint arXiv:1912.00873}, 2019.

\bibitem{kharazmi2021hp}
E.~Kharazmi, Z.~Zhang, and G.~E. Karniadakis.
\newblock {hp-VPINNs: Variational physics-informed neural networks with domain
  decomposition}.
\newblock {\em Computer Methods in Applied Mechanics and Engineering},
  374:113547, 2021.

\bibitem{khoo2021solving}
Y.~Khoo, J.~Lu, and L.~Ying.
\newblock {Solving parametric PDE problems with artificial neural networks}.
\newblock {\em European Journal of Applied Mathematics}, 32(3):421--435, 2021.

\bibitem{kidger2020universal}
P.~Kidger and T.~Lyons.
\newblock {Universal approximation with deep narrow networks}.
\newblock In {\em Conference on learning theory}, pages 2306--2327. PMLR, 2020.

\bibitem{kingma2014adam}
D.~P. Kingma and J.~Ba.
\newblock {Adam: A method for stochastic optimization}.
\newblock {\em arXiv preprint arXiv:1412.6980}, 2014.

\bibitem{koroglu2010comparison}
S.~Koroglu, P.~Sergeant, and N.~Umurkan.
\newblock {Comparison of analytical, finite element and neural network methods
  to study magnetic shielding}.
\newblock {\em Simulation Modelling Practice and Theory}, 18(2):206--216, 2010.

\bibitem{kumar2011multilayer}
M.~Kumar and N.~Yadav.
\newblock {Multilayer perceptrons and radial basis function neural network
  methods for the solution of differential equations: A survey}.
\newblock {\em Computers \& Mathematics with Applications}, 62(10):3796--3811,
  2011.

\bibitem{kutyniok2022theoretical}
G.~Kutyniok, P.~Petersen, M.~Raslan, and R.~Schneider.
\newblock {A theoretical analysis of deep neural networks and parametric PDEs}.
\newblock {\em Constructive Approximation}, 55(1):73--125, 2022.

\bibitem{LaTorre2015inverse}
D.~{La Torre}, H.~Kunze, F.~Mendivil, M.~{Ruiz Galan}, and R.~Zaki.
\newblock {Inverse Problems: Theory and Application to Science and Engineering
  2015}.
\newblock {\em Mathematical Problems in Engineering}, 2015:1--3, 2015.

\bibitem{lagaris1998artificial}
I.~E. Lagaris, A.~Likas, and D.~I. Fotiadis.
\newblock {Artificial neural networks for solving ordinary and partial
  differential equations}.
\newblock {\em IEEE Transactions on Neural Networks}, 9(5):987--1000, 1998.

\bibitem{lagaris2000neural}
I.~E. Lagaris, A.~C. Likas, and D.~G. Papageorgiou.
\newblock {Neural-network methods for boundary value problems with irregular
  boundaries}.
\newblock {\em IEEE Transactions on Neural Networks}, 11(5):1041--1049, 2000.

\bibitem{lawal2022physics}
Z.~K. Lawal, H.~Yassin, D.~T.~C. Lai, and A.~Che~Idris.
\newblock {Physics-informed neural network (PINN) evolution and beyond: a
  systematic literature review and bibliometric analysis}.
\newblock {\em Big Data and Cognitive Computing}, 6(4):140, 2022.

\bibitem{lee1990neural}
H.~Lee and I.~S. Kang.
\newblock {Neural algorithm for solving differential equations}.
\newblock {\em Journal of Computational Physics}, 91(1):110--131, 1990.

\bibitem{leobacher2014introduction}
G.~Leobacher and F.~Pillichshammer.
\newblock {\em {Introduction to quasi-Monte Carlo integration and
  applications}}.
\newblock Springer, 2014.

\bibitem{leshno1993multilayer}
M.~Leshno, V.~Y. Lin, A.~Pinkus, and S.~Schocken.
\newblock {Multilayer feedforward networks with a nonpolynomial activation
  function can approximate any function}.
\newblock {\em Neural networks}, 6(6):861--867, 1993.

\bibitem{leveque2007finite}
R.~J. LeVeque.
\newblock {\em {Finite Difference Methods for Ordinary and Partial Differential
  Equations: Steady-State and Time-Dependent Problems}}.
\newblock SIAM, 2007.

\bibitem{li2022deep}
J.~Li, J.~Yue, W.~Zhang, and W.~Duan.
\newblock {The deep learning Galerkin method for the general Stokes equations}.
\newblock {\em Journal of Scientific Computing}, 93(1):5, 2022.

\bibitem{li2015constructing}
X.~Li and X.~Wu.
\newblock {Constructing long short-term memory based deep recurrent neural
  networks for large vocabulary speech recognition}.
\newblock In {\em 2015 ieee international conference on acoustics, speech and
  signal processing (icassp)}, pages 4520--4524. IEEE, 2015.

\bibitem{li2021fourier}
Z.~Li, N.~Kovachki, K.~Azizzadenesheli, B.~Liu, K.~Bhattacharya, A.~Stuart, and
  A.~Anandkumar.
\newblock {Fourier Neural Operator for Parametric Partial Differential
  Equations}, 2021.

\bibitem{li2020multipole}
Z.~Li, N.~Kovachki, K.~Azizzadenesheli, B.~Liu, A.~Stuart, K.~Bhattacharya, and
  A.~Anandkumar.
\newblock {Multipole graph neural operator for parametric partial differential
  equations}.
\newblock {\em Advances in Neural Information Processing Systems},
  33:6755--6766, 2020.

\bibitem{ming2021deep}
Y.~Liao and P.~Ming.
\newblock {Deep Nitsche Method: Deep Ritz Method with Essential Boundary
  Conditions}.
\newblock {\em Communications in Computational Physics}, 29(5):1365--1384,
  2021.

\bibitem{lin2018resnet}
H.~Lin and S.~Jegelka.
\newblock {ResNet with one-neuron hidden layers is a universal approximator}.
\newblock {\em Advances in neural information processing systems}, 31, 2018.

\bibitem{liu2022adaptive}
M.~Liu, Z.~Cai, and J.~Chen.
\newblock {Adaptive two-layer ReLU neural network: I. Best least-squares
  approximation}.
\newblock {\em Computers {\&} Mathematics with Applications}, 113:34--44, 2022.

\bibitem{logan2010first}
D.~L. Logan.
\newblock {\em {A First Course in the Finite Element Method}}.
\newblock Cengage Learning, 4th edition, 2010.

\bibitem{lu2021priori}
Y.~Lu, J.~Lu, and M.~Wang.
\newblock {A priori generalization analysis of the deep Ritz method for solving
  high dimensional elliptic partial differential equations}.
\newblock In {\em Conference on learning theory}, pages 3196--3241. PMLR, 2021.

\bibitem{lu2017expressive}
Z.~Lu, H.~Pu, F.~Wang, Z.~Hu, and L.~Wang.
\newblock {The expressive power of neural networks: A view from the width}.
\newblock {\em Advances in neural information processing systems}, 30, 2017.

\bibitem{lysaker2003noise}
M.~Lysaker, A.~Lundervold, and X.-C. Tai.
\newblock {Noise removal using fourth-order partial differential equation with
  applications to medical magnetic resonance images in space and time}.
\newblock {\em IEEE Transactions on Image Processing}, 12(12):1579--1590, 2003.

\bibitem{mai2001numerical}
N.~Mai-Duy and T.~Tran-Cong.
\newblock {Numerical solution of differential equations using multiquadric
  radial basis function networks}.
\newblock {\em Neural networks}, 14(2):185--199, 2001.

\bibitem{mai2002mesh}
N.~Mai-Duy and T.~Tran-Cong.
\newblock {Mesh-free radial basis function network methods with domain
  decomposition for approximation of functions and numerical solution of
  Poisson's equations}.
\newblock {\em Engineering Analysis with Boundary Elements}, 26(2):133--156,
  2002.

\bibitem{mai2003approximation}
N.~Mai-Duy and T.~Tran-Cong.
\newblock {Approximation of function and its derivatives using radial basis
  function networks}.
\newblock {\em Applied Mathematical Modelling}, 27(3):197--220, 2003.

\bibitem{mai2005solving}
N.~Mai-Duy and T.~Tran-Cong.
\newblock {Solving high-order partial differential equations with indirect
  radial basis function networks}.
\newblock {\em International Journal for Numerical Methods in Engineering},
  63(11):1636--1654, 2005.

\bibitem{maiorov1999lower}
V.~Maiorov and A.~Pinkus.
\newblock {Lower bounds for approximation by MLP neural networks}.
\newblock {\em Neurocomputing}, 25(1-3):81--91, 1999.

\bibitem{malek2006numerical}
A.~Malek and R.~S. Beidokhti.
\newblock {Numerical solution for high order differential equations using a
  hybrid neural network—optimization method}.
\newblock {\em Applied Mathematics and Computation}, 183(1):260--271, 2006.

\bibitem{manevitz2005neural}
L.~Manevitz, A.~Bitar, and D.~Givoli.
\newblock {Neural network time series forecasting of finite-element mesh
  adaptation}.
\newblock {\em Neurocomputing}, 63:447--463, 2005.

\bibitem{mao2020physics}
Z.~Mao, A.~D. Jagtap, and G.~E. Karniadakis.
\newblock {Physics-informed neural networks for high-speed flows}.
\newblock {\em Computer Methods in Applied Mechanics and Engineering},
  360:112789, 2020.

\bibitem{margossian2019review}
C.~C. Margossian.
\newblock {A review of automatic differentiation and its efficient
  implementation}.
\newblock {\em Wiley interdisciplinary reviews: data mining and knowledge
  discovery}, 9(4):e1305, 2019.

\bibitem{mcfall2009artificial}
K.~S. McFall and J.~R. Mahan.
\newblock {Artificial neural network method for solution of boundary value
  problems with exact satisfaction of arbitrary boundary conditions}.
\newblock {\em IEEE Transactions on Neural Networks}, 20(8):1221--1233, 2009.

\bibitem{mcgee1994applications}
D.~L. McGee~Jr.
\newblock {\em {Applications of neural networks to Partial Differential
  Equations}}.
\newblock PhD thesis, The University of Arizona, 1994.

\bibitem{meade1994solution}
A.~J. Meade~Jr and A.~A. Fernandez.
\newblock {Solution of nonlinear ordinary differential equations by feedforward
  neural networks}.
\newblock {\em Mathematical and Computer Modelling}, 20(9):19--44, 1994.

\bibitem{meade1994numerical}
A.~J. Meade~Jr and A.~A. Fernandez.
\newblock {The numerical solution of linear ordinary differential equations by
  feedforward neural networks}.
\newblock {\em Mathematical and Computer Modelling}, 19(12):1--25, 1994.

\bibitem{mikhlin1964variational}
S.~G. Mikhlin and L.~Chambers.
\newblock {\em {Variational methods in mathematical physics}}, volume~50.
\newblock Pergamon Press Oxford, 1964.

\bibitem{milne2000calculus}
L.~M. Milne-Thomson.
\newblock {\em {The calculus of finite differences}}.
\newblock American Mathematical Society, 2000.

\bibitem{misyris2020physics}
G.~S. Misyris, A.~Venzke, and S.~Chatzivasileiadis.
\newblock {Physics-informed neural networks for power systems}.
\newblock In {\em 2020 IEEE Power \& Energy Society General Meeting (PESGM)},
  pages 1--5. IEEE, 2020.

\bibitem{moshfegh2017direct}
J.~Moshfegh and M.~N. Vouvakis.
\newblock {Direct solution of FEM models: Are sparse direct solvers the best
  strategy?}
\newblock In {\em 2017 International Conference on Electromagnetics in Advanced
  Applications (ICEAA)}, pages 1636--1638. IEEE, 2017.

\bibitem{munoz2021equivalence}
J.~Mu{\~n}oz-Matute, D.~Pardo, and L.~Demkowicz.
\newblock {Equivalence between the DPG method and the exponential integrators
  for linear parabolic problems}.
\newblock {\em Journal of Computational Physics}, 429:110016, 2021.

\bibitem{murray2002mathematical}
J.~D. Murray.
\newblock {\em {Mathematical Biology I.: An Introduction}}.
\newblock Springer, 2002.

\bibitem{nareklishvili2023deep}
M.~Nareklishvili, N.~Polson, and V.~Sokolov.
\newblock {Deep Partial Least Squares for Instrumental Variable Regression},
  2023.

\bibitem{nesterov1983method}
Y.~Nesterov.
\newblock {A method for unconstrained convex minimization problem with the rate
  of convergence $\mathcal{O}(1/k^2)$}.
\newblock In {\em Dokl. Akad. Nauk. SSSR}, volume 269, page 543, 1983.

\bibitem{newman1999monte}
M.~E. Newman and G.~T. Barkema.
\newblock {\em {Monte Carlo methods in statistical physics}}.
\newblock Clarendon Press, 1999.

\bibitem{nguyen1999approximation}
T.~Nguyen-Thien and T.~Tran-Cong.
\newblock {Approximation of functions and their derivatives: A neural network
  implementation with applications}.
\newblock {\em Applied Mathematical Modelling}, 23(9):687--704, 1999.

\bibitem{oden2010general}
J.~T. Oden, A.~Hawkins, and S.~Prudhomme.
\newblock {General diffuse-interface theories and an approach to predictive
  tumor growth modeling}.
\newblock {\em Mathematical Models and Methods in Applied Sciences},
  20(03):477--517, 2010.

\bibitem{omella2021sensitivity}
{\'{A}}.~J. Omella, R.~Celorrio, and D.~Pardo.
\newblock {Sensitivity and uncertainty analysis by discontinuous Galerkin of
  lock-in thermography for crack characterization}.
\newblock {\em Computer Methods in Applied Mechanics and Engineering},
  373:113523, 2021.

\bibitem{ongie2020deep}
G.~Ongie, A.~Jalal, C.~A. Metzler, R.~G. Baraniuk, A.~G. Dimakis, and
  R.~Willett.
\newblock {Deep Learning Techniques for Inverse Problems in Imaging}.
\newblock {\em IEEE Journal on Selected Areas in Information Theory},
  1(1):39--56, 2020.

\bibitem{pang2019fpinns}
G.~Pang, L.~Lu, and G.~E. Karniadakis.
\newblock {fPINNs: Fractional physics-informed neural networks}.
\newblock {\em SIAM Journal on Scientific Computing}, 41(4):A2603--A2626, 2019.

\bibitem{pardo2004integration}
D.~Pardo.
\newblock {\em {Integration of hp-adaptivity with a two grid solver:
  applications to electromagnetics}}.
\newblock PhD thesis, University of Texas at Austin, 2004.

\bibitem{park1993approximation}
J.~Park and I.~W. Sandberg.
\newblock {Approximation and radial-basis-function networks}.
\newblock {\em Neural computation}, 5(2):305--316, 1993.

\bibitem{paszke2017automatic}
A.~Paszke, S.~Gross, S.~Chintala, G.~Chanan, E.~Yang, Z.~DeVito, Z.~Lin,
  A.~Desmaison, L.~Antiga, and A.~Lerer.
\newblock {Automatic differentiation in PyTorch}.
\newblock In {\em NIPS-W}, 2017.

\bibitem{paszke2019pytorch}
A.~Paszke, S.~Gross, F.~Massa, A.~Lerer, J.~Bradbury, G.~Chanan, T.~Killeen,
  Z.~Lin, N.~Gimelshein, L.~Antiga, A.~Desmaison, A.~Kopf, E.~Yang, Z.~DeVito,
  M.~Raison, A.~Tejani, S.~Chilamkurthy, B.~Steiner, L.~Fang, J.~Bai, and
  S.~Chintala.
\newblock {PyTorch: An Imperative Style, High-Performance Deep Learning
  Library}.
\newblock In H.~Wallach, H.~Larochelle, A.~Beygelzimer, F.~d\textquotesingle
  Alch\'{e}-Buc, E.~Fox, and R.~Garnett, editors, {\em Advances in Neural
  Information Processing Systems 32}, pages 8024--8035. Curran Associates,
  Inc., 2019.

\bibitem{penwarden2023metalearning}
M.~Penwarden, S.~Zhe, A.~Narayan, and R.~M. Kirby.
\newblock {A metalearning approach for physics-informed neural networks
  (PINNs): Application to parameterized PDEs}.
\newblock {\em Journal of Computational Physics}, 477:111912, 2023.

\bibitem{petersen2021topological}
P.~Petersen, M.~Raslan, and F.~Voigtlaender.
\newblock {Topological properties of the set of functions generated by neural
  networks of fixed size}.
\newblock {\em Foundations of computational mathematics}, 21:375--444, 2021.

\bibitem{pinkus1999approximation}
A.~Pinkus.
\newblock {Approximation theory of the MLP model in neural networks}.
\newblock {\em Acta numerica}, 8:143--195, 1999.

\bibitem{poggio2017and}
T.~Poggio, H.~Mhaskar, L.~Rosasco, B.~Miranda, and Q.~Liao.
\newblock {Why and when can deep---but not shallow---networks avoid the curse
  of dimensionality: A review}.
\newblock {\em International Journal of Automation and Computing},
  14(5):503--519, 2017.

\bibitem{polyak1964some}
B.~T. Polyak.
\newblock {Some methods of speeding up the convergence of iteration methods}.
\newblock {\em USSR Computational Mathematics and Mathematical Physics},
  4(5):1--17, 1964.

\bibitem{press2007numerical}
W.~H. Press.
\newblock {\em {Numerical recipes 3rd edition: The art of scientific
  computing}}.
\newblock Cambridge university press, 2007.

\bibitem{qian1999momentum}
N.~Qian.
\newblock {On the momentum term in gradient descent learning algorithms}.
\newblock {\em Neural networks}, 12(1):145--151, 1999.

\bibitem{raissi2017physics1}
M.~Raissi, P.~Perdikaris, and G.~E. Karniadakis.
\newblock {Physics Informed Deep Learning (Part I): Data-driven Solutions of
  Nonlinear Partial Differential Equations}, 2017.

\bibitem{raissi2017physics2}
M.~Raissi, P.~Perdikaris, and G.~E. Karniadakis.
\newblock {Physics Informed Deep Learning (Part II): Data-driven Discovery of
  Nonlinear Partial Differential Equations}, 2017.

\bibitem{raissi2019physics}
M.~Raissi, P.~Perdikaris, and G.~E. Karniadakis.
\newblock {Physics-informed neural networks: A deep learning framework for
  solving forward and inverse problems involving nonlinear partial differential
  equations}.
\newblock {\em Journal of Computational physics}, 378:686--707, 2019.

\bibitem{ramuhalli2005finite}
P.~Ramuhalli, L.~Udpa, and S.~S. Udpa.
\newblock {Finite-element neural networks for solving differential equations}.
\newblock {\em IEEE Transactions on Neural Networks}, 16(6):1381--1392, 2005.

\bibitem{reddy2019introduction}
J.~N. Reddy.
\newblock {\em {Introduction to the Finite Element Method}}.
\newblock McGraw-Hill Education, 2019.

\bibitem{ritz1909neue}
W.~Ritz.
\newblock {{\"U}ber eine neue Methode zur L{\"o}sung gewisser
  Variationsprobleme der mathematischen Physik.}
\newblock 1909.

\bibitem{rivera2022quadrature}
J.~A. Rivera, J.~M. Taylor, {\'A}.~J. Omella, and D.~Pardo.
\newblock {On quadrature rules for solving Partial Differential Equations using
  Neural Networks}.
\newblock {\em Computer Methods in Applied Mechanics and Engineering},
  393:114710, 2022.

\bibitem{robbins1951stochastic}
H.~Robbins and S.~Monro.
\newblock {A stochastic approximation method}.
\newblock {\em Annals of Mathematical Statistics}, 22:400--407, 1951.

\bibitem{rojas2023robust}
S.~Rojas, P.~Maczuga, J.~Mu{\~n}oz-Matute, D.~Pardo, and M.~Paszynski.
\newblock {Robust Variational Physics-Informed Neural Networks}.
\newblock {\em arXiv preprint arXiv:2308.16910}, 2023.

\bibitem{rojas2016quadrature}
S.~Rojas, I.~Muga, and D.~Pardo.
\newblock {A quadrature-free method for simulation and inversion of 1.5D direct
  current (DC) borehole measurements}.
\newblock {\em Computational Geosciences}, 20(6):1301--1318, 2016.

\bibitem{rosenblatt1958perceptron}
F.~Rosenblatt.
\newblock {The perceptron: a probabilistic model for information storage and
  organization in the brain.}
\newblock {\em Psychological review}, 65(6):386, 1958.

\bibitem{ruder2016overview}
S.~Ruder.
\newblock {An overview of gradient descent optimization algorithms}.
\newblock {\em arXiv preprint arXiv:1609.04747}, 2016.

\bibitem{salsa2016partial}
S.~Salsa.
\newblock {\em {Partial Differential Equations in Action: From Modelling to
  Theory}}, volume~99.
\newblock Springer, 2016.

\bibitem{saltzman2009biomedical}
W.~M. Saltzman.
\newblock {\em {Biomedical engineering: bridging medicine and technology}}.
\newblock Cambridge University Press, 2009.

\bibitem{samek2017explainable}
W.~Samek, T.~Wiegand, and K.-R. M{\"u}ller.
\newblock {Explainable artificial intelligence: Understanding, visualizing and
  interpreting deep learning models}.
\newblock {\em arXiv preprint arXiv:1708.08296}, 2017.

\bibitem{schmidhuber2022annotated}
J.~Schmidhuber.
\newblock {Annotated history of modern AI and Deep Learning}.
\newblock {\em arXiv preprint arXiv:2212.11279}, 2022.

\bibitem{shahriari2020deep}
M.~Shahriari, D.~Pardo, A.~Picon, A.~Galdran, J.~{Del Ser}, and
  C.~Torres-Verd{\'{i}}n.
\newblock {A deep learning approach to the inversion of borehole resistivity
  measurements}.
\newblock {\em Computational Geosciences}, 24(3):971--994, 2020.

\bibitem{shahriari2021error}
M.~Shahriari, D.~Pardo, J.~A. Rivera, C.~Torres‐Verd{\'{i}}n, A.~Picon,
  J.~{Del Ser}, S.~Ossand{\'{o}}n, and V.~M. Calo.
\newblock {Error control and loss functions for the deep learning inversion of
  borehole resistivity measurements}.
\newblock {\em International Journal for Numerical Methods in Engineering},
  122(6):1629--1657, 2021.

\bibitem{shang2022deep}
Y.~Shang, F.~Wang, and J.~Sun.
\newblock {Deep Petrov-Galerkin method for solving partial differential
  equations}.
\newblock {\em arXiv preprint arXiv:2201.12995}, 2022.

\bibitem{shirvany2008numerical}
Y.~Shirvany, M.~Hayati, and R.~Moradian.
\newblock {Numerical solution of the nonlinear Schrodinger equation by
  feedforward neural networks}.
\newblock {\em Communications in Nonlinear Science and Numerical Simulation},
  13(10):2132--2145, 2008.

\bibitem{shukla2020physics}
K.~Shukla, P.~C. Di~Leoni, J.~Blackshire, D.~Sparkman, and G.~E. Karniadakis.
\newblock {Physics-informed neural network for ultrasound nondestructive
  quantification of surface breaking cracks}.
\newblock {\em Journal of Nondestructive Evaluation}, 39:1--20, 2020.

\bibitem{sirignano2018dgm}
J.~Sirignano and K.~Spiliopoulos.
\newblock {DGM: A deep learning algorithm for solving partial differential
  equations}.
\newblock {\em Journal of computational physics}, 375:1339--1364, 2018.

\bibitem{smaoui2004modelling}
N.~Smaoui and S.~Al-Enezi.
\newblock {Modelling the dynamics of nonlinear partial differential equations
  using neural networks}.
\newblock {\em Journal of Computational and Applied Mathematics},
  170(1):27--58, 2004.

\bibitem{smith1985numerical}
G.~D. Smith.
\newblock {\em {Numerical solution of Partial Differential Equations: Finite
  Difference Methods}}.
\newblock Oxford University press, 1985.

\bibitem{smyrlis2004local}
G.~Smyrlis and V.~Zisis.
\newblock {Local convergence of the steepest descent method in Hilbert spaces}.
\newblock {\em Journal of mathematical analysis and applications},
  300(2):436--453, 2004.

\bibitem{sneddon2006elements}
I.~N. Sneddon.
\newblock {\em {Elements of Partial Differential Equations}}.
\newblock Courier Corporation, 2006.

\bibitem{srivastava2015highway}
R.~K. Srivastava, K.~Greff, and J.~Schmidhuber.
\newblock {Highway networks}.
\newblock {\em arXiv preprint arXiv:1505.00387}, 2015.

\bibitem{strauss2007partial}
W.~A. Strauss.
\newblock {\em {Partial Differential Equations: An introduction}}.
\newblock John Wiley \& Sons, 2007.

\bibitem{strouboulis2000design}
T.~Strouboulis, I.~Babu{\v{s}}ka, and K.~Copps.
\newblock {The design and analysis of the generalized finite element method}.
\newblock {\em Computer methods in applied mechanics and engineering},
  181(1-3):43--69, 2000.

\bibitem{sutskever2013momentum}
I.~Sutskever, J.~Martens, G.~Dahl, and G.~Hinton.
\newblock {On the Importance of Initialization and Momentum in Deep Learning}.
\newblock In {\em Proceedings of the 30th International Conference on Machine
  Learning - Volume 28}, ICML'13, page III–1139–III–1147. JMLR.org, 2013.

\bibitem{tarantola2005inverse}
A.~Tarantola.
\newblock {\em {Inverse Problem Theory and Methods for Model Parameter
  Estimation}}.
\newblock Society for Industrial and Applied Mathematics, Philadelphia, PA,
  Jan. 2005.

\bibitem{taylor2023deep2}
J.~M. Taylor, M.~Bastidas, D.~Pardo, and I.~Muga.
\newblock {Deep Fourier Residual method for solving time-harmonic Maxwell's
  equations}, 2023.

\bibitem{taylor2023deep1}
J.~M. Taylor, D.~Pardo, and I.~Muga.
\newblock {A Deep Fourier Residual method for solving PDEs using Neural
  Networks}.
\newblock {\em Computer Methods in Applied Mechanics and Engineering},
  405:115850, 2023.

\bibitem{tikhonov2013equations}
A.~N. Tikhonov and A.~A. Samarskii.
\newblock {\em {Equations of Mathematical Physics}}.
\newblock Courier Corporation, 2013.

\bibitem{tsoulos2008neural}
I.~Tsoulos, D.~Gavrilis, and E.~Glavas.
\newblock {Neural network construction and training using grammatical
  evolution}.
\newblock {\em Neurocomputing}, 72(1-3):269--277, 2008.

\bibitem{tsoulos2009solving}
I.~G. Tsoulos, D.~Gavrilis, and E.~Glavas.
\newblock {Solving differential equations with constructed neural networks}.
\newblock {\em Neurocomputing}, 72(10-12):2385--2391, 2009.

\bibitem{uchiyama1993solving}
T.~Uchiyama and N.~Sonehara.
\newblock {Solving inverse problems in nonlinear PDEs by recurrent neural
  networks}.
\newblock In {\em IEEE International Conference on Neural Networks}, pages
  99--102. IEEE, 1993.

\bibitem{uriarte2023deep}
C.~Uriarte, D.~Pardo, I.~Muga, and J.~Mu{\~n}oz-Matute.
\newblock {A Deep Double Ritz Method (D$^2$RM) for solving Partial Differential
  Equations using Neural Networks}.
\newblock {\em Computer Methods in Applied Mechanics and Engineering},
  405:115892, 2023.

\bibitem{uriarte2022finite}
C.~Uriarte, D.~Pardo, and {\'A}.~J. Omella.
\newblock {A Finite Element based Deep Learning solver for parametric PDEs}.
\newblock {\em Computer Methods in Applied Mechanics and Engineering},
  391:114562, 2022.

\bibitem{uriarte2023memory}
C.~Uriarte, J.~M. Taylor, D.~Pardo, O.~A. Rodr{\'\i}guez, and P.~Vega.
\newblock {Memory-Based Monte Carlo Integration for Solving Partial
  Differential Equations Using Neural Networks}.
\newblock In {\em International Conference on Computational Science}, pages
  509--516. Springer, 2023.

\bibitem{valueva2020application}
M.~V. Valueva, N.~Nagornov, P.~A. Lyakhov, G.~V. Valuev, and N.~I. Chervyakov.
\newblock {Application of the residue number system to reduce hardware costs of
  the convolutional neural network implementation}.
\newblock {\em Mathematics and computers in simulation}, 177:232--243, 2020.

\bibitem{van2007electromagnetic}
J.~G. Van~Bladel.
\newblock {\em {Electromagnetic fields}}, volume~19.
\newblock John Wiley \& Sons, 2007.

\bibitem{vaswani2017attention}
A.~Vaswani, N.~Shazeer, N.~Parmar, J.~Uszkoreit, L.~Jones, A.~N. Gomez,
  {\L}.~Kaiser, and I.~Polosukhin.
\newblock {Attention is all you need}.
\newblock {\em Advances in neural information processing systems}, 30, 2017.

\bibitem{versteeg2007introduction}
H.~K. Versteeg and W.~Malalasekera.
\newblock {\em {An Introduction to Computational Fluid Dynamics: The Finite
  Volume Method}}.
\newblock Pearson Education, 2007.

\bibitem{virtanen2020scipy}
P.~Virtanen, R.~Gommers, T.~E. Oliphant, M.~Haberland, T.~Reddy, D.~Cournapeau,
  E.~Burovski, P.~Peterson, W.~Weckesser, J.~Bright, et~al.
\newblock {SciPy 1.0: fundamental algorithms for scientific computing in
  Python}.
\newblock {\em Nature methods}, 17(3):261--272, 2020.

\bibitem{weinzierl2000introduction}
S.~Weinzierl.
\newblock {Introduction to Monte Carlo methods}.
\newblock {\em arXiv preprint hep-ph/0006269}, 2000.

\bibitem{wojtowytsch2020can}
S.~Wojtowytsch and E.~Weinan.
\newblock {Can shallow neural networks beat the curse of dimensionality? A mean
  field training perspective}.
\newblock {\em IEEE Transactions on Artificial Intelligence}, 1(2):121--129,
  2020.

\bibitem{yadav2015introduction}
N.~Yadav, A.~Yadav, M.~Kumar, et~al.
\newblock {\em {An introduction to neural network methods for differential
  equations}}, volume~1.
\newblock Springer, 2015.

\bibitem{yaman2013recent}
F.~Yaman, V.~G. Yakhno, and R.~Potthast.
\newblock {Recent Theory and Applications on Inverse Problems}.
\newblock {\em Mathematical Problems in Engineering}, 2013:303154, 2013.

\bibitem{yang2021b}
L.~Yang, X.~Meng, and G.~E. Karniadakis.
\newblock {B-PINNs: Bayesian physics-informed neural networks for forward and
  inverse PDE problems with noisy data}.
\newblock {\em Journal of Computational Physics}, 425:109913, 2021.

\bibitem{yarotsky2017error}
D.~Yarotsky.
\newblock {Error bounds for approximations with deep ReLU networks}.
\newblock {\em Neural Networks}, 94:103--114, 2017.

\bibitem{ye2018deep}
J.~C. Ye, Y.~Han, and E.~Cha.
\newblock {Deep Convolutional Framelets: A General Deep Learning Framework for
  Inverse Problems}.
\newblock {\em SIAM Journal on Imaging Sciences}, 11(2):991--1048, 2018.

\bibitem{zang2020weak}
Y.~Zang, G.~Bao, X.~Ye, and H.~Zhou.
\newblock {Weak adversarial networks for high-dimensional partial differential
  equations}.
\newblock {\em Journal of Computational Physics}, 411:109409, 2020.

\bibitem{zeiler2012adadelta}
M.~D. Zeiler.
\newblock {Adadelta: an adaptive learning rate method}.
\newblock {\em arXiv preprint arXiv:1212.5701}, 2012.

\bibitem{ziemianski2003hybrid}
L.~Ziemia{\'n}ski.
\newblock {Hybrid neural network/finite element modelling of wave propagation
  in infinite domains}.
\newblock {\em Computers \& structures}, 81(8-11):1099--1109, 2003.

\end{thebibliography}

\clearpage
\thispagestyle{empty}
\hfill
\clearpage

\end{document}